%% file: asymptotics_fr.tex
\documentclass{amsbook}

\usepackage{amssymb}
\usepackage{amsfonts}
\usepackage{graphicx}

\newtheorem{theorem}{Th\'{e}or\`{e}me}
\theoremstyle{plain}

\newtheorem{definition}{Definition}

\newtheorem{lemma}{Lemme}

\newtheorem{proposition}{Proposition}

\numberwithin{equation}{section}

\begin{document}
\Large
\pagenumbering{roman}
\begin{center}

\huge \textbf{Gane Samb LO}\\
\vskip 6cm
\Huge \textbf{Convergence Vague (IA)} \\
\bigskip
-\\
\bigskip \textbf{Suites de Vecteurs Al\'eatoires}
\vskip 6cm

\huge \textit{\textbf{Soci\'et\'e Africaine de Probabilit\'es et de Statistiques (SPAS) Editions Series}.\\
 \textbf{Saint-Louis, SENEGAL - Calgary, Alberta. 2016}}.\\

\bigskip \Large  \textbf{DOI} : http://dx.doi.org/10.16929/sbs/2016.0002\\
\bigskip \textbf{ISBN} 978-2-9559183-2-6
\end{center}

\newpage

\huge \textbf{SPAS Series Books}\\
\bigskip \bigskip

\Large \textbf{Advisers}\\

\bigskip \bigskip

\Large \textbf{List of published books}\\

\newpage

\noindent \textbf{Library of Congress Cataloging-in-Publication Data}\\

\noindent Gane Samb LO, 1958-\\

\noindent Convergence Vague (IA). Suites de Vecteurs Al\'eatoires.\\

\noindent SPAS Editions, 2016.\\

\noindent Copyright \copyright Société Africaine de Probabilité et de Statistiques (SPAS).\\

\noindent \textit{DOI} : 10.16929/sbs/2016.0002\\

\noindent \textit{ISBN}  978-2-9559183-2-6

\newpage

\noindent \textbf{Auteur : Gane Samb LO}\\
\bigskip

\bigskip
\noindent \textbf{Emails}:\\
\noindent gane-samb.lo@ugb.edu.sn, ganesamblo@ganesamblo.net.\\

\bigskip
\noindent \textbf{Url's}:\\
\noindent www.ganesamblo@ganesamblo.net\\
\noindent www.statpas.net/cva.php?email.ganesamblo@yahoo.com.\\

\bigskip \noindent \textbf{Affiliations}.\\
Affiliation principale : Université Gaston Berger, UGB, SENEGAL.\\
African University of Sciences and Technology, AUST, ABuja, Nigeria.\\
Chercheur associ\'e au : LSTA, Université Pierre et Marie, Paris VI, France.\\

\noindent \textbf{L'auteur enseigne ou a enseigné} au niveau Master dans les universités suivantes:\\
Saint-Louis, Sénégal (UGB)\\
Banjul, Gambié (TUG)\\
Bamako, Mali (USTTB)\\
Ouagadougou - Burkina Faso (UJK)\\
African Institute of Mathematical Sciences, Mbour, SENEGAL, AIMS.\\
Franceville, Gabon\\

\bigskip \noindent \textbf{Dédicaces}.\\

\noindent \textbf{ A mon épouse Mbaye Ndaw Fall, ma compagne depuis des décades}

\bigskip \noindent \textbf{Manifestation de reconnaissance de Soutien Financier et Matériel}.\\

\noindent L'auteur manifeste sa gratitude au Centre d'Excellence de la Banque Mondiale pour les Mathematiques, l'Informatique et les TIC 
(CEA-MITIC) pour un support continu de ses projets dans les années 2014, 2015 and 2016. Il remercie aussi les autorités de l'Université Gaston Berger pour le soutien permenant sous toutes ses formes.\\

\title{Convergence Vague (IA). Suites de Vecteurs Al\'eatoires}

\begin{abstract} \large (\textbf{English}) This monograph aims at presenting the core weak convergence theory for sequences of random vectors with values in $\mathbb{R}^k$. In some places, a more general formulation in metric spaces is provided. It lays out the necessary foundation that paves the way to applications in particular subfields of the theory. In particular, the needs of Asymptotic Statistics are addressed. A whole chapter is devoted to weak convergence in $\mathbb{R}$ where specific tools, for example for handling weak convergence of sequences using independent and indentically distributed random variables such that the Renyi's representations by means of standard uniform or exponential random variables, are stated. The function empirical process is presented as a powerful tool for solving a  considerable number of asymptotic problems in Statistics. The text is written in a self-contained approach whith the proofs of all used results at the exception of the general Skorohod-Wichura Theorem.\\

\noindent (\textbf{Fran\c{c}ais}) \noindent Cet ouvrage a l'ambition de pr\'esenter le noyau dur de la th\'eorie de la convergence vague de suite de vecteurs al\'eatoires dans $\mathbb{R}^k$. Autant que possible, dans certaines situations, la th\'eorie g\'en\'erale dans des espaces m\'etriques est donn\'ee. Il pr\'epare la voie \`a une sp\'ecialization dans certains sous-domaines de la convergence vague. En particulier, les besoins de la statistique asymptotique ont \'et\'e satisfaits. Un chapitre de l'ouvrage concerne la convergence vague dans $\mathbb{R}$ avec des outils sp\'ecifiques, par exemple, pour \'etudier les suites de variables al\'eatoires ind\'ependantes et identiquement distribu\'ees tels que la repr\'esentation de Renyi au moyen de variables al\'eatoires uniformes ou exponentielles standard. Le processus empirique function est introduit comme un outil puissant pour \'etudier des probl\`emes asymptotiques en Statistiques. Le texte est r\'edig\'e dans une approche auto-citante avec toutes les preuves des r\'esultats utilis\'es, \`a l'exception du Th\'eor\`eme de Skorohod-Wichura.\\

\noindent \textbf{Keywords.} Convergence vague; convergence en distribution; Théorème Portmanteau; Caractérisation d'une loi de probabilités; fonction de distribution, fonction de répartition; fonction caratéristiques; densité de probabilité; Marches aléatoires; Processus empirique; Loi Multinomiale; Compacité Relative; Tension Asymptotique; Tension uniform; Théorème de la transformation Continue; Représentation de Renyi et de Malmquist; Statistiques d'ordre; Méthodes Delta Multivariariées; Processus empirique fonctionnel.\\

\noindent \textbf{AMS 2010 Classification Subjects :} 60XXX; 62G30
\end{abstract}

\maketitle

\frontmatter
\tableofcontents
\mainmatter
\Large

\include{preface_gen_fr}
\include{asymptotics_cv_00_fr} 
\include{asymptotics_cv_01_fr} 
\include{asymptotics_cv_02_fr} 
\include{asymptotics_cv_03_fr} 
\include{asymptotics_cv_04_fr} 
\include{asymptotics_math_01_fr} 

\include{asymptotics_biblio_fr}
\end{document}

%% file: preface_gen_fr.tex
\noindent \textbf{Ce texte} fait partie d'une s\'{e}rie dont l'ambition est de brasser une
grande partie des probabilit\'{e}s \`{a} des fins p\'{e}dagogiques. Ces
textes permettront aux apprenant de se former tous seuls.

\noindent Ils pourront constituer pour les professeurs de documents de cours et
d'exercices. Pour les plus ambitieux, ils seront une base de d\'{e}part pour
des textes plus avanc\'{e}s et personnalis\'{e}s. Ils sont mis gracieusement 
\`{a} la disposition des apprenants et des ma\^{\i}tres. 

\bigskip 

\noindent \textbf{Nos ouvrages sont rang\'{e}s dans trois cat\'{e}gories :}

\bigskip 

\noindent \textbf{Une cat\'{e}gorie d'initiation pour les d\'{e}butants}. Il s'agit d'ouvrages
souvent accessibles d\`{e}s la premi\`{e}re ann\'{e}e d'universit\'{e} et ne
demandant pas de pr\'{e}-requis en math\'{e}matiques avanc\'{e}es. Les
ouvrages de probabilit\'{e}s \'{e}l\'{e}mentaires et de \ Statistiques \'{e}l\'{e}mentaires sont \`{a} ranger dansc cette cat\'{e}gorie. Cette initiation
pr\'{e}c\`{e}de les versions math\'{e}matiques de ces th\'{e}ories et pr\'{e}parent les applications de ces derni\`{e}res.

\bigskip 

\noindent \textbf{Une cat\'{e}gorie d'ouvrages d'applications et d'outils}. Des \'{e}%
tudiants ou des chercheurs dans des disciplines annexes comme l'\'{e}%
conomie, la m\'{e}decine, l'hydrologie, la finance, etc. peuvent besoin
d'outils assez avanc\'{e}s de la th\'{e}orie des probabilit\'{e}s ou des
statistiques. Ils sont plus int\'{e}ress\'{e}s plus par l'application des
outils que leur fondements ou leur d\'{e}veoppement. Des ouvrages adapt\'{e}%
s \`{a} ce besoin peuvent \^{e}tre compos\'{e}s. Un parfait exemple est un
ouvrage de statistiques math\'{e}matiques destin\'{e}s \`{a} des \'{e}%
conomistes n'ayant pas n\'{e}cessairement fait la th\'{e}orie de la mesure.

\bigskip 

\noindent \textbf{Une cat\'{e}gorie d'ouvrages sp\'{e}cialis\'{e}s}. Il s'agit d'ouvrages de niveau
internation rigoureusement \'{e}crits avec tous les argumentaires n\'{e}%
cessaires. Ces ouvrages sont con\c{c}us pour \^{e}tre consult\'{e}s partout
sur le globe dans leurs versions fran\c{c}aises et anglaises. Ils sont bas%
\'{e}s sur les ouvrages de base : Theorie de la mesure, Fondements math\'{e}%
matiques des probabilit\'{e}s. A partir de cette base, les ouvrages
seront-auto cit\'{e}s, ce qui veut dire que que toutes les math\'{e}matiques
utilis\'{e}es dans un ouvrage de cette s\'{e}rie, en dehors des notions de
premier cylce, seront d\'{e}montr\'{e}es quelque part dans un ouvrage de la
cat\'{e}gorie.

\bigskip 

\noindent La lecture de ces ouvrage ne requiert alors que le niveau de Licence. Un
lecteur capable de lire l'ouvrage de mesure et int\'{e}gration aura les
moyens de se sp\'{e}cialiser gr\^{a}ce \`{a} cette cat\'{e}gorie, dans
beaucoup de branches des probabilit\'{e}s et des statistiques.

\bigskip 

\noindent Sans insister, nous dirons que nous n'incluons la \textbf{cat\'{e}gorie des ouvrages
de recherche} qui sont et seront disponibles. Ces ouvrages partagent la m%
\^{e}me orientation, \`{a} savoir qu'ils sont \'{e}crits pour \^{e}tre lus
apr\_s les cours fondamentaux de mesure et des probabilit\'{e}s math\'{e}%
matiques.

\bigskip

\noindent Nos collaborateurs et anciens \'{e}l\`{e}ves sont pri\'{e}s de faire vivre
la chaine de sorte que le centre de Saint-Louis, donc le S\'{e}n\'{e}gal et
l'Afrique, soit un v\'{e}ritable creuset et une grande \'{e}cole de math\'{e}%
matiques, de probabilit\'{e}s et de statistique.

%% file: asymptotics_cv_00_fr.tex
\chapter{Revue de convergence vague dans $\mathbb{R}^k$} \label{ChapRevCvRk}

\section{Introduction} \label{cv.review.sec1}

Dans ce chapitre, nous allons voir que les lecteurs, pour la plupart d'entre
eux, ont d\'{e}j\`{a} \'{e}tudi\'{e} plusieurs exemples de convergence
vague. Ce qui leur manque peut \^{e}tre, c'est la coh\'{e}rence math\'{e}%
matique de cette th\'{e}orie et sa place dans le cadre g\'{e}n\'{e}ral des
convergences. Ces manquements seront combl\'{e}s dans ce cours.\\

\noindent Nous allons exhumer des r\'{e}sultats connus de convergence vague rencontr%
\'{e}s au cours du parcourt d\'{e}j\`{a} effectu\'{e}. Pour commencer, nous
allons emettre une assertion que nous ne serons en mesure de d\'{e}montrer
que dans le chapitre \ref{cv}, in particular dans le th\'eor\`eme \ref{cv.theo.portmanteau.rk} de ce chapitre.\\

\section{Convergence vague sur $\mathbb{R}^{k}$} \label{cv.review.sec2}

\noindent Commen\c{c}ons par rappeler que d'apr\`{e}s l'ouvrage \textbf{Base math\'{e}%
matiques des probabilit\'{e}s} \cite{bmtp} que la loi de probabilit\'{e}
d'un vecteur al\'{e}atoire $X:(\Omega ,A,\mathbb{P})\mapsto \mathbb{R}^{k}$
est caract\'{e}ris\'{e}e par

\bigskip \noindent \textbf{(a)} sa fonction de r\'{e}partition: \newline
\begin{equation*}
\mathbb{R}^{k}\ni x\hookrightarrow F_{X}(x)=\mathbb{P}(X\leq x),
\end{equation*}

\bigskip \noindent \textbf{(b)} sa fonction caract\'{e}ristique (ici, $i$
est le nombre imaginaire pur v\'{e}rifiant $i^{2}=-1$, et $<.,.>$ d\'{e}%
signe le produit scalaire classique sur $\mathbb{R}^{k}$) 
\begin{equation*}
\mathbb{R}^{k}\ni u\hookrightarrow \Phi (u)=E(exp(i<u,X>)),
\end{equation*}

\bigskip \noindent \textbf{(c)} La fonction des moments (si elle existe dans
un voisinage du vecteur nul) 
\begin{equation*}
\mathbb{R}^{k}\ni u\hookrightarrow \Psi _{X}(x)=E(exp(<u,X>)).
\end{equation*}

\bigskip \noindent et

\bigskip \noindent \textbf{(d)} par sa d\'{e}riv\'{e}e de Radon-Nikodym, (si elle
existe), par rapport \`{a} une mesure $\sigma $-finie $\nu $ sur $\mathbb{R}%
^{k}$,  
\begin{equation*}
d\mathbb{P}/d\nu =f_{X}.
\end{equation*}

\noindent Il est int\'{e}ressant que ces caract\'{e}risations s'\'{e}tendent \`{a} la convergence vague ainsi, dans un th\'{e}or\`{e}me que nous d\'{e}montrerons plus tard dans le th\'eor\`eme \ref{cv.theo.portmanteau.rk} du chapitre \ref{cv}.\\ 

\begin{theorem} (\textbf{THEOREME - DEFINITION - LEMME}) \label{cv.review.rk} Soit une suite vecteurs
al\'{e}atoires $X_{n} :(\Omega _{n},\mathcal{A}_{n},\mathbb{P}_{n})\mapsto (\mathbb{R}^{k},\mathbb{B}(\mathbb{R}^{k}))$ et un autre vecteur al\'{e}atoire $X : (\Omega_{\infty} ,\mathcal{A}_{\infty},\mathbb{P}_{\infty})\mapsto (\mathbb{R}^{k}, \mathbb{B}(\mathbb{R}^{k})).$ Alors les propri\'{e}t\'{e}s (a) et (b) suivantes sont \'{e}quivalentes \newline

\bigskip \noindent \textbf{(a)} Pour tout $u\in \mathbb{R}^{k}$,\newline

\begin{equation*}
\Phi _{X_{n}}(u)\rightarrow \Phi _{X}(u)\text{ quand }n\rightarrow +\infty 
\end{equation*}

\bigskip \noindent \textbf{(b)} Pour tout point de continuit\'{e} $x\in 
\mathbb{R}^{k}$, 
\begin{equation*}
F_{X_{n}}(x)\rightarrow F_{X}(x)\text{ quand }n\rightarrow +\infty .
\end{equation*}

\bigskip \noindent Si l'un de ces deux points a lieu, alors nous disons que $%
X_{n}$ converge vaguement vers $X,$ ou $X_{n}$ converge en distribution vers 
$X$ ou $X_{n}$ converge en loi $X,$ not\'{e}e%
\begin{equation*}
X_{n}\rightsquigarrow X\text{ ou }X_{n}\longrightarrow _{d}X\text{ ou }X_{n}%
\overset{\mathcal{L}}{\longrightarrow }X\text{ ou }X_{n}\overset{w}{%
\longrightarrow }X\text{ ou }X_{n}\longrightarrow _{w}X
\end{equation*}%

\noindent (ici $w$ fait r\'{e}ference au term anglais weak convergence).\\

\bigskip \noindent Nous avons des conditions suffisantes de convergence vague :\\

\bigskip \noindent \textbf{(c)} Si de plus, les fonction des moment  $\Psi
_{X_{n}}$ sont definies sur $B_{n}$, $n\geq 1$ et $\Psi _{X}$ est d\'{e}%
finie sur $B$, o\`{u} les $B_{n}$ et $B$ sont des voisinages de $0$, et si
pour tout $x\in B$, 
\begin{equation*}
\Psi _{X_{n}}(x)\rightarrow \Psi _{X}(x)\text{ quand }n\rightarrow +\infty 
\end{equation*}

\noindent alors $X_{n}$ converge vaguement $X$.\newline

\bigskip \noindent \textbf{(d)} Enfin, si les distributions admettent des d%
\'{e}riv\'{e}es de Radon-Nikodym, c'est-\`{a}-dire des densit\'{e}s de
probabilit\'{e}s par rapport \`{a} la m\^{e}me mesure $\sigma $-finie $\nu $%
, not\'{e}es $d\mathbb{P}_{n}/d\nu =f_{X_{n}}$ et $d\mathbb{P}/d\nu =f_{X}$
existent et si pour tout $x\in D_{X}=\{x,f_{X}(x)>0\}$, 
\begin{equation*}
f_{X_{n}}(x)\rightarrow f_{X}(x)\text{ as }n\rightarrow +\infty ,
\end{equation*}

\noindent alors $X_{n}\rightsquigarrow X$.\newline

\noindent  Nous avons ce dernier point.\\

\noindent \textbf{(e)} Supposons que $\{X_n, \ \ n\geq 1\} \subset \mathbb{R}^k$ converge vaguement vers $X \in \mathbb{R}^k$, quand $n \rightarrow +\infty$ et soit $A$ une matrice r\'elle de $m$ lignes et $k$ colonnes avec $m\geq 1$. Alors $\{AX_n, \ \ n\geq 1\} \subset \mathbb{R}^m$ converge vaguement vers $AX \in \mathbb{R}^m$.
\end{theorem}

\noindent \textbf{Remarque}. Le point (e) ci-dessus est une cons\'equence du th\'eor\`eme \ref{cv.mappingTh} de la transformation continue 
 \'etablie et prouv\'ee dans le chapitre \ref{cv}.\\

\bigskip

\bigskip \noindent \textbf{En r\'{e}sum\'{e}}, la convergence en loi sur $%
\mathbb{R}^{k}$ a lieu lorsque les fonctions de r\'{e}partition, les
fonctions caract\'{e}ristiques, les densit\'{e}s de probabilit\'{e}s ou les
fonctions des moments convergent vers celle d'une loi de probabilit\'{e},
pourvu que dans les deux derniers cas, les fonctions en question existent.\\

\bigskip \noindent  Tout cela est \'{e}norme. Cela nous sonne l'occasion d'aller directement sur
les exemples connus \ de convergence de ces fonctions et d'en ajouter de
nouveaux.\\

\noindent Avant d'aller plus loin, nous aurons souvent besoin de cet outil int\'eressant pour obtenir la convergence vague \`a partir de la 
convergence de fonction caract\'eeristique.\\

\begin{proposition} \label{cv.wold} \textbf{Cri\`ere de Wold}. La suite de vecteurs al\'eatoires $\{X_n, \ \ n\geq 1\} \subset \mathbb{R}^k$ converge vaguement vers to $X \in \mathbb{R}^k$, quand $n \rightarrow +\infty$ si et seulement si pour $a \in \mathbb{R}^k$, la suite $\{<a,X_n>, \ \ n\geq 1\} \subset \mathbb{R}$ converge veguement vers $X \in \mathbb{R}$ quand $n \rightarrow +\infty$.
\end{proposition}

\noindent \textbf{Preuve}. La preuve est rapide. Elle utilise les notations ant\'erieures. Supposons que $X_n$ converge faiblement vers $X$ 
dans $\mathbb{R}^k$ quand $n \rightarrow +\infty$. En utilisant la convergence des fonctions carat\'eristiques, nous avons pour tout $u\in \mathbb{R}^k$,
$$
\mathbb{E}(exp(i<X_n,u>) \rightarrow \mathbb{E}(exp(i<X,u>) \ \ quand \ \ n \rightarrow +\infty.
$$  

\noindent Il s'en suit que $a\in \mathbb{R}^k$ et pour tout $t \in \mathbb{R}$, nous avons

\begin{equation}
\mathbb{E}(exp(it<X_n,a>) \rightarrow \mathbb{E}(exp(it<X,a>) \ \ quand \ \ n \rightarrow +\infty. \label{cv.proj}
\end{equation}

\noindent c-a-d que, en prenant $u=ta$ dans la formule pr\'ec\'edente, et en notant $Z_n=<X_n,a>$ et $Z=<X,a>$, 

$$
\mathbb{E}(exp(itZ_n) \rightarrow \mathbb{E}(exp(itZ) \ \ quand \ \ n \rightarrow +\infty.
$$

\noindent Cela signifie que $Z_n \rightsquigarrow Z$, qui est \'egal à $<a,X_n>$, converges vaguement vers $<a,X>$.\\

\noindent Inversement, supposons que pour tout $a \in \mathbb{R}^k$, la suite $\{<a,X_n>, \ \ n\geq 1\} \subset \mathbb{R}$ converge vaguement vers $X \in \mathbb{R}$ as $n \rightarrow +\infty$. Alors, en prenant $t=1$ dans (\ref{cv.proj}), nous obtenons pour tout $a=u \in \mathbb{R}^k$,

$$
\mathbb{E}(exp(i<X,u>) \rightarrow \mathbb{E}(exp(i<X,u>) \ \ quand \ \ n \rightarrow +\infty.
$$

\noindent ce qui signifie que $X_n \rightsquigarrow +\infty$ quand $n\rightarrow +\infty$.\\

\section{Examples de convergence vague dans $\mathbb{R}$} \label{cv.review.sec3}

\subsection{Convergence de la loi hyperg\'{e}om\'{e}trique vers la loi binomiale} \label{cv.review.subsec.HypBin}

Soit $X_{N}$ suivant une loi hyperg\'{e}om\'{e}trique $\mathcal{H}(N,M,n)$
avec $M/N\rightarrow p,$ $N\rightarrow \infty ,$ n restant fixe$.$ Alors $%
X_{N}$ tend en loi vers une variable al\'{e}atoire X suivant une loi
binomiale $\mathcal{B}(n,p)$.\\

\noindent \textbf{PREUVE}. Pour prouver cela, utilisons les densit\'{e}s par rapport \`{a} la mesure de
comptage $\nu $ sur $\mathbb{N}$. Nous avons :%
\begin{equation*}
f_{X_{n}}(k)=\frac{\left( 
\begin{tabular}{l}
$M$ \\ 
$k$%
\end{tabular}%
\right) \left( 
\begin{tabular}{l}
$N-M$ \\ 
$n-k$%
\end{tabular}%
\right) C_{M}^{k}C_{M-N}^{n-k}}{\left( 
\begin{tabular}{l}
$N$ \\ 
$n$%
\end{tabular}%
\right) },0\leq k\leq \min (n,M).
\end{equation*}

\noindent Supposons que $M/N\rightarrow p,$ $N\rightarrow \infty$. Nous aurons%
\begin{eqnarray*}
f_{X_{n}}(k) &=&\frac{M!}{k!(M-k)!}\frac{(M-N)!}{(n-k)!(N-M-(n-k))!}\frac{%
n!(N-n)!}{N!} \\
&=&\frac{n!}{k!(n-k)!}\times \frac{M!}{(M-k)!}\times \frac{!(N-M-(n-k))!}{%
(N-M-(n-k))!}\times \frac{(M-n)!}{N!} \\
&=&\left( 
\begin{tabular}{l}
$n$ \\ 
$k$%
\end{tabular}%
\right) \times \left\{ \frac{M!}{(M-k)!}\right\} \left\{ \frac{(N-M)!}{%
(N-M-(n-k))!}\right\} \left\{ \frac{(M-n)!}{N!}\right\} .
\end{eqnarray*}%

\noindent Mais%
\begin{eqnarray*}
\left\{ \frac{M!}{(M-k)!}\right\}  &=&(M-k+1)(M-k+2)...(M+1)M \\
&=&M^{k}(1-\frac{k-1}{M})(1-\frac{k-2}{M})\times ...\times (1-\frac{1}{M}) \\
&=&M^{k}(1+o(1))
\end{eqnarray*}

\noindent puisque\ $M\rightarrow \infty $ et $k$ est fixe. Ensuite%
\begin{eqnarray*}
\left\{ \frac{(N-M)!}{(N-M-(n-k))!}\right\}  &=&(N-M-(n-k)+1)\times
...\times (N-M-1)(N-M) \\
&=&(N-M)^{n-k}(1+\frac{n-k-1}{N-M})(1+\frac{n-k-2}{N-M})\times ...\times ((1+%
\frac{1}{N-M}) \\
&=&(N-M)^{n-k}-1+o(1))
\end{eqnarray*}

\noindent puisque, aussi, $N-M=N(1-M/N)\sim N(1-p)\rightarrow \infty $ et $n-k$ est fix%
\'{e}. Enfin%
\begin{eqnarray*}
\left\{ \frac{(M-n)!}{N!}\right\}  &=&\frac{1}{(N-n+1)(N-n+2)...(N-1)N} \\
&=&\frac{1}{N^{n}(1-\frac{n-1}{N})(1-\frac{n-2}{N})...(1-\frac{1}{N})} \\
&=&\frac{1}{N^{n}(1=o(1))}.
\end{eqnarray*}

\noindent pour des raisons similaires. Au total, pour tout $0\leq k\leq n$%
\begin{equation*}
f_{X_{n}}(k)=\left( 
\begin{tabular}{l}
$n$ \\ 
$k$%
\end{tabular}%
\right) \left( \frac{M}{N}\right) ^{k}\left( \frac{N-M}{N}\right)
^{k}(1+o(1))\rightarrow \left( 
\begin{tabular}{l}
$n$ \\ 
$k$%
\end{tabular}%
\right) p^{k}(1-p)^{n-k}.
\end{equation*}

\noindent Ainsi pour tout point $k$ du domaine de la densit\'{e} de la loi binomiale
par rapport \`{a} la mesure de comptage $\nu ,$ not\'{e}e 
\begin{equation*}
f_{X}(k)=\left( 
\begin{tabular}{l}
$n$ \\ 
$k$%
\end{tabular}%
\right) p^{k}(1-p)^{n-k},
\end{equation*}

\noindent nous avons%
\begin{equation*}
\forall (1\leq k\leq n),f_{X_{n}}(k)\rightarrow f_{X}(k).
\end{equation*}

\bigskip \noindent Nous avons donc la convergence en loi cherch\'{e}e.\\

\noindent \textbf{Remarque utile en th\'eorie des sondages}. Cette approximation permet de consid\'{e}rer que le tirage avec remise et
celui sans remise sont \'{e}quivalente dans une enqu\^{e}te avec une population tr\`{e}s large. C'est un peu la notion suivante : Lorsque la
population est tr\`{e}s grande, dans un tirage avec remise d'un nombre d'individu assez modeste, il est presque impossible qu'un m\^{e}me individu
sorte plus d'une fois.\\

\bigskip \bigskip

\subsection{Convergence de la loi binomiale vers la loi de Poisson} \label{cv.review.subsec.BinPois}

Soit $X_{n}$ suivant une loi $\mathcal{B}(n,p)$ \ aec $p=p_{n}\rightarrow 0$
et $np_{n}\rightarrow \lambda ,$ $0<\lambda ,$ quand $n\rightarrow \infty .$
Alors $X_{n}$ tend en loi vers une variable al\'{e}atoire $X$ suivant une
loi de Poisson de param\`{e}tre $\lambda$.\\

\noindent \textbf{PREUVE}. Pour attester cela, utilisons les fonctions des moments. Soit $X$ une
variable suivant la loi de Poisson de param\`{e}tre $\lambda >0.$ Nous avons 
\begin{equation*}
\Psi _{X_{n}}(t)=(p_{n}+(1-p_{n})e^{t})^{n}\text{ pour }n\geq 1;\Psi
_{X}(t)=\exp (\lambda (e^{t}-1)),t\in \mathbb{R}.
\end{equation*}

\noindent Notons $\lambda _{n}=np_{n}\rightarrow \lambda .$ Pour tout $t$ fix\'{e},%
\begin{equation*}
\Psi _{X_{n}}(t)=(\frac{\lambda _{n}}{n}+(1-\frac{\lambda _{n}}{n}%
)e^{t})^{n}=\left( 1-\frac{\lambda _{n}(e^{t}-1)}{n}\right) ^{n}\rightarrow
\exp (\lambda (e^{t}-1))=\Psi _{X}(t)
\end{equation*}

\noindent par le r\'{e}sultat classique d'analyse qui affirme que 
\begin{equation*}
(1+\frac{x_{n}}{n})^{n}\rightarrow e^{x}\text{ quand }n\rightarrow +\infty 
\text{ pourvu que }x_{n}\rightarrow x\in \mathbb{R}\text{ lorsque }%
n\rightarrow +\infty .
\end{equation*}

\subsection{Convergence de la loi de Poisson vers la loi normale} \label{cv.review.subsec.PoisGauss}

Soit $Z_{\lambda }$ une variable al\'{e}atoire suivant la loi de Poisson de
param\'etre $\lambda $ : $Z_{\lambda }\sim \mathcal{P}(\lambda ).$ Alors la
variable 
\begin{equation*}
\frac{Z_{\lambda }-\lambda }{\sqrt{\lambda }}
\end{equation*}

\noindent converge en loi vers une variable al\'eatoire $X$ suivant une loi normale standard, i.e $X \sim \mathcal{N}(0,1)$, quand 
$\lambda \rightarrow +\infty$.\\

\bigskip \noindent \textbf{PREUVE}. Nous utilisons les fonctions g\'{e}n\'{e}ratrices des moments.\ Rappelons la
fonction des moments de $Z_{\lambda }\sim \mathcal{P}(\lambda )$ : 
\begin{equation*}
\Psi _{Z_{\lambda }}(t)=\exp (\lambda (e^{t}-1)).
\end{equation*}

\noindent Soit 
\begin{equation*}
Y(\lambda )=\frac{Z-\lambda }{\sqrt{\lambda }}=\frac{Z-\mathbb{E}(X)}{\sigma
_{Z}}.
\end{equation*}

\noindent Nous avons 
\begin{equation*}
\Psi _{Y(\lambda )}(u)=e^{-\sqrt{\lambda }}\times \varphi _{Z}(u/\sqrt{%
\lambda })=e^{-\sqrt{\lambda }}\times \exp (\lambda (e^{u/\sqrt{\lambda }%
}-1)).
\end{equation*}

\noindent Quand $\lambda \rightarrow \infty ,$ nous pouvons d\'{e}velopper 
\begin{equation*}
\lambda (e^{u/\sqrt{\lambda }}-1)=\lambda (1+\frac{u}{\sqrt{\lambda }}+\frac{%
u^{2}}{2\lambda }+O(\lambda ^{-3/2})-1
\end{equation*}%
\begin{equation*}
=u\sqrt{\lambda }+\frac{u^{2}}{2}+O(\lambda ^{-1/2}).
\end{equation*}

\noindent Donc 
\begin{equation*}
\Psi _{Y(\lambda )}(u)=\exp (\frac{u^{2}}{2}+O(\lambda ^{-1/2}))\rightarrow
\exp (u^{2}/2).
\end{equation*}

\noindent On conclut aussi que 
\begin{equation*}
\frac{Z-\lambda }{\sqrt{\lambda }}\rightarrow \mathcal{N}(0,1)
\end{equation*}%
quand $\lambda \rightarrow \infty$.

\bigskip \bigskip

\subsection{Convergence de la loi binomiale vers la loi normale} \label{cv.review.subsec.Bin}

Soit $X_{n}$ suivant une loi $\mathcal{B}(n,p)$ \ avec $p\in ]0,1[$ fix\'{e}%
. Alors quand $n\rightarrow \infty $ 
\begin{equation}
Z_{n}=\frac{X_{n}-np}{\sqrt{npq}}\rightsquigarrow N(0,1) \ \ a \ \ as n\rightarrow +\infty. \label{cv.cltBin}
\end{equation}

\noindent \textbf{Preuve}. Utilisons les fonctions des moments. $Soit$ $X\sim \mathcal{B}(n,p).$ Nous
avons 
\begin{equation*}
\Psi _{X_{n}}(u)=(q+pe^{u})^{n}.
\end{equation*}%
D'o\`{u} 
\begin{equation}
\Psi _{(X_{n}-np)/\sqrt{npq}}(u)=e^{-\sqrt{np/q}}\times \Psi _{X_{n}}(u/%
\sqrt{npq})  \label{binom00}
\end{equation}

\noindent avec 
\begin{equation*}
\Psi _{X}(u/\sqrt{npq})=(q+pe^{u/\sqrt{npq}})^{n}.
\end{equation*}

\noindent L'id\'{e}e de la suite des calculs est d'utiliser un d\'{e}veloppement
d'ordre 2 de $e^{u/\sqrt{npq}}$ au voisinage de 0 quand $n\rightarrow \infty 
$ et $u$ fix\'{e}. Ensuite l'expression obtenue sera de la forme $1+v_{n}$, o%
\`{u} $v_{n}$ tend vers z\'{e}ro. Ensuite, un developpement de $\log $($%
1+v_{n})$ d'ordre 2 est op\'{e}r\'{e}.\\

\noindent Ainsi, quand $n\rightarrow \infty $ et u fix\'{e}, 
\begin{equation*}
e^{u/\sqrt{npq}}=1+\frac{u}{\sqrt{npq}}+\frac{u^{2}}{2npq}+O(n^{-3/2}).
\end{equation*}

\noindent D'o\`{u} 
\begin{equation*}
(q+pe^{u/\sqrt{npq}})=1+u\sqrt{p/nq}+\frac{u^{2}}{2nq}+O(n^{-3/2})=1+v_{n}
\end{equation*}

\noindent avec 
\begin{equation*}
v_{n}=u\sqrt{p/nq}+\frac{u^{2}}{2nq}+O(n^{-3/2})\rightarrow 0.
\end{equation*}

\noindent Si bien que 
\begin{eqnarray*}
\log (1+u\sqrt{p/nq}+\frac{u^{2}}{2nq}+O(n^{-3/2}))&=&\log (1+v_{n})\\
&=&v_{n}-\frac{1}{2}v_{n}^{2}+O(v_{n}^{3})\\
&=&u\sqrt{p/nq}+\frac{u^{2}}{2nq}-\frac{pu^{2}}{2nq}+O(n^{-3/2}).
\end{eqnarray*}

\noindent Alors
\begin{eqnarray*}
\Psi _{X_{n}}(u/\sqrt{npq})&=&(q+pe^{u/\sqrt{npq}})^{n}=\exp (n\log (q+pe^{u/\sqrt{npq}}))\\
&=&\exp (n(u\sqrt{p/nq}+\frac{u^{2}}{2nq}-\frac{pu^{2}}{2nq}+O(n^{-3/2})))\\
&=&\exp (u\sqrt{np/q}+\frac{u^{2}}{2q}-\frac{pu^{2}}{2q}+O(n^{-1/2})))\\
&=&e^{u\sqrt{np/q}}e^{u^{2}/2+O(n^{-1/2})}.
\end{eqnarray*}

\noindent En retournant \`{a} (\ref{binom00}), on arrive \`{a} 
\begin{equation*}
\Psi _{(X_{n}-np)/\sqrt{npq}}(u)\rightarrow \exp (u^{2}/2).
\end{equation*}

\noindent D'o\`{u} l'approximation 
\begin{equation*}
(\beta (n,p)-np)/\sqrt{npq}\rightarrow \mathcal{N}(0,1)
\end{equation*}

\noindent QED.\\

\bigskip \noindent \textbf{Remarque}. Nous reviendrons sur une deuxi\`eme preuve directe de ce r\'esultat en utilisant le th\'eor\`eme central limite standard ci-dessous.\\

\bigskip \bigskip

\subsection{Th\`{e}or\`{e}me Central Limite Standard dans $\mathbb{R}$.} \label{cv.review.subsec.cltR}

Les deux cas d\'{e}j\`{a} vus sont des cas sp\'{e}ciaux d'un cas g\'{e}n\'{e}%
ral. En effet, si $(X_{n})_{n\geq 1}$ est une suive de variables al\'{e}%
atoires r\'{e}elles ayant des moments de second ordre finis, on peut
s'attendre \`{a} ce que 
\begin{equation*}
\frac{X_{n}-E(X_{n})}{\sigma _{X_{n}}}
\end{equation*}

\noindent converge en loi vers la loi normale centr\'{e}e r\'{e}duite. Ceci n'est pas
toujours vraie. Mais si cela est le cas, nous dirons qu'on a un th\'{e}or%
\`{e}me central limite (Central limit th\'{e}or\`{e}me, CLT). Cela est vrai
pour l'\'{e}chantillon dans le cas suivant.\\

\noindent Soit $X_{1},X_{2},...$ une suite de variables al\'{e}toires ind\'{e}pendantes, identiquement distribu\'{e}es selon la fonction de r\'{e}partition $F$ avec 
\begin{equation*}
E(X_{i})=\mu =\int xdF(x)=0,\sigma _{X_{i}}^{2}=\sigma ^{2}=\int (x-\mu
)^{2}dF(x).
\end{equation*}

\noindent Posons, pour $n\geq 1,$%
\begin{equation*}
S_{n}=X_{1}+...+X_{n}.
\end{equation*}

\noindent Alors quand $n\rightarrow \infty$%
\begin{equation*}
\frac{S_{n}}{\sqrt{n}}\rightarrow \mathcal{N}(0,1).
\end{equation*}

\noindent  \textbf{PREUVE}. Consid\'{e}rons les fonctions caract\'{e}ristiques 
\begin{equation*}
\mathbb{R}\ni u\hookrightarrow \Phi _{X_{i}}(u)=E(e^{iuX_{i}})=\Psi (u).
\end{equation*}%

\noindent Par l'existence des moments \`{a} l'ordre $2$, on a le d\'{e}veloppement
limit\'{e} \`{a} l'ordre $2$, 
\begin{equation*}
\Phi (u)=1+iu\Phi ^{\prime }(0)+\frac{1}{2}u^{2}\Phi ^{\prime \prime
}(0)+O(u^{3})
\end{equation*}%
\begin{equation*}
=1-\frac{1}{2}u^{2}+O(u^{2})
\end{equation*}%

\noindent puisque 
\begin{equation*}
\Phi ^{\prime }(0)=i\text{ }\mathbb{E}(X)=0,\text{ }\Phi ^{\prime \prime
}(0)=-\mathbb{E}(X^{2})=-1.
\end{equation*}

\noindent D\`{e}s lors 
\begin{equation*}
\Phi _{S_{n}/\sqrt{n}}(u)=(\Phi (u/\sqrt{n}))^{n}.
\end{equation*}

\noindent Pour $u$ fix\'{e} et $n\rightarrow \infty ,$%
\begin{eqnarray*}
\Phi _{S_{n}/\sqrt{n}}(u) &=&(\Phi (u/\sqrt{n}))^{n}=\exp (n\log (1-\frac{%
u^{2}}{2n}+O(n^{-3/2})) \\
&=&\exp (n(-\frac{u^{2}}{n}+O(n^{-3/2})) \\
&=&\exp (-u^{2}/2+O(n^{-1/2})) \\
&\rightarrow &\exp (-u^{2}/2)
\end{eqnarray*}

\noindent Nous venons d'\'{e}tablir que 
\begin{equation*}
\frac{S_{n}}{\sqrt{n}}\rightarrow \mathcal{N}(0,1).
\end{equation*}

\noindent Dans le cas g\'{e}n\'{e}ral non centr\'e et non normalis\'e, nous avons le Th\'{e}or\`{e}me central limite

\begin{equation*}
\frac{1}{\sigma \sqrt{n}} (S_{n}-n\mu) \rightarrow \mathcal{N}(0,1).
\end{equation*}

\bigskip \noindent Donnons deux exemples importants du th\'eor\`eme central limite standard relatifs aux exp\'eriences de Bernouilli.\\

\noindent \textbf{Example 1 : Convergence vague de la loi binomiale}.\\

\noindent Nous allons donner une autre preuve du r\'esultat (\ref{cv.cltBin}) de la sous-section \ref{cv.review.subsec.Bin} relatif à la loi limite d'une suite de variables al\'eatoires suivant une loi binomiale lorsque le nombre d'essais $n$ cro\^it ind\'efin\'ement pendant que la probabilit\'e de succ\`es $\in ]0,1[$ reste fix\'e. Nous gardons les notations de cette sous-section.\\

\noindent En nous fondant sur les cours de probabilit\'es \'elementaires que nous pouvons trouver dans un grand nombre d'ouvrages, en particulier dans \cite{ept-en2016}, de la pr\'esente s\'erie de Th\'eorie de Probabilit\'e et de Statistiques, au chapitre, Lemme 1, que si $X_n \sim \mathcal{B}(n,p)$, alors  $X_{n}$ est la somme de $n$ variables al\'eatoires $Y_{1},...,Y_{n}$, identiquement distribu\'ees selon un loi de Bernouilli $\mathcal{B}(p)$ random variables, i.e.,

\begin{equation*}
X_{n}=Y_{1}+...+Y_{n}.
\end{equation*}

\noindent Pour chaque $Y_{i}$, $1\leq i \leq n$, nous avons

\begin{equation*}
\mathbb{E}(Y_{i})=p\text{ and }\sigma ^{2}=\mathbb{V}ar(Y_{i})=pq\text{ where }q=1-p.
\end{equation*}

\noindent D\`es lors, la variable $Z_n$ de la formule (\ref{cv.cltBin}) devient

\begin{equation*}
Z_{n}=\frac{X_{n}-np}{\sqrt{npq}}=\frac{1}{\sigma \sqrt{n}}
\sum_{i=1}^{n}(Y_{i}-\mathbb{E}(Y_{i}).
\end{equation*}

\noindent Ainsi, la convergence vague de $Z_{n}$ to $\mathcal{N}(0,1)$  quand  $n \rightarrow +\infty$ d\'ecoule de l'application du th\'eor\`eme central limite sur $\mathbb{R}$.\\

\bigskip \noindent \textbf{Remarque}. Cette preuve est simple et belle. La premi\`ere est toujours utile. Simplement parce que nous pouvons être appel\'es à utiliser ou \`a enseigner ce r\'esultat \`a un niveau o\`u le th\'eor\`eme central limite n'est pas disponible. Au del\`a de cette raison, cette preuve fait partie de l'histoire des probabilit\'es. Dans le m\^eme esprit, les premi\`eres d\'ecouvertes de cette loi remontent en 
1732 avec \textit{de} Moivre et en 1801 avec Laplace (voir Lo\`{e}ve \cite{loeve}, page 23). Ces m\'ethides historiques peuvent \^etre revisit\'es dans dans \cite{ept-en2016} ou dans \cite{ept-fr2016} avec une r\'edaction adapt\'ee au niveau de la premi\`ere ann\'ee universitaire.\\

\bigskip \noindent \textbf{Example 2 : Loi Binomiale N\'egative}.\\

\noindent Pour un entier $k\geq 1$ fix\'ee, une variable suivant la loi Binomiale N\'egative $X_{k}$ peut \^etre d\'efinie relativement aux essais de  Bernouilli avec une probabilit\'e de succ\`es $p \in ]0,1[$. Le nombre d'essais ind\'ependants de l'exp\'erience de Bernouilli de $p$ n\'ecessaires pir avoir $k$ succ\`es, suit par d\'efinition la loi Binomiale N\'egative de param\`etres $k$ et $p$, not\'ee  $X_{k} \sim \mathcal{NB}(k,p)$. Pour $k=1$, la  $X_{1}$ suit une loi g\'eom\'etrique de param\`etre $p$, not\'ee $X_1 \sim \mathcal{G}(p)$.\\

\noindent De m\^eme que pour la suite de variables al\'eatoires binomiales, nous pouvons appliquer le   à une suite de variables al\'eatoires ninomiales n\'egatives  $X_k$, $k\geq 1$ pour avoir le r\'esultat

\begin{equation}
Z_{n}=\frac{p(X_{k}-\frac{k}{p})}{\sqrt{nq}} \rightsquigarrow \mathcal{N}(0,1) \text{ quand } k \rightarrow +\infty. \label{cv.cltNegBin01}
\end{equation} 

\noindent A cet effet, le lecteur peut trouver dans les cours de probabilit\'es \'elementaires de son choix,  en particulier dans \cite{ept-en2016}, de la pr\'esente s\'erie de Th\'eorie de Probabilit\'e et de Statistiques, au chapitre 2, que le lemme 2 assure que la variable al\'eatoire $X_{k}$  suivant la loi $\mathcal{NB}(k,p)$ est la somme de $k$ variables al\'eatoires $Y_{1},...,Y_{k}$, ind\'ependantes et identiquement distribu\'ees selon la loi g\'eom\'etrique $\mathcal{G}(p)$, i.e.,

\begin{equation*}
X_{k}=Y_{1}+...+Y_{n},
\end{equation*}

\noindent et que pour chacune des variables $Z_i$, nous avons

\begin{equation*}
\mathbb{E}(Y_{i})=\frac{1}{p}\text{ et }\sigma ^{2}=\mathbb{V}ar(Y_{i})=\frac{q}{p^{2}}\text{ o\`u }q=1-p.
\end{equation*}

\noindent Alors, en appliquant le th\'eor\`eme central limite, nous obtenons

\begin{equation*}
Z_{n}=\frac{p(X_{k}-\frac{k}{p})}{\sqrt{nq}}=\frac{1}{\sigma \sqrt{n}}\sum_{i=1}^{n}(Y_{i}-\mathbb{E}(Y_{i}))\rightsquigarrow \mathcal{N}(0,1) \text{ as } k \rightarrow +\infty. \label{cv.cltNegBin02}
\end{equation*}

\noindent Ceci prouve (\ref{cv.cltNegBin01}).\\

\bigskip \bigskip

\subsection{Lois limites des valeurs extr\^{e}mes} \label{cv.review.subsec.evt}

\bigskip \noindent Consid\'{e}rons $X_{1},X_{2},$ etc...,\ une suite de variables al\'{e}%
atoires r\'{e}elles ind\'{e}pendantes de m\^{e}me fonction de r\'{e}%
partition $F.$ Consid\'{e}rons pour chaque $n\geq 1$

\begin{equation*}
M_{n}=\max (X_{1},...,X_{n}).
\end{equation*}

\bigskip \noindent Rappelons tout de suite que nous avons 
\begin{equation*}
P(M_{n}\leq x)=F(x)^{n},x\in \mathbb{R}.
\end{equation*}

\noindent La th\'{e}orie des valeurs extr\^{e}mes a commenc\'{e} ses beaux jours par
la d\'{e}couverte des lois limites de la suite $M_{n}$ en type. On dira que $%
M_{n}$\ converge en type vers Z si et seulement il existe des suites $\left(
a_{n}>0\right) _{n\geq 1}$ et $\left( b_{n}\right) _{n\geq 1}$\ telles que
la suite de variables al\'{e}atoires 
\begin{equation*}
\frac{M_{n}-b_{n}}{a_{n}}
\end{equation*}

\noindent converge vaguement vers $Z$, 
\begin{equation*}
\frac{M_{n}-b_{n}}{a_{n}}\rightsquigarrow Z.
\end{equation*}

\noindent D\'{e}couvrons les trois limites trois types de lois extremales d\'{e}finies
par leur fonctions de r\'{e}partition ($f.r$):\newline

\noindent \textbf{Type de Gumbel :} $\Lambda (x)=\exp (-\exp (-x)),x\in R.$
(On note $\Lambda $ une $v.a$ de $f.r$ $\Lambda )$\newline

\bigskip \noindent \textbf{Type de Frechet} de param\`{e}tre $\alpha >0:$ 
\begin{equation*}
\varphi _{\alpha }(x)=\exp (-x^{-\alpha })1_{(x\geq 0)}.
\end{equation*}

\noindent  On note $FR(\alpha )$ une $v.a$ de $f.r$ $\varphi _{\alpha })$%
\newline

\bigskip \noindent \textbf{Type de Weibull} de param\`{e}tre $\beta >0:$ 
\begin{equation*}
\psi _{\beta }(x)=\exp (-(-x)^{\beta })1_{(x<0)}+1_{(x\geq 0)}.
\end{equation*}

\noindent (On note $W(\beta )$ une $v.a$ de $f.r$ $\psi _{\beta })$.\newline

\noindent \textbf{Donnons quelques exemples de convergence en type de maxima}.\\

\noindent Nous allons utiliser la convergence des fonction de r\'{e}partition. Puisque les limites
vagues ont des fonction de r\'{e}partition continues, nous allons voir la
limite en tout point de $\mathbb{R}.$\\

\noindent \textbf{(a)} Loi de exponentielle : Supposons que 
\begin{equation*}
F(x)=(1-\exp (-x))1_{(x\geq 0)}
\end{equation*}

\noindent  est celle d'une loi exponentielle standard. En utilisant les
fonctions de r\'{e}partition, \noindent montrons que

\begin{equation*}
M_{n}-\log n\overset{d}{\rightarrow }\Lambda .
\end{equation*}

\noindent En effet
\begin{equation*}
P(M_{n}-\log n\leq x)=P(M_{n}\leq x+\log n)=F(M_{n}\leq x+\log n)^{n}.
\end{equation*}

\noindent Pour tout $x\in \mathbb{R},$ $x+\log n\geq 0$ pour $n\geq \exp (-x).$ Donc
pour $n$ assez grand $F(M_{n}\leq x+\log n)=(1-\exp (-x-\log n))$ et donc
pour tout $x\in \mathbb{R}$
\begin{equation*}
P(M_{n}-\log n\leq x)=(1-\frac{e^{-x}}{n})\rightarrow e^{-e^{-x}}=\Lambda(x).
\end{equation*}

\noindent  \textbf{(b)} Loi de par\'{e}to de param\`{e}tre $\alpha >0:$

\begin{equation*}
F(x)=(1-x^{-\alpha }))1_{(x\geq 1)}.
\end{equation*}

\noindent En utilisant les fonctions de r\'{e}partition, \noindent montrons que

\begin{equation*}
n^{-1/\alpha }M_{n}\overset{d}{\rightarrow }C(\alpha ).
\end{equation*}

\noindent Remarquons que les variables de Pareto sont positives et donc $M_{n}$ est
positif pour tout $n\geq 1.$ Etudions deux cas.\\

\noindent Cas $x\leq 0.$ Dans ce cas%
\begin{equation*}
P(n^{-1/\alpha }M_{n}\leq 0)=0=\varphi _{\alpha }(x),
\end{equation*}

\noindent et la limite des fonctins de r\'{e}partitions a lieu.\\

\noindent Cas $x>0$. Dans ce cas
\begin{equation*}
P(n^{-1/\alpha }M_{n}\leq x)=P(M_{n}\leq n^{1/\alpha }x).
\end{equation*}

\noindent Donc pour $n$ assez grand, on aura $n^{1/\alpha }x>1$ (par exemple prendre $%
n\geq (1/x)^{-\alpha })$ et pour ces $n,$%
\begin{eqnarray*}
P(n^{-1/\alpha }M_{n} &\leq &x)=F(n^{1/\alpha }x)^{n}=(1-(n^{1/\alpha
}x)^{-\alpha })^{n} \\
&=&(1-\frac{x^{-\alpha }}{n})^{n}\rightarrow \exp (-x^{-\alpha }) \\
&=&\varphi _{\alpha }(x).
\end{eqnarray*}

\noindent Nous avons bien que pour tout $x\in \mathbb{R}$,%
\begin{equation*}
P(n^{-1/\alpha }M_{n}\leq x)\longrightarrow \varphi _{\alpha }(x).
\end{equation*}

\noindent D\`{e}s lors
\begin{equation*}
n^{-1/\alpha }M_{n}\rightsquigarrow FR(\alpha ).
\end{equation*}

\noindent \textbf{(c)} Loi uniforme sur $(0,1)$ : 
\begin{equation*}
F(x)=x1_{(0\leq x\leq 1)}+1_{(x\geq 1)}.
\end{equation*}

\noindent En utilisant les fonctions de r\'{e}partition, montrons que

\begin{equation*}
n(M_{n}-1)\overset{d}{\rightarrow }W(1).
\end{equation*}

\noindent Nous avons
\begin{equation*}
P(n(M_{n}-1)\leq x)=F(1+\frac{x}{n})^{n}.
\end{equation*}

\noindent Etudions deux cas.\\

\noindent Cas $x\geq 0.$ Dans ce cas $1+x/n$ est positif pour tout $n\geq 1$ et 
\begin{equation*}
P(n(M_{n}-1)\leq x)=F(1+\frac{x}{n})^{n}=1=\psi _{1}(x)
\end{equation*}

\noindent et nous avons la convergence des fonctions de r\'{e}partition.\\

\noindent Cas $x<0.$ Pour $n$ suffisamment grand, nous aurons $0\leq 1+x/n\leq 1$
(prendre par exemple $x\geq -n$ $i.e.$ $n\geq -(x)\geq 0).$ Pour ces valeurs
de $n$,

\begin{eqnarray*}
P(n(M_{n}-1) &\leq &x)=F(1+\frac{x}{n})^{n} \\
&=&(1+\frac{x}{n})^{n}\rightarrow e^{x}=\psi _{1}(x).
\end{eqnarray*}

\noindent Nous avons bien que pour tout $x\in \mathbb{R}$,%
\begin{equation*}
P(n(M_{n}-1)\leq x)\longrightarrow \psi _{1}(x).
\end{equation*}

\noindent D\`{e}s lors%
\begin{equation*}
n(M_{n}-1)\rightsquigarrow W(1).
\end{equation*}

\bigskip \noindent \textbf{R\'{e}sum\'{e}} : En th\'{e}orie des valeurs extr\^{e}mes, il est d\'{e}montr%
\'{e} que les seules lois limites non d\'{e}g\'{e}n\'{e}r\'{e}es possibles
sont bien celles-la. Dans le chapitre r\'{e}serv\'{e} \`{a} une \'{e}tude sp%
\'{e}cifique de la convergence dans $\mathbb{R}$ utilisant les inverses g%
\'{e}n\'{e}ralis\'{e}es, les formes g\'{e}n\'{e}rales des distributions dont
les maxima convergent chacun des types seront donn\'{e}es.

\section{Exemples de convergence dans $\mathbb{R}^{k}$} \label{cv.review.sec4}

\subsection{Th\'{e}or\`{e}me Central Limite Standard sur $\mathbb{R}^{k}$} \label{cv.review.subsec.clt}

Passons au th\'{e}or\`{e}me central limite dans $\mathbb{R}^{k}.$ Soit $%
X_{1},X_{2},....$\ une suite de variables al\'{e}atoires centr\'ees, ind\'{e}pendantes
et identiquement distribu\'{e}es (i.i.d) centr\'{e}es de matrice de
covariance $\Sigma =(\sigma _{ij})_{1\leq i\leq k,1\leq j\leq k}$,  c'est 
\`{a} dire

\begin{equation*}
\sigma _{ij}=Cov(X_{i},X_{j}).
\end{equation*}

\bigskip \noindent Consid\'{e}rons les sommes partielles 
\begin{equation*}
S_{n}=X_{1}+X_{2}+...+X_{n}.
\end{equation*}

\bigskip \noindent Nous avons le th\'{e}or\`{e}me central limit sur $\mathbb{R}^{k}$, 
\begin{equation*}
S_{n}/\sqrt{n}\rightsquigarrow \mathcal{N}(0,\Sigma )
\end{equation*}

\bigskip \noindent \textbf{PREUVE}. Remarquons que la matrice $\Sigma $\ est sym\'{e}trique et semi-positive
puisque, pour tout $u\in \mathbb{R}^{k}$%
\begin{equation*}
^{t}u\Sigma u=\text{ }^{t}u\mathbb{E}(XX^{\prime })u=\mathbb{E}%
((^{t}Xu)(^{t}Xu))=\mathbb{E}((^{t}Xu)^{2})\geq 0.
\end{equation*}

\noindent  D'apr\`{e}s la th\'{e}orie des matrices (voir cours d'alg\`{e}bre de deuxi%
\`{e}me ann\'{e}e), $\Sigma $\ poss\`{e}de des valeurs propres non n\'{e}%
gatives et elle est diagonalisable au moyen d'une matrice orthogonale $T$,
telle que 
\begin{equation*}
^{t}T\Sigma T=diag(\lambda _{1},\lambda _{2},...,\lambda _{n})=\Lambda .
\end{equation*}%

\noindent Posons 
\begin{equation*}
Y_{i}=\text{ }^{t}TX_{i}.
\end{equation*}

\noindent Les variables $Y_{i}$\ sont centr\'{e}es, iid et de matrice de covariance 
\begin{equation*}
\Sigma _{Y}=\text{ }^{t}T\Sigma T=\Lambda .
\end{equation*}

\noindent Cela veut dire que les composantes de $Y_{i}$\ sont non corrol\'{e}es et ont
pour variances respectives les nombres $\lambda _{1},\lambda
_{2},...,\lambda _{n}.$\ Posons 
\begin{equation}
M_{n}=\frac{1}{\sqrt{n}}(Y_{1}+Y_{2}+...+Y_{n})=\text{ }^{t}T(\frac{S_{n}}{%
\sqrt{n}}).  \label{rk00}
\end{equation}

\bigskip \noindent Pour tout $A=$ $^{t}(a_{1},a_{2},...,a_{k})\in \mathbb{R}^{k},$%
\begin{equation*}
<A,M_{n}>=\frac{1}{\sqrt{n}}\sum_{i=1}^{i=n}<A,Y_{i}>.
\end{equation*}

\noindent Or les variables $<A,Y_{i}>$\ sont centr\'{e}e, i.i.d, de variance

\begin{equation*}
\mathbb{E}<A,Y_{i}>^{2}=\sum_{i=1}^{i=n}a_{i}^{2}\lambda _{i}=\text{ }%
^{t}A\Lambda A,
\end{equation*}

\noindent en vertu de la non corr\'{e}lation des composantes de chaque Y$_{i}.$\ Le th%
\'{e}or\`{e}me central limite dans $\mathbb{R}$ implique que 
\begin{equation*}
<A,M_{n}>\rightarrow \mathcal{N}(0,\sum_{i=1}^{i=n}a_{i}^{2}\lambda _{i})=%
\mathcal{N}(0,^{t}A\Lambda A)
\end{equation*}

\noindent Or $\mathcal{N}(0,^{t}A\Lambda A)$\ est la loi d'un vecteur gaussien r\'{e}%
sultant de la transformation lin\'{e}aire $^{t}AZ=<A,Z>,$\ o\`{u} $Z$ suit
la loi $\mathcal{N}(0,\Lambda ).$\ D'o\`{u} 
\begin{equation*}
\forall A\in \mathbb{R}^{k},<A,M_{n}>\rightsquigarrow <A,Z>.
\end{equation*}

\noindent En terme de fonction caract\'{e}ristique, cela veut dire qur pour $t\in 
\mathbb{R}$ et pour tout $A\in \mathbb{R}^{k},$

\begin{equation*}
\mathbb{E}\exp (it<A,M_{n}>)\rightarrow \mathbb{E}\exp (it<A,Z>).
\end{equation*}%

\noindent Pour $t=1,$ nous avons pour tout  $A\in \mathbb{R}^{k}$%
\begin{equation*}
\Phi _{M_{n}}(A)=\mathbb{E}\exp (it<A,M_{n}>)\rightarrow \Phi _{Z}(A)=\mathbb{E}\exp (it<A,Z>).
\end{equation*}%

\noindent Ceci veut bien dire que 
\begin{equation*}
M_{n}\rightsquigarrow Z.
\end{equation*}

\noindent Ceci,la formule (\ref{rk00}) et le point (e) du th\'eor\`eme \ref{cv.review.rk} ensemble impliquent 

\begin{equation*}
S_{n}/\sqrt{n}=TM_{n}\rightarrow TZ
\end{equation*}%
et 
\begin{equation*}
^{t}TZ\sim \mathcal{N}(0,T\Lambda ^{t}T)=\mathcal{N}(0,\Sigma ).
\end{equation*}%

\noindent D'o\`{u}, enfin, 
\begin{equation*}
S_{n}/\sqrt{n}\rightsquigarrow \mathcal{N}(0,\Sigma )
\end{equation*}

\noindent Ceci exprime la version simple du th\`{e}or\`{e}me central limite dans $%
\mathbb{R}^{k}.$

\subsection{Convergence de la loi multinomiale} \label{cv.review.subsec.multinomial}

\bigskip Un k-uplet $X_n=(X_{1,n},...,X_{k,n})$\ suit une loi multimoniale de param%
\`{e}tres $n\geq 1$\ et $p=(p_{1},p_{2},...p_{k})$\ avec 
\begin{equation*}
\forall (1\leq i\leq k),p_{i}>0\text{ \ }et\text{ \ }\sum_{1\leq i\leq
i1}p_{i}=1,
\end{equation*}

not\'{e}e $\mathcal{M}_{k}(n,p)$, ssi sa loi de probabilit\'{e} est: 

\begin{equation*}
\mathbb{P}(X_{1}=n_{1},...,X_{k}=n_{k})=\frac{n!}{n_{1}!\times ...\times
n_{k}!}p_{1}^{n_{1}}\times p_{2}^{n_{2}}\times ...\times p_{k}^{n_{k}}
\end{equation*}

pour $(n_{1},...,n_{k})$\ v\'{e}rifiant 
\begin{equation*}
\forall (1\leq i\leq k),\text{ }n_{i}\geq 0\text{ \ et \ }\sum_{1\leq i\leq
k}n_{i}=n.
\end{equation*}

\noindent Elle est g\'{e}n\'{e}r\'{e}e de la mani\`{e}re suivante. Soit une exp\'{e}%
rience \`{a} $k$ issues $E_{i},1$\ $\leq i\leq k,$\ chacune se r\'{e}alisant
avec une probabilit\'{e} $p_{i}>0.$\ On la r\'{e}p\`{e}te $n$ fois de mani%
\`{e}re ind\'{e}pendante. A l'issue de ces $n$ essais, soit $X_{i,n}$\ le
nombre de r\'{e}alisations de l'issue $E_{i}.$\ Le vecteur ainsi obtenu suit
une loi $\mathcal{M}_{k}(n,p)$. Bien s\^{u}r chaque $X_{i,n}$\ suit une loi
binomiale $\mathcal{B}(n,p_{i}).$\\

\noindent Nous avons le r\'esultat de convergence vague suivant.\\

\begin{equation}
Z_{n}=^{t}(\frac{X_{1,n}-np_{1}}{\sqrt{np_{1}}},...,\frac{X_{k,n}-np_{k}}{\sqrt{np_{k}}}) \rightsquigarrow \mathcal{N}_{k}(0,\Sigma) \ \ as \ \ n\rightarrow +\infty, \label{cv.tclMultiNom01}
\end{equation}

\noindent o\`u $\Sigma$ est matrice carr\'ee d'ordre $k$ dont les \'el\'ements sont $\Sigma_{i,i}=1-p_{i}$ et $\Sigma_{i,j}=\sqrt{p_{i}p_{j}}$, $1\leq i,j\leq k$.\\ 

\noindent \textbf{Remarque importante}. Ce r\'esultat a d'importantes applications. Nous pouvons mentionner son utilisation
pour trouver la loi des distributions finies du processus empirique en probabilit\'es. Il est aussi à la base des tests du khi-deux en Statistiques. Nous verrons ces tests plus tard dans la partie r\'eserv\'ee à cet effet dans cette s\'erie.\\

\bigskip \noindent \textbf{Preuve}. Nous allons pr\'esenter deux preuves. La premi\`ere utilise la convergence de la fonction des moments et du development de la fonction logarithmique. Elle est plus adapt\'ee à un enseignement de premier cycle. La deuxi\`eme exploite le th\'eor\`eme central limit dans $\mathbb{R}^k$, que nous avons vu pr\'ec\'edemment. Elle est plus adapt\'ee à expos\'e de niveau sup\'erieur.\\ 

\noindent \textbf{Premi\`ere preuve}.\\

\noindent Nous connaissons d\'{e}j\`{a} sa fonction g\'{e}n\'{e}ratrice des moments
qui est   
\begin{equation*}
\phi _{X_n}(u)=(\sum_{1\leq i\leq k}p_{i}e^{u_{i}})^{n}.
\end{equation*}

\noindent Posons 
\begin{eqnarray*}
Z_{n}&=&(\frac{X_{1,n}-np_{1}}{\sqrt{np_{1}}},...,\frac{X_{k,n}-np_{k}}{\sqrt{%
np_{k}}})\\
&=&AX+B
\end{eqnarray*}

\noindent avec 
\begin{equation*}
A=\left( 
\begin{array}{cccc}
1/\sqrt{np_{1}} &  &  &  \\ 
& \sqrt{np_{2}} &  &  \\ 
&  & ... &  \\ 
&  &  & \sqrt{np_{k}}%
\end{array}%
\right)
\end{equation*}%

\noindent et 
\begin{equation*}
B=\left( 
\begin{array}{c}
-\sqrt{np_{1}} \\ 
-\sqrt{np_{2}} \\ 
... \\ 
-\sqrt{np_{k}}%
\end{array}%
\right) .
\end{equation*}

\noindent D\`{e}s lors 
\begin{equation*}
\phi _{Z_{n}}(u)=\exp (<B,u>)\times \phi _{X}(^{t}Au)
\end{equation*}%
\begin{equation*}
=\exp (\sum_{1\leq i\leq k}-\sqrt{np_{i}}u_{i})\times (\sum_{1\leq i\leq
k}p_{i}e^{u_{i}/\sqrt{np_{i}}})^{n}
\end{equation*}

\noindent Notons que $u$ est fix\'{e}. Pour tout $i$ fix\'{e}, \ $u_{i}/\sqrt{np_{i}}%
\rightarrow \rightarrow 0$\ quand $n\rightarrow \infty $\ car chaque $%
p_{i}>0.$\ D'o\`{u} 
\begin{equation*}
e^{u_{i}/\sqrt{np_{i}}}=1+u_{i}/\sqrt{np_{i}}+\frac{1}{2}\frac{u_{i}^{2}}{%
np_{i}}+O(n^{-3/2}).
\end{equation*}

\noindent D'o\`{u}\begin{eqnarray*}
A&=&(\sum_{1\leq i\leq k}p_{i}e^{u_{i}/\sqrt{np_{i}}})^{n}=\exp (n\log
(\sum_{1\leq i\leq k}p_{i}e^{u_{i}/\sqrt{np_{i}}})).\\
&=&\exp (n\log (1+\sum_{1\leq i\leq k}u_{i}\sqrt{p_{i/}n_{i}}+\sum_{1\leq
i\leq k}\frac{1}{2}\frac{u_{i}^{2}}{n_{i}}+O(n^{-3/2}))).
\end{eqnarray*}

\noindent Posons aussi 
\begin{equation*}
a=\sum_{1\leq i\leq k}u_{i}\sqrt{p_{i/}n}+\sum_{1\leq i\leq k}\frac{1}{2}%
\frac{u_{i}^{2}}{n}\rightarrow 0\text{ }quand\text{ }n\rightarrow \infty .
\end{equation*}

\noindent Nous aurons 
\begin{equation*}
A=\exp (n\log (1+a)).
\end{equation*}

\noindent D\'{e}veloppons $log(1+a)$ \`{a} l'ordre $2$ en mettant dans $O(n^{-3/2})$\
tous les autres termes \ tendant vers z\'{e}ro:

\begin{eqnarray*}
A&=&\exp (n(a-\frac{1}{2}a^{2}+O(a^{3})).\\
&=&\exp (n(\sum_{1\leq i\leq k}u_{i}\sqrt{p_{i/}n}+\sum_{1\leq i\leq k}\frac{1%
}{2}\frac{u_{i}^{2}}{n}-\frac{1}{2}(\sum_{1\leq i\leq k}u_{i}\sqrt{p_{i/}n}%
)^{2}+O(n^{-3/2})))\\
&=&\exp (\sum_{1\leq i\leq k}u_{i}\sqrt{np_{i}}+\sum_{1\leq i\leq k}\frac{1}{2}%
u_{i}^{2}-\frac{1}{2}(\sum_{1\leq i\leq k}u_{i}\sqrt{p_{i}})^{2}+O(n^{-1/2}))\\
&=&\exp (\sum_{1\leq i\leq k}u_{i}\sqrt{np_{i}})\times \exp (\sum_{1\leq i\leq
k}\frac{1}{2}u_{i}^{2}-\frac{1}{2}(\sum_{1\leq i\leq k}u_{i}\sqrt{p_{i}}%
)^{2}+O(n^{-1/2})).
\end{eqnarray*}

\noindent En mettant tout cela ensemble, nous obtenons
\begin{equation*}
\phi _{Z_{n}}(u)=\exp (\sum_{1\leq i\leq k}\frac{1}{2}u_{i}^{2}-\frac{1}{2}%
(\sum_{1\leq i\leq k}u_{i}\sqrt{p_{i}})^{2}+O(n^{-1/2}))
\end{equation*}%
\begin{equation*}
\rightarrow \phi _{Z}(u)=\exp (\sum_{1\leq i\leq k}\frac{1}{2}u_{i}^{2}-%
\frac{1}{2}(\sum_{1\leq i\leq k}u_{i}\sqrt{p_{i}})^{2}).
\end{equation*}

\noindent Et 
\begin{equation}
\phi _{Z}(u)=\exp (\left\{ \sum_{1\leq i\leq k}\frac{1}{2}%
(1-p_{i})u_{i}^{2}-\sum_{1\leq i,j\leq k}u_{i}u_{j}\sqrt{p_{i}p_{j}}\right\} 
\label{cv7001}
\end{equation}

\noindent est la fonction des moments d'un vecteur gaussien $Z$ centr\'{e} dont la
matrice de variances-covariances $\Sigma $\ v\'{e}rifie 
\begin{equation}
\Sigma _{ii}=(1-p_{i})  \label{cv7002}
\end{equation}

\noindent et 
\begin{equation}
\Sigma _{ij}=-\sqrt{p_{i}p_{j}}.  \label{cv7003}
\end{equation}

\noindent Donc 
\begin{equation*}
Z_{n}\rightsquigarrow \mathcal{N}(O,\Sigma ).
\end{equation*}

\noindent La premi\`ere preuve finit ici.\\

\bigskip \noindent \textbf{Deuxi\`eme preuve}. Au $i$-i\`eme, $i\in \{1,...,n\}$, nous obtenons le vecteur variable al\'eatoire  
\begin{equation*}
Z^{(i)}=\left( 
\begin{tabular}{l}
$Z_{1}^{(i)}$ \\ 
... \\ 
$Z_{k}^{(i)}$%
\end{tabular}%
\right) 
\end{equation*}

\noindent d\'efini ainsi : pour chaque $1\leq r \leq k$, nous avons 

\begin{equation*}
Z_{r}^{(i)}=\left\{ 
\begin{tabular}{lll}
$1$ & si  & l'issue $E_{r}$ se realise au  $i$-i\`eme essai et auncune autre ne se r\'ealise \\ 
$0$ & si  & un autre issue que $E_r$ se r\'ealise\  
\end{tabular}%
\right. 
\end{equation*}

\noindent Il est \'evident que $Z^{(i)}$ a un distribution multinomiale $\mathcal{M}_{k}(1,k)$ et que les vecteurs al\'eatoires $Z^{(i)}$ sont independents.\\

\noindent De plus, pour tout $i\in \{1,...,n\}$, chaque  $Z_{r}^{(i)}$, $1\leq r \leq k$, suit une loi de Bernouilli de param\`etre 
$p$ et un seul des $Z_{r}^{(i)}$ ($1\leq r\leq k$) prend la valeur un $(1)$, les autres \'etant nuls. Cela implique que  
\begin{equation*}
Z_{r}^{(i)}Z_{s}^{(i)}=0\text{ for }1\leq r\neq s\leq k,1\leq i\leq n.
\end{equation*}

\noindent Nous avons aussi 
\begin{equation*}
Z_{1}^{(i)}+...+Z_{n}^{(i)}=1.
\end{equation*}

\noindent Alors, pour tout $i\in \{1,...,n\},$

\begin{equation*}
\mathbb{E}(Z_{r}^{(i)})=p_{i}\text{ et }\mathbb{V}ar(Z_{r}^{(i)})=p_{i}(1-p_{i}),\text{ }%
1\leq r\leq k
\end{equation*}

\noindent et pour tout $1\leq r\neq s\leq k$
\begin{equation*}
cov(Z_{r}^{(i)},Z_{s}^{(i)})=\mathbb{E}(Z_{r}^{(i)}Z_{s}^{(i)})-\mathbb{E}(Z_{r}^{(i)})\mathbb{E}(Z_{s}^{(i)})=-p_{r}p_{s},
\end{equation*}

\noindent puisque $Z_{r}^{(i)}Z_{s}^{(i)}=0.$ D\`es lors, chaque $Z^{(i)}$ poss\`ede la matrice de
variance-covariance 

\begin{equation*}
\Sigma ^{0}=\left( 
\begin{tabular}{lllll}
$p_{1}(1-p_{1})$ & $-p_{1}p_{2}$ & ... & $-p_{1}p_{k-1}$ & $-p_{1}p_{k}$ \\ 
$-p_{2}p_{1}$ & $p_{2}(1-p_{2})$ & ... & $-p_{2}p_{k-1}$ & $-p_{2}p_{k1}$ \\ 
... & ... & ... & ... & ... \\ 
$-p_{k-1}p_{1}$ & $-p_{k-1}p_{2}$ & ... & $p_{k-1}(1-p_{k-1})$ & $%
-p_{k-1}p_{k}$ \\ 
$-p_{k}p_{1}$ & $-p_{k}p_{2}$ & ... & $-p_{k}p_{k-1}$ & $-p_{k}(1-p_{k})$%
\end{tabular}%
\right) 
\end{equation*}

\noindent ou, par une autre notation,

\begin{equation*}
\Sigma _{0}=(\sigma _{ij}^{0})_{1\leq i,j\leq k}\text{ with }\sigma
_{ij}^{0}=\left\{ 
\begin{tabular}{lll}
$p_{i}(1-p_{i})$ & if & $i=j$ \\ 
$-p_{i}p_{j}$ & if & $i\neq j$%
\end{tabular}
\right. .
\end{equation*}

\noindent Apr\`es $n$ essais, la somme des vecteurs al\'eatoires $Z^{(1)},...,Z^{(n)}$, qui sont id\'ependants et suivent identiquement une loi
$\mathcal{M}_{k}(1,k)$, donne $X_{n}$, autrement 
\begin{equation*}
X_{n}=Z^{(1)}+...+Z^{(n)}.
\end{equation*}

\noindent Par le th\'eor\`eme central limit standard multivari\'e, nous avons, quand $n\rightarrow +\infty$,

\begin{equation*}
S_{n}=\frac{1}{\sqrt{n}}\sum_{i=1}^{n}\left( Z^{(i)}-\mathbb{E}(Z^{(i)}\right)
\rightsquigarrow Z_{0}\sim \mathcal{N}_{k}(0,\Sigma _{0}).
\end{equation*}

\noindent Mais nous pouvons ais\'ement v\'erifier que

\begin{equation*}
S_{n}=\frac{1}{\sqrt{n}}\sum_{i=1}^{n}\left( Z^{(i)}-\mathbb{E}(Z^{(i)}\right)
=^{t}\left( \frac{X_{1,n}-np_{1}}{\sqrt{n}},\frac{X_{2,n}-np_{2}}{\sqrt{n}}%
,...,\frac{X_{k,n}-np_{k}}{\sqrt{n}}\right) .
\end{equation*}

\noindent Et nous obtenons la relation matricielle 

\begin{equation*}
DS_{n}=Z_{n},
\end{equation*}

\noindent o\`u $D$ est la matrice diagonale

\begin{equation*}
D=diag(1/\sqrt{p_{1}},...,1/\sqrt{p_{k}}).
\end{equation*}

\noindent Par the th\'eor\`eme de transformation continue (Point (e) du th\'eor\`eme \ref{cv.review.rk}), nous avons

\begin{equation*}
Z_{n}=DS_{n}\rightsquigarrow DZ_{0}\sim \mathcal{N}_{k}(0,D\Sigma _{0}D),
\end{equation*}

\noindent puisque $D$ est sym\'etrique. Il reste à calculer 

\begin{equation*}
\Sigma =D\Sigma _{0}D=(\sigma _{ij})_{1\leq i,j\leq k}.
\end{equation*}

\noindent For $1\leq h,j\leq k,$ $(\Sigma _{0}D)_{hj}$ est la produit matriciel de la $h$-i\`eme ligne de $\Sigma _{0}$ par la
$j$-i\`eme colonne de $D$. En exploitant la diagonalit\'e de $D$, nous obtenons pour $1\leq h,j\leq k,$

\begin{equation*}
(\Sigma _{0}D)_{hj}=\sigma _{hj}^{0}/\sqrt{p_{j}}.
\end{equation*}

\noindent Ensuite,  $\sigma _{ij}=(D\Sigma _{0}D)_{ij}$ est le produit de la  $i$-i\`eme
ligne de $D$ par la $j$-i\`eme colonne de $(\Sigma _{0}D)^{(j)}=$ $^{t}((\Sigma
_{0}D)_{1j},(\Sigma _{0}D)_{2j},...,(\Sigma _{0}D)_{kj})$. En utilisant encore le fait que $D$ est diagonale, nous arrivons \`a
 
\begin{equation*}
(D\Sigma _{0}D)_{ij}=\frac{1}{\sqrt{p_{i}}}(\Sigma _{0}D)_{ij},
\end{equation*}

\noindent ce qui aboutit \`a
\begin{equation*}
\sigma _{ij}=(D\Sigma _{0}D)_{ij}=\frac{1}{\sqrt{p_{i}p_{j}}}\sigma
_{ij}^{0}=\left\{ 
\begin{tabular}{lll}
$\sigma _{ii}^{0}/p_{i}=1-p_{i}$ & if & $i=j$ \\ 
-$\sqrt{p_{i}p_{j}}$ & if & $i\neq j$%
\end{tabular}%
\right. ..
\end{equation*}

\noindent Nous obtenons encore
\begin{equation*}
Z_{n}\rightsquigarrow \mathcal{N}_{k}(0,\Sigma ), \ \ as  \ \ n\rightarrow +\infty.
\end{equation*}

\noindent o\`u $\Sigma$ is d\'efini dans la ligne qui suit la formule (\ref{cv.tclMultiNom01}) dans le partie haut de cette sous-section. Ceci finit la deuxi\`eme preuve.\\

\noindent 

En conclusion, nous avons le r\'{e}sultat suivant.

\begin{proposition} \label{cv.multinomial}
Soit une suite de vecteurs al\'{e}atoires $X(n)=(X_{1}(n),...,X_{k}(n))$\
suivant une loi multimoniale de param\`{e}tres $n\geq 1$\ et $%
p=(p_{1},p_{2},...p_{k})$\ avec 
\begin{equation*}
\forall (1\leq i\leq k),p_{i}>0\text{ \ et \ }\sum_{1\leq i\leq i1}p_{i}=1.
\end{equation*}

\noindent Alors la suite de vecteurs 
\begin{equation*}
Z_{n}=(\frac{X_{1}-np_{1}}{\sqrt{np_{1}}},...,\frac{X_{1}-np_{k}}{\sqrt{%
np_{k}}})
\end{equation*}

\noindent converge vers une loi normale k-dimensionn\'{e}e, centr\'{e}e, de matrice de variances-covariances $\Sigma $\ avec 
\begin{equation*}
\Sigma _{ii}=(1-p_{i})
\end{equation*}

\noindent et 
\begin{equation*}
\Sigma _{ij}=-\sqrt{p_{i}p_{j}}.
\end{equation*}
\end{proposition}

\subsection{Limites des dimensions finies du processus empirique uniforme} \label{cv.review.subsec.ep}

\bigskip \noindent Consid\'{e}rons $U_{1},$\ $U_{2},...$\ une suite de variables al\'{e}atoires
ind\'{e}pendantes et identiquement distribu\'{e}es selon la loi uniforme sur
(0,1), de fonction de r\'{e}partition commune $F(s)=s1_{(0\leq s\leq
1)}+1_{(s\geq 1)}$. Pour $n\geq 1,$\ on d\'{e}finit la  fonction r\'{e}%
partition empirique uniforme bas\'{e}e sur l'\'{e}chantillon $U_{1},$\ $%
U_{2},...,U_{n}$\ la fonction 
\begin{equation*}
\begin{array}{ccc}
\mathbb{R}\ni x & \mapsto  & U_{n}(s)=\frac{1}{n}Card\{i,\text{ }1\leq i\leq
n\text{, \ }U_{i}\leq s\}%
\end{array}%
\end{equation*}

\noindent Et d\'{e}finissons le processus empirique uniforme%
\begin{equation*}
\alpha _{n}(s)=\sqrt{n}(U_{n}(s)-s),0\leq s\leq 1.
\end{equation*}

\bigskip \noindent Consid\'{e}rons $0=t_{0}<t_{1}<...<t_{k}<t_{k+1}=1$\ et posons  
\begin{equation*}
Y_{n}=(\alpha _{n}(t_{1}),...,\alpha _{n}(t_{k+1})).
\end{equation*}

\noindent Alors nous avons le r\'{e}sultat :

\begin{proposition} \label{cv.ep.df}
Les distributions finies du processus empirique uniformes de la forme \ $%
(\alpha _{n}(t_{1}),...,\alpha _{n}(t_{k+1})),$ avec $%
0=t_{0}<t_{1}<...<t_{k}<t_{k+1}=1,$ v\'{e}rifient%
\begin{equation*}
(\alpha _{n}(t_{1}),...,\alpha _{n}(t_{k+1}))\rightarrow \mathcal{N}%
_{k}(0,\left( \min (t_{i},t_{j})-t_{i}t_{j}\right) _{1\leq i,j\leq k})
\end{equation*}
\end{proposition}

\bigskip \noindent \textbf{PREUVE}. Posons 

\begin{equation}
Z_{n}=(\frac{\alpha _{n}(t_{1})}{\sqrt{t_{1}}},\frac{\alpha
_{n}(t_{2})-\alpha _{n}(t_{1})}{\sqrt{t_{2}-t_{1}}}...,\frac{\alpha
_{n}(t_{k+1})-\alpha _{n}(t_{k})}{\sqrt{t_{k+1}-t_{k}}}).  \label{cv730}
\end{equation}

\noindent Remarquons que 
\begin{equation*}
N_{n}=(nF_{n}(t_{1}),nF_{n}(t_{2})-nF_{n}(t_{1}),...,nF_{n}(t_{k+1})-nF_{n}(t_{k}))
\end{equation*}

\noindent suit une loi multinomiale dont les composantes du vecteur des probabilit\'es des issues sont : $t_{1}$, $t_{2}-t_{1}$,...,$t_{k+1}-t_{k}$. En effet,  
\begin{equation*}
nF_{n}(t_{j})-nF_{n}(t_{j-1})
\end{equation*}

\noindent est le nombre d'observations tombant dans $]t_{j-1}, t_{j}]$ et pour tout $j$, la probabilit\'e qu'une observation tombe
 dans $]t_{j-1}, t_{j}]$ est bien $p_{j}=t_{j}-t_{j-1}$.\\

\noindent Nous pouvons appliquer le th\'eor\`eme de convergence faible de la loi multinomiale \'etablie dans la sous-section 
\ref{cv.review.subsec.multinomial}.\\

\noindent D\'efinisson $Z_{n}$ en centrant chaque  $j$-i\`eme composante de $N_{n}$ par $n(t_{j}-t_{j-1})$ en en la normalisant par
$\sqrt{n(t_{j}-t_{j-1})}$.\\

\noindent Rappelons nous que nous avons : $Y_{n}^{T}=(\alpha _{n}(t_{1}),...,\alpha _{n}(t_{k+1}))$. Nous obtenons la relation matricielle

$$
Z_n=AY_n \Leftrightarrow Y_n=BZ_n
$$

\noindent o\`u la relation $y=Bz$ traduit la correspondance suivante :  

$$
y_i=\sqrt{t_1}x_{1} + \sqrt{t_{2}-t_{1}}x_{2} + ... + \sqrt{t_{2}-t_{1}}x_{i}, \ \ 1\leq i \leq k+1.
$$

\noindent D'apr\`es la convergence de la loi multinomiale sus-mentionn\'ee, $Z_n$ converges faiblement vers un vecteur gaussian $Z=(Z_{1},Z_{2},...,Z_{k+1})$ v\'erifiant 

\begin{equation*}
\mathbb{E}(Z_{j}^{2})=1-(t_{j}-t_{j-1})
\end{equation*}

\noindent et
\begin{equation*}
\mathbb{E}(Z_{i}Z_{j})=-\sqrt{(t_{i}-t_{i-1})(t_{j}-t_{j-1})}.
\end{equation*}

\noindent Par le th\'eor\`eme de la transformation continue, (Point (e) du th\'eor\`eme \ref{cv.review.rk}), $Y_n=BZ_n$ concerge vaguement vers $Y=BZ$, avec 

$$
Y_i=\sqrt{t_1}Z_{1} + \sqrt{t_{2}-t_{1}}Z_{2} + ... + \sqrt{t_{2}-t_{1}}Z_{i}, \ \ 1\leq i \leq k+1.
$$

\bigskip \noindent Soit le vecteur $T=(T_{1},...,T_{k+1})$ d\'efini par $(t_{j}-t_{j-1})Z_{j}=T_{j}$, $1\leq i \leq k$, i.e., 
\begin{equation*}
Z=(\frac{T_{1}}{\sqrt{t_{1}}},\frac{T_{2}}{\sqrt{(t_{2}-t_{1})}},...,\frac{T_{j}}{\sqrt{(t_{j}-t_{j-1})}},...,\frac{T_{k+1}}{\sqrt{(t_{k+1}-t_{k})}})
\end{equation*}

\bigskip \noindent Nous avons
\begin{equation*}
\mathbb{E}(T_{j}^{2})=\mathbb{E}((Z_{j}\sqrt{(t_{j}-t_{j-1})}%
)^{2}=(t_{j}-t_{j-1})(1-(t_{j}-t_{j-1})).
\end{equation*}

\bigskip \noindent et
\begin{equation*}
\mathbb{E}(T_{i}T_{j})=\sqrt{(t_{j}-t_{j-1})(t_{i}-t_{i-1})}\mathbb{E}%
(Z_{i}Z_{j})=-(t_{j}-t_{j-1})(t_{i}-t_{i-1}).
\end{equation*}

\noindent Avant de calculer la covariance de $Y_{i}$ et de $Y_{j}$, v\'erifions que pour $t_i \leq t_j$, nous avons
\begin{eqnarray*}
t_i t_j&=&(\sum_{h=1}^{h=i} (t_h-t_{h-1}))(\sum_{r=1}^{r=j} (t_r-t_{r-1}))\\
&=& (\sum_{h=1}^{h=i} (t_h-t_{h-1}))(\sum_{r=1}^{r=i} (t_r-t_{r-1})+\sum_{r=i+1}^{r=j} (t_r-t_{r-1}))\\
&=&(\sum_{h=1}^{h=i} (t_h-t_{h-1}))^2+ \sum_{h=1}^{h=i} \sum_{r=i+1}^{r=j} (t_h-t_{h-1})(t_r-t_{r-1})\\
&=&\sum_{h=1}^{h=i} (t_h-t_{h-1})^2+\sum_{1\leq h\neq r \leq i} (t_h-t_{h-1})(t_r-t_{r-1})\\
&-& \sum_{h=1}^{h=i} \sum_{r=i+1}^{r=j} (t_h-t_{h-1})(t_r-t_{r-1})
\end{eqnarray*}

\noindent En mettant ensemble les points pr\'ec\'edents, nous sommes en mesure de calculer la matrice de variance-covariance de $Y$. Pour $1\leq i\leq j\leq 1$, nous avons
\begin{eqnarray*}
Y_{i}Y_{j}&=&(\sum_{h=1}^{h=i}T_{h})^{2}+\sum_{h=1}^{h=i}\sum_{r=i+1}^{r=j}T_{h}T_{k}\\
&=& \sum_{h=1}^{h=i}T_{h}^2+\sum_{1\leq h\neq r \leq i}T_{h}T_{r} +\sum_{h=1}^{h=i}\sum_{r=i+1}^{r=j}T_{h}T_{r}.
\end{eqnarray*}

\noindent Finalement, nous obtenons
\begin{eqnarray*}
\mathbb{E}(Y_{i}Y_{j})&=& \sum_{h=1}^{h=i}(1-(t_h-t_{h-1})) -\sum_{1\leq h\neq r \leq i} (t_h-t_{h-1})(t_r-t_{r-1})\\
&-&\sum_{h=1}^{h=i}\sum_{r=i+1}^{r=j} (t_h-t_{h-1})(t_r-t_{r-1})\\
&=& \sum_{h=1}^{h=i} (t_h-t_{h-1})-\sum_{h=1}^{h=i}(t_h-t_{h-1})^2\\
&-& \sum_{1\leq h\neq r \leq i} (t_h-t_{h-1})(t_r-t_{r-1}) -\sum_{h=1}^{h=i}\sum_{r=i+1}^{r=j} (t_h-t_{h-1})(t_r-t_{r-1})\\
&=& t_i -\sum_{h=1}^{h=i}(t_h-t_{h-1})^2 -\sum_{1\leq h\neq r \leq i} (t_h-t_{h-1})(t_r-t_{r-1})\\
&-& \sum_{h=1}^{h=i}\sum_{r=i+1}^{r=j} (t_h-t_{h-1})(t_r-t_{r-1})\\
&=& t_i - t_i t_j=\min(t_i,t_j)-t_i t_j.
\end{eqnarray*}

\noindent Ceci ach\`eve la preuve.\\

\section{Principe d'invariance} \label{cv.review.sec5}

Soit $X_{1},X,.....$une suite iid de variables al\'{e}atoires centr\'{e}es
avec $E\left\vert X_{i}\right\vert ^{2}<\infty $\ et soit pour$n\geq 1,$

\begin{equation*}
S_{n}=X_{1}+......+X_{n},n\geq 1.
\end{equation*}%

\noindent Posons pour $0\leq t\leq 1,$\  
\begin{equation*}
S_{n}(t)=\frac{S_{\left[ nt\right] }}{\sqrt{n}}
\end{equation*}

\noindent Nous allons \'{e}tudier la loi limite du processus \{$S_{n}(t),$\ $0\leq
t\leq 1\}$ en distribution finie.\\

\noindent Pour cela, soit $0=t_{0}<t_{2}<...<t_{k}=1,$\ $\ k\geq 1.$ Nous avons

\begin{proposition} \label{cv.review.ip}
The sequence of  finite distributions
\begin{equation*}
\left( \frac{S_{\left[ nt_{j}\right] }}{\sqrt{n}},1\leq j\leq k\right), \ \ n\geq 1,
\end{equation*}

\noindent weakly converges to k-dimensional centered Gaussian vector of
variance-covariance matric%
\begin{equation*}
\left( \min (t_{i},t_{j})\right) _{1\leq i,j\leq k}.
\end{equation*}
\end{proposition}

\bigskip \noindent \textbf{Proof}.   Nous avons
\begin{equation*}
\left\{ 
\begin{array}{c}
Y_{n}(t_{1})=X_{n}(t_{1})-X_{n}(t_{0})=\frac{1}{\sqrt{n}}\sum_{[nt_{0}]<j%
\leq \lbrack nt_{1}\}}X_{j} \\ 
.... \\ 
Y_{n}(t_{i})=X_{n}(t_{i})-X_{n}(t_{i-1})=\frac{1}{\sqrt{n}}%
\sum_{[nt_{i-1}]<j\leq \lbrack nt_{i}\}}X_{j} \\ 
.... \\ 
Y_{n}(t_{k})=X_{n}(t_{k})-X_{n}(t_{k-1})=\frac{1}{\sqrt{n}}%
\sum_{[nt_{k-1}]<j\leq \lbrack nt_{k}\}}X_{j}%
\end{array}%
\right\} .
\end{equation*}

\noindent Nous constatons que les variables $Y_{n}(t_{i})$\ sont ind\'{e}pendantes et
pour chaque $t_{i}$, nous pouvons appliquer le th\'{e}or\`{e}me central limite dans $\mathbb{R}$,  
\begin{equation*}
Y_{n}(t_{i})=\frac{1}{\sqrt{n}}\sum_{[nt_{i-1}]<j\leq \lbrack
nt_{i}\}}X_{j}\rightarrow N(0,t_{i}-t_{i-1})
\end{equation*}

\noindent D\`{e}s lors, pour tout $u=(u_{1},...,u_{k})\in R^{k}$\  
\begin{equation*}
\mathbb{E}(\exp (\sum_{1\leq i\leq 1}Y_{n}(t_{i})u_{i})=\prod_{1\leq i\leq 1}%
\mathbb{E}(\exp (Y_{n}(t_{i})u_{i})\rightarrow \prod_{1\leq i\leq 1}e^{\frac{%
1}{2}u_{i}^{2}/(t_{i}-t_{i-1})}.
\end{equation*}

\noindent Donc le vecteur $Y_{n}=(Y_{n}(t_{i}),$\ $1\leq i\leq k)$\ tend vers un
vecteur gaussien $Z$ \`{a} composantes ind\'{e}pendantes et dont le i$^{i%
\grave{e}me}$\ composante a la variance $t_{i}-t_{i-1}$.\\

\noindent  Mais le vecteur $X_{n}=(X_{n}(t_{i}),$\ $1\leq i\leq k)$\ est une transformation lin\'{e}aire
de $Y_{n}$\ de la forme  
\begin{equation*}
X_{n}=AY_{n}=\left( 
\begin{array}{cccc}
1 & 0 & ... & 0 \\ 
1 & 1 & ... & 0 \\ 
1 & ... & 1 & 0 \\ 
1 & ... & ... & 1%
\end{array}%
\right) Y_{n}
\end{equation*}

\noindent avec  
\begin{equation*}
A_{ij}=1_{(i\leq j)}.
\end{equation*}

\noindent Alors $X_{n}$ converges en loi vers le vecteur gaussien $V=AZ$, de sorte que  
\begin{equation*}
V_{i}=Z_{1}+...+Z_{i}
\end{equation*}

\noindent et  
\begin{equation*}
Z_{i}=V_{i}-V_{i-1}.
\end{equation*}

\noindent Donc pour tout $1\leq i\leq k,$%
\begin{equation*}
\mathbb{E}(V_{i}^{2})=\sum_{1\leq j\leq i}\mathbb{E}(Z_{j}^{2})=\sum_{1\leq
j\leq i}(t_{j}-t_{j-1})=t_{j}.
\end{equation*}

\noindent Et pour tout 1$\leq i\leq j\leq k,$%
\begin{equation*}
\mathbb{E}(V_{i}V_{j})=\mathbb{E}(V_{i}(V_{i}+(V_{j}-V_{i})
\end{equation*}

\noindent et

\begin{equation*}
=\mathbb{E}(V_{i}^{2})+\mathbb{E}(V_{i}(V_{j}-V_{i})).
\end{equation*}

\noindent Puisque  
\begin{equation*}
V_{i}=Z_{1}+...+Z_{i}
\end{equation*}

\noindent et
  
\begin{equation*}
Z_{i}=V_{i}-V_{i-1}
\end{equation*}%

\noindent sont ind\'{e}pendants et centr\'{e}s, il vient que   
\begin{equation*}
\mathbb{E}(V_{i}V_{j})=\mathbb{E}(V_{i})=t_{i}=t_{i}\wedge t_{j}.
\end{equation*}

\noindent Cela suffit pour d\'emontrer le r\'esultat.\\

 

%% file: asymptotics_cv_01_fr.tex
\chapter{Th\'eorie de la Convergence Vague} \label{chapTmCv} \label{cv}

\section{Introduction}

Dans ce chapitre, nous traitons de la th\'{e}orie unifi\'{e}e de la
convergence vague par sa caract\'{e}risaton fonctionnelle. Nous souhaitons nous limiter dans ce texte \`a l'\'etude de cette th\'eorie sur des mesures de probabilit\'e sur $\mathbb{R}^k$, $k\geq 1$.\\

\noindent Cependant, il arrive que les preuves soient exactement les m\^emes que le cas g\'en\'eral pour des mesures de probabilit\'e sur un espace espace m\'{e}trique $(S,d)$ muni de sa $\sigma $-alg\`{e}bre bor\'{e}lienne $\mathcal{B}(S)$, ne faisant intervenir que les propri\'et\'es g\'en\'erales de la m\'etrique. Dans de telles situations, nous \'enon\c{c}ons la propri\'et\'e et les preuves dans le cas g\'en\'eral.\\

\noindent Pour traiter de la tension de suites de mesures de probabilit\'e, c'est-\`a-dire, de l'existence pour ces suites, de sous-suites convergentes au sens vague, le traitement sera fait sp\'{e}cifiquement \`{a} $\mathbb{R}^{k}$ par le biais du th\'eor\'eme de Helly-Bray.\\

\noindent Comme dans toute th\'{e}orie limite, il faut necessairement traiter de l'unicit\'{e} de la limite, de crit\`{e}res de convergence, de compacit\'{e} relative et de transformations \'{e}ventuelles. Pour compacit\'{e} relative, nous parlerons plut\^{o}t de tension.

\section{D\'{e}finition, Unicit\'{e} et Th\'{e}or\`{e}me Portmanteau}

\begin{definition} \label{cv.cvDEF}
La suite d'applications mesurables $X_{n}:(\Omega _{n},\mathcal{A}_{n},\mathbb{P}_{n})\mapsto (S,B(S))$ converge vaguement vers 
$X :\ (\Omega_{\infty},\mathcal{A}_{\infty},\mathbb{P}_{\infty})\mapsto (S,\mathcal{B}(S))$\ mesurable ssi pour toute fonction 
$f:S\mapsto \mathbb{R},$\ continue et born\'{e}e (not\'{e} $f\in \mathcal{C}_{b}(S)$\ ), 
\begin{equation}
\mathbb{E}f(X_{n}))\rightarrow \mathbb{E}f(X) \text{ quand } n \rightarrow +\infty. \label{cv21}
\end{equation}
\end{definition}

\bigskip \noindent Nous remarquons que les espaces de d\'{e}part n'ont aucune importance dans cette th\'{e}orie, d'o\`{u} le nom de convergence
vague. Notons $L=\mathbb{P}_{X}=\mathbb{P}_{\infty}\circ X^{-1}$, la loi de $X$ d\'{e}finie 

\begin{equation*}
\forall \text{ }B\in \mathcal{B}_{\infty}(S),\text{ }L(B)=\mathbb{P}_{\infty}(X^{-1}(B))=%
\mathbb{P}(X\in B),
\end{equation*}

\bigskip \noindent et pour tout $n\geq 1$, $\mathbb{P}^{(n)}$ la loi de probabilit\'e de $X_n$ d\'efinie par   

\begin{equation*}
\forall \text{ } B\in \mathcal{B}(S),\text{ } \mathbb{P}^{(n)}(B)=\mathbb{P}_{n}(X_{n}^{-1}(B))=\mathbb{P}_n(X_{n}\in B).
\end{equation*}

\bigskip \noindent La d\'efinition dit que $X_n$ converge vaguement vers $X$ si et seulement si pour tout $f\in \mathcal{C}_{b}(S)$,
\begin{equation*}
\int_{S}\text{ }f(x)\text{ }d\mathbb{P}^{(n)}(x) \rightarrow \int_{S}\text{ }f(x)\text{ }dL(x)  \text{ quand } n \rightarrow +\infty.
\end{equation*}

\noindent On pourrait ainsi remplacer (\ref{cv21}) par 
\begin{equation}
\mathbb{E}f(X_{n }))\rightarrow \int_{S}\text{ }f\text{ }dL \text{ quand } n \rightarrow +\infty,  \label{cv22}
\end{equation}

\noindent et dire que la suite ($X_{n})_{n\geq 1\ }$converge vaguement vers une
probabilit\'{e} $L$. Dans la suite, nous utiliserons les deux terminologies.\\

\bigskip \noindent \textbf{Faites attention \`a ce point important}. Il est aussi important de voir que les symboles d'\'esp\'erance math\'ematique dans (\ref{cv21}) d\'ependent des mesures de probabilit\'e qu'ils utilisent. En cons\'equence, ils doivent \^etre labellis\'es par exemple sous la forme,

$$
\mathbb{E}_{\infty}(f(X))=\int f(X) d\mathbb{P}_{\infty}, \text{  }\mathbb{E}_{n}(f(X_n))=\int f(X_n) d\mathbb{P}_{n}, \text{  }n\geq 1.
$$
 
\noindent avec $\mathbb{E}_{\infty}=\mathbb{E}_{\mathbb{P}_{\infty}}$ et $\mathbb{E}_{n}=\mathbb{E}_{\mathbb{P}_{n}}$. Mais, nous avons choisi d'omettre ces labels en vue de garder les notations simples. Nous les utiliserons lorsque cela est n\'ecessaire. Ainsi, on ne doit utiliser la lin\'earit\'e de l'esp\'erance ou effecter des op\'erations sur les esp\'erances que lorsque les espaces de probabilit\'e de d\'epart sont les 
m\^emes, comme dans la section \ref{cv.CvCp}.\\

\noindent Nous allons montrer que la limite est unique en distribution dans le sens
suivant.\\

\bigskip 

\begin{proposition} \label{tm_cv_prop1}
Soit une suite d'applications mesurables $X_{n}:(\Omega _{n},\mathcal{A}_{n},\mathbb{P}_{n})\mapsto (S,B(S))$, et soit  $\mathbb{Q}_{1}$ et $\mathbb{Q}_{2}$ deux probabilit\'{e}s sur $(S,\mathcal{B}(S)).$ Supposons que $X_{n}$
converge vaguement \`{a} la fois vers $\mathbb{Q}_{1}$ et $\mathbb{Q}_{2}.$ 
Alors, n\'ecessairement, nous avons  
\begin{equation*}
\mathbb{Q}_{1}=\mathbb{Q}_{2}.
\end{equation*}
\end{proposition}

\noindent \textbf{PREUVE}. Supposons que  $X_{n}$ converge vaguement \`{a} la fois
vers $\mathbb{Q}_{1}$ et $\mathbb{Q}_{2}.$ Pour montrer que $\mathbb{P}_{1}=%
\mathbb{P}_{2}$, il suffit de montrer qu'elles ont \'{e}gales sur la classe
des ouverts $\Theta $ de $(S,d)$. \ En effet $\Theta $ est un $\pi $-system
qui engendre $\mathcal{B}(S)$. Alors deux probabilit\'{e}s sur $(S,\mathcal{B%
}(S))$\ qui sont \'{e}gales  sur $\Theta ,$\ sont par ailleurs \'{e}gales
sur $\mathcal{B}(S).$\\

\noindent Soit $G$ un ouvert de $S$. Pour tout entier $m\geq 1,$\ posons $f_{m}(x)=\min
(m$\ $d(x,G^{c}),1).$\ Pour tout $m$, la fonction $f_{m}$\ est \`{a} valeurs
dans $[0,1]$, donc born\'{e}e. Puisque $G^{c}$\ est ferm\'{e}e, on a 
\begin{equation*}
d(x,G^{c})=\left\{ 
\begin{array}{c}
>0\text{ si }x\in G \\ 
0\text{ si }x\in G^{c}\text{ }%
\end{array}%
\right. .
\end{equation*}

\bigskip \noindent Montrons que $f_{m}$\ est lipschitzienne. Evaluons $%
\left\vert f_{m}(x)-f_{m}(y)\right\vert $\ selon trois cas.

\bigskip \noindent Cas 1. $(x,y)\in (G^{c})^{2}$. Donc 
\begin{equation*}
\left\vert f_{m}(x)-f_{m}(y)\right\vert =0\leq m\text{ }d(x,y).
\end{equation*}

\bigskip \noindent Cas 2. $x\in G$\ et $y\in G^{c}$\ (ou en permutant les r%
\^{o}les de x et y). On a 
\begin{equation*}
\left\vert f_{m}(x)-f_{m}(y)\right\vert =\left\vert \min
(md(x,G^{c}),1)\right\vert \leq m\text{ }d(x,G^{c})\leq m\text{ }d(x,y),
\end{equation*}

\bigskip \noindent par d\'{e}finition m\^{e}me de $d(x,G^{c})=\inf \{d(x,z),$%
\ z$\in G^{c}\}.$\newline

\bigskip \noindent Cas 3. Soit $(x,y)\in G^{2}.$\ On a, en utilisant la
propri\'{e}t\'{e} (\ref{annexe2}) de la section appendice \ref{cv.annexe}
\begin{equation*}
\left\vert f_{m}(x)-f_{m}(y)\right\vert =\left\vert \min
(md(x,G^{c}),1)-\min (md(y,G^{c}),1)\right\vert \leq \left\vert
md(x,G^{c})-md(y,G^{c})\right\vert ,
\end{equation*}%
\begin{equation*}
\leq m\text{ }d(x,y)
\end{equation*}

\bigskip \noindent par l'in\'{e}galit\'{e} triangulaire bis. Donc $f_{m}$\
est m-lipschitzienne. De plus 
\begin{equation*}
f_{m}\uparrow 1_{G}\text{ quand m}\uparrow \infty .
\end{equation*}

\bigskip \noindent En effet, si x$\in G^{c},$\ $f_{m}(x)=0\uparrow
0=1_{G}(x).$\ Si $x\in G,\ d(x,G^{c})>0$\ et $md(x,G^{c})\uparrow \infty .$\
Pour $m$ assez grand, 
\begin{equation}
f_{m}(x)=1\uparrow 1_{G}(x)  \label{limfm}
\end{equation}

\bigskip \noindent En r\'{e}sum\'{e}, chaque fonction  $f_{m}$\ est
lipschitzienne, born\'{e}e, positive.\\ 

\noindent Une fois que nous avons les propri\'{e}t\'{e}s de $f_{m},$ revenons \`{a}
notre hypoth\`{e}se \`{a} savoir que $X_{n}$ converge vaquement \`{a} la
fois vers $\mathbb{P}_{1}$ et $\mathbb{P}_{2}$, c'est-\`{a}-dire que que
pour toute fonction $f:S\mapsto \mathbb{R},$\ continue et born\'{e}e, nous
avons, quand $n \rightarrow +\infty$,
\begin{equation}
\mathbb{E}f(X_{n}))\rightarrow \int f\text{ }d\mathbb{Q}_{1}\text{ \ and \ }%
\mathbb{E}f(X_{n}))\rightarrow \int f\text{ }d\mathbb{Q}_{2}.
\end{equation}

\noindent Par l'unicit\'{e} de la limite de suites r\'{e}elles, nous avons%
\begin{equation*}
\forall (f\in C_{b}(S)),\int f\text{ }d\mathbb{Q}_{1}=\int f\text{ }d\mathbb{%
Q}_{2.}
\end{equation*}%

\noindent On peut donc appliquer cette \'{e}galit\'{e} pour toutes les fonctions $f_{m}
$ d\'{e}finies ci-haut. Cela nous donnera

\begin{equation*}
\forall (m\geq 1),\text{ }\int f_{m}\text{ }d\mathbb{Q}_{1}=\int f_{m}\text{ 
}d\mathbb{Q}_{2.}
\end{equation*}

\noindent En faisant croitre $m$ vers l'infini, en utilisant (\ref{limfm})\ et en
appliquant le Th\'{e}or\`{e}me de convergence monotone, nous obtenons%
\begin{equation*}
\int 1_{G}\text{ }d\mathbb{Q}_{1}=\int 1_{G}\text{ }d\mathbb{Q}_{2},
\end{equation*}

\noindent signifiant  
\begin{equation*}
\mathbb{Q}_{1}(G)=\mathbb{Q}_{2}(G).
\end{equation*}

\noindent Puisque $G$ a \'{e}t\'{e} arbitrairement fix\'{e}, il s'ensuit que l'\'{e}%
galit\'{e} est vraie pour tous les ouverts de $S$. Alors, par nos remarques
liminaires, nous avons bien $\mathbb{Q}_1=\mathbb{Q}_{2}.$\\

\bigskip \noindent \textbf{Notation}. Lorsque $(X_{n})_{n\geq 1}$ converge vaguement $X$, nous utilisons la notation principale
\begin{equation*}
X_{n}\rightsquigarrow X \text{ quand } \rightarrow +\infty,
\end{equation*}

\noindent mais il nous arrivera aussi d'utiliser ces deux autres : $X_{n}\rightarrow _{\mathcal{L}}X$ (pour convergence en loi) or $X_{n}\rightarrow _{d}X$ (pour convergence en distribution) ou $X_{n}\rightarrow _{w}X$ (pour convergence vague, \textit{weakly} en Anglais).\\ 

\bigskip \noindent Nous avons maintenant besoin de la caract\'{e}risation de
cette convergence. Selon les besoins, on peut avoir besoin d'autres angles
d'attaque, pour l'\'{e}tablir.\\

\begin{theorem} \label{cv.theo.portmanteau} (\textbf{Portmanteau}). La suite d'applications $X_{n }:(\Omega _{n },\mathcal{A}%
_{n },P_{n })\mapsto (S,B(S))$ converge vaguement vers une probabilit\'{e} $L$ ssi

\bigskip \noindent (ii) Pour tout ouvert $G$ de S , 
\begin{equation*}
\liminf_{n\rightarrow +\infty} \mathbb{P}_n(X_{n }\in G)\geq L(G).
\end{equation*}

\bigskip \noindent (iii) Pour tout ferm\'{e} $F$ de $S$, 
\begin{equation*}
\limsup_{n\rightarrow +\infty} \mathbb{P}_n(X_{n }\in F)\leq L(F).
\end{equation*}

\bigskip \noindent (iv) Pour toute fonction $f$ semi-continue inf\'{e}rieurement (\textit{s.c.i}) et minor\'{e}e, 
\begin{equation*}
\liminf_{n\rightarrow +\infty} \mathbb{E}f(X_n)\geq \int f \text{ }dL.
\end{equation*}

\bigskip \noindent (v) Pour toute fonction $f$ semi-continue sup\'{e}rieurement (\textit{s.c.s}) et major\'{e}e,
 
\begin{equation*}
\limsup_{n\rightarrow +\infty} \mathbb{E}f(X_n)\leq \int f dL.
\end{equation*}

\bigskip \noindent (vi) Pour tout bor\'{e}lien $B$ de S tel que $L(\partial
B)=0,$%
\begin{equation*}
\lim_{n\rightarrow +\infty} \mathbb{P}_n(X_{n }\in B)=\lim_{n\rightarrow +\infty} \mathbb{P}_n(X_{n }\in B)=L(B).
\end{equation*}

\bigskip \noindent (vii) Pour toute fonction $f$ positive, born\'{e}e et
lipschitzienne. 
\begin{equation*}
\liminf_{n\rightarrow +\infty} \mathbb{E}f(X_n)\geq \int f \\
dL.
\end{equation*}
\end{theorem}

\bigskip

\bigskip \noindent Avant de commencer la preuve, rappelons que $\partial B$ est la fronti\`{e}re de $B$. Si $L(\partial B)=0,$\ on dit que $B$ est $L$ -continue. Quant aux fonctions semi-continues, nous ferons un rappel dans la sous-section \ref{cv.subsec.annexe.semic} de la section annexe \ref{cv.annexe}.\\

\bigskip \noindent \textbf{PREUVE}.\\

\noindent Les points (ii) et (iii) sont \'{e}quivalents par compl\'{e}mentation. De m\^{e}me pour les points (iv) et (v)
en passant de $f$ \`{a} $-f$ et en utilisant les propri\'{e}t\'{e}s vues dans la sous-section \ref{cv.subsec.annexe.semic} de la section annexe \ref{cv.annexe} ci-bas. Maintenant notons (i) la formule (\ref{cv21}) de la d\'{e}finition de la convergence vague. L'implication $(i)\Rightarrow (vii)$\ est 
\'{e}vidente car une fonction lipschitzienne (de param\`{e}tre $k$), c'est 
\`{a} dire, telle que 
\begin{equation*}
\forall (x,y)\in S^{2},\left\vert f(x)-f(y)\right\vert \leq k\text{ }d(x,y),
\end{equation*}

\noindent est continue.\\

\noindent Prouvons $(vii)\Rightarrow (ii)$. Soit $G$ un ouvert de $S$. Pour tout entier $m\geq 1,$\ posons $f_{m}(x)=\min (m$\ $d(x,G^{c}),1).$\ Pour tout $m$, la
fonction f$_{m}$\ est \`{a} valeurs dans [0,1], donc born\'{e}e. Nous avons d%
\'{e}j\`{a} recontr\'{e} cette finction dans la preuve de la proposition \ref{tm_cv_prop1}. Nous savons que chaque $m$, $f_{m}$\ est lipschitzienne, born\'{e}e, positive et
\begin{equation*}
f_{m}\uparrow 1_{G}\text{ quand m}\uparrow \infty .
\end{equation*}

\bigskip \noindent Nous avons pour tout $n\geq 1$ et pour tout $m\geq 1,$ 
\begin{equation*}
\mathbb{E}(1_{G}(X_{n}))\geq \mathbb{E}f_{m}(X_{n}).
\end{equation*}

\bigskip \noindent Appliquons (vii) pour avoir 
\begin{equation}
\liminf_{n\rightarrow +\infty}\text{ }\mathbb{E}(1_{G}(X_{n}))\geq \liminf_{n\rightarrow +\infty}\mathbb{E}f_{m}(X_{n})\geq \int f_{m}\text{ }dL.
\label{cv23}
\end{equation}

\noindent Or pour toute mesure de probabilit\'e et pour toute partie mesurable

\begin{equation*}
\mathbb{E}_{\mathbb{Q}}(1_{B})=\mathbb{Q}(B)
\end{equation*}

\bigskip \noindent Pour $B=1_{X_{n}^{-1}(G)}=1_{(X_{n}\in G)},$\ et en
passant \`{a} la limite sur m et en utilisant le th\'{e}or\`{e}me de
convergence monotone \`{a} (\ref{cv23}), on obtient 
\begin{equation*}
\liminf_{n\rightarrow +\infty}\mathbb{P}_n(X_{n}\in G)\geq \int 1_{G}\text{ 
}dL=L(G)
\end{equation*}

\bigskip \noindent Donc (ii) est d\'{e}montr\'{e}.\newline

\bigskip \noindent Prouvons que : $(ii)\Rightarrow (iv).$\newline

\noindent Soit (ii) vraie. Soit $f$ une fonction semi-continue inf%
\'{e}rieurement minor\'{e}e par M. Nous pouvons prouver (iv) pour $f-M=g$
positive, qui est encore s.c.i. Alors les ensembles ($g\leq c)$\ sont ferm\'{e}s selon la proposition
 \ref{cv.annexe.SC} de la section annexe \ref{cv.annexe}. Posons pour $m\geq 1$\ $%
fix\acute{e}.$\ 
\begin{equation*}
G_{i}=\{g>i/m\},\text{ i}\geq 1.
\end{equation*}

\noindent et 
\begin{equation*}
g_{m}=\frac{1}{m}\sum_{i=1}^{m^{2}}1_{G_{i}}.
\end{equation*}

\noindent Les ensembles $G_{i}$\ sont ouverts car $g$ est s.c.i. Remarquons que 
\begin{equation}
g_{m}(x)=\frac{i}{m}\text{ }pour\text{ }\frac{i}{m}<g(x)\leq \frac{i+1}{m},%
\text{ }pour\text{ }i=1,...,m^{2}-1  \label{cv24}
\end{equation}

\noindent et
 
\begin{equation*}
g_{m}(x)=m\text{ }pour\text{ }g(x)>m.
\end{equation*}

\noindent Donc, nous avons 
 
\begin{equation*}
g_{m}\leq g.
\end{equation*}

\noindent De plus, d'apr\`{e}s (\ref{cv24}), il est vrai que 
\begin{equation*}
\left| g_{m}(x)-g(m)\right| \leq 1/m\text{ }pour\text{ }g(x)\leq m.
\end{equation*}

\noindent On a 
\begin{equation*}
g(X_{n })\geq g_{m}(X_{n })=\frac{1}{m}\sum_{i=1}^{m^{2}}1_{G_{i}}(X_{n })=%
\frac{1}{m}\sum_{i=1}^{m^{2}}1_{(X_{n }\in G_{i})},
\end{equation*}

\noindent signifiant

\noindent 
\begin{equation}
\mathbb{E}g(X_{n })\geq \mathbb{E}g_{m}(X_{n })=\frac{1}{m}\mathbb{E}.
\sum_{i=1}^{m^{2}}1_{(X_{n }\in G_{i})}.  \label{cv25}
\end{equation}

\bigskip \noindent Donc (\ref{cv25}) donne 
\begin{equation*}
\mathbb{E}g(X_{n })\geq \mathbb{E}g_{m}(X_{n })\geq \frac{1}{m}%
\sum_{i=1}^{m^{2}}\mathbb{E}1_{(X_{n }\in G_{i})} = \frac{1}{m}%
\sum_{i=1}^{m^{2}}\mathbb{P}_n(X_{n }\in G_{i}).
\end{equation*}

\noindent En passant \`{a} la limite sur $n ,$\ et en appliqunant (ii), on
aura 
\begin{equation*}
\liminf_{n\rightarrow +\infty} \mathbb{E}g(X_{n })\geq \liminf_{n\rightarrow +\infty} \text{ }\mathbb{E}g_{m}(X_{n
})\geq \frac{1}{m}\sum_{i=1}^{m^{2}}L(G_{i})=\int g_{m}\text{ }dL\geq
\int_{(g\leq m)}g_{m}\text{ }dL,
\end{equation*}

\noindent et donc

\begin{equation*}
\geq \int_{(g\leq m)}g\text{ }dL+\int_{(g\leq m)}(g_{m}-g)\text{ }dL.
\end{equation*}

\noindent Quand $m\rightarrow \infty$

\begin{equation*}
\int_{(g\leq m)}g\text{ }dL\rightarrow \int g\text{ }dL
\end{equation*}

\noindent et
 
\begin{equation*}
|\int_{(g\leq m)}(g_{m}-g)\text{ }dL|\leq L(S)/m\rightarrow 0.
\end{equation*}

\noindent D'o\`{u} 

\begin{equation*}
\liminf_{n\rightarrow +\infty} \text{ }\mathbb{E}g(X_{n })\geq \int g\text{ }dL.
\end{equation*}

\bigskip \noindent Maintenant, en rempla\c{c}ant $g$ par $f-M$, la m\^{e}me
formule reste vraie, par simplification de M. Donc (iv) vraie.

\bigskip \noindent Prouvons $(ii)\Rightarrow (vi).$\ Rappelons que $\partial
B=\overline{B}-int(B),$\ autrement dit, $B$ est la diff\'{e}rence entre la
fermeture de B et de son int\'{e}rieur. Donc, puisque $int(B)\subseteq
B\subseteq \overline{B}$%
\begin{equation}
L(\partial B)=L(int(B))-L(\overline{B})=0 \Rightarrow L(int(B))=L(\overline{B%
})=L(B)  \label{cv26a}
\end{equation}

\bigskip \noindent Puisque $int(B)$ est ouvert et que $\overline{B}$ est ferm\'e, on peut utiliser (ii) et
(iii) \`{a} la fois pour avoir 
\begin{equation}
L(int(B)) \leq \liminf_{n\rightarrow +\infty} \mathbb{P}_n(X_{n}\in int(B))\leq \limsup_{n\rightarrow +\infty} \mathbb{P}_n(X_{n }\in B),  \label{cv26b}
\end{equation}

\begin{equation}
\leq \mathbb{P}_n(X_{n}\in \overline{B}) \leq L(int(B),  \label{cv26X}
\end{equation}

\bigskip \noindent D'o\`{u}, par (\ref{cv26a})

\begin{equation*}
L(B)=\lim_{n\rightarrow +\infty} \mathbb{P}_n(X_{n }\in B)=\lim_{n\rightarrow +\infty} \mathbb{P}_n(X_{n}\in B).
\end{equation*}

\bigskip \noindent Ce qui \'{e}tait \`{a} d\'{e}montrer.\newline

\bigskip \noindent Prouvons que $(vi)\Rightarrow (iii).$

\bigskip \noindent Soit (vi) vraie et soit $F$ un ferm\'{e} de $S$. Posons ($%
F_{\varepsilon})=\{x,$\ $d(x,F)\leq \varepsilon \}$\ pour $\varepsilon \geq 0.$\ On 
\begin{equation*}
F\subseteq F(\varepsilon )
\end{equation*}

\noindent et
 
\begin{equation*}
F(\varepsilon )\downarrow F\text{ }pour\text{ }\varepsilon \downarrow 0
\end{equation*}

\bigskip \noindent Maintenant $\partial F(\varepsilon )\subseteq \{x,$\ d(x,F)=$%
\varepsilon \}.$\ Donc les ensembles $\partial F(\varepsilon )$\ sont disjoints,
donc au plus un nombre d\'{e}nombrable d'ensembles parmi eux, ont une
probabilit\'{e} non nulle (voir Proposition \ref{cv.subsec.annexe.Disjoint} de la section annexe \ref{cv.annexe} ci-bas). Donc on peut trouver
une suite $\varepsilon _{n}\downarrow 0$\ telle que pour tout $n$, 
\begin{equation*}
L(\partial F(\varepsilon _{n}))=0.
\end{equation*}

\noindent Pour $n$ fix\'{e}, $F\subseteq F(\varepsilon _{n})$,
\begin{equation*}
\limsup_{n\rightarrow +\infty} \mathbb{P}_n(X_{n }\in F)\leq \limsup_{n\rightarrow +\infty} \mathbb{P}_n(X_{n }\in
F(\varepsilon _{n}))
\end{equation*}

\noindent et par application de (vi) 
\begin{equation*}
\limsup_{n\rightarrow +\infty} \mathbb{P}_n(X_{n }\in F)\leq \limsup_{n\rightarrow +\infty} \mathbb{P}_n(X_{n }\in
F(\varepsilon _{n}))=L(F(\varepsilon _{n})).
\end{equation*}

\bigskip \noindent Maintenant, en passant \`{a} la limite quand $n\uparrow\infty$, 
\begin{equation*}
\limsup_{n\rightarrow +\infty} \mathbb{P}_n(X_{n }\in F)\leq L(F),
\end{equation*}

\noindent ce qui est bien (iii).\newline

\bigskip \noindent Prouvons $(iv)\Rightarrow (i)$\newline

\bigskip \noindent Si (iv) est vraie, alors (v) est vraie. Donc une fonction
$f$ continue et born\'{e}e, est \`{a} la fois s.c.i. et minor\'{e}e, et s.c.s.
et major\'{e}e, on aura 
\begin{equation*}
\int fdL\leq \liminf_{n\rightarrow +\infty} \mathbb{E}f(X_{n })\leq \liminf_{n\rightarrow +\infty} \mathbb{E}f(X_{n
})\leq \limsup_{n\rightarrow +\infty} \sup \mathbb{E}f(X_{n })\leq \int fdL
\end{equation*}

\bigskip \noindent D'o\`{u} 
\begin{equation*}
\int fdL=\lim_{n\rightarrow +\infty} \mathbb{E}f(X_{n })=\lim_{n\rightarrow +\infty} \mathbb{E}f(X_{n}).
\end{equation*}

\bigskip \noindent En r\'{e}sum\'{e} nous avons montr\'{e} le th\'{e}or\`{e}%
me par ce sch\'{e}ma 
\begin{equation*}
\begin{array}{ccccccc}
(i) & \Rightarrow & (vii) & \Rightarrow & (ii) & \Leftrightarrow & (iii) \\ 
\Uparrow &  &  &  & \Downarrow &  & \Uparrow \\ 
(v) & \Leftrightarrow & (iv) & = & 
\begin{array}{cc}
(iv) & (vi)%
\end{array}
& = & (vi)%
\end{array}%
\end{equation*}

\noindent qui montre que les six assertions sont \'{e}quivalentes
entre elles.\\

\section{Tansformations continues} \label{cv.mappingCv}

Soit une suite $(X_{n })_{n\geq 1}$ \`{a} valeurs dans l'espace m\'{e}trique $(S, d)$ qui converge vaguement vers $X$ et soit une application 
$g$ de $S$ dans un autre espace m\'{e}trique $(D, r)$. Alors, la suite  $g(X_{n})_{n\geq 1}$\ converge-t-elle vaguement vers $g(X)$?\\

\noindent Pour r\'epondre partiellement \`a cette question, supposons que $g$ soit continue. Alors, il est \'{e}vident que si $f\in C_{b}(D)$, alors $f\circ g\in C_{b}(S)$.\\

\noindent Donc, nous avons
 
\begin{equation*}
\forall \text{ }f\in C_{b}(S),\text{ }\mathbb{E}f(X_{n }))\rightarrow 
\mathbb{E}f(X)
\end{equation*}

 \noindent et puisque $f\in C_{b}(D)\Rightarrow f\circ g\in C_{b}(S)$,  cela m\`ene \`a :
\begin{equation*}
\forall \text{ }f\in C_{b}(D),\text{ }\mathbb{E}f\circ g(X_{n }))\rightarrow 
\mathbb{E}f\circ g(X).
\end{equation*}

\noindent D'o\`{u} la proposition :

\begin{proposition} Soit une suite d'applications $X_{n }:(\Omega _{n}, \mathcal{A}_{n},\mathbb{P}_{n
})\mapsto (S,B(S))$\ convergeant vaguement vers la variable al\'{e}atoire $X :(\Omega_{\infty},\mathcal{A}_{\infty},\mathbb{P}_{\infty})\mapsto (S,B(S))$\ (ou vers la probabilit\'{e} $L$ sur $S$), et $g : S \mapsto D$\ une application continue. Alors 
\begin{equation*}
g(X_{n })\rightarrow _{w}g(X).
\end{equation*}

\noindent ou par une autre \'{e}criture,
 
\begin{equation*}
g(X_{n })\rightarrow _{w}L\circ g^{-1}
\end{equation*}
\end{proposition}

\bigskip \noindent Cette proposition est tr\`{e}s importante. Mais, on a
plus que cela. En effet, on n'a pas besoin de la continuit\'{e} sur tout
l'ensemble $S$. Il suffira que l'ensemble des points de discontinuit\'{e} de $g$
soit de mesure nulle par rapport \`{a} $L$. Soit $discont(g)$ l'ensemble
mesurable des points de discontinuit\'{e} de $g$. D'apr\`es le lemme \ref{cv.annexe.Discont} de la Section Appendice \ref{cv.annexe} ci-bas, cet ensemble est mesurable.\\

\noindent Nous avons le r\'esultat g\'en\'eral suivant.\\

\begin{proposition} \label{cv.mappingTh} Soit une suite d'applications $X_{n }:(\Omega _{n }, \mathcal{A}_{n },\mathbb{P}_{n
})\mapsto (S,B(S)))$\ convergeant vaguement vers la variable al\'{e}atoire $X:(\Omega_{\infty},\mathcal{A}_{\infty},\mathbb{P}_{\infty})\mapsto (S,B(S))$\ (ou vers la probabilit\'{e} $L$ sur $S$), et $g : S \mapsto D$ une application telle que $L(discont(g))=\mathbb{P}(X\in discont(g))=0$. Alors 
\begin{equation*}
g(X_{n })\rightarrow _{w}g(X)
\end{equation*}

\noindent ou par une autre \'{e}criture, 

\begin{equation*}
g(X_{n })\rightarrow _{w}L\circ g^{-1}.
\end{equation*}
\end{proposition}

\bigskip \noindent \textbf{PREUVE}. Soit $X_{n }\rightarrow _{w}L$\ avec $L(discont(g))=0$. Soit $F$ une partie ferm%
\'{e}e de D. Montrons que nous avons
\begin{equation}
\overline{g^{-1}(F))}\text{ }\subseteq g^{-1}(F)\cup discont(g).  \label{cv28}
\end{equation}

\bigskip \noindent En effet soit x$\in \overline{g^{-1}(F)}.$\ Donc il
existe une suite $(y_{n})_{n\geq 1}\in $\ $g^{-1}(F)$\ telle que y$%
_{n}\rightarrow x$\ et pour tout n$\geq 1,$\ $g(y_{n})\in F.$\ Alors \textit{ou bien} 

$$x\in discont(g),
$$

\noindent \textit{ou bien} $x$ est un point de continuit\'{e}. Dans ce dernier cas, $(y_{n}\rightarrow x,$\ $g(y_{n})\in F$\ $)\Rightarrow $\ $g(y_{n})\rightarrow
g(x)\in F$\ puisque $F$ est ferm\'{e} et donc :
\begin{equation*}
x\in g^{-1}(F)
\end{equation*}

\noindent Nous voyons donc que (\ref{cv28}) est \'etablie. Combinons cette formule avec le point (iii) du th\'{e}
or\`{e}me Portmanteau \ref{cv.theo.portmanteau}. Nous avons 

\begin{equation*}
\limsup_{n\rightarrow +\infty} \mathbb{P}_n(g(X_{n })\in F)=\limsup_{n\rightarrow +\infty} \mathbb{P}_n(X_{n }\in
g^{-1}(F))\leq \limsup_{n\rightarrow +\infty} \mathbb{P}_n(X_{n }\in \overline{g^{-1}(F))}\text{ })
\end{equation*}

\noindent et
 
\begin{equation*}
\limsup_{n\rightarrow +\infty} \mathbb{P}_n(X_{n }\in \overline{g^{-1}(F))}\text{ })\leq L(\overline{g^{-1}(F))}\text{ )}\leq L(g^{-1}(F))+L(discont(g)).
\end{equation*}

\noindent Ce qui donne
 
\begin{equation*}
\limsup_{n\rightarrow +\infty} \mathbb{P}_n(g(X_{n })\in F\text{ })\leq L\circ g^{-1}(F).
\end{equation*}

\noindent D'o\`{u} le r\'{e}sultat cherch\'{e}.\\

\section{Cas particulier de $\mathbb{R}^{k}$, $k\geq 1$}

Int\'{e}ressons nous au cas particulier $S=\mathbb{R}^{k}$. Avant d'aller plus loin, rappelons les caract\'erisations de lois de probabilit\'es sur 
$S=\mathbb{R}^{k}$. \\

\noindent Pour commencer, faisons cette pr\'ecision.\\

\noindent \textbf{Terminilogie} Dans toute la partie traitant de $\mathbb{R}^{k}$, un vecteur al\'atoire de dimension $k\geq 1$ est simplement une application mesurable d\'efinie sur un espace mesurable donn\'e \`a valeurs dans $\mathbb{R}^{k}$.\\

\noindent Ensuite, adoptons ces notations relatives aux coordonnées d'un vecteur $X$ et d'une suite de vecteurs $(X_n)_{n\geq 1}$ :\\

$$
X=\left[ 
\begin{array}{c}
X_{1} \\ 
\cdot \cdot \cdot \\ 
X_{k}%
\end{array}
\right],
$$

\noindent 

$$
X_{n}=\left[ 
\begin{array}{c}
X^{(n)}_{1} \\ 
\cdot \cdot \cdot \\ 
X^{(n)}_{k}
\end{array}
\right].
$$

\bigskip \noindent Nous arppelons (voir \cite{bmtp} dans notre s\'erie) que la loi de probabilit\'e d'un vecteur al\'eatoire $X$ de $\mathbb{R}^{k}$ est caract\'eris\'ee par sa \textit{fonction de r\'epartition}, d\'efinie par  
\begin{equation*}
^{t}(t_{1},t_{2},...,t_{k})\mapsto F_{X}(t_{1},t_{2},...,t_{k})=\mathbb{P}(X_{1}\leq t_{1},X_{2}\leq t_{2},...,X_{k}\leq t_{k}).
\end{equation*}

\bigskip \noindent Elle est aussi caract\'eris\'ee par sa \textit{fonction caract\'eristique} donn\'ee par : 
\begin{equation*}
^{t}(u_{1},u_{2},...,u_{k})\mapsto \Phi _{X}(u_{1},u_{2},...,u_{k})=\mathbb{E}(\exp (\sum_{j}^{k}i\text{ u}_{j}X_{j}))
\end{equation*}

\noindent Si sa \textit{fonction des moments} donn\'ee 
\begin{equation*}
^{t}(u_{1},u_{2},...,u_{k})\mapsto \Psi _{X}(u_{1},u_{2},...,u_{k})=\mathbb{E}(\exp (\sum_{j}^{k}\text{ u}_{j}X_{j})),
\end{equation*}

\noindent existe autour d'un voisinage du vecteur nul, alors elle d\'efinit la loi de probabilit\'e de $X$ de mani\`ere unique.\\

\noindent Il en est de m\^eme pour \textit{densité de probabilit\'e} par rapport \`a la mesure de Lebesgues sur $\mathbb{R}^{k}$. Celle-ci existe si et seulement si 
\begin{equation*}
^{t}(t_{1},t_{2},...,t_{k})\mapsto f_{X}(t_{1},t_{2},...,t_{k})=\frac{%
\partial ^{(k)}F_{X}(t_{1},t_{2},...,t_{k})}{\partial t_{1}\partial
t_{2}\cdot \cdot \cdot \partial t_{k}},
\end{equation*}

\noindent \textit{presque partout} - en $(t_{1},t_{2},...,t_{k}))$ - par rapport \`a la mesure de Lebesgues sur $\mathbb{R}^{k}$.\\

\noindent Il est remarquable que ces m\^emes d\'eterminants jouent aussi les grands r\^oles en convergence vague\\

\noindent Avec les notations d\'ej\`a introduites, nous avons les propositions suivantes.\\

\begin{proposition} \label{cv.FRDIR}
Une suite de variables al\'eatoirs $X_{n }: (\Omega _{n }, \mathcal{A}_{n },\mathbb{P}_{n})\mapsto (\mathbb{R}^{k},\mathcal{B}(\mathbb{R}^{k}))$\ converge vaguement vers la variable al\'{e}atoire $X : (\Omega_{\infty},\mathcal{A}_{\infty},\mathbb{P}_{\infty})\mapsto (\mathbb{R}^{k},\mathcal{B}(\mathbb{R}^{k}))$, alors pour tout point $t=(t_{1},t_{2},...,t_{k})$ de continuit\'{e} de $F_{X}$\ ou de $L$,

\begin{equation}
\mathbb{P}_n(X_{n }\in \prod_{i=1}^{k}\left] -\infty ,t_{i}\right]
)\rightarrow F_{X}(t_{1},t_{2},...,t_{k}).  \label{cv32}
\end{equation}%
\end{proposition} 

\bigskip \noindent \textbf{PREUVE}. Soit la fonction de r\'{e}partition de X 
\begin{equation*}
F_{X}(t_{1},t_{2},...,t_{k})=\mathbb{P}_{\infty}(X_{1}\leq t_{1},X_{2}\leq
t_{2},...,X_{k}\leq t_{k})
\end{equation*}%
\begin{equation*}
=\mathbb{P}_{\infty}(X\in \prod_{i=1}^{k}\left] -\infty ,t_{i}\right] )
\end{equation*}

\noindent Notons $t=(t_{1},...,t_{k})$\ et $t(n)=(t_{1}(n),t_{2}(n),...,t_{k}(n)).$\
On dira que $t(n)\uparrow t$\ (resp $t(n)\downarrow t)$\ ssi 

\begin{equation*}
\forall (1\leq i\leq k),\text{ }t_{i}(n)\uparrow t_{i}\text{ }(resp.\text{ }%
t_{i}(n)\downarrow t_{i})
\end{equation*}%

\noindent Posons A(t)= $\prod_{i=1}^{k}\left] -\infty ,t_{i}\right] .$\ Nous avons
quand $n\uparrow \infty $,

\begin{equation*}
A(t(n))\downarrow A(t)
\end{equation*}%

\noindent et donc, par la limite monotone des probabilit\'{e}s,
 
\begin{equation*}
F_{X}(t)=\mathbb{P}_{\infty}(X\in A(t(n))\downarrow \mathbb{P}_{\infty}(X\in A(t))=F_{X}(t)
\end{equation*}%
Par suite, F$_{X}$\ est continue \`{a} droite en tout t. Mais 
\begin{equation*}
A(t(n))\uparrow A^{+}(t)=\prod_{i=1}^{k}\left] -\infty ,t_{i}\right[ 
\end{equation*}%

\noindent et par suite
 
\begin{equation*}
F_{X}(t)=\mathbb{P}_{\infty}(X\in A(t(n))\uparrow \mathbb{P}_{\infty}(X\in A^{+}(t))
\end{equation*}

\noindent Mais nous avons 
\begin{eqnarray}
D(t)&=&A(t)\text{ }\backslash \text{ }A^{+}(t)\\
&=&\{x=(x_{1},...,x_{k})\in A(t),\exists 1\leq i\leq k,\text{ }x_{i}=t_{i}\} \notag
\end{eqnarray}

\noindent Pour mieux comprendre cette formule, regardez la pour k=1 
\begin{equation*}
D(a)=]-\infty ,\text{ }a]\text{ }\backslash \text{ }]-\infty ,\text{ }a[=\{a\}
\end{equation*}

\noindent et pour k=2 
\begin{equation*}
D(a,b)=]-\infty ,\text{ }a]\text{ }\times \text{ }]-\infty ,\text{ }b]\text{ }%
\backslash \text{ }]-\infty ,\text{ }a[\text{ }\times \text{ }]-\infty ,%
\text{ }b[
\end{equation*}

\begin{equation*}
=\{(x,y)\in ]-\infty ,\text{ }a]\text{ }\times \text{ }]-\infty ,\text{ }b],%
\text{ }x=a\text{ ou }y=b\}
\end{equation*}%

\noindent D'o\`{u}, si 
\begin{equation}
\mathbb{P}_{\infty}(X\in D(t))=L(D(t))=0  \label{cv30}
\end{equation}

\noindent alors, quand n$\rightarrow \infty$,
\begin{eqnarray*}
F_{X}(t)&=&\mathbb{P}_{\infty}(X\in A(t(n))\uparrow \mathbb{P}_{\infty}(X\in A^{+}(t))\\
&=&\mathbb{P}_{\infty}(X\in A(t))-\mathbb{P}_{\infty}(X\in D(t))=F_{X}(t)
\end{eqnarray*}

\noindent Donc (\ref{cv30}) est la condition de continuit\'{e} de F$_{X}$\ en t. Mais
la fronti\`{e}re de A(t) est exactement D(t), i.e.,
\begin{equation}
\partial A(t)=D(t).
\end{equation}

\bigskip \noindent Car l'int\'{e}rieur de A(t) est s\^{u}rement A$^{+}(t).$\
Donc d'apr\`{e}s le point (vi) du th\'{e}or\`{e}me Portmanteau, nous avons
la partie directe de la proposition : quand $n \rightarrow +\infty$,

\begin{equation*}
(X_{n} \rightarrow_{w} X) \Longrightarrow (\mathbb{P}_n(X_n \in ]-\infty ,t]) \rightarrow F_{X}(t)\text{ pour }
F_{X} \text{ continue en } t.
\end{equation*}

\bigskip \noindent Ce qui finit la preuve.\\

\begin{proposition} \label{cv.FRINV}
Une suite d'applications $X_{n}:(\Omega _{n},\mathcal{A}_{n},P_{n})\mapsto (\mathbb{R}^{k},\mathbb{B}(\mathbb{R}^{k}))$\ mesurables
converge vaguement vers la variable al\'{e}atoire $X:\ (\Omega_{\infty},\mathcal{A}_{\infty},\mathbb{P}_{\infty})\mapsto (\mathbb{R}^{k},\mathcal{B}(\mathbb{R}^{k}))$, ssi pour tout point $t=(t_{1},t_{2},...,t_{k})$\ de continuit\'{e} de $F_{X}$, 
\begin{equation}
F_{X_{n}}(t)\rightarrow F_{X}(t_{1},t_{2},...,t_{k})  \label{CFR}
\end{equation}%
\end{proposition}

\noindent \textbf{PREUVE}. Supposons que pour tout point $%
t=(t_{1},t_{2},...,t_{k})$\ de continuit\'{e} de $F_{X}$, $\
F_{X_{n}}(t)\rightarrow F_{X}(t_{1},t_{2},...,t_{k}).$ Pour montrer que $%
X_{n}$ converge vaquement, nous allons montrer le point (ii) du Th\'{e}or%
\`{e}me Portmanteau, c'est-\`{a}-dire, pour tout ouvert $G$ de $\mathbb{R}%
^{k},$%
\begin{equation*}
\liminf_{n\rightarrow +\infty} P(X_{n}\in G)\geq P(X\in G).
\end{equation*}

\noindent Soit un ouvert $G$ de $\mathbb{R}^{k}.$ D'apr\`{e}s la proposition \ref{cv.GFcontinuous} de la section appendice \ref{cv.annexe} ci-bas, $G$ peut s'\'{e}crire sous la forme 
\begin{equation*}
G=\bigcup_{j\geq 1}]a^{j},b^{j}]
\end{equation*}

\bigskip \noindent avec $]a^{j},b^{j}]$ $F_{X}$-continu, c'est-\`{a}-dire
que pour tout point $c$ telle que%
\begin{equation*}
c_{i}=a_{i}^{(j)}\text{ ou }c_{i}=b_{i}^{(j)},
\end{equation*}

\noindent $c$ est un point de continuit\'{e} de $F_{X}.$ Gr\^{a}ce \`{a} la continuit%
\'{e} de la probabilit\'{e} $\mathbb{P}_{X}$, on peut trouver pour tout $%
\eta >0,$\ un rang $m$ tel que 
\begin{equation}
\mathbb{P}_{X}(G)-\eta \leq \mathbb{P}_{X}(\bigcup_{j=1}^{m}]a^{j},b^{j}])
\label{cv31b}
\end{equation}

\noindent Posons $A_{j}=]a^{j},b^{j}]$ a formule de Poincarr\'{e} donne 

\begin{eqnarray}
\mathbb{P}_{X}(\bigcup_{j=1}^{m}A_{j})&=&\sum \mathbb{P}_{X}(A_{j})-\sum 
\mathbb{P}_{X}(A_{i}A_{j})+\sum \mathbb{P}_{X}(A_{i}A_{j}A_{k}) \label{FP1}\\
&+&...+(-1)^{n+1}\mathbb{P}_{X}(A_{1}A_{2}...A_{n}) \notag
\end{eqnarray}

\noindent et 

\begin{eqnarray}
\mathbb{P}_{X_{n}}(\bigcup_{j=1}^{m}A_{j})&=&\sum \mathbb{P}_{X_{n}}(A_{j})-\sum \mathbb{P}_{X_{n}}(A_{i}A_{j})+\sum \mathbb{P}
_{X_{n}}(A_{i}A_{j}A_{k}) \label{FP2}\\
&+&...+(-1)^{n+1}\mathbb{P}_{X_{n}}(A_{1}A_{2}...A_{n})  \notag
\end{eqnarray}

\noindent Maintenant traitons chacun des termes de ces expressions. Prenons un terme
quelconque
\begin{equation*}
\mathbb{P}_{X}(A_{i_{1}}A_{i_{2}}...A_{i_{k}}).
\end{equation*}

\noindent Puisque la classe \ $\mathcal{U}$ des intervalles  $F_{X}$-continus est
stable par intersection finie (Sous-section \ref{subsubsecUstable} de la section \ref{cv.annexe}), cette intersection $%
A_{i_{1}}A_{i_{2}}...A_{i_{k}}$ est de type $]a,b]\in \mathcal{U}$ et la
formule de la mesure de Lebesgues-Stieljes donne 

\begin{equation*}
\mathbb{P}_{X}(A_{i_{1}}A_{i_{2}}...A_{i_{k}})=\sum\limits_{\varepsilon \in
\{0,1\}^{k}}(-1)^{\sum_{1\leq i\leq k}\varepsilon _{i}}F_{X}(b+\varepsilon
\ast (a-b))
\end{equation*}

\noindent et

\begin{equation*}
\mathbb{P}_{X_{n}}(A_{i_{1}}A_{i_{2}}...A_{i_{k}})=\sum\limits_{\varepsilon
\in \{0,1\}^{k}}(-1)^{\sum_{1\leq i\leq k}\varepsilon
_{i}}F_{X_{n}}(b+\varepsilon \ast (a-b)).
\end{equation*}%

\noindent Ici les points $c=b+\varepsilon \ast (a-b)$ sont tous des points de continuit%
\'{e} puisque $]a,b]\in \mathcal{U}$. Il suit de cela et de la convergence (%
\ref{CFR}) que

\begin{equation*}
\mathbb{P}_{X_{n}}(A_{i_{1}}A_{i_{2}}...A_{i_{k}})\rightarrow \mathbb{P}%
_{X}(A_{i_{1}}A_{i_{2}}...A_{i_{k}}).
\end{equation*}%

\noindent En op\'{e}rant par terme dans (\ref{FP1}) et (\ref{FP2}), nous concluons que

\begin{equation*}
\mathbb{P}_{X_{n}}(\bigcup_{j=1}^{m}A_{j})\rightarrow \mathbb{P}_{X}(\bigcup_{j=1}^{m}A_{j}),
\end{equation*}

\noindent et par la suite,

\begin{eqnarray*}
\liminf_{n\rightarrow +\infty} \mathbb{P}_n(X_{n}\in G)&=&\liminf_{n\rightarrow +\infty} \mathbb{P}_n(X_{n}\in
\bigcup_{j\geq 1}]a^{j},b^{j}])\\
&&\geq \lim_{n\rightarrow +\infty}\mathbb{P}_n(X_{n}\in \bigcup_{j=1}^{m}]a^{j},b^{j}])\geq \mathbb{P}_n(X\in G)-\eta 
\end{eqnarray*}

\noindent pour tout $\eta >0$. D'o\`{u}
 
\begin{equation*}
\liminf_{n\rightarrow +\infty} \mathbb{P}_n(X_{n}\in G)\geq \mathbb{P}_{\infty}(X\in G).
\end{equation*}%

\noindent pour tout $G$ ouvert. Il s'en suit que 
\begin{equation*}
X_{n}\rightarrow _{d}X \text{ quand } n \rightarrow +\infty.
\end{equation*}

\noindent Passons aux fonctions caract\'eristiques. Nous avons :

\begin{proposition} \label{cv.FC}
Une suite d'applications $X_{n}:\ (\Omega _{n }, \mathcal{A}_{n },\mathbb{P}_{n})\mapsto (\mathbb{R}^{k},\mathcal{B}(\mathbb{R}^{k}))$\ mesurables converge vaguement vers la variable al\'{e}atoire $X:\ (\Omega_{\infty},\mathcal{A}_{\infty},\mathbb{P}_{\infty})\mapsto (\mathbb{R}^{k},\mathcal{B}(\mathbb{R}^{k}))$, ssi pour tout point $%
^{t}(u_{1},u_{2},...,u_{k})\in R^{k},$%
\begin{equation*}
\Phi _{X_{n }}(u_{1},u_{2},...,u_{k})\mapsto \Phi _{X}(u_{1},u_{2},...,u_{k}).
\end{equation*}
\end{proposition}

\noindent \textbf{Remarque}. La preuve que nous pr\'{e}sentons ici, est bas\'{e}e sur le th\'{e}or\`{e}me de  Stone-Weirstrass, qui est un \'{e}l\'{e}ment important de la topologie des espaces de fonctions continues d\'{e}finies sur un compact. Ce th\'{e}or\`{e}me est rappel\'{e}
dans la proposition \ref{cv.subsec.annexe.SW} de la section \ref{cv.annexe}, voir ci-bas. Cepandant, une autre preuve beaucoup plus jolie 
\`{a} nos yeux, est donn\'{e}e dans le the th\'{e}or\`{e}me \ref{cv.tension.ConvFC} du chapitre \ref{cv.tensRk}. Ce dernier est bas\'{e}
sur le concept de tension et le th\'{e}or\`{e}me de continuit\'{e} de L\'{e}vy.\\

\noindent \textbf{Proof}. Rappelons la d\'{e}finition de la fonction charact%
\'{e}ristique ainsi qu'il suit 
\begin{equation*}
^{t}{\large (u}_{1}{\large ,u}_{2}{\large ,...,u}_{k}{\large )\mapsto \Phi }%
_{X}{\large (u}_{1}{\large ,u}_{2}{\large ,...,u}_{k}{\large )=\mathbb{E}(}%
\exp {\large (}\sum_{j}^{k}{\large i}\text{ u}_{j}{\large X}_{j}{\large )),}
\end{equation*}

\bigskip \noindent qui peut \^{e}tre r\'{e}-\'{e}crite ainsi 
\begin{equation*}
^{t}(u_{1},u_{2},...,u_{k})\mapsto \exp (\sum_{j}^{k}i\text{ u}%
_{j}X_{j})=\cos (\sum_{j}^{k}i\text{ u}_{j}X_{j})+i\sin (\text{ }\sum_{j}^{k}%
\text{u}_{j}X_{j}).
\end{equation*}

\bigskip \noindent Cette fonction complexe a des composantes qui sont des
fonctions born\'{e}es et continues de  $X$ et par d\'{e}finition, 
\begin{equation*}
\mathbb{E}\exp (\sum_{j}^{k}i\text{ u}_{j}X_{j})=\mathbb{E}\cos
(\sum_{j}^{k}i\text{ u}_{j}X_{j})+i\text{ }\mathbb{E}\sin (\text{ }%
\sum_{j}^{k}\text{u}_{j}X_{j}).
\end{equation*}

\noindent Donc par la simple d\'{e}finition de la convergence vague, nous
avons en tout $^{t}(u_{1},u_{2},...,u_{k})\in \mathbb{R}^{k},$%
\begin{equation}
\Phi _{X_{n}}(u_{1},u_{2},...,u_{k})\mapsto \Phi _{X}(u_{1},u_{2},...,u_{k}). \label{cv33}
\end{equation}

\bigskip \noindent Ceci \'{e}tablit le sens direct de la preuve. Pour
prouver le sens indirect one, nous en appelons au th\'{e}or\`{e}me de
Stone-Weirstrass (Voir Proposition \ref{cv.subsec.annexe.SW}
de la sous-section \ref{cv.subsec.annexe.SW} de la section \ref{cv.annexe},
voir ci-bas).\newline

\noindent A cet effet, posons  $\Delta (a)=[-a,a]^{k}$, for $0<a\in \mathbb{R}$.
\noindent Notons par  $H$ la classe des fonctions continue born\'{e}es sur $%
\mathbb{R}^{k}$ dont les \'{e}l\'{e}ments sont de la forme 
\begin{equation}
^{t}(x_{1},x_{2},...,x_{k})\mapsto \sum_{r=1}^{m}a_{r}\exp (\sum_{j}^{k}\pi
n_{j,r}\text{ }i\text{ }x_{j}/a),  \label{cv.DefHFR}
\end{equation}

\noindent o\`{u} les coefficients $a_{r}$ sont des nombres complexes et les $n_{j,r}$
sont des entiers. En d'autres termes, les \'{e}l\'{e}ments de $H$ sont des
combinaisons lin\'{e}aires de fonctions exponentielles de combinaisons lin%
\'{e}aires de  $x_{1},x_{2,}...,x_{k}$. \ Remarquons que chaque fonction  
\begin{equation*}
x_{j}\mapsto \exp (\pi n_{j,r}\text{ }i\text{ }x_{j}/a)
\end{equation*}

\noindent est p\'{e}riodique de p\'{e}riode $2a,$ si bien que $h\in H$ attenit ses
valeurs dans $\Delta (a$, et alors $\{h(x),x\in \mathbb{R}^{k}\}\subset
\{h(x),x\in \Delta (a)\}$ et donc 
\begin{equation*}
\left\Vert h\right\Vert =\left\Vert h\right\Vert _{\Delta (a)}.
\end{equation*}

\noindent Pour tout $h\in H,$ notons $h_{\Delta (b)}$ la restriction de $h$ sur $%
\Delta (b),$ pour  $b>0,$ et soit  $H_{b}=\{h_{\Delta (b)},h\in H\}.$ V\'{e}%
rifions que pour tout $0<b<a$, les conditions du th\'{e}or\`{e}me de Stone-Weirstrass sont v\'{e}rifi\'{e}es pour  $H_{b}$.\\

\noindent (1) $H$ contient les constantes. Soit $d$ un nombre complexe quelconque.
Dans (\ref{cv.DefHFR}), prenons $m=1,a_{1}=d$ et $n_{1,1}=n_{2,1}=...=n_{k,1}=0.$ Nous obtenons $h(x)=a_{1}=d$.\\

\noindent (2) $H$ est stable par somme et produits finis. Cela est \'{e}vident.\\

\noindent (3) Les conjugu\'{e}s $\overline{h}$ \ des \'{e}l\'{e}ments $h$ de $H$
restent dans $H$.\\ 

\noindent (4) $H$ s\'{e}pare les points de $\mathbb{R}^{k}$. \ Pour le voir, prenons $%
^{t}(z_{1},z_{2},...,z_{k})\neq ^{t}(y_{1},y_{2},...,y_{k})$ dans $\Delta
(b),$. Donc au moins pour un indice  $j_{0}\in \{1,...,k\},$ nous avons $%
z_{j_{0}}\neq y_{j_{0}}$. Mais $\left\vert z_{j_{0}}-y_{j_{0}}\right\vert
\leq 2b<2a.$ D\'{e}finissons $h_{0}\in H$ par 
\begin{equation*}
h_{0}(x)=e^{i\pi x_{j_{0}}/a}.
\end{equation*}%
L'\'{e}galit\'{e}  $e^{i\pi z_{j_{0}}/a}=e^{i\pi y_{j_{0}}/a}$ aurait entra%
\^{\i}n\'{e} qu'il existe un entier $k\in \mathbb{Z}$ tel que $%
(z_{j_{0}}-y_{j_{0}})=2ka,$ ce qui est impossible. Donc, nous avons 
\begin{equation*}
h_{0}(z)=e^{i\pi z_{j_{0}}/a}\neq e^{i\pi y_{j_{0}}/a}=h_{0}(y).
\end{equation*}

\noindent Ce qui implique que $h_{0}$ s\'{e}pare $^{t}(z_{1},z_{2},...,z_{k})$ et $%
^{t}(y_{1},y_{2},...,y_{k})$.\\

\bigskip \noindent Les conditions du th\'eor\`eme de Stone-Weirstrass sont v\'erifi\'es. Maintenant, supposons que (\ref{cv33}) est vraie. Soit  $f\in C_{b}(\mathbb{R}^{k},\mathbb{R})\subseteq C_{b}(\mathbb{R}^{k},\mathbb{C})$. Soit  $f_{b}$
la restriction de $f$ sur  $\Delta (b)$ si bien que  $f_{b}\in C(\Delta (b),%
\mathbb{C})$ o\`{u} $b=b(a)<a,$ est une fonction croissante non born\'{e}e
de  $a$ ($b(a)=a/2$ for instance)$.$ A cet \'{e}tape, nous appliquons le th%
\'{e}or\`{e}me de Stone-Weirstrass pour trouver, pour tout  $\varepsilon >0,$
un \'{e}l\'{e}ment $h$ de $H$, 
\begin{equation*}
h(x)=\sum_{r=1}^{m}a_{r}\exp (\sum_{j}^{k}i\text{ }u_{r,j}x_{j})
\end{equation*}

\noindent tel que 
\begin{equation*}
\underset{x\in \Delta (b)}{\sup }\left\vert f(x)-h(x)\right\vert =\left\Vert
f-h\right\Vert _{\Delta (b)}\leq \varepsilon /3.
\end{equation*}

\noindent En appliquant (\ref{cv33}), nous avons 
\begin{equation*}
\mathbb{E}(h(X_{n}))\rightarrow \mathbb{E}(h(X))\text{ as }n\rightarrow
+\infty .
\end{equation*}

\noindent Soit donc un entier positif $n_{0}$ tel que pour tout $n\geq n_{0}$,
 
\begin{equation}
\normalsize
\left\vert \mathbb{E}(h(X_{n}))\rightarrow \mathbb{E}(h(X))\right\vert
=\left\vert \int h\text{ }d\mathbb{P}_{n}\circ X_{n}^{-1}-\int h\text{ }d%
\mathbb{P}\circ X^{-1}\right\vert \leq \varepsilon /3  \label{cv34g}
\end{equation}

\bigskip \noindent Nous avons
\begin{eqnarray*}
\mathbb{E}(f(X_{n}))-\mathbb{E}(f(X)) &=&(\int f\text{ }d\mathbb{P}_{n}\circ
X_{n}^{-1}-\int h\text{ }d\mathbb{P}_{n}\circ X_{n}^{-1})  \notag \\
&&+(\int h\text{ }d\mathbb{P}_{n}\circ X_{n}^{-1}-\int_{^{c}}h\text{ }d%
\mathbb{P}\circ X^{-1}) \notag \\
&&+(\int h\text{ }d\mathbb{P}\circ X^{-1}-\int_{^{c}}f\text{ }d\mathbb{P}%
\circ X^{-1}). \notag
\end{eqnarray*}

\bigskip \noindent Le premier terme satisfait 

\begin{eqnarray}
\mathbb{E}\left\vert \int f\text{ }d\mathbb{P}_{n}\circ X_{n}^{-1}-\int h%
\text{ }d\mathbb{P}_{n}\circ X_{n}^{-1}\right\vert  &\leq &\int_{\Delta
(b)}\left\vert f-h\right\vert \text{ }d\mathbb{P}_{n}\circ X_{n}^{-1} \notag \\
&&+\int_{\Delta ^{c}(b)}\left\vert f-h\right\vert \text{ }d\mathbb{P}%
_{n}\circ X_{n}^{-1} \notag \\
&\leq &\varepsilon /3+(\left\Vert f\right\Vert +\left\Vert h\right\Vert )%
\mathbb{P}_{n}(X_{n}\in \Delta ^{c}(b)). \label{cv34d}
\end{eqnarray}

\bigskip \noindent Par la m\^{e}me m\'{e}thode, le troisi\`{e}me terme aussi
v\'{e}rifie 
\begin{equation}
\normalsize
\mathbb{E}\left\vert \int f\text{ }d\mathbb{P}_{\infty }\circ X^{-1}-\int h%
\text{ }d\mathbb{P}_{\infty }\circ X^{-1}\right\vert \leq \varepsilon
/3+(\left\Vert f\right\Vert +\left\Vert h\right\Vert )\text{ }\mathbb{P}%
_{\infty }(X\in \Delta ^{c}(b)).  \label{cv34e}
\end{equation}

\bigskip \noindent En mettant ensemble les formules (\ref{cv34g}), (\ref{cv34d}) and (\ref{cv34e}), nous obtenons pour tout $n\geq n_{0}$\ 
\begin{equation*}
\left\vert \mathbb{E}(f(X_{n}))-\mathbb{E}(f(X))\right\vert \leq \varepsilon
+(\left\Vert f\right\Vert +\left\Vert h\right\Vert )(\mathbb{P}_{n}(X_{n}\in
\Delta ^{c}(b))+\mathbb{P}_{\infty }(X\in \Delta ^{c}(b))).
\end{equation*}

\noindent Alors, pour tout $n\geq n_{0}$ fix\'{e}$,$ nous faisons cro\^{\i}tre $a$ vers $+\infty .$ En cons\'{e}quence $b(a)\uparrow +\infty $ et par la
suite, $\mathbb{P}_{n}(X_{n}\in \Delta ^{c}(b))+\mathbb{P}_{\infty}(X\in
\Delta ^{c}(b))\downarrow 0.$ Alors pour tout $n\geq n_{0}$

\bigskip\ 
\begin{equation*}
\left\vert \mathbb{E}(f(X_{n}))-\mathbb{E}(f(X))\right\vert \leq \varepsilon 
\end{equation*}

\noindent Ce qui prouve que  
\begin{equation*}
\mathbb{E}(f(X_{n}))\rightarrow \mathbb{E}(f(X))\text{ as }n\rightarrow
+\infty 
\end{equation*}

\noindent et donc

\begin{equation*}
X_{n}\rightsquigarrow X\text{ as }n\rightarrow +\infty .
\end{equation*}

\bigskip \bigskip

\noindent En mettant ensemble les propositions (\ref{cv.FRDIR}), (\ref{cv.FRINV}) et (\ref{cv.FC}), nous obtenons le théorème complet Pormanteau dans $\mathbb{R}^k$. 

\begin{theorem} \label{cv.theo.portmanteau.rk} Let $k$ be a positive integer. The sequence of random vectors $X_{n}:(\Omega _{n },\mathcal{A}_{n },
\mathbb{P}_{n })\mapsto (\mathbb{R}^k,\mathcal{B}(\mathbb{R}^k))$, $\geq 1$, weakly converges to the random vector $X:(\Omega _{\infty},\mathcal{A}_{\infty},\mathbb{P}_{\infty})\mapsto (\mathbb{R}^k,\mathcal{B}(\mathbb{R}^k))$ if and only if one of these assertions holds.\\

\bigskip \noindent (i) For any real-valued continuous and bounded function $f$ defined on $\mathbb{R}^k)$,
\begin{equation*}
\lim_{n\rightarrow +\infty} \mathbb{E}f(X_n) =\mathbb{E}f(X).
\end{equation*}

\bigskip \noindent (ii) For any open set $G$ in $\mathbb{R}^k$, 
\begin{equation*}
\liminf_{n\rightarrow +\infty} \mathbb{P}_n(X_{n}\in G)\geq \mathbb{P}_{\infty}(X \in G).
\end{equation*}

\bigskip \noindent (iii) For any closed set $G$ of  $S$, we have
\begin{equation*}
\limsup_{n\rightarrow +\infty} \mathbb{P_{n}}(X_{n}\in F)\leq \mathbb{P_{\infty}}(X \in F).
\end{equation*}

\bigskip \noindent (iv) For any inferior semi-continuous and bounded below function $f$, we have
\begin{equation*}
\liminf_{n\rightarrow +\infty} \mathbb{E}f(X_n)\geq \mathbb{E}f(X).
\end{equation*}

\bigskip \noindent (v) For any superior semi-continuous and bounded above function $f$, we have
\begin{equation*}
\limsup_{n\rightarrow +\infty} \mathbb{E}f(X_n)\leq \mathbb{E}f(X).
\end{equation*}

\bigskip \noindent (vi) For any Borel set $B$ of $S$ that is $\mathbb{P}_{X}$-continuous, that is $\mathbb{P}_{\infty}(X \in \partial B)=0$, we have
\begin{equation*}
\lim_{n\rightarrow +\infty} \mathbb{P}_n(X_{n }\in B)=\mathbb{P}_{X}(B)=\mathbb{P}_{\infty}(X \in B).
\end{equation*}

\bigskip \noindent (vii) For any nonnegative and bounded Lipschitz function $f$, we have 
\begin{equation*}
\liminf_{n\rightarrow +\infty} \mathbb{E}f(X_n)\geq \mathbb{E}f(X).\\
\end{equation*}

\bigskip \noindent (viii) For any continuity point $t=(t_{1},t_{2},...,t_{k})$ of $F_{X}$, we have,

\begin{equation*}
F_{X_n}(t) \rightarrow F_{X}(t) \text{ as } n\rightarrow +\infty.
\end{equation*}

\noindent where for each $n\geq 1$, $F_{X_n}$ is the distribution function of $X_n$ and $F_{X}$ that of $X$.\\

\noindent (ix) For any point $u=(u_{1},u_{2},...,u_{k})\in \mathbb{R}^{k}$, 
\begin{equation*}
\Phi _{X_{n}}(u)\mapsto \Phi _{X}(u) \text{ as }n\rightarrow +\infty .
\end{equation*}

\noindent where for each $n\geq 1$, $\Phi_{X_n}$ is the characteristic function of $X_n$ and $\Phi_{X}$ is that of $X$
\end{theorem}

\noindent \textbf{Cri\`ere de Wold}. La suite de vecteurs aléatoires $\{X_n, \ \ n\geq 1\} \subset \mathbb{R}^k$ converge vaguement
vers to $X \in \mathbb{R}^k$, quand $n \rightarrow +\infty$ si et seulement si pour $a \in \mathbb{R}^k$, la suite $\{<a,X_n>, \ \ n\geq 1\} \subset \mathbb{R}$ converge veguement vers $X \in \mathbb{R}$ quand $n \rightarrow +\infty$.

\noindent \textbf{Preuve}. La preuve est rapide. Elle utilise les notations ant\'erieures. Supposons que $X_n$ converge faiblement vers $X$ 
dans $\mathbb{R}^k$ quand $n \rightarrow +\infty$. En utilisant la convergence des fonctions caratéristiques, nous avons pour tout $u\in \mathbb{R}^k$,
$$
\mathbb{E}(exp(i<X_n,u>) \rightarrow \mathbb{E}(exp(i<X,u>) \ \ quand \ \ n \rightarrow +\infty.
$$  

Il s'en suit que $a\in \mathbb{R}^k$ et pour tout $t \in \mathbb{R}$, nous avons

\begin{equation}
\mathbb{E}(exp(it<X_n,a>) \rightarrow \mathbb{E}(exp(it<X,a>) \ \ quand \ \ n \rightarrow +\infty. \label{cv.proj}
\end{equation}

\noindent c-a-d que, en prenant $u=ta$ dans la formule précédente, et en notant $Z_n=<X_n,a>$ et $Z=<X,a>$, 

$$
\mathbb{E}(exp(itZ_n) \rightarrow \mathbb{E}(exp(itZ) \ \ quand \ \ n \rightarrow +\infty.
$$

\noindent Cela signifie que $Z_n \rightsquigarrow Z$, qui est égal à $<a,X_n>$, converges vaguement vers $<a,X>$.\\

\noindent Inversement, supposons que pour tout $a \in \mathbb{R}^k$, la suite $\{<a,X_n>, \ \ n\geq 1\} \subset \mathbb{R}$ converge vaguement vers $X \in \mathbb{R}$ as $n \rightarrow +\infty$. Alors, en prenant $t=1$ dans (\ref{cv.proj}), nous obtenons pour tout $a=u \in \mathbb{R}^k$,

$$
\mathbb{E}(exp(i<X,u>) \rightarrow \mathbb{E}(exp(i<X,u>) \ \ quand \ \ n \rightarrow +\infty.
$$

\noindent ce qui signifie que $X_n \rightsquigarrow +\infty$ quand $n\rightarrow +\infty$.\\

\section{Th\'{e}or\`{e}me de Scheff\'{e}}

Dans la section pr\'{e}c\'{e}dente, nous avons li\'{e} la convergence vague
et quelques caract\'{e}ristiques de variables al\'{e}atoires dans $\mathbb{R}%
^{k}$\ telles que la fonction de r\'{e}partition et la fonction caract\'{e}%
ristique. On peut se demander ce qu'il en est par rapport aux densit\'{e}s
de probabilit\'{e}s par rapport \`{a} la mesure de Lebesgues dans $\mathbb{R}%
^{k}.$\ Le th\'{e}or\`{e}me de Scheff\'{e}(1947) r\'{e}pond \`{a} cette pr%
\'{e}occupation dans le cadre g\'{e}n\'{e}ral. Enon\c{c}ons-le d'abord.
\newline

\begin{theorem} \label{cv.scheffe} (Th\'eor\`eme de Scheff\'e) Soit $\lambda $\ une mesure sur un espace mesurable $(E,B)$.
Et soit p, $(p_{n})_{n\geq 1}$\ des densit\'{e}s de probabilit\'{e}s par
rapport \`{a} $\lambda $, c'est-\`{a}-dire des applications num\'{e}riques d%
\'{e}finies sur E, postitives, mesurables telles que 
\begin{equation}
\forall n\geq 1,\int p_{n}\text{ }d\lambda =\int p\text{ }d\lambda =1.
\label{scheffe1}
\end{equation}

\bigskip \noindent Si 
\begin{equation*}
p_{n}\rightarrow p,\text{ }\lambda -pp
\end{equation*}

\bigskip \noindent alors 
\begin{equation}
\sup_{B\in \mathcal{B}}\left\vert \int_{B}p_{n}\text{ }d\lambda -\int_{B}p%
\text{ }d\lambda \right\vert =\frac{1}{2}\int \left\vert p_{n}-p\right\vert 
\text{ }d\lambda \rightarrow 0  \label{scheffe2}
\end{equation}
\end{theorem}

\noindent \textbf{PREUVE}. Supposons que $p_{n}\rightarrow p,$\ $\lambda -pp.$\ Posons $\Delta
_{n}=p-p_{n}.$\ Alors, (\ref{scheffe1}) implique 
\begin{equation*}
\int \Delta _{n}\text{ }d\lambda =0.
\end{equation*}

\noindent Donc, pour $B\in \mathcal{B}$,
\begin{equation*}
\int_{B^{c}}\Delta _{n}\text{ }d\lambda =\int \Delta _{n}\text{ }d\lambda
-\int_{B}\Delta _{n}\text{ }d\lambda =-\int_{B}\Delta _{n}\text{ }d\lambda .
\end{equation*}

\noindent D'o\`{u}
 
\begin{equation}
2\left| \int_{B}\Delta _{n}\text{ }d\lambda \right| =\left| \int_{B}\Delta
_{n}\text{ }d\lambda \right| +\left| \int_{B^{c}}\Delta _{n}\text{ }d\lambda
\right| \leq \int_{B}\left| \Delta _{n}\right| \text{ }d\lambda
+\int_{B^{c}}\left| \Delta _{n}\right| \text{ }d\lambda \leq \int \left|
\Delta _{n}\right| \text{ }d\lambda ,  \label{scheffe3}
\end{equation}

\noindent c'est-\`{a}-dire 
\begin{equation}
\left| \int_{B}\Delta _{n}\text{ }d\lambda \right| \leq \frac{1}{2}\int
\left| \Delta _{n}\right| \text{ }d\lambda .  \label{scheffe4}
\end{equation}

\bigskip \noindent En prenant $B=(\Delta _{n}\geq 0)$\ dans (\ref{scheffe3}%
), nous avons 
\begin{equation*}
2\left| \int_{B}\Delta _{n}\text{ }d\lambda \right| =\left| \int_{B}\Delta
_{n}^{+}\text{ }d\lambda \right| +\left| \int_{B^{c}}-\Delta _{n}^{-}\text{ }%
d\lambda \right| =\int \Delta _{n}^{+}d\lambda +\int \Delta _{n}^{-}d\lambda
=\int \left| \Delta _{n}\right| d\lambda .
\end{equation*}

\bigskip \noindent En mettant ensemble les deux derni\`{e}res formules, nous
avons 
\begin{equation}
\sup_{B\in \mathcal{B}}\left| \int_{B}p_{n}\text{ }d\lambda -\int_{B}p\text{ 
}d\lambda \right| =\frac{1}{2}\int \left| p_{n}-p\right| \text{ }d\lambda .
\label{scheffe6}
\end{equation}

\bigskip \noindent Maintenant, 
\begin{equation*}
0\leq \Delta _{n}^{+}=\max (0,p-p_{n})\leq p.
\end{equation*}

\bigskip \noindent De plus, 
\begin{equation*}
\int \Delta _{n}^{+}\text{ }d\lambda =\int_{(\Delta _{n}\geq 0)}\Delta _{n}%
\text{ }d\lambda =\int \Delta _{n}\text{ }d\lambda -\int_{(\Delta _{n}\leq
0)}\Delta _{n}\text{ }d\lambda =\int_{(\Delta _{n}\leq 0)}-\Delta _{n}\text{ 
}d\lambda =\int \Delta _{n}^{-}\text{ }d\lambda ,
\end{equation*}

\bigskip \noindent de sorte que 
\begin{equation}
\int \left| \Delta _{n}\right| \text{ }d\lambda =2\int \Delta _{n}^{+}\text{ 
}d\lambda  \label{scheffe7}
\end{equation}

\bigskip \noindent Appliquons le th\'{e}or\`{e}me de convergence domin\'{e}e
de Lebesgues \`{a} $0\leq \Delta _{n}^{+}\leq \left| \Delta _{n}\right|
\rightarrow 0$\ $\lambda -pp$, $0\leq \Delta _{n}^{+}\leq p.$\ Nous aurons 
\begin{equation*}
\int \Delta _{n}^{+}\text{ }d\lambda \rightarrow 0,
\end{equation*}

\bigskip \noindent en vertu de (\ref{scheffe6}), 
\begin{equation*}
\sup_{B\in \mathcal{B}}\left| \int_{B}p_{n}\text{ }d\lambda -\int_{B}p\text{ 
}d\lambda \right| =\frac{1}{2}\int \left| p_{n}-p\right| \text{ }d\lambda
=\int \Delta _{n}^{+}\text{ }d\lambda \rightarrow 0
\end{equation*}

\bigskip \noindent Le th\'{e}or\`{e}me de Scheff\'{e} peut alors s'appliquer
aux densit\'{e}s de probabilit\'{e} dans $\mathbb{R}^{k}$\ en particulier.
Nous aurons :

\begin{proposition} \label{cv.scheffeExt}
\noindent (A) Soit $X_{n}: (\Omega _{n },\mathcal{A}_{n},\mathbb{P}_{n}) \mapsto (\mathbb{R}^{k},\mathcal{B}(\mathbb{R}^{k}))$ une suite de vecteurs al\'eatoires et $X : (\Omega_{\infty},\mathcal{A}_{\infty},\mathbb{P}_{\infty})\mapsto (\mathbb{R}^{k},\mathcal{B}(\mathbb{R}^{k}))$ un autre vecteur al\'eatoire, tous de lois absol\^ument continues par rapport \`a la mesure de Lebesgues sur  $\mathbb{R}^{k}$ not\'ee $\lambda_k$. Soit, pour tout $n\geq 1$,  $f_{X_{n}}$ la densit\'e de probabilit\'e de $X_n$, et $f_{X}$ celle de $X$. Supposons que : 

\begin{equation*}
f_{X_{n}}  \rightarrow f_{X}, \text{ }\lambda_k-p.s., \text{ as } n\rightarrow +\infty.
\end{equation*}

\noindent Alors $X_{n}$ converge vaguement vers $X$ quand $n \rightarrow +\infty$.\\

\noindent (B) Soit $X_{n }: (\Omega _{n },\mathcal{A}_{n},\mathbb{P}_{n })\mapsto (\mathbb{R}^{k},\mathcal{B}(\mathbb{R}^{k}))$ une suite vecteurs al\'eroires à valeurs discrètes $X : (\Omega_{\infty},\mathcal{A}_{\infty},\mathbb{P}_{\infty})\mapsto (\mathbb{R}^{k},\mathcal{B}(\mathbb{R}^{k}))$ un autre vecteur al\'eatoires discret. Notons pour chaque $n$, le support d\'enombrable de $X_n$ par $D_n$, c'est-\`a-dire : 
$$
\mathbb{P}_n(X_n \in D_n)=1 \text{ et pour chaque  } x \in D_n, \text{ } \mathbb{P}_n(X_n=x)\neq 0
$$

\noindent et par $D_{\infty}$ le support d\'enombrable de $X$. Soit  $D=D_{\infty} \cap \cap_{n \geq 1} D_n$ et soit $\nu$ la mesurede comptage sur $D$. Alors les densit\'es de probabibilit\'e des vecteurs $X_n$ et de $X$ par rapport \`a $\nu$ sont d\'efinies sur $D$ par
$$
f_{X_n}(x)=\mathbb{P}_n(X_n=x), \text{ } n\geq 1, \text{ } f_{X}(x)=\mathbb{P}_{\infty}(X=x), \text{ }x\in D.
$$

\noindent De plus, si 
\begin{equation*}
(\forall x\in D), f_{X_{n }}(x))  \rightarrow f_{X}(x).
\end{equation*}
\noindent Alors $X_{n}$ converge vaguement vers $X$.
\end{proposition}

\section{Convergence vague et convergence en probabilit\'e sur le m\^eme espace de probabili\'e}\label{cv.CvCp}

\noindent Cette section met la convergence vague \`{a} sa place dans le sch%
\'{e}ma de l'\'{e}tude de convergence de suites de variables al\'{e}atoires d%
\'{e}finies sur le m\^{e}me espace de probabilit\'{e} $(\Omega ,\mathcal{A},\mathbb{P})$ et \`{a} valeur dans un espace m\'{e}trique $(S,d).$\newline

\noindent Nous avons d\'{e}j\`{a} vu que la convergence vague n'exige pas la
connaissance des espaces de d\'{e}part. Lorsque les suites \'{e}tudi\'{e}es
sont sur le m\^{e}me espace, nous avons alors des relations int\'{e}%
ressantes avec les autres types de convergence.\newline

\bigskip \noindent A l'inverse, il existe des r\'{e}sultats regroup\'{e}s
sous le manteau de Skorohod-Wichura-Dudley qui permettent de donner une
version de convergence vague en convergence presque-s\^{u}re dans un espace
convenablement choisi. La d\'{e}monstration ne sera pas abord\'ee dans cet ouvrage. Elle sera donn\'ee dans l'ouvrage avanc\'e de 
convergence vague. Ici, Il sera expos\'{e} et illustr\'{e} seulement
pour le cas de l'espace $\mathbb{R}$ dans le chapitre \ref{cv.R}.\\

\noindent Commen\c{c}ons par les d\'{e}finitions.

\subsection{D\'efinitions}

Dans tout ce document,\`a l'exception de la section sur le th\'eor\`eme de
Skorohod-Wichura, les suites de variables al\'{e}atoires $(X_{n})_{n\geq 0},$
$(Y_{n})_{n\geq 0}$, etc., et les variables al\'{e}atoires $X,$ $Y,$ $etc.$
sont d\'{e}finies sur le m\^{e}me espace de probabilit\'{e} $(\Omega,\mathcal{A},\mathbb{P})$ \ et sont \`{a} valeurs dans un espace m\'{e}troque 
$(S,d)$. Nous aurons aussi \`{a} utiliser des constantes $c$ dans $S$.\\

\bigskip \noindent \textbf{(a) Convergence presque-s\^ure}.\\

\noindent La suite $(X_{n})_{n\geq 0}$ converge presque-surement vers $X$,
not\'{e} $X_{n}\longrightarrow X,$ $a.s.,$ si et seulement l'ensemble sur
lequel \ $(X_{n})_{n\geq 0}$\ ne converge pas vers $X$\ est n\'{e}gligeable,
c'est-\`{a}-dire que%
\begin{equation*}
\mathbb{P}(\{\omega \in \Omega ,X_{n}\nrightarrow X\})=\mathbb{P}(\{\omega
\in \Omega ,d(X_{n},X)\nrightarrow 0\})=0.
\end{equation*}

\noindent Nous caract\'{e}risons cet ensemble ais\'{e}ment par 
\begin{equation*}
(X_{n}\nrightarrow X)=\bigcup\limits_{k\geq 1}\bigcap\limits_{n\geq
0}\bigcup\limits_{p\geq n}(d(X_{p},X)>k^{-1}).
\end{equation*}

\noindent Ceci m\`{e}ne \`{a} la definition suivante : $(X_{n})_{n\geq 0}$
converge presque-s\^urement vers $X$ si et seulement si : 
\begin{equation}
\forall k\geq 1,\text{ }\mathbb{P}\left( \bigcap\limits_{n\geq
0}\bigcup\limits_{p\geq n}(d(X_{p},X)>k^{-1})\right) =0.  \label{cps}
\end{equation}

\bigskip \noindent \textbf{(b) Convergencen probabilit\'{e}}.\\

\noindent La suite $(X_{n})_{n\geq 0}$ converge en probabilit\'{e} vers $X,$ not\'{e} $%
X_{n}\longrightarrow _{\mathbb{\mathbb{P}}}X,$ $a.s.$\ si et seulement si 
\begin{equation*}
\forall \varepsilon >0,\text{ }\lim_{n\rightarrow +\infty} \mathbb{P}%
\left( d(X_{n},Y)>\varepsilon \right) =0.
\end{equation*}

\noindent Nous allons faire une br\`{e}ve comparaison entre ces deux types de
convergence que vous connaissez d\'{e}j\`{a} :

\bigskip

\begin{proposition} \label{cv.CvCp.01}
Si $X_{n}\longrightarrow X,$ $a.s.,$ alors $X_{n}\longrightarrow _{\mathbb{P}%
}X$.
\end{proposition}

\noindent \textbf{Preuve}. La preuve est classique. Il suffit de supposer
qu'il y a convergence presque-sure, c'est-\`{a}-dire que (\ref{cps}) a lieu
et de voir que pour tout $k\geq 1,$%
\begin{equation*}
(d(X_{n},X)>k^{-1})\subset \bigcup\limits_{p\geq
n}(d(X_{p},X)>k^{-1})=:B_{n,k}.
\end{equation*}

\noindent Mais la suite $B_{n,k}$ (en $n$) est d\'{e}croissante vers%
\begin{equation*}
\bigcap\limits_{n\geq 0}\bigcup\limits_{p\geq n}(d(X_{p},X)>k^{-1})=:B_{k}.
\end{equation*}

\noindent Donc pour tout $n\geq 0$ et pour tout $k\geq 1$, 
\begin{equation*}
\mathbb{P}(d(X_{n},X)>k^{-1})\leq \mathbb{P}\left( B_{n,k}\right)
\end{equation*}

\noindent et par suite, par la continuit\'{e} de la probabilit\'{e},

\begin{equation*}
\limsup_{n\rightarrow +\infty}\mathbb{P}(d(X_{n},X)>k^{-1})\leq
\lim_{n\rightarrow +\infty}\mathbb{P}\left( B_{n,k}\right) =\mathbb{P}\left(
B_{k}\right) =0,
\end{equation*}

\noindent en appliquant (\ref{cps}).\\

\subsection{Convergence vague et convergence en probabilit\'{e}}

Avant d'\'{e}noncer les propri\'{e}t\'{e}s, nous pouvons enrichir le Th\'{e}%
or\`eme Pormanteau par ce point suppl\'{e}mentaire.

\begin{lemma}
\label{cvcp.lem01} La suite $(X_{n})_{n\geq 0}$ converge vaguement vers $X$
si et seulement si : \newline

\noindent (viia) Pour toute fonction $f:S\longrightarrow \mathbb{R}$ born%
\'{e}e et Lipschitzienne,%
\begin{equation*}
\mathbb{E}f(X_{n})\rightarrow f(X).
\end{equation*}
\end{lemma}

\noindent \textbf{Preuve}. Nous nous mettons dans le cadre de la d\'{e}%
monstration du Th\'{e}or\`{e}me Partmanteau. Ce point $(viia)$ entra\^{\i}ne
le point $(i)$ puisque $f$ est continue born\'{e}e. De plus, $(vii)$ est un
cas particulier de $(viiia),$ ce qui fait que $(viia)\Longrightarrow (vii).$
En retour, si $(vii)$ a lieu, on peut consid\'{e}rer l'infimum $A$ de la
fonction $f$ et son supremum B et appliquer le point $(vii)$ \`{a} la
fonction $f-A$ puis \`{a} $-f+B,$ ce qui entra\^{\i}nera le point $(viia).$
On aura ainsi $(i)\Longleftrightarrow (vii)\Longleftrightarrow (viia)$. Ce
qui termine la preuve.\\

\noindent Nous allons maintenants \'enoncer un ensemble de propri\'et\'es.\\

\bigskip \noindent \textbf{(a) La convergence en probabilit\'{e} entra\^ine la
convergenve vague}

\begin{proposition}
\label{cvcp.prop2} Si\ $X_{n}\longrightarrow _{\mathbb{P}}X$ \ implique que $%
X_{n}\rightsquigarrow X$\newline
\end{proposition}

\noindent \textbf{Preuve}. Supposons que $X_{n}\longrightarrow _{\mathbb{P}%
}X.$ Montrons que $X_{n}\rightsquigarrow X$ en utilisant le point $(viia)$
du Lemma \ref{cvcp.lem01} ci-dessus. Soit donc une fonction $f$
lipschitzienne de rapport $\ell $ et born\'{e}e par $M.$ Nous avons pour
tout $n\geq 0,$%
\begin{equation*}
\left\vert f(X_{n})-f(X)\right\vert \leq \ell d(X_{n},X).
\end{equation*}

\noindent Nous avons pour tout $n\geq 0$ et pour $\varepsilon >0,$

\begin{eqnarray*}
\left\vert \mathbb{E}f(X_{n})-\mathbb{E}f(X)\right\vert &\leq &\mathbb{E}%
\left\vert f(X_{n})-f(X)\right\vert \\
&\leq &\int_{(d(X_{n},X)\leq \varepsilon )}\left\vert
f(X_{n})-f(X)\right\vert d\mathbb{P}\\
&+&\int_{(d(X_{n},X)>\varepsilon
)}\left\vert f(X_{n})-f(X)\right\vert d\mathbb{P}.
\end{eqnarray*}

\noindent Mais, pour tout $n\geq 0$ et pour $\varepsilon >0,$%
\begin{equation*}
\int_{(d(X_{n},X)\leq \varepsilon )}\ell d(X_{n},X)d\mathbb{P}\leq \ell
\varepsilon .
\end{equation*}

\noindent De plus, pour tout $n\geq 0$ et pour $\varepsilon >0,$%
\begin{equation*}
\int_{(d(X_{n},X)>\varepsilon )}\left\vert f(X_{n})-f(X)\right\vert d\mathbb{%
P}\leq \int_{(d(X_{n},X)>\varepsilon )}2Md\mathbb{P}\leq 2M\text{ }\mathbb{P}%
(d(X_{n},X)>\varepsilon ).
\end{equation*}

\noindent D'o\`{u}, pour tout $n\geq 0$ et pour $\varepsilon >0,$%
\begin{equation*}
\left\vert \mathbb{E}f(X_{n})-\mathbb{E}f(X)\right\vert \leq \ell
\varepsilon +2M\text{ }\mathbb{P}(d(X_{n},X)>\varepsilon ).
\end{equation*}

\noindent D'o\`{u} pour tout $\varepsilon >0,$%
\begin{equation*}
\limsup_{n\rightarrow +\infty}\left\vert \mathbb{E}f(X_{n})-\mathbb{E}%
f(X)\right\vert \leq \ell \varepsilon .
\end{equation*}

\noindent Nous obtenons, lorsque $\varepsilon \downarrow 0.$ 
\begin{equation*}
\mathbb{E}f(X_{n})\longrightarrow \mathbb{E}f(X).
\end{equation*}

\noindent Ce qui finit la preuve.\newline

\bigskip \noindent \textbf{(b) La convergence vague et la convergence en probabilit\'{e}
vers une constante sont \'{e}quivalentes}

\begin{proposition} \label{cv.CvECp}
\bigskip\ $X_{n}\longrightarrow _{\mathbb{P}}c$ \ si et seulement si $%
X_{n}\rightsquigarrow c$\newline
\end{proposition}

\noindent \textbf{Preuve}. Le sens ($X_{n} \rightarrow _{\mathbb{P}}c)
\Rightarrow \left( X_{n}\rightsquigarrow c\right)$. Prouvons maintenant que $%
\left( X_{n}\rightsquigarrow c\newline
\right) \Rightarrow (X_{n}\rightarrow _{\mathbb{P}}c).$ Pour cela, supposons
que$\left( X_{n}\rightsquigarrow c\right) $ et soit $\varepsilon >0.$ Le
point ($ii)$ du Th\'{e}or\`{e}me Portmanteau donne%
\begin{eqnarray*}
\liminf_{n\rightarrow +\infty}\mathbb{P}(d(X_{n},c) &<&\varepsilon )=\liminf_{n\rightarrow +\infty}\mathbb{P}(X_{n}\in B(c,\varepsilon ))\leq 
\mathbb{P}(c\in B(c,\varepsilon )) \\
&=&\mathbb{P}(d(c,c)>\varepsilon ) \\
&=&\mathbb{P}(\Omega )=1.
\end{eqnarray*}

\noindent D'o\`{u} 
\begin{equation*}
\limsup_{n\rightarrow +\infty}\mathbb{P}(d(X_{n},c)\geq \varepsilon
)=1-\liminf_{n\rightarrow +\infty}\mathbb{P}(d(X_{n},X)\leq \varepsilon
)\leq 1-1=0.
\end{equation*}

\noindent D'o\`{u}, pour tout $0<\varepsilon$,

\begin{equation*}
\limsup_{n\rightarrow +\infty}\mathbb{P}(d(X_{n},c)>\varepsilon )\leq
\limsup_{n\rightarrow +\infty}\mathbb{P}(d(X_{n},c)\geq \varepsilon )=0.
\end{equation*}

\bigskip \noindent \textbf{(c) Deux suites \'{e}quivalentes en probabilit\'{e} converge
vaguement vers la m\^{e}me limite, s'il y a lieu.}

\begin{proposition}
\label{cvcp.prop4} \ Si \ $X_{n}\rightsquigarrow X$ \ et $%
d(X_{n},Y_{n})\longrightarrow _{\mathbb{P}}0$, alors $Y_{n}\rightsquigarrow X
$\newline
\end{proposition}

\noindent \textbf{Preuve}. Supposons que $X_{n}\rightsquigarrow X$ \ et $%
d(X_{n},Y_{n})\longrightarrow _{\mathbb{P}}0.$ Montrons que $%
Y_{n}\rightsquigarrow X$ en utilisant le point $(viia)$ du Lemma \ref%
{cvcp.lem01} ci-dessus. Soit une fonction $f$ lipschitzienne de rapport $%
\ell $ et born\'{e}e par $M.$ Nous aurons$,$Nous avons pour tout $n\geq 0$
et pour $\varepsilon >0,$ \newline
\ 
\begin{eqnarray*}
\left\vert \mathbb{E}f(Y_{n})-\mathbb{E}f(X)\right\vert &\leq &\mathbb{E}%
\left\vert f(Y_{n})-f(X)\right\vert \\
&\leq &\mathbb{E}\left\vert f(X_{n})-f(X)\right\vert +\mathbb{E}\left\vert
f(Y_{n})-f(X_{n})\right\vert .
\end{eqnarray*}

\noindent Alors, en appliquant (viia) et se fondant sur le fait que $%
X_{n}\rightsquigarrow X$, on a%
\begin{equation*}
\limsup_{n\rightarrow +\infty}\left\vert \mathbb{E}f(Y_{n})-\mathbb{E}%
f(X)\right\vert \leq \limsup_{n\rightarrow +\infty}\mathbb{E}\left\vert
f(Y_{n})-f(X_{n})\right\vert .
\end{equation*}

\noindent Avec la m\^{e}me m\'{e}thode utilis\'{e}e dans la preuve de la
proposition \ref{cvcp.prop2}, nous aurons%
\begin{eqnarray*}
\mathbb{E}\left\vert f(Y_{n})-f(X_{n})\right\vert &\leq
&\int_{(d(Y_{n},X_{n})\leq \varepsilon )}\left\vert
f(Y_{n})-f(X_{n})\right\vert d\mathbb{P}\\
&+&\int_{(d(Y_{n},X_{n})>\varepsilon
)}\left\vert f(Y_{n})-f(X_{n})\right\vert d\mathbb{P} \\
&\leq &\ell \varepsilon +2M\text{ }d(Y_{n},X_{n}),
\end{eqnarray*}

\noindent ce qui tend vers zero quand $n\rightarrow \infty $ puis $%
\varepsilon \downarrow 0.$ Et par suite%
\begin{equation*}
\limsup_{n\rightarrow +\infty}\left\vert \mathbb{E}f(Y_{n})-\mathbb{E}%
f(X)\right\vert =0.
\end{equation*}

\bigskip \noindent \textbf{(d) Theor\`{e}me de Slutsky}.\\

\noindent Nous avons cet important outil de convergence vague.

\begin{proposition} \label{cv.slutsky}
\bigskip\ Si \ $X_{n}\rightsquigarrow X$ \ et $Y_{n}\rightsquigarrow c$,
alors $(X_{n},Y_{n})\rightsquigarrow (X,c)$
\end{proposition}

\noindent \textbf{Preuve}. Soit $X_{n}\rightsquigarrow X$ \ et $%
Y_{n}\longrightarrow _{\mathbb{P}}c.$ On veut montrer que $%
(X_{n},Y_{n})\rightsquigarrow (X,c).$ Remarquons d'abord que \ $%
Y_{n}\longrightarrow _{\mathbb{P}}c$ puisque $Y_{n}\rightsquigarrow c. $\ De
plus, sur $S^{2}$ munie de la m\'{e}trique euclidienne%
\begin{equation*}
d_{e}((x^{\prime },y^{\prime }),(x^{\prime \prime },y^{\prime \prime }))=%
\sqrt{d(x^{\prime },x^{\prime \prime })^{2}+d(y^{\prime },y^{\prime \prime
})^{2}},
\end{equation*}

\noindent nous avons 
\begin{equation*}
d_{e}((X_{n},Y_{n}),(X_{n},c))=d(Y_{n},c).
\end{equation*}

\noindent Il s'en suit que pour $\varepsilon >0$%
\begin{equation*}
\limsup_{n\rightarrow +\infty}\mathbb{P}(d_{e}((X_{n},Y_{n}),(X_{n},c))>%
\varepsilon )=\limsup_{n\rightarrow +\infty}\mathbb{P}(d(Y_{n},c)>%
\varepsilon )=0,
\end{equation*}

\noindent puisque $Y_{n}\longrightarrow _{\mathbb{P}}c$. Ainsi $%
d_{e}((X_{n},Y_{n}),(X_{n},c))\longrightarrow _{\mathbb{P}}0.$ D'apr\`{e}s
la proposition de la sous-section XX, il suffit de montrer la convergence
vague de $(X_{n},c)$ pour avoir celle de $(X_{n},Y_{n})$.\newline

\noindent Maintenant, montrer la convergence vague de $(X_{n},c)$ vers $(X,c)
$, revient \`{a} montrer que pour fonction $g(\cdot ,\cdot )$ r\'{e}elle
continue et born\'{e}e d\'{e}finie sur $S^{2}$, on a $\mathbb{E}%
g(X_{n},c)\rightarrow \mathbb{E}g(X,c).$ Mais ceci vient du fait que, $c$ 
\'{e}tant fix\'{e}e, la fonction $f(x)=g(x,c)$ est continue et born\'{e}e et 
$\mathbb{E}f(X_{n})\rightarrow \mathbb{E}f(X)$ puisque $X_{n}%
\rightsquigarrow X$. Mais $\mathbb{E}f(X_{n})\rightarrow \mathbb{E}f(X)$ 
\'{e}quivaut \`{a} $\mathbb{E}g(X_{n},c)\rightarrow \mathbb{E}g(X,c).$ Ce
qui finit la preuve.\newline

\bigskip \noindent \textbf{(e) Les convergence par coordonn\'{e}es vague et en probabilit\'{e} ne sont pas \'{e}quivalentes}.\\

\begin{proposition} \label{cv.CPcoordinates}
\bigskip Si \ $X_{n}\longrightarrow _{\mathbb{P}}X$ \ et $%
Y_{n}\longrightarrow _{\mathbb{P}}Y$, alors $(X_{n},Y_{n})\longrightarrow _{%
\mathbb{P}}(X,Y)$\newline
\end{proposition}

\noindent \textbf{Preuve}. Utilisons la distance de Manhattan sur $S^{2}:$%
\begin{equation*}
d_{m}((x^{\prime },y^{\prime })x^{\prime },y^{\prime }),(x^{\prime \prime
},y^{\prime \prime }))=d(x^{\prime },x^{\prime \prime })+d(y^{\prime
},y^{\prime \prime }).
\end{equation*}

\noindent Nous aurons pour $\varepsilon >0,$ $\limsup_{n\rightarrow +\infty}\mathbb{P}(d_{m}((X_{n},Y_{n}),(X_{n},c))>\varepsilon )$ est 
\begin{eqnarray*}
&=&\limsup_{n\rightarrow +\infty}\mathbb{P}(d(X_{n},Y_{n})+d(X,Y))>%
\varepsilon ) \\
&\leq &\limsup_{n\rightarrow +\infty}\mathbb{P}(d(X_{n},X)>\varepsilon
/2)+\mathbb{P}(d(Y_{n},Y)>\varepsilon /2) \\
&\leq &\limsup_{n\rightarrow +\infty}\mathbb{P}(d(X_{n},X)>\varepsilon
/2)+\limsup_{n\rightarrow +\infty}\mathbb{P}(d(Y_{n},Y)>\varepsilon /2) \\
&=&0.
\end{eqnarray*}

\bigskip

\subsection{Th\'{e}or\`{e}me de Skorohod-Wichura} \label{cv.subsec.skorohod}

\noindent Pour simplifier, mettons nous dans un espace m\'{e}trique complet s%
\'{e}parable. Enon\c{c}ons le th\'{e}or\`{e}me.

\begin{theorem} \label{cv.skorohodWichura} Soit une suite de variables al\'{e}atoires $(X_{n})_{n\geq 0}$ et $X$ une
autre variable al\'{e}atoire, toutes \`{a} valeurs dans $(S,d)$ complet s%
\'{e}parable, pas n\'{e}cessairement d\'{e}finies sur le m\^{e}me espace de
probabilit\'{e}.\newline

\noindent Si $X_{n}\rightsquigarrow X,$ alors on peut construire un espace
de probabilit\'{e} $(\Omega ,\mathcal{A},\mathbb{P})$ portant des variables
al\'{e}atoires $\ (Y_{n})_{n\geq 0}$ et $Y$ , $\ $telles que 
\begin{equation*}
\mathbb{P}_{X}=\mathbb{P}_{Y}\text{ et }\left( \forall n\geq 0,\mathbb{P}%
_{X_{n}}=\mathbb{P}_{Y_{n}}\right)
\end{equation*}%
et%
\begin{equation*}
Y_{n}\rightarrow Y,\text{ }p.s.
\end{equation*}
\end{theorem}

\bigskip

\noindent Ce th\'{e}or\`{e}me est puissant et peut se r\'{e}v\'{e}ler utile
dans certaines situations.\\

\section{Annexe} \label{cv.annexe}
\subsection{Intervalles $F$-continus o\`u $F$ est une fonction de distribution} \label{cv.subsFcontinuous}

\bigskip \noindent Soit $\mathbb{P}$ une probabilit\'{e} $\mathbb{P}$ sur $(\mathbb{R}^{k},%
\mathcal{B}(\mathbb{R}^{k})).$ Consid\'{e}rons sa fonction de r\'{e}partion 
\begin{equation*}
(x_{1},...,x_{k})\hookrightarrow F(x_{1},...,x_{k})=P\left(
\prod\limits_{i=1}^{k}]-\infty ,x_{i}]\right) .
\end{equation*}

\subsubsection{Intervalles de continuit\'{e} de F}.

\bigskip \noindent Soit un intervalle born\'{e} de $\mathbb{R}^{k}$ de la forme 
\begin{equation*}
]a,b]=\prod\limits_{i=1}^{k}]a_{i},b_{i}].
\end{equation*}

\bigskip \noindent D\'{e}finissons 
\begin{equation*}
E(a,b)=\{c=(c_{1},...,c_{k})\in \mathbb{R}^{k},\text{ }\forall 1\leq i\leq
k,(c_{i}=a_{i}\text{ ou }c_{i}=b_{i})\}.
\end{equation*}

\bigskip \noindent  Nous dirons que l'invervalle $]a,b]$ est F-continue ssi $]a,b]$ born\'{e} et
nous avons
\begin{equation*}
\forall c\in E(a,b),\mathbb{P}(\partial ]-\infty ,c])=0.
\end{equation*}

\bigskip \noindent Soit $\mathcal{U}$ la classe des intervalles $F$-continus. Par convention, nous mettons
l'ensemble vide dans $\mathcal{U}$ puisque la condition d'appartenance ne peut \^etre v\'{e}rifi\'{e}e. Nous allons voir quelques propri\'{e}t\'{e}s de la classe $\mathcal{U}$.

\subsubsection{$\mathcal{U}$ est stable par intersection finies} \label{subsubsecUstable}

\bigskip \noindent  Soit $]a,b]=\prod\limits_{i=1}^{k}]a_{i},b_{i}]\in U$ et $%
]c,d]=\prod\limits_{i=1}^{k}]c_{i},d_{i}]$. Nous avons%
\begin{equation*}
]a,b]\cap ]c,d]=\prod\limits_{i=1}^{k}]a_{i}\vee c_{i},b_{i}\wedge
d_{i}]=]\alpha ,\beta ]
\end{equation*}

\bigskip \noindent o\`u $x\vee y$ et $x\wedge y$ d\'{e}signent respectivement le maximum et
minimum de $x$ et $y$ et, \ \ $\alpha =(a_{1}\vee c_{1},..,a_{k}\vee c_{k})$
et $\beta =(b_{1}\wedge d_{1},...,b_{k}\wedge d_{k}).$ Si $]a,b]\cap ]c,d]$
est vide, il est dans $\mathcal{U}$. Sinon, aucun des facteurs $]a_{i}\vee
c_{i},b_{i}\wedge d_{i}]$ n'est vide. Nous allons montrer que%
\begin{equation}
\forall e\in E(\alpha ,\beta ),\partial ]-\infty ,e]\subset
\bigcup\limits_{z\in E(a,b)\cup E(c,d)}\partial ]-\infty ,z].
\label{unionEab}
\end{equation}

\bigskip \noindent En effet, Soit $e\in E(\alpha ,\beta ).$ On a donc%
\begin{equation*}
e_{i}=a_{i}\vee c_{i}\text{ ou }b_{i}\wedge d_{i},\text{ }1\leq i\leq k.
\end{equation*}%
Soit $t\in \partial ]-\infty ,e].$ Cela veut dire que%
\begin{equation*}
(t_{i}\leq c_{i},\text{ }1\leq i\leq k)\text{ et }(\exists i_{0}\text{,}%
t_{i_{0}}=e_{i_{0}})
\end{equation*}

\bigskip \noindent Puisque $]\alpha ,\beta ]$ est dans $]a,b]$ et $]c,d],$ $t$ v\'{e}rifie%
\begin{equation*}
t_{i}\leq b_{i}\text{ et }t_{i}\leq d_{i},1\leq i\leq k.
\end{equation*}

\bigskip \noindent Maintenant, on consid\`{e}re un $i_{0}$ tel que $t_{i_{0}}=c_{i_{0}}.$ Nous
avons quatre cas%
\begin{equation*}
\left\{ 
\begin{tabular}{lll}
$t_{i_{0}}=e_{i_{0}}=a_{i_{0}}\vee c_{i_{0}}=a_{i_{0}}$ & $\Longrightarrow $
& $t_{i_{0}}=a_{i_{0}}$ et $t_{i}\leq b_{i},1\leq i\leq k$ \\ 
$t_{i_{0}}=e_{i_{0}}=a_{i_{0}}\vee c_{i_{0}}=c_{i_{0}}$ & $\Longrightarrow $
& $t_{i_{0}}=c_{i_{0}}$ et $t_{i}\leq d_{i},1\leq i\leq k$ \\ 
$t_{i_{0}}=e_{i_{0}}=b_{i_{0}}\wedge d_{i_{0}}=b_{i_{0}}$ & $\Longrightarrow 
$ & $t_{i_{0}}=b_{i_{0}}$ et $t_{i}\leq b_{i},1\leq i\leq k$ \\ 
$t_{i_{0}}=e_{i_{0}}=b_{i_{0}}\wedge d_{i_{0}}=d_{i_{0}}$ & $\Longrightarrow 
$ & $t_{i_{0}}=d_{i_{0}}$ et $t_{i}\leq d_{i},1\leq i\leq k$%
\end{tabular}%
.\right. 
\end{equation*}

\bigskip \noindent We conclude the following conclusions for the four lines. First line : $t\in
\partial ]-\infty ,z_{1}]$ where $%
z_{1}=(b_{1},...,b_{i_{0}-1},a_{i_{0},}b_{i_{0}+1},b_{k})\in E(a,b).$ Second
line : $t\in \partial ]-\infty ,z_{2}]$ where $%
z_{2}=(d_{1},...,d_{i_{0}-1},c_{i_{0},}d_{i_{0}+1},d_{k})\in E(c,d).$ Third
line : $t\in \partial ]-\infty ,b]$ and of course $b\in E(a,b).$ Fourth line
: $t\in \partial ]-\infty ,d]$ and of course $d\in E(c,d).$ So $t$ is one of
the  $\partial ]-\infty ,z]$ with $z\in E(a,b)\cup E(c,d).$ So \ref{unionEab}
est vraie et puisque cette union est finie et, est compos\'{e}e de parties
mesurables de probabilit\'{e} nulle nous avons

\begin{equation*}
\forall e\in E(\alpha ,\beta ),P(\partial ]-\infty ,e])=0.
\end{equation*}

\noindent D\`{e}s lors, $\mathcal{U}$ est stable par intersection finie. Nous avons :\\

\begin{lemma} \label{cv.annexe.Inc}
Tout voisinage d'un point $x$ contient un intervalle $]a,b]$
F-continue contenant $x$.
\end{lemma}

\bigskip \noindent Soit $V$ un voisinage de $x.$ Il existe un intervalle $]a,b[$ tel
que
\begin{equation*}
x\in \prod\limits_{i=1}^{k}]a_{i},b_{i}[
\end{equation*}

\bigskip \noindent Soit
\begin{equation*}
\varepsilon _{0}=\min (x_{i}-a_{i},1\leq i\leq k)\wedge \min
(b_{i}-x_{i},1\leq i\leq k).
\end{equation*}

\bigskip \noindent Notons $\delta =(1,...,1)$ le vecteur de R$^{k}$ dont toutes les composantes
sont nulles. Nous avons $0<\varepsilon <$ $\varepsilon _{0},$%
\begin{equation*}
]a+\varepsilon \delta ,x+\varepsilon \delta ]\subset ]a,b[.
\end{equation*}

\bigskip \noindent Tout point $e$ de $E(a+\varepsilon \delta ,x+\varepsilon \delta )$ est de la
forme
\begin{equation*}
t(\varepsilon )=(t_{1}+\varepsilon ,t_{2}+\varepsilon ,...,t_{k}+\varepsilon
)
\end{equation*}

\bigskip \noindent avec bien sur $t_{i}=a_{i}$ ou $t_{i}=x_{i}.$ Pour un choix de $%
t=(t_{1},...,t_{k}),$ les ensembles $\partial ]-\infty ,t(\varepsilon )]$
sont disjoints. Donc, en dehors d'un partie d\'{e}nombrable $D(t)$ de $%
]0,\varepsilon _{0}[$ on aura 
\begin{equation*}
P(\partial ]-\infty ,t(\varepsilon )])=0
\end{equation*}

\bigskip \noindent En dehors de l'ensemble d\'{e}nombrable $D=\cup _{t}D(t)$ $\subset
]0,\varepsilon _{0}[$, (puisque $D$ est union de $2^{k}$ ensembles d\'{e}%
nombrables), nous pouvons choisir un $\varepsilon$ de $]0,\varepsilon _{0}[$ tel que pour tout $e$ de composantes
\begin{equation*}
e_{i}=a_{i}+\varepsilon \text{ or }x_{i}+\varepsilon,
\end{equation*}

\bigskip \noindent nous avons 
\begin{equation*}
\mathbb{P}(\partial ]-\infty ,t(\varepsilon )])=0
\end{equation*}

\bigskip \noindent et donc
\begin{equation*}
x\in ]a+\varepsilon \delta ,x+\varepsilon \delta ]\subset ]a,b[.
\end{equation*}

\bigskip \noindent Appelons $\varepsilon(x)$ la valeur trouv\'ee de $\varepsilon$. Nous venons de montrer que pour tout voisinage V de $x$, il existe $]A_{x},B_{x}[=]a+\varepsilon(x) \delta
,x+\varepsilon(x) \delta /2[$ et $]a_{x},a_{x}]=]a+\varepsilon(x) \delta
,x+\varepsilon(x) \delta ] \in \mathcal{U}$ tels que 

\begin{equation}
x\in ]A_{x},B_{x}[\subset ]a_{x},b_{x}]\subset V.  \label{cleEab}
\end{equation}

\noindent A partir de l\`a, montreons que \textit{Tout ouvert $G$ de $\mathbb{R}^{k}$ est union d\'{e}nombrable d'invervalles $F$-continue}.\\

\bigskip \noindent  En effet, d'apr\`{e}s (\ref{cleEab}), tout ouvert $G$ peut s'\'{e}crire%
\begin{equation*}
G=\bigcup\limits_{x\in G}]A_{x},B_{x}[.
\end{equation*}

\bigskip \noindent Puisque $\mathbb{R}^{k}$ un espace m\'{e}trique s\'{e}parable, ce
recouvrement ouvert se r\'{e}duit \`{a} un ses sous-recouvrements d\'{e}nombrables,
i.e., il existe une suite $(x_{j})_{j\geq 0}\subset G$ telle que%
\begin{equation*}
G=\bigcup\limits_{j\geq 0}]A_{x_{j}},B_{x_{j}}[.
\end{equation*}

\bigskip \noindent On en d\'{e}duit que%
\begin{equation*}
G=\bigcup\limits_{j\geq 0}]a_{x_{j}},b_{x_{j}}].
\end{equation*}

\bigskip \noindent o\`{u} les $]a_{x_{j}},b_{x_{j}}]$ sont des intervalles $F$-continus. Nous concluons par :

\begin{proposition} \label{cv.GFcontinuous}
Soit $F$ une fonction de distribution sur $\mathbb{R}^{k}$, $k\geq 1$. Alors, tout ouvert $G$ dans $\mathbb{R}^{k}$ est une union d\'enombrable 
d'intervalles $F$-continus de la form $]a,b]$ or $]a,b[$, o\`u par d\'efinition, un intervalle $(a,b)$ est $F$-continu si et seulement si pour tout  
$$
\varepsilon=(\varepsilon_1, \varepsilon_2, ..., \varepsilon_k) \in \{0,1\}^k,
$$

\noindent le point
$$
b+\varepsilon*(a-b)=(b_1+\varepsilon_1 (a_1-b_1), b_2+\varepsilon_2 (a_2-b_2), ..., b_k+\varepsilon_k (a_k-b_k))
$$ 

\noindent est un point de continuit\'e de $F$.
\end{proposition}

\subsection{Fonctions semi-continues} \label{cv.subsec.annexe.semic}

\bigskip \noindent Une fonction $f : S \mapsto \overline{\mathbb{R}}$\ est
continue en tout $x$ ssi

\bigskip \noindent (i) pour tout x$\in \mathbb{R}$, pour tout $\varepsilon >0,$%
\ il existe un voisinage $V$ de $x$ tel que 
\begin{equation*}
y\in V\Rightarrow f(y)\in ]f(x)-\varepsilon ,f(x)+\varepsilon \lbrack .
\end{equation*}

\bigskip \noindent Dans cette formule, on s'int\'{e}resse \`{a} tout
l'intervalle $]f(x)-\varepsilon ,f(x)+\varepsilon \lbrack .$\ Mais on peut s'int%
\'{e}resser uniquement \`a l'une des bornes de l'intervalle. Cela nous donne les
fonctions semi-continues. Pr\'{e}cis\'{e}ment, $f$ est dite semi-continue sup%
\'{e}rieurement (not\'{e} s.c.s ) ssi\\

\bigskip \noindent (ii) pour tout $x\in \mathbb{R}$, pour tout $\varepsilon >0,$%
\ il existe un voisinage $V$ de $x$ tel que 
\begin{equation*}
y\in V\Rightarrow f(y)<f(x)+\varepsilon
\end{equation*}

\bigskip \noindent Elle est dite semi-continue inf\'{e}rieurement (not\'{e}
s.c.i) ssi\\

\bigskip \noindent (iii) Pour tout $x\in \mathbb{R}$, pour tout $\varepsilon
>0, $\ il existe un voisinage $V$ de $x$ tel que 
\begin{equation*}
y\in V\Rightarrow f(y)>f(x)-\varepsilon
\end{equation*}

\bigskip \noindent Nous avons ces propri\'et\'es imm\'ediates.\\

\noindent (a) Une function num\'erique est continue si et seulement si elle est \`a la fois \textit{s.c.s} and  \textit{s.c.i}.\\

\noindent (b) Une function numérique $f$ est \textbf{s.c.s} si et seulement si son opposé $-f$ is \textit{s.c.i}.\\

\noindent Voici une caract\'erisation des fonctions semi-continues.\\

\begin{proposition}\label{cv.annexe.SC} Nous avons les propr\'et\'es suivantes :\\

\noindent (1) Une fonction $f : S \mapsto \mathbb{R}$ est \textit{s.c.s} si et seulement si l'ensemble $(f\geq c)$ est ferm\'e pour tout nombre r\'eel $c$.\\

\noindent (2) Une fonction $f : S \mapsto \mathbb{R}$ est \textit{s.c.s} si et seulement si l'ensemble $(f\geq c)$ est ferm\'e pour tout nombre r\'eel $c$.\\

\noindent (3) Si $f$ est \textit{s.c.s} or \textit{s.c.i}, alors elle est mesurable.\\
\end{proposition}

\bigskip \noindent \textbf{Preuve}. Preuve du point (1). Commen\c{c}ons par le sens direct. Sout $f$ une fonction \textit{s.c.s} from $S$ to $\mathbb{R}$. Montrons que pour tout r\'eel $c$ fix\'e, l'ensemble $(f\geq c)$ est ferm\'e en montrant que l'ensemble $(f < c)$ est ouvert. Soit $x \in G=(f < c)^{c}$, c'est-\`a-dire que, $f(x)<c$. Posons $\varepsilon=c-f(x) >0$. Puisque $f$ est \textit{s.c.s}, il existe un voisinage de $V$ de $x$ tel que  

\begin{equation*}
y\in V\Rightarrow f(y) < f(x)+\varepsilon=c,
\end{equation*}

\noindent ce qui peut \^etre \'ecrit sous la forme
\begin{equation*}
y\in V \Rightarrow f(y) < c,
\end{equation*}

\noindent ce qui signifie que $V\subseteq G^c$. Nous avons prouv\'e que $G^c$ contient chacun de ses points avec un de ses voisinages. Donc $G^c$ est ouvert. Le sens direct est prouv\'e.\\

\bigskip \noindent Inversement, supposons que pour tout nombre r\'eel $c$, l'ensemble $(f\geq c)$ est ferm\'e. Soit $x\in S$\ quelquonque. Donc pour tout $\varepsilon >0$, $G=(f<f(x)+\varepsilon )$ est ouvert. Mais $x$ appartient \`{a} $G$, donc $G$ le contient avec un de ses
voisinages $V\in V(x)$\ et donc 
\begin{equation*}
y\in V\Rightarrow y\in (f<f(x)+\varepsilon )\Rightarrow f(y)\leq f(x)+\varepsilon.
\end{equation*}

\bigskip \noindent Donc f est semi-continue sup\'{e}rieurement.\\

\noindent Le point (2) se prouve directement \`a partir du point 1, en utilisant la transformation $-f$.\\

\noindent Le point (3) d\'ecoule des crit\`eres classiques de mesurabilit\'e pour les fonctions r\'eelles.\\

\subsection{Propri\'et\'e caract\'eristique d'une famille parties mesurables disjointes} \label{cv.subsec.annexe.Disjoint}

\begin{proposition} \label{cv.annexe.FDS}
Soit une famille (B$_{\lambda })_{\lambda \in \Gamma }$\ de parties
mesurables deux \`{a} deux disjointes d'un espace probabilis\'{e} $(\Omega,\mathcal{A},L)$. Alors un nombre au plus d\'{e}nombrable d'entre elles ont une probabilit\'{e} non nulle.
\end{proposition}

\bigskip \noindent \textbf{Preuve}.\newline
Soit l'ensemble des indices $\lambda $\ pour lesquels B$_{\lambda }$\ est de
probabilit\'{e} non nulle, 
\begin{equation*}
D=\{\lambda \in \Gamma ,\text{ }L(B_{\lambda })>0\}.
\end{equation*}
On a surement 
\begin{equation*}
D=\cup _{k\geq 1}D_{k},
\end{equation*}
avec 
\begin{equation*}
D_{k}=\{\lambda \in \Gamma ,\text{ }L(B_{\lambda })>1/k\}.
\end{equation*}

\bigskip \noindent Maintenant soit r \'{e}l\'{e}ments de D$_{k}$\ not\'{e}s $%
\lambda _{1},\lambda _{2},...,\lambda _{r},$\ on a, par le fait que les
ensembles sont disjoints deux \`{a} deux, 
\begin{equation*}
1\geq L(\bigcup_{1}^{r}B_{\lambda _{j}})=\sum_{1}^{r}L(B_{\lambda _{j}})\geq
r/k.
\end{equation*}
D'o\`{u} 
\begin{equation*}
r\leq k.
\end{equation*}
Donc $D_{k}$\ contient au plus $k$ \'{e}l\'{e}ments. Par suite, $D$ est une
union d\'{e}nombrable d'ensembles finis. D est donc au plus d\'{e}nombrable.

\subsection{Mesurabilit\'{e} de l'ensemble des points de discontinuit\'{e}.} \label{cv.subsec.annexe.Discontinuity}

Voil\`{a} un r\'{e}sultat surprenant, \`{a} savoir que l'ensemble des points
de discontinuit\'{e} d'une application quelconque g, not\'{e} $discont(g),$\
d'un espace m\'{e}trique (S, d) dans un autre (D, r) est mesurable. En effet, nous avons

\begin{lemma} \label{cv.annexe.Discont}
Soit $g$ une fonction de l'espace m\'etrique $(S, d)$ vers l'espace m\'etrique $(D, r)$. Soit $discont(g)$, l'ensemble des points continuit\'e de $g$.  Nous avons

\begin{equation}
discont(g)=\bigcup_{s=1}^{\infty }\bigcap_{t=1}^{\infty }B_{s,t} \label{annexe4}
\end{equation}%

\noindent o\`u pour chaque couple d'entiers $(s,t)$, l'ensemble 
\begin{equation*}
B_{s,t}=\left\{ x\in S,\exists (y,z)\in S^{2},\text{ }d(x,y)<1/t,d(z,x)<1/t,%
\text{ }r(g(y),g(z))\geq 1/s\right\} .
\end{equation*}

\noindent est ouvert.
\end{lemma}

\bigskip \noindent La cons\'equence de ce lemme est que $discont(g)$ est mesurable.\\

\bigskip \noindent \textbf{Preuve du lemme}. Montrons d'abord que
 
\begin{equation*}
\bigcup_{s=1}^{\infty }\bigcap_{t=1}^{\infty }B_{s,t}\subseteq discont(g)
\end{equation*}

\bigskip \noindent Soit $x\in \bigcup_{s=1}^{\infty }\bigcap_{t=1}^{\infty
}B_{s,t}.$\ Il existe un entier s$\geq 1\geq fix\acute{e}$\ tel que pour
chaque entier t$\geq 1,$\ il existe y$_{t}$\ et z$_{t}$\ tel tel que 
\begin{equation*}
d(x,y_{t})<1/t,
\end{equation*}

\bigskip \noindent et 
\begin{equation*}
d(x,z_{t})<1/t,
\end{equation*}

\bigskip \noindent et 
\begin{equation}
\forall \text{ }t\geq 1,\text{ }r(g(y_{t}),g(z_{t}))\geq 1/s  \label{cv29}
\end{equation}

\bigskip \noindent Et si $g$ est continue en $x$, alors par continuit\'{e},\\
 
\begin{equation*}
r(g(y_{t}),g(z_{t}))\leq r(g(y_{t}),g(x))+r(g(x_{t}),g(z_{t}))\rightarrow 0
\end{equation*}

\bigskip \noindent ce qui est en contradiction avec (\ref{cv29}). Donc x
n'est pas un point de continuit\'{e}, d'o\`{u} $x\in discont(g)$.\newline

\bigskip \noindent Montrons 
\begin{equation*}
discont(g)\subseteq \bigcup_{s=1}^{\infty }\bigcap_{t=1}^{\infty }B_{s,t}.
\end{equation*}

\bigskip \noindent Soit $x$ un point de discontinuit\'{e} de g. Par n\'{e}%
gation de la continuit\'{e}, 
\begin{equation*}
\exists \text{ }\varepsilon >0,\forall \text{ }\eta >0,\exists \text{ }y\in S,%
\text{ }d(x,y)<\eta ,\text{ \ }r(g(y),g(x))\geq \varepsilon .
\end{equation*}

\bigskip \noindent Soit s un entier tel que $\varepsilon \geq 1/s,$\ donc pour
tout $1/t$\ o\`{u} t en entier positif non nul, 
\begin{equation*}
\exists \text{ }y\in S,\text{ }d(x,y)<1/t,\text{ \ }r(g(y),g(x))\geq 1/s.
\end{equation*}

\bigskip \noindent En posant $z=x$, on a bien 
\begin{equation*}
d(x,z)<1/t,\text{ }d(x,y)<1/t,\text{ \ }r(g(y),g(x))\geq 1/s.
\end{equation*}

\bigskip \noindent Donc $x\in \bigcup_{s=1}^{\infty }\bigcap_{t=1}^{\infty
}B_{s,t}.$\ D'o\`{u} l'\'{e}galit\'{e}.\newline

\bigskip \noindent Montrons enfin que chaque $B_{s,t}$\ est ouvert. Posons $%
a=1/s>0$\ $et$\ $b=1/t>0$. Soit x$\in B_{s,t},$\ Donc 
\begin{equation*}
\exists (y,z)\in S^{2},\text{ }d(x,y)<b,d(z,x)<b,\text{ }r(g(y),g(z))\geq a
\end{equation*}

\bigskip \noindent Soit $c=min(b-d(x,y),b-d(z,x))>0$. Soit $x^{\prime }\in
B(x,c),$\ donc 
\begin{equation*}
d(x^{\prime },y)\leq d(x^{\prime },x)+d(x,y)<c+d(x,y)\leq b
\end{equation*}

\bigskip \noindent et 
\begin{equation*}
d(x^{\prime },z)<d(x^{\prime },x)+d(x,z)\leq c+d(x,z)\leq b
\end{equation*}

\bigskip \noindent et 
\begin{equation*}
r(g(y),g(z))\geq a
\end{equation*}

\bigskip \noindent D'o\`{u} $x^{\prime }\in B_{s,t}.$\ D'o\`{u} 
\begin{equation*}
x\in B(x,c)\subseteq B_{s,t}
\end{equation*}

\bigskip \noindent Ainsi chaque $B_{s,t}$\ contient ses points avec des
boules ouvertes. Donc chaque $B_{s,t}$\ est ouvert. Donc $discont(g)$\ est
mesurable.

\subsection{Th\'{e}or\`{e}me de Stone-Weieirstrass} {cv.subsec.annexe.SW1}

\begin{proposition} {cv.subsec.annexe.SW.prop}
Soit (S, d) un espace m\'{e}trique compact et H une partie non vide de $%
\mathcal{C}(S, \mathbb{R})$ l'ensemble des fonctions continues de S dans $R$%
\ verifiant les propri\'{e}t\'{e}s\newline

\bigskip \noindent (i) H est r\'{e}ticul\'{e}e, i.e, si f et g sont deux 
\'{e}l\'{e}ments de $H$, alors $f\wedge g$\ et $f\vee g$\ appartient \`{a} $g$%
\newline

\bigskip \noindent (ii) Si x et y sont \'{e}l\'{e}ments de S, et (a,b) un
couple de r\'{e}els (avec a=b si x=y), alors il existe deux \'{e}lements h
et k de H tels que\newline
\begin{equation*}
h(x)=a\text{ et }k(y)=b.
\end{equation*}

\bigskip \noindent Alors $H$ est dense de $C(S, \mathbb{R})$ munie de sa
topologie uniforme, c'est-\`{a}-dire, que toute fonction continue de $S$
dans $\mathbb{R}$ est limite uniforme d'une suite d'\'{e}l\'{e}ments de $H$.%
\end{proposition}

\begin{theorem} \label{cv.subsec.annexe.SW}
Soit $(S, d)$ un espace m\'{e}trique compact et $H$ une partie non vide de $%
C(S, \mathbb{R})$ l'ensemble des fonctions continues de $S$ dans $C$\
verifiant les propri\'{e}t\'{e}s :\\

\noindent (i) $H$ contient les fonctions constantes.\\

\noindent (ii) $Si(h,k)\in H^{2}$, $h+k\in H,h\times k\in H,\overline{u}\in H$.\\

\noindent (iii) $H$ s\'{e}pare les points de $S$, i.e., pour tous \'{e}l\'{e}ments x et
y distinctes de S, $x\neq y,$\ alors il existe $h\in H$%
\begin{equation*}
h(x)\neq h(y).
\end{equation*}

\noindent Alors $H$ est dense de $C(S, \mathbb{R})$ muni de sa topologie uniforme,
c'est-\`{a}-dire, que toute fonction continue de S dans $C$\ est limite
uniforme d'une suite d'\'{e}l\'{e}ments de $H$.
\end{theorem}

\bigskip \noindent \textbf{Remarque}.\newline

\noindent Si on travaille sur $\mathbb{R}$, le th\'{e}or\`{e}me est vraie et la
condition $\overline{u}\in H,$\ n'a pas de sens.\\

\subsection{Divers} \label{cv.subsec.annexe.divers}

\bigskip \noindent \textbf{Une relation utile}. A prouver pour des r\'{e}els
x, y, X, et Y, 
\begin{equation}
\left\vert \min (x,y)-\min (X,Y)\right\vert \leq \left\vert x-X\right\vert
+\left\vert y-Y\right\vert  \label{annexe2}
\end{equation}%
En effet si min(x,y)=x et min(X,Y)=X, 
\begin{equation*}
\left\vert \min (x,y)-\min (X,Y)\right\vert \leq \left\vert x-X\right\vert
\end{equation*}%
\ si min(x,y)=y et min(X,Y)=Y, 
\begin{equation*}
\left\vert \min (x,y)-\min (X,Y)\right\vert \leq \left\vert y-Y\right\vert
\end{equation*}%
Maintenant soit min(x,y)=x et min(X,Y)=Y. On peut supposer que x$\leq Y.$\
On aura 
\begin{equation*}
\min (x,y)-\min (X,Y)=Y-x\leq X-x
\end{equation*}%
puisque $X\geq Y.$\ Le cas min(x,y)=y et min(X,Y)=X se traite comme le cas pr%
\'{e}c\'{e}dent. Donc (\ref{annexe2}) est vraie.
 

%% file: asymptotics_cv_02_fr.tex
\chapter{Tension uniforme et tension asymptotique} \label{cv.tensRk}

\section{Introduction}

Toute th\'{e}orie limite poss\`{e}de une partie qui traite de la notion de
compacit\'{e}, c'est-\`{a}-dire de l'existence de sous-suites convergentes,
l'\'{e}quivalent du th\'{e}or\`{e}me de Bolzano-Weierstrass pour les suites
num\'{e}riques. Il s'agira pour la convergence vague de la partie sur la
tension des suites de variables al\'{e}atoires. Au plan g\'{e}n\'{e}ral, le th%
\'{e}or\`{e}me de Prohorov pr\'{e}vaut et affirme que toute suite
asymptotiquement tendue de variables al\'{e}atoires poss\`{e}de une
sous-suite vaguement convergente.\newline

\noindent Dans ce chapitre, nous allons nous restreindre au cas o\`{u} $S$
est $\mathbb{R}^{k}$, puisque le traitement entre le cas g\'{e}n\'{e}ral et
le cas de $\mathbb{R}^{k}$ sont un peu \'{e}loign\'{e}s l'un de l'autre.
Dans ce cas particulier, le th\'{e}or\`{e}me de Helly-Bray jour le grand r%
\^{o}le.\newline

\noindent Ce chapitre est donc centr\'{e} sur $\mathbb{R}^{k}$, o\`{u} $%
k\geq 1$. Deux points sont derri\`{e}re toutes les propri\'{e}t\'{e}s expos%
\'{e}es ici. Le premier est que les parties compactes de $\mathbb{R}^{k}$,
sont les ensembles de ferm\'{e}s et born\'{e}s. Le second est que $\mathbb{R}%
^{k}$ est un espace m\'{e}trique complet s\'{e}p\'{e}rable.\newline

\noindent Dans ce chapitre, sauf indication contraire, la norme $max$ ainsi d\'efinie

\begin{equation*}
\left\Vert x\right\Vert =\max_{1\leq i\leq k}\left\vert x_{i}\right\vert .
\end{equation*}

\noindent est utilis\'{e}e. Ainsi les boules ouvertes $B(x,r)$ et ferm\'{e}s $B^{f}(x,r)$ sont d\'{e}%
finies par
\begin{equation*}
B(x,r)=\{x\in \mathbb{R}^{k},\left\Vert x\right\Vert
<r\}=\prod\limits_{i=1}^{k}]x_{i}-r,x_{i}+r[
\end{equation*}

\noindent 
pour $x=(x_{1},...,x_{k})$ et $r>0$ et 
\begin{equation*}
B^{f}(x,r)=\{x\in \mathbb{R}^{k},\left\Vert x\right\Vert \leq
r\}=\prod\limits_{i=1}^{k}[x_{i}-r,x_{i}+r]
\end{equation*}

\noindent pour $x=(x_{1},...,x_{k})$ et $r\geq 0$.

\bigskip \noindent Avant de commencer, faisons quelques notations.\newline

\noindent Soit $a=(a_{1},...,a_{k})$ et $b=(b_{1},...,b_{k}).$ D\'{e}%
finissons ces relations d'ordre

\begin{equation*}
(a\leq b)\Longleftrightarrow (\forall (1\leq i\leq k),\text{ }a_{i}\leq
b_{i}),
\end{equation*}

\noindent ensuite, 
\begin{equation*}
(a < b)\Longleftrightarrow (\forall (1\leq i\leq k),\text{ }a_{i}\leq
b_{i},\text{ }\exists (1\leq i_{0}\leq k),\text{ }a_{i_{0}}<b_{i_{0}})
\end{equation*}

\noindent et enfin, 
\begin{equation*}
(a \prec b)\Longleftrightarrow (\forall (1\leq i\leq k),\text{ }a_{i}<b_{i},)
\end{equation*}%
et son sym\'{e}trique 
\begin{equation*}
(a \succ b)\Longleftrightarrow (\forall (1\leq i\leq k),\text{ }a_{i}>b_{i},)
\end{equation*}

\bigskip \noindent D\'{e}finissons ces classes d'ensembles compacts.\newline

\bigskip \noindent Pour $A=(A_{1},...,A_{k}) \prec V=(B_{1},...,B_{k}),$ notons%
\begin{equation*}
K_{A,B}=\prod\limits_{i=1}^{k}[A_{i},B_{i}].
\end{equation*}

\noindent Pour $A=(A_{1},...,A_{k}) \succ 0,$ notons%
\begin{equation*}
K_{A}=\prod\limits_{i=1}^{k}[-A_{i},A_{i}].
\end{equation*}

\noindent Pour $M\in \mathbb{R},$ $M>0,$ notons%
\begin{equation*}
K_{c,M}=[-M,M]^{k}.
\end{equation*}

\noindent Ces ensembles $K_{A,B},$ $K_{A}$ et $K_{c,M}$, tous compacts,
servent \`{a} d\'{e}finir la tension de suites de vecteurs al\'eatoires. Ils jouent des r\^{o}les \'{e}%
quivalents dans la d\'{e}finition de la tension. Pour cela \'{e}non\c{c}ons
ici ce r\'{e}sultat d'\'{e}quivalence.\newline

\begin{proposition}
\label{tensTensprop1} Soit $\{\mathbb{P}_{n},n\geq 1\}$ une suite de
probabilit\'{e}s sur $(\mathbb{R}^{k},\mathcal{B}(\mathbb{R}^{k})).$ Les
propositions suivantes sont \'{e}quivalentes.\newline

\noindent \textbf{(1a)} Pour tout $\varepsilon >0,$ il exist un compact $K$
de \ $\mathbb{R}^{k}$ tel que 
\begin{equation*}
\inf_{n\geq 1}\mathbb{P}_{n}(K)\geq 1-\varepsilon .
\end{equation*}

\noindent \textbf{(2a)} Pour tout $\varepsilon >0,$ il existe un r\'{e}el $%
M>0$ tel que%
\begin{equation*}
\inf_{n\geq 1}\mathbb{P}_{n}(K_{c,M})\geq 1-\varepsilon .
\end{equation*}

\noindent \textbf{(3a)} Pour tout $\varepsilon >0,$ il existe un vecteur $%
A=(A_{1},...,A_{k}) \succ 0, $ de \ $\mathbb{R}^{k}$ tels que

\begin{equation*}
\inf_{n\geq 1}\mathbb{P}_{n}(K_{A})\geq 1-\varepsilon .
\end{equation*}

\noindent \textbf{(4a)} Pour tout $\varepsilon >0,$ il existe deux vecteurs $%
A=(A_{1},...,A_{k}) \prec V=(B_{1},...,B_{k})$ de \ $\mathbb{R}^{k}$ tels que%
\begin{equation*}
\inf_{n\geq 1}\mathbb{P}_{n}(K_{A,B})\geq 1-\varepsilon .
\end{equation*}

\noindent \textbf{(1b)} Pour tout $\varepsilon >0,$ il exist un compact $K$
de \ $\mathbb{R}^{k} $ tel que 
\begin{equation*}
\liminf_{n\rightarrow \infty }\mathbb{P}_{n}(K)\geq 1-\varepsilon .
\end{equation*}

\noindent \textbf{(2b)} Pour tout $\varepsilon >0,$ il existe un r\'{e}el $%
M>0$ tel que%
\begin{equation*}
\liminf_{n\rightarrow \infty }\mathbb{P}_{n}(K_{c,M})\geq 1-\varepsilon .
\end{equation*}

\noindent \textbf{(3b)} Pour tout $\varepsilon >0,$ il existe un vecteur $%
A=(A_{1},...,A_{k}) \succ 0,$ de \ $\mathbb{R}^{k}$ tels que%
\begin{equation*}
\liminf_{n\rightarrow \infty }\mathbb{P}_{n}(K_{A})\geq 1-\varepsilon .
\end{equation*}

\noindent \textbf{(4b)} Pour tout $\varepsilon >0,$ il existe deux vecteurs $%
A=(A_{1},...,A_{k}) \prec V=(B_{1},...,B_{k})$ de \ $\mathbb{R}^{k}$ tels que%
\begin{equation*}
\liminf_{n\rightarrow \infty }\mathbb{P}_{n}(K_{A,B})\geq 1-\varepsilon .
\end{equation*}
\end{proposition}

\noindent \textbf{PREUVE}. Il faut noter que nous avons deux groupes de
formules : $(1a)-(4a)$ et $(1b)-(4b)$. En fait, nous allons d\'{e}montrer
les \'{e}quivalences des diff\'erents points de chaque groupe, puis entre
les deux premiers points des deux groupes.\newline

\noindent \textbf{Equivalence entre les points du groupe (1a)-(4b)}: Soit $%
\varepsilon >0$ fix\'{e}. Montrons :\newline

\noindent $(1a)\Longrightarrow (2a)$. Soit K tel $\sup_{n\geq 1}\mathbb{P}_{n}(K)\geq 1-\varepsilon .$ Puisque $K$ est compact, il est born\'{e}. Il
est dans un ensemble $\{x,\left\Vert x\right\Vert \leq M\}=K_{c,M}$ et%
\begin{equation*}
\inf_{n\geq 1}\mathbb{P}_{n}(K_{c,M})\geq \inf_{n\geq 1}\mathbb{P}_{n}(K)\geq 1-\varepsilon .
\end{equation*}

\noindent $(2a)\Longrightarrow (3a)$. Ceci est \'evident, puisque un $K_{c,M} $ est un $K_{A}$ avec $A=(M,M,...,M).$\newline

\noindent $(3a)\Longrightarrow (4a)$. Ceci est aussi \'evident, puisque un $K_{A}$, pour $A \succ 0$, est exactement $K_{-A,A}$.\newline

\noindent $(4a)\Longrightarrow (1a)$. Ceci est encore \'evident puisque $%
K_{A,B}$ est un compact de $\mathbb{R}^k$.\newline

\noindent \textbf{Equivalence entre les points du groupe (1b)-(4b)}. La
preuve est exactement la m\^{e}me que pour le premier groupe.\newline

\noindent \textbf{Equivalence entre les deux groupes}. Il suffit de prouver
ceci : $(1a)\Longleftrightarrow (1b)$\newline

\noindent \textbf{Si (1a) est vraie}, alors pour tout $\varepsilon \geq 1,$
il existe un compact $K$ de $\mathbb{R}^{k}$ tel que, pour tout $n\geq 1$, 
\begin{equation*}
\mathbb{P}_{n}(K)\geq 1-\varepsilon .
\end{equation*}

\noindent D\`{e}s lors%
\begin{equation*}
\liminf_{n\rightarrow \infty }\mathbb{P}_{n}(K)\geq 1-\varepsilon ,
\end{equation*}

\noindent ce qui donne (1b).\newline

\noindent \textbf{Si (1b) est vraie}, alors, pour tout $\varepsilon \geq 1,$
il existe un compact $K$ tel que 
\begin{equation*}
\left\{ \sup_{n\geq 1}\inf_{p\geq n}\mathbb{P}_{p}(K)\right\} \geq
1-\varepsilon /2.
\end{equation*}

\noindent Alors, il existe $N\geq 1,$ tel que 
\begin{equation*}
\inf_{p\geq N+1}\mathbb{P}_{p}(K)\geq 1-\varepsilon ,
\end{equation*}

\noindent c'est-\`{a}-dire que pour tout $n>N,$%
\begin{equation*}
\mathbb{P}_{n}(K)\geq 1-\varepsilon .
\end{equation*}

\noindent Puisque $K$ est un ensemble compact, il est dans un $K_{c,M_{\infty }}$ et
alors, pour tout $n>N$, 
\begin{equation*}
\mathbb{P}_{n}(K_{c,M_{\infty }})\geq 1-\varepsilon .
\end{equation*}

\noindent Maintenant, pour chaque $1\leq j\leq N$ fix\'{e}, l'ensemble $%
(\left\Vert x\right\Vert \leq M)=K_{c,M}$ cro\^{\i}t, avec $M$, vers $\mathbb{R}^{k}$ et donc,  $\mathbb{P}_j(\left\Vert X_{j}\right\Vert \leq M)\uparrow 1$. Ainsi, il existe $M_{j}>0$ r\'{e}el tel que pour tout $1\leq j\leq N,$%
\begin{equation*}
\mathbb{P}_{j}(K_{c,M_{j}})\geq 1-\varepsilon .
\end{equation*}

\noindent Nous venons de montrer que toute probabilit\'{e} $\mathbb{P}^{(0)}$
sur $(\mathbb{R}^{k},\mathcal{B}(\mathbb{R}^{k}))$\ est \textbf{tendue},
c'est-\`{a}-dire que pour tout $\varepsilon >0,$ il existe un compact $%
K^{(0)}=K_{c,M^{(0)}}$\ de $\mathbb{R}^{k}$\ tel que%
\begin{equation}
\mathbb{P}^{0}(K^{(0)})\geq 1-\varepsilon .  \label{tensTindiv}
\end{equation}

\bigskip \noindent Maintenant, si on prend
\begin{equation*}
M=\max (M_{1},...,M_{N},M_{\infty }),
\end{equation*}

\noindent nous avons que les $K_{c,M_{j}},$ $1\leq j\leq M$ et $%
K_{c,M_{\infty }}$ sont tous dans $K_{c,M}$, un compact de $\mathbb{R}^{k}$
et enfin pour tout $n\geq 1,$%
\begin{equation*}
\mathbb{P}_{n}(K_{c,M})\geq 1-\varepsilon 
\end{equation*}

\noindent et enfin%
\begin{equation*}
\inf_{n\geq 1}\mathbb{P}_{n}(K_{c,M})\geq 1-\varepsilon ,
\end{equation*}

\noindent ce qui est (1a).\newline

\bigskip \noindent Voici un autre lien entre de telles formules et les
fonctions de r\'{e}partition associ\'{e}s. A toute probabilit\'{e}, nous
associons sa fonction de r\'{e}partition
\begin{equation*}
F_{\mathbb{P}}(x)=\mathbb{P}(]-\infty ,x])=P(\prod\limits_{i=1}^{k}]-\infty
,x_{i}]),x=(x_{1},...,x_{k})\in \mathbb{R}^{k}.
\end{equation*}

\noindent Cette fonction de r\'{e}partition, \`{a} son tour, d\'{e}termine
la probabilit\'{e} de Lebesgues-Stieljes ainsi : pour tout $(a,b)\in \mathbb{%
R}^{k}\times \mathbb{R}^{k},$ $a\leq b,$

\begin{equation*}
\mathbb{P}(]a,b])=\Delta _{a,b}F=\sum\limits_{\epsilon \in
\{0,1\}^{k}}(-1)^{s(\epsilon )}F(b+\epsilon \ast (a-b))\geq 0,
\end{equation*}

\noindent o\`{u} pour $\epsilon =(\epsilon _{1},...,\epsilon _{k})\in \{0,1\}^{k},$ $%
s(\epsilon )=\epsilon _{1}+...+\epsilon _{k},$ pour $x=(x_{1},...,x_{k})\in 
\mathbb{R}^{k}$, $y=(y_{1},...,y_{k}),$ $x\ast y=(x_{1}y_{1},...,x_{k}y_{k}).
$.\\
 
\noindent Nous allons utiliser les mesures de Lebesques-Stieljes que le lecteur peut r%
\'{e}viser dans les ouvrages \cite{bmtp} et surtout dans le chapitre 1 de 
\cite{ips}. Nous avons besoin de cette notation. Soit $M>0.$ Notons%
\begin{equation*}
L_{M}=\{x,\exists (1\leq i\leq k),x_{i}\leq -c\}
\end{equation*}

\noindent Nous avons :

\begin{proposition}
\label{tensTensprop2} \bigskip Soit $\{\mathbb{P}_{n},n\geq 1\}$ une suite
de probabilit\'{e}s sur $(\mathbb{R}^{k},\mathcal{B}(\mathbb{R}^{k}))$ et
soit la suite de leurs fonctions de r\'{e}partition $\{F_{n}\geq 1\}$ avec $%
F_{\mathbb{P}_{n}}=F_{n}$ for $n\geq 1.$ Alors les trois points suivants
sont \'{e}quivalentes entre eux.\newline

\noindent (1c) Pour tout $\varepsilon >0,$ il existe $0<C\in \mathbb{R}^{k}$
et il existe $c>0$\ such that%
\begin{equation*}
\inf_{n\geq 1}F_{n}(C)\geq 1-\varepsilon 
\end{equation*}

\noindent and

\begin{equation*}
\inf_{n\geq 0}\mathbb{P}_{n}(L_{c})\leq \varepsilon .
\end{equation*}

\noindent (2c) Pour tout $\varepsilon >0,$ il existe $0<c$ tel que pour $%
c^{(k)}=(c,...c)$ et il existe $M>0$ such that 
\begin{equation*}
\inf_{n\geq 1}F_{n}(c^{(k)})\geq 1-\varepsilon 
\end{equation*}%
and%
\begin{equation*}
\inf_{n\geq 0}\mathbb{P}_{n}(L_{M})\leq \varepsilon .
\end{equation*}

\noindent (3c) Pour tout $\varepsilon >0,$ il existe $M>0$ tel que%
\begin{equation*}
\inf_{n\geq 1}P_{n}(K_{c,M})\geq 1-\varepsilon .
\end{equation*}

\noindent Puisque le point (3c) est le point (2c) de la proposition \ref%
{tensTensprop1}, alors les points (3a) et (3b) sont \'{e}quivalents \`{a}
tous les points de cette proposition.
\end{proposition}

\noindent \textbf{Preuve}. Proc\'{e}dons aux preuves des diff\'{e}rentes 
\'{e}quivalences.\\

\noindent \textbf{Preuve de }$\mathbf{(1c)\Longrightarrow (2c)}$. Pour tout $%
\varepsilon >0,$ il existe $0<C\in \mathbb{R}^{k}$ such that%
\begin{equation*}
\inf_{n\geq 1}F_{n}(C)\geq 1-\varepsilon .
\end{equation*}

\noindent Soit $c=max\{C_{i},1\leq i\leq k\}$ . we have ]-$\infty ,C]\subset
]-\infty ,c^{(k)}]$ and $F_{n}(c^{(k)})\geq F_{n}(C),$ 
\begin{equation*}
\inf_{n\geq 1}F_{n}(c^{(k)})\geq 1-\varepsilon .
\end{equation*}%
La preuve est termin\'{e}e puisque la deuxi\`{e}me formule est la m\^{e}me
pour les deux points.\\

\noindent \textbf{Preuve de }$\mathbf{(2c)\Longrightarrow (3c)}$. De $(2c)$,
nous tirons un $d^{(k)}=(d,...,d),$ avec $d>0$, tel que 
\begin{equation*}
\inf_{n\geq 1}F_{n}(d^{(k)})\geq 1-\varepsilon /2,
\end{equation*}

\noindent et un nombre $e>0$ tel que

\begin{equation*}
\inf_{n\geq 1}\mathbb{P}_{n}(L_{e})\leq \varepsilon/2.
\end{equation*}

\noindent En faisant $M=\max (d,e)$, nous avons

\begin{equation*}
\inf_{n\geq 1}F_{n}(M^{(k)})\geq 1-\varepsilon /2.
\end{equation*}

\bigskip \noindent et un nombre $e>0$ tel que

\begin{equation*}
\inf_{n\geq 1}\mathbb{P}_{n}(L_{M})\leq \varepsilon /2.
\end{equation*}

\noindent Maintenant, d\'{e}composons $\mathbb{R}^{k}$ en \ $\mathbb{R}^{k}=L_{M}+L_{M}^{c},$ avec
\begin{equation*}
L_{M}^{c}=\{x,\forall (1\leq i\leq k),x_{i}\geq -M\}
\end{equation*}

\noindent qui se d\'{e}compose lui-m\^{e}me en
 
\begin{eqnarray*}
L_{M}^{c} &=&\{x,\forall (1\leq i\leq k),-M\leq x_{i}\leq M\}\\
&+&\{x,\forall (1\leq i\leq k),x_{i}\geq -M\text{ et }\exists (1\leq i\leq k),x_{i}>M\} \\
&=&K_{c,M}+B,
\end{eqnarray*}

\noindent o\`{u}, de mani\`{e}re \'{e}vidente,

\begin{equation*}
B\subset ]-\infty ,M^{(k)}]^{c}.
\end{equation*}%

\noindent D\`{e}s lors, on peut tirer du faut que $\mathbb{R}^{k}=L_{M}+K_{c,M}+B,$ que

\begin{equation}
K_{c,M}^{c}=L_{M}+B.  \label{decomKcM}
\end{equation}

\noindent Donc, pour tout $n\geq 1,$

\begin{equation*}
\mathbb{P}_{n}(K_{c,M})=\mathbb{P}_{n}(L_{M})+\mathbb{P}_{n}(B)\leq
\varepsilon /2+\varepsilon /2=\varepsilon ,
\end{equation*}

\bigskip \noindent puisque $B\subset ]-\infty ,M^{(k)}]^{c}$ et donc pour tout $n\geq 1,$ 
\begin{eqnarray*}
\mathbb{P}_{n}(B) &\leq &\mathbb{P}_{n}(]-\infty ,M^{(k)}]^{c}) \\
&\leq &1-\mathbb{P}_{n}(]-\infty ,M^{(k)}]) \\
&\leq &1-F_{n}(M^{(k)})\leq \varepsilon /2.
\end{eqnarray*}

\noindent Ce qui termine la preuve de ce point.\\

\noindent \textbf{Preuve de }$\mathbf{(3c)\Longrightarrow (1c)}$. Supposons (3c) : \
pour tout $\varepsilon >0,$ il existe $M>0$ tel que

\begin{equation*}
\inf_{n\geq 1}\mathbb{P}_{n}(K_{c,M})\geq 1-\varepsilon .
\end{equation*}

\noindent Nous avons 
\begin{equation*}
\inf_{n\geq 1}F_{n}(M^{(k)})=\inf_{n\geq 1}\mathbb{P}_{n}(]-\infty
,M^{(k)}])\geq \inf_{n\geq 1}\mathbb{P}_{n}(K_{c,M})\geq 1-\varepsilon .
\end{equation*}

\noindent Ensuite, \`{a} cause de (\ref{decomKcM}),%
\begin{equation*}
\mathbb{P}_{n}(L_{M})\leq \mathbb{P}_{n}(K_{c,M}^{c})\leq \varepsilon .
\end{equation*}

\noindent Donc (1c) a lieu. Et la boucle est boucl\'{e}e.\newline

\noindent Maintenant, exposons la notion de tension.\newline

\section{Tension}

\subsection{Tension individuelle}

Dans ce cas particulier, toute probabilit\'{e} sur \ $(\mathbb{R}^{k},%
\mathcal{B}(\mathbb{R}^{k}))$ est tendue dans le sens suivant

\begin{proposition}
\label{tensTensprop3} Pour toute probabilit\'{e} $\mathbb{P}$ sur $(\mathbb{R%
}^{k},\mathcal{B}(\mathbb{R}^{k}))$\ est tendue, c'est-\`{a}-dire que pour
tout $\varepsilon >0,$ il existe un compact de $\mathbb{R}^{k}$\ tel que

\begin{equation*}
\mathbb{P}(K)\geq 1-\varepsilon .
\end{equation*}
\end{proposition}

\noindent Ce r\'{e}sultat est d\'{e}j\`{a} montr\'{e} dans (\ref{tensTindiv}%
). L'enjeu se trouve dans la tension uniforme ou tension asymptotique.

\subsection{Tension asymptotique. Tension uniforme}

\begin{definition} \label{tensTensDef1}

\noindent (a) Une suite de probabilit\'{e}s $\{\mathbb{P}_{n},n\geq 1\}$\
sur $(\mathbb{R}^{k},\mathcal{B}(\mathbb{R}^{k}))$ est dite asymptotiquement
tendue ou uniform\'{e}ment tendue ssi pour tout $\varepsilon >0,$ il existe
un compact $K$ de \ $\mathbb{R}^{k}$\ \ tel que 
\begin{equation}
\inf_{n\geq 1}\mathbb{P}_{n}(K)\geq 1-\varepsilon   \label{tensTensUnif}
\end{equation}

\noindent ou de mani\`{e}re \'{e}quivalente%
\begin{equation}
\liminf_{n\rightarrow \infty }\mathbb{P}_{n}(K)\geq 1-\varepsilon .
\label{tensTensAsymp}
\end{equation}

\noindent (b) Une suite de variables al\'{e}atoires $\{X_{n},n\geq 1\}$\ 
\`{a} valeurs dans \ $(\mathbb{R}^{k},\mathcal{B}(\mathbb{R}^{k}))$ est
asymptotiquement tendue ou uniform\'{e}ment tendue ssi la famille de
probabilit\'{e}s $\{\mathbb{P}_{X_{n}},n\geq 1\}$ est asymptotiquement
tendue ou uniform\'{e}ment tendue, c'est-\`{a}-dire, pour tout $\varepsilon
>0,$ il existe un compact $K$ de \ $\mathbb{R}^{k}$\ \ tel que 
\begin{equation*}
\inf_{n\geq 1}\mathbb{P}(X_{n}\in K)\geq 1-\varepsilon 
\end{equation*}

\noindent ou de mani\`{e}re \'{e}quivalente%
\begin{equation*}
\liminf_{n\rightarrow \infty }\mathbb{P}(X_{n}\in K)\geq 1-\varepsilon .
\end{equation*}

\noindent (c) Une suite $\{F_{n},n\geq 1\}$\ de fonctions de r\'{e}partition
sur $(\mathbb{R}^{k},\mathcal{B}(\mathbb{R}^{k}))$ est asymptotiquement
tendue ou uniform\'{e}ment tendue ssi la suite des mesures de
Lebesgues-Stieljes $\{\mathbb{P}_{n},n\geq 1\}$ d\'{e}finie par ses \'{e}l%
\'{e}ments est asymptotiquement tendue ou uniform\'{e}ment tendue, c'est-%
\`{a}-dire, pour tout $\varepsilon >0,$ il existe un compact $K$ de $\mathbb{R}^{k}$\ \ tel que 
\begin{equation*}
\inf_{n\geq 1}\mathbb{P}(X_{n}\in K)\geq 1-\varepsilon 
\end{equation*}

\noindent ou de mani\`{e}re \'{e}quivalente, si et seulement si, pour tout $%
\varepsilon >0,$ il existe $0<C\in \mathbb{R}^{k}$ et il existe $c>0$\ tels que
\begin{equation*}
\inf_{n\geq 1}F_{n}(C)\geq 1-\varepsilon 
\end{equation*}

\noindent et
\begin{equation*}
\inf_{n\geq 0}\mathbb{P}_{n}(L_{c})\leq \varepsilon,
\end{equation*}

\noindent si et seulement si, pour tout $\varepsilon >0,$ il existe $0<c$
tel que pour $c^{(k)}=(c,...c),$ 
\begin{equation*}
\inf_{n\geq 1}F_{n}(c^{(k)})\geq 1-\varepsilon
\end{equation*}

\noindent et
\begin{equation*}
\inf_{n\geq 0}\mathbb{P}_{n}(L_{c})\leq \varepsilon,
\end{equation*}
\end{definition}

\bigskip \noindent Il vient de la proposition \ref{tensTensprop1} que l'\'{e}quivalence entre la tension uniforme (\ref{tensTensUnif}) et la tension (\ref{tensTensAsymp}) sont identiques dans le cas sp\'ecifique de $\mathbb{R}^k$. Pour cette raison, nous parlerons uniquement de
famille tendue de mesure de probabilit\'{e}s ou de variables al\'{e}atoires.\newline

\bigskip \noindent Avant d'en venir au th\'{e}or\`{e}me de Helly-Bray,
donnons trois importantes propri\'{e}t\'{e}s de la tension.

\bigskip

\subsection{Tension et transformation continue.}

La tension d'une famille de variables al\'{e}atoires est conserv\'{e}e par
transformation continue. Nous avons :

\begin{proposition} \label{tensTensprop5} Soit une suite de variables al\'{e}atoires $%
\{X_{n},n\geq 1\}$\ \`{a} valeurs dans \ $(\mathbb{R}^{k},\mathbb{B}(\mathbb{%
R}^{k}))$ et $g:\mathbb{R}^{k}\longmapsto \mathbb{R}^{m},$ $m\geq 1,$ une
fonction continue. Alors, si $\{X_{n},n\geq 1\}$ est tendue, alors $%
\{g(X_{n}),n\geq 1\}$ est tendue.
\end{proposition}

\noindent \textbf{Preuve}. Soit \bigskip $\{X_{n},n\geq 1\}$ tendue et $g:%
\mathbb{R}^{k}\longmapsto \mathbb{R}^{m}$. Pour tout $\varepsilon >0,$ il
existe un compact $K$ de \ $\mathbb{R}^{k}$ tel que 
\begin{equation}
\inf_{n\geq 1}\mathbb{P}(X_{n}\in K)\geq 1-\varepsilon .
\end{equation}%
Mais $(X_{n}\in K)\subset (g(X_{n})\in g(K)),$ o\`{u} 
\begin{equation*}
K_{0}=g(K)=\{g(x),x\in K\}
\end{equation*}

\noindent est l'image directe de $K$ par $g$, qui est un ensemble compact.
En effet, toute suite de $K_{0}$ poss\`{e}de une sous-suite convergente dans 
$K_{0}.$ Pour le voir soit $\{g(x_{n}),x_{n}\in K,n\geq 1\}$ une suite dans $%
K_{0}.$ Puisque $K$ est compact, la suite $(x_{n})_{n\geq 1},$ qui est dans $%
K,$ poss\`{e}de une sous-suite convergente dans $K$, soit $%
x_{n(k)}\rightarrow x\in K.$ Puisque $g$ est continue, $g(x_{n(k)})%
\rightarrow g(x)\in K_{0}.$ Il s'en suit que $K_{0}$ est un compact de $%
\mathbb{R}^{m}$\ et%
\begin{equation*}
\inf_{n\geq 1}\mathbb{P}(g(X_{n})\in K_{0})\geq \inf_{n\geq 1}\mathbb{P}%
(X_{n}\in K)\geq 1-\varepsilon .
\end{equation*}

\noindent Ce qui finit la preuve.\newline

\subsection{Caract\'{e}risation de la tension de vecteurs par celles des
composantes.}

\noindent Dans le cas particulier de $\mathbb{R}^k$, nous avons cette caract\'{e}risation :

\begin{proposition}
\label{tensTensprop6} Soit une suite de variables al\'{e}atoires $%
\{X_{n},n\geq 1\}$\ \`{a} valeurs dans \ $(\mathbb{R}^{k},\mathbb{B}(\mathbb{%
R}^{k})).$ Alors $\{X_{n},n\geq 1\}$ est tendue \ si et seulement les suites
des composantes $\{X_{n}^{(i)},n\geq 1\}$, $1\leq i\leq k$, sont tendues.
\end{proposition}

\noindent \textbf{Preuve}. Soit une suite de variables al\'{e}atoires $%
\{X_{n},n\geq 1\}$ \`{a} valeurs dans \ $(\mathbb{R}^{k},\mathbb{B}(\mathbb{R%
}^{k}))$.\newline

\noindent Supposons qu'elle est tendue. Par la proposition \ref{tensTensprop5}, chaque suite de composantes $\{X_{n}^{(i)},n\geq 1\}=\{\pi
_{i}(X_{n}),n\geq 1\}$, $1\leq i\leq k$, est une transformation continue par
la $i$-i\`{e}me projection $\pi _{i}$. Elle est donc tendue.\newline

\noindent Supposons que chaque de composante $\{X_{n}^{(i)},n\geq 1\},$ $%
1\leq i\leq k$, est tendue. Donc, pour tout \ $\varepsilon >0$, pour tout $%
1\leq i\leq k$, il existe $A_{i}>0$ such%
\begin{equation*}
\inf_{n\geq 1}\mathbb{P}(X_{n}^{(i)}\in \lbrack -A_{i},A])\geq 1-\varepsilon
/k.
\end{equation*}

\noindent En posant $A=(A_{1},...,A_{k}),$ nous avons $A>0$ et puisque%
\begin{equation*}
\bigcap\limits_{i=1}^{k}\left( X_{n}^{(i)}\in \lbrack -A_{i},A]\right)
=\left( X_{n}\in \prod\limits_{i=1}^{k}[-A_{i},A_{i}]\right) ,
\end{equation*}

\noindent il en r\'{e}sulte que, pour tout $n\geq 1,$%
\begin{eqnarray*}
\mathbb{P}\left( X_{n}\notin \prod\limits_{i=1}^{k}[-A_{i},A_{i}]\right) &=&%
\mathbb{P}\left( \bigcup\limits_{i=1}^{k}\left( X_{n}^{(i)}\notin \lbrack
-A_{i},A]\right) \right) \\
&\leq &\sum\limits_{i=1}^{k}P\left( X_{n}^{(i)}\notin \lbrack
-A_{i},A]\right) \leq \varepsilon ,
\end{eqnarray*}

\noindent si bien que pour tout $n\geq 1,$%
\begin{equation*}
\mathbb{P}\left( X_{n}\in K_{A}\right) \geq 1-\varepsilon .
\end{equation*}

\noindent La suite $\{X_{n},n\geq 1\}$ est tendue. La preuve est achev\'{e}e.

\subsection{Une suite convergente vaguement est tendue}

\noindent Nous avons ce r\'{e}sultat int\'{e}ressant qui fonde la th\'{e}%
orie de la tension. Il repose sur la tension de la limite, qui est une
probabilit\'{e} sur $\mathbb{R}^{k}$, donc tendue par la proposition \ref%
{tensTensprop3}.

\begin{proposition}
\label{tensTensprop7} Soit une suite de variables al\'{e}atoires $%
\{X_{n},n\geq 1\}$\ \`{a} valeurs dans \ $(\mathbb{R}^{k},\mathbb{B}(\mathbb{%
R}^{k}))$ convergente vaguement vers $X$. Alors elle est tendue.
\end{proposition}

\noindent \textbf{Preuve}. Utilisons la tension de $\mathbb{P}_{X}$. Pour tout $\varepsilon >0,$ il existe un compact $K_{A}=[-A,A]$ de $\mathbb{R}^k$ \ tel que \begin{equation*}
\mathbb{P}(X\in K)\geq 1-\varepsilon .
\end{equation*}

\noindent Soit $0<\delta <1$ et pose $A+\delta =(A_{1}+\delta
,...,A_{k}+\delta )$ et 
\begin{equation*}
\overset{o}{K}_{A+\delta }=\prod\limits_{i=1}^{k}]-A_{i}-\delta
,A_{i}+\delta \lbrack .
\end{equation*}

\noindent Puisque $X_{n}\rightsquigarrow X$ et que $\overset{o}{K}_{A+\delta
}$ est ouvert, nous pouvons utiliser le point $(ii)$ du th\'{e}or\`{e}me
Portmanteau pour tout $0<\delta <1$%
\begin{equation*}
\liminf_{n\rightarrow \infty }\mathbb{P}(X_{n}\in K_{A+1})\geq \liminf_{n\rightarrow \infty }P(X_{n}\in \overset{o}{K}_{A+\delta })\geq 
\mathbb{P}(X\in \overset{o}{K}_{A+\delta }),
\end{equation*}

\noindent donc pour tout $0<\delta <1,$%
\begin{equation*}
\liminf_{n\rightarrow \infty }\mathbb{P}(X_{n}\in K_{A+1})\geq \mathbb{P}%
(X\in \overset{o}{K}_{A+\delta }),
\end{equation*}

\bigskip \noindent Puisque $K$ est ferm\'{e}, $\overset{o}{K}_{A+\delta }\downarrow 
\overline{K}=K$ quand $\delta \downarrow 0$ et alors%
\begin{equation*}
\liminf_{n\rightarrow \infty }\mathbb{P}(X_{n}\in K_{A+1})\geq \mathbb{P}%
(X\in K)\geq 1-\varepsilon .
\end{equation*}

\noindent Il s'en suit que la suite $\{X_{n},n\geq 1\}$ est asymptotiquement
tendue. Elle a en fait h\'{e}rit\'{e} la tension de $X$.\newline

\bigskip \noindent Enfin, nous pouvons maintenant aborder le th\'{e}or\`{e}%
me fondamental de la tension

\bigskip

\section{Th\'{e}or\`{e}me de compacit\'{e} de Prohorov dans $\mathbb{R}^{k}.$%
}

\noindent Le th\'{e}or\`{e}me suivant est l'inverse de la proposition \ref%
{tensTensprop7}, du point de vue de la convergence de sous-suites.

\begin{theorem} \label{tensTheo2}
(Prohorov - Helly-Bray) Soit une suite tendue de variables al\'{e}atoires $%
\{X_{n},n\geq 1\}$\ \`{a} valeurs dans \ $(\mathbb{R}^{k},\mathbb{B}(\mathbb{%
R}^{k})).$ Alors la suite $\{X_{n},n\geq 1\}$ contient une sous-suite $%
\{X_{n(k)},k\geq 1\}$ convergente vers une variable al\'{e}atoire de
probabilit\'{e} $L=\mathbb{P}_{X}.$
\end{theorem}

\bigskip \noindent Ce th\'{e}or\`{e}me peut \^{e}tre directement d\'{e}montr%
\'{e} comme fait dans les ouvrages de Billingsley \cite{billingsley}et van
der Vaart et Wellner \cite{vaart}. Dans ce texte, nous allons passer par le
th\'{e}or\`{e}me de Helly-Bray comme dans les ouvrages de van dervaart \cite%
{vaart_asymp} et Lo\`{e}ve \cite{loeve}.\newline

\noindent Mais, nous passons par le th\'{e}or\`{e}me de Helly-Bray
ci-dessous.\newline

\begin{theorem} \label{tensTheo1}
(Helly-Bray) Toute suite tendue $\{F_{n},n\geq 1\}$\ de fonctions de r\'{e}%
partition \ d\'{e}finies sur $(\mathbb{R}^{k},\mathcal{B}(\mathbb{R}^{k}))$
contient une sous-suite $\{F_{n(k)},k\geq 1\}$\ convergente vaguement vers
une fonction de distribution F qui n'est pas n\'{e}cessairement une function
de r\'{e}partition.
\end{theorem}

\bigskip \noindent \textbf{Preuve}. Soit une suite $\{F_{n},n\geq 1\}$\ de
fonctions de r\'{e}partition \ d\'{e}finies sur $(\mathbb{R}^{k},\mathbb{B}(%
\mathbb{R}^{k})).$ Notons $\mathbb{Q}^{k}$ l'ensemble des \'{e}lements de $%
\mathbb{R}^{k}$\ \`{a} composantes rationnelles. $\mathbb{Q}^{k}$ est d\'{e}%
nombrable et dense dans $\mathbb{R}^{k}$. Enum\'{e}rons $\mathbb{Q}^{k}$
sous la forme $\mathbb{Q}^{k}=\{q_{1},q_{2},...\}.$ Nous allons maintenant
proc\'{e}der \`{a} plusieurs \'{e}tapes.\newline

\noindent \textbf{Etape 1}. Trouvons $F$ sur $\mathbb{Q}^{k}$. Nous allons
mettre en place la construction classique de la suite diagonale. Nous avons
que $(F_{n}(q_{1}))_{n\geq 1}$ $\subset \lbrack 0,1].$ D'apr\`{e}s la propri%
\'{e}t\'{e} de Bolzano-Weierstrass sur $\mathbb{R}$, cette suite admet une
sous-suite $(F_{1,n}(q_{1}))_{n\geq 1}$ convergente vers $F(q_{1}).$\newline

\noindent Ensuite, nous appliquons la sous-suite de fonctions de r\'{e}%
partition $(F_{1,n})_{n\geq 1}$ \`{a} $q_{2}$ ainsi : $(F_{1,n}(q_{2}))_{n%
\geq 1}$ $\subset \lbrack 0,1].$ Nous en d\'{e}duisons une sous-suite $%
(F_{2,n}(q_{2}))_{n\geq 1}$ convergente vers $F(q_{2}).$ Nous pr\'{e}c\'{e}%
dons ainsi de proche en proche. Nous obtenons des sous-suites $%
(F_{j,n})_{n\geq 1},$ $j=1,2,...$ v\'{e}rifiant \newline

\noindent (a) Pour tout $j\geq 1,$ $(F_{j+1,n})_{n\geq 1}$ est une
sous-suite de chaque $(F_{i,n})_{n\geq 1}1\leq i\leq j$.\newline

\noindent \bigskip (b) Pour tout $j\geq 1,$ pour tout $1\leq j\leq i,$ $%
F_{j,n}(q_{i})\rightarrow F(q_{i})$.

\noindent Maintenant, consid\'{e}rons la suite diagonale $(F_{j,j})_{j\geq
1}.$ Vous pouvez voir avec l'aide d'un tableau simple, comme celui
ci-dessous que : pour tout $i\geq 1$ fix\'{e}, la suite $\{F_{j,j},j\geq i\}$
est une sous-suite de $(F_{i,n})_{n\geq i}$ et donc%
\begin{equation*}
F_{j,j}(q_{i})\rightarrow F(q_{i}).
\end{equation*}

\noindent Pour lire ce tableau, il faut retenir que chaque ligne est une
sous-suite des celles qui la devancent. On voit ainsi que chaque terme
diagonal $\mathbf{F}_{j,j}$ est membre de toutes les lignes de $1$ \`{a} $j$.\\

\begin{equation*}
\begin{tabular}{lllllllllll}
$\mathbf{F}_{1,1}$ & $F_{1,2}$ & $F_{1,3}$ & $F_{1,4}$ & $F_{1,5}$ & $%
F_{1,6} $ & $F_{1,7}$ & $F_{1,8}$ & $F_{1,9}$ & $F_{1,10}$ & ... \\ 
& $\mathbf{F}_{2,2}$ & $F_{2,3}$ & $F_{2,4}$ & $F_{2,5}$ & $F_{2,6}$ & $%
F_{1,7}$ & $F_{1,8}$ & $F_{1,9}$ & $F_{1,10}$ & ... \\ 
&  & $\mathbf{F}_{3,3}$ & $F_{3,4}$ & $F_{3,5}$ & $F_{3,6}$ & $F_{3,7}$ & $%
F_{3,8}$ & $F_{3,9}$ & $F_{3,10}$ & ... \\ 
&  &  & $\mathbf{F}_{4,4}$ & $F_{4,5}$ & $F_{4,6}$ & $F_{4,7}$ & $F_{4,8}$ & 
$F_{4,9}$ & $F_{4,10}$ & ... \\ 
&  &  &  & ... & ... & ... & ... & ... & .. & ... \\ 
&  &  &  &  & $\mathbf{F}_{j,j}$ & $F_{j,j+1}$ & $F_{j,j+2}$ & $F_{j,j+3}$ & 
$F_{j,j+5}$ & ...%
\end{tabular}%
\end{equation*}

\bigskip \noindent Nous concluons que la sous-suite diagonale $\left( F_{j,j}\right)
_{j\geq 1}$ renomm\'{e}e $(F_{n(j)})_{j\geq 1}$ v\'{e}rifie%
\begin{equation*}
\forall q\in \mathbb{Q}^{k},\text{ } F_{n(j)}(q)\rightarrow F(q) \text{
quand }j\rightarrow +\infty.
\end{equation*}

\bigskip \noindent \textbf{Etape 2}. Propri\'{e}t\'{e}s de $F$ sur $\mathbb{Q%
}^{k}.$ Pour tout $(a,b)\in \mathbb{Q}^{k}\times \mathbb{Q}^{k},$ quand $%
j\rightarrow +\infty ,$%
\begin{equation*}
0\leq \Delta _{a,b}F_{n(j)}=\sum\limits_{\epsilon \in
\{0,1\}^{k}}(-1)^{s(\epsilon )}F_{n(j)}(b+\epsilon \ast (a-b))
\end{equation*}
\begin{equation*}
\rightarrow \Delta _{a,b}F=\sum\limits_{\epsilon \in
\{0,1\}^{k}}(-1)^{s(\epsilon )}F(b+\epsilon \ast (a-b))\geq 0,
\end{equation*}

\noindent puisque tous les points $b+\epsilon \ast (a-b)$ sont dans $\mathbb{%
Q}^{k}$. Il s'en suit que $F$ est \`{a} volume positif sur $\mathbb{Q}^{k}$.%
\newline

\noindent $F$ est aussi croissant sur $\mathbb{Q}^{k}$ en h\'{e}ritant la
croissance des $F_{n(j)}$, $j\geq 1$, on $\mathbb{Q}^{k}$.\newline

\noindent \textbf{Etape 3}. D\'{e}finissons $G$ \`{a} sur $\mathbb{J}^{k}=%
\mathbb{R}^{k}\backslash \mathbb{Q}^{k}$ par la formule 
\begin{equation*}
G(x)=\inf \{F(q),q\in \mathbb{Q}^{k},x<q\}\in \lbrack 0,1].
\end{equation*}

\noindent pour $x\in \mathbb{J}^{k}$. Il est \'{e}vident que $G$ est bien d%
\'{e}finie maintenant sur $\mathbb{J}^{k}$. Elle est aussi \'{e}videmment
croissante au sens large.\newline

\noindent \textbf{(a) Montrons que }$G$\textbf{\ est continue \`{a} droite}.
Soit $x\in \mathbb{J}^{k}$ et soit $\varepsilon >0.$ Par d\'{e}finition de
l'infimum, il existe $q\in $ tel que $F(q)<G(x)+\varepsilon .$ Pour tout $%
y\in \mathbb{J}^{k},$ $x<y<q,$ $F(y)\leq G(q)$ and $\varepsilon
>F(q)-G(x)\geq G(y)-G(x).$ Donc%
\begin{equation}
\forall \varepsilon >0,\exists q,x<y<q\Longrightarrow 0\leq
G(y)-G(x)<\varepsilon .  \label{cadG}
\end{equation}

\noindent Alors $G$ est continue \`{a} droite.\newline

\noindent \textbf{(c) Montrons que $F_{n(j)}(x)\rightarrow G(x)$ aux points
de continuit\'{e} }$x$\textbf{\ de }$G.$ \newline

\noindent Soit $x$ un point de continuit\'{e} de $G.$ Pour tout $\varepsilon
>0,$ nous pouvons trouver $(y^{\prime },y^{\prime \prime })$ $\in \mathbb{Q}%
^{k}$ tel que $y^{\prime }<x<y^{\prime \prime }$ et $G(y^{\prime \prime
})-G(y^{\prime })<\varepsilon /2.$ Soit $(y^{\prime },y^{\prime \prime })$ $%
\in \mathbb{Q}^{k}$ tel que  $y^{\prime }<q^{\prime }<x<q^{\prime \prime
}<y^{\prime \prime }.$  Then $\ F(q^{\prime \prime })-F(q^{\prime })\leq
G(y^{\prime \prime })-G(y^{\prime })\leq \varepsilon .$ Next

\begin{eqnarray*}
G(y^{\prime })\leq F(q^{\prime })&=&F_{n(j)}(q^{\prime })\leq \liminf_{j\rightarrow +\infty}F_{n(j)}(x)\\
&\leq& \limsup_{j\rightarrow +\infty} F_{n(j)}(x)\\
&\leq& F_{n(j)}(q^{\prime \prime })=F(q^{\prime \prime })\leq G(y^{\prime \prime }).
\end{eqnarray*}

\noindent Ainsi les nombres 
$$
\liminf_{j\rightarrow +\infty} F_{n(j)}(x)
$$ 

\noindent et 

$$
\limsup_{j\rightarrow +\infty} F_{n(j)}(x)
$$ 

\noindent et $F(x)$ sont tous dans l'intervalle $[G(y^{\prime }),G(y^{\prime \prime })]$ de longueur \'{e}gale au plus \`{a} $\varepsilon$. Cela
qui implique 
\begin{equation}
\max (\left\vert G(x)-\liminf_{j\rightarrow +\infty} F_{n(j)}(x)\right\vert ,\left\vert G(x)-\lim\sup_{j\rightarrow +\infty} F_{n(j)}(x)\right\vert )\leq \varepsilon ,  \label{cadGG}
\end{equation}

\noindent pour tout $\varepsilon >0$. D\`{e}s lors, 
\begin{equation*}
F_{n(j)}(x)\rightarrow G(x)\text{ as }j\rightarrow +\infty .
\end{equation*}

\bigskip \noindent \textbf{(c) Montrons que }$G$\textbf{\ est \`{a} volume positif}.\\

\noindent Pour tout $(a,b)\in \mathbb{J}^{k}\times \mathbb{J}^{k},$ soit $q^{\prime
}\downarrow a$ et $q^{\prime \prime }\downarrow b$ avec $q^{\prime }>a$ and $%
q^{\prime \prime }>b$ and $(q^{\prime },q^{\prime \prime })\in \mathbb{Q}%
^{k}\times \mathbb{Q}^{k}.$ Par limite d\'{e}croissante et par d\'{e}%
finition de $G$,%
\begin{equation}
0\leq \Delta _{q^{\prime },q^{\prime \prime }}F\rightarrow \Delta
_{a,b}G\geq 0.  \label{vpG}
\end{equation}

\noindent \textbf{R\'{e}sum\'{e}}. $G$ est une fonction de distribution sur $%
\mathbb{J}^{k}$\ et $F_{n(j)}(x)\rightarrow G(x)$\\

\bigskip 

\noindent \textbf{Etape 4 (finale)}. \ Maintenant, on peut prolonger $G$ sur $\mathbb{Q}^{k}\ $%
par

\begin{equation*}
G(q)=\inf \{F(x),x\in \mathbb{J}^{k},q<x\}\in \lbrack 0,1],\text{ }q\in 
\mathbb{Q}^{k}.
\end{equation*}

\noindent En retour montrons que $G$ est continue \`{a} droite. Que $G$ soit continue 
\`{a} droite en $x\in \mathbb{J}^{k}$ vient de (\ref{cadG})%
\begin{equation*}
\forall \varepsilon >0,\exists q,x<y<q\Longrightarrow 0\leq
G(y)-G(x)<\varepsilon .
\end{equation*}

\noindent Cette formule est vraie, \`{a} l'origine, pour $y\in \mathbb{J}^{k}.$ Pour l'%
\'{e}tendre, choisissons un $y_{0}\in \mathbb{J}^{k}$\ \ et soit  $q_{0}\in 
\mathbb{Q}^{k}$ tel que  
\begin{equation*}
x<y_{0}<q_{0}<q.
\end{equation*}

\noindent Alors pour tout $z\in \mathbb{R}^{k}$

\begin{equation*}
x<z<y_{0}\Longrightarrow G(z)-G(x)<\varepsilon .
\end{equation*}

\noindent Puisque si  $z\in \mathbb{Q}^{k},$ nous avons%
\begin{equation*}
x<z<y_{0}\Longrightarrow G(z)-G(x)\leq G(y)-G(x)\leq \varepsilon .
\end{equation*}

\noindent Montrons que  $G$ soit continue \`{a} droite en $x\in \mathbb{Q}^{k}.$
Re-utilisons la technique qui a donn\'{e} (\ref{cadG}). Par d\'{e}finition,
pour tout $\varepsilon >0,$ il existe  $y\in \mathbb{J}^{k}$ tel que  $%
G(x)\leq G(y)<G(x)+\varepsilon $ avec  $x<y.$ Nous pouvons trouver $q_{0}\in 
\mathbb{J}^{k}$ tel que  $x<q_{0}<y.$ A partir de l\`{a}, nous reconduisons
es m\^{e}mes m\'{e}thodes pour avoir
  
\begin{equation*}
x<z<q_{0}\Longrightarrow G(z)-G(x)\leq \varepsilon .
\end{equation*}

\noindent Maintenant, en utilisant la continuit\'{e} \`{a} droite et le fait que $G$
est \`{a} volume positif sur $\mathbb{Q}^{k}$\ pour l'\'{e}tendre \`{a} \ $%
\mathbb{R}^{k}$, en proc\'{e}dant de la mani\`{e}re qui a conduit \`{a} (\ref%
{vpG}).\\

\noindent Enfin, pour d\'{e}montrer que $F_{n(j)}(x)\rightarrow G(x)$ aux
points de continuit\'{e}, il suffit de le faire pour $x\in \mathbb{Q}^{k}.$
Il suffit de reconsid\'{e}rer la technique qui a conduit \`{a} (\ref{cadGG})
et de remplacer $(y^{\prime },y^{\prime \prime })$ $\in \mathbb{Q}^{k}$ par $%
(z^{\prime },z^{\prime \prime })$ $\in \mathbb{J}^{k}$ points de continuit%
\'{e} de $G$ avec $z^{\prime }\uparrow x$ and $z^{\prime \prime }\downarrow
x$.\\

\noindent Pouvoir choisir de tels points de continuit\'{e} vient du fait que les
frontieres des ensembles $]-\infty ,z^{\prime }]$ et $]-\infty ,z^{\prime
\prime }]$ sont distincts lorsque $z^{\prime }$ croit strictement et $%
z^{\prime \prime }$ d\'{e}croit strictement. On peut dont les choisir tels
que leur valeur par la mesure $m_{G}(]a,b])=\Delta _{a,b}G$ soient nulles,
ce qui fait des $z^{\prime }$ and $z^{\prime \prime }$ des points de
continuit\'{e}. Nous pensons que le lecteur est d\'{e}j\`{a} familier avec
ce type de raisonnement depuis l'expos\'{e} sur la mesure de l'exposition.\\

\bigskip \noindent \textbf{Remarque sur la preuve}. Nous avons voulu avoir une preuve aussi compl\`{e}te
que possible. L'\'{e}tape 4 est superflue si on peut montrer d\`{e}s l'\'{e}%
tape 2 que la fonction $F$ est continue \`{a} droite sur $\mathbb{Q}^{k}.$
En effet la preuve dans van der Vaart \cite{vaart_asymp} ne permet pas \`{a}
notre avis, de dire que $F$ est continue en $x\in $ $\mathbb{Q}^{k}$.\\

\bigskip \noindent Maintenant passons \`{a} la preuve du th\'{e}or\`{e}me de
Prohorov.\newline

\noindent \textbf{Preuve du th\'{e}or\`{e}me \ref{tensTheo2} de Prohorov}.
Supposons que la suite de fonctions de r\'{e}partition $\{F_{n},n\geq 1\}$
soit tendue, c'est-\`{a}-dire que la suite des mesures de Lebesgues-Stieljes 
$\{\mathbb{P}_{n}(]a,b])=\Delta _{a,b}F_{n},n\geq 1\}$ est tendue. D'apr\`{e}%
s la proposition \ref{tensTensprop2}, pour tout $\varepsilon >0,$ nous avons
un $C>0,$ $C\in \mathbb{R}^{k}$ tel que pour tout $n\geq 1,$%
\begin{equation*}
F_{n}(C)\geq 1-\varepsilon .
\end{equation*}

\noindent D'apr\`{e}s le th\'{e}or\`{e}me \ref{tensTheo1}, il existe une
sous-suite $\left( F_{n(j)}\right) _{j\geq 1}$ de $\left( F_{n}\right)
_{n\geq 1}$ convergente vaguement vers une fonction de distribution $F$
associ\'{e}e \`{a} une mesure $L$ d\'{e}finie par $L(]a,b])=\Delta _{a,b}F$
et born\'{e}e par l'unit\'{e}.\\

\noindent Consid\'{e}rons la  famille $\{C_{h}=C+h^{(k)},$ avec $h>0\}.$ Ces points sont tels que les fronti\`{e}res $\partial ]-\infty, C_{h}]$ sont disjointes. Nous pouvons alors
choisir une suite $C_{h_{p}}$ de sorte que $L(\partial ]-\infty
,C_{h_{p}}])=0$ pour tout $p\geq 1$\ and  $C_{h_{p}}\uparrow (+\infty )^{(k)}
$ quand $p\uparrow +\infty$. Ces points sont donc des points de continuit%
\'{e} de $F$ et d\'{e}passent $C.$ Donc pour tout $p\geq 1$ fix\'e,

\begin{equation*}
F_{n(j)}(C_{h_{p}})\geq 1-\varepsilon .
\end{equation*}

\noindent En faisant $j\rightarrow \infty ,$ nous obtenons
\begin{equation*}
F(C_{h_{p}})=L(]-\infty ,C_{h_{p}})\geq 1-\varepsilon .
\end{equation*}

\noindent On conclut en faisant $p\uparrow +\infty $ puis $\varepsilon
\downarrow 0$ pour obtenir

\begin{equation}
F((+\infty )^{(k)})=1.  \label{condFRPINF}
\end{equation}

\bigskip \noindent Sur un autre plan, pour tout $\varepsilon >0,$ il existe $M>0$,
tel que  
\begin{equation*}
\sup P_{n}(L_{M})\leq \varepsilon .
\end{equation*}

\noindent Nous devons montrer que

\begin{equation}
\lim_{\exists (1\leq i\leq k),x_{i}\rightarrow -\infty }F(x)=0,
\label{condFRMINF}
\end{equation}%

\noindent ce qui revient \`{a} dire que pour tout $\varepsilon >,$ il existe $M>0$ tel
que

\begin{equation*}
\exists (1\leq i\leq k),x_{i}<-M\Longrightarrow F(x)\leq \varepsilon .
\end{equation*}%

\noindent Mais

\begin{equation*}
\exists (1\leq i\leq k),x_{i}<-M,]-\infty ,x]\subset L_{M},
\end{equation*}

\noindent donc, pour tout $n\geq 1$,%
\begin{equation*}
\exists (1\leq i\leq k),x_{i}<-M,\text{ \ }F_{n}(x)\leq \varepsilon .
\end{equation*}

\noindent Maintenant, soit un point $x$ fix\'{e} de sorte que :  $\exists (1\leq i\leq
k),x_{i}<-M.$ Soit $x(h)=x+h^{(k)},$ avec  $0<h<-(M+x_{i}).$ Par la m\'{e}%
thode, que nous croyons maintenant classique pour nos lecteurs, il existe
une suite de points $x(h_{p})$, $p\geq 1,$ qui sont des points de continuit%
\'{e} de  $F$ avec  $h_{p}\downarrow 0.$ Nous avons donc pour tout $p\geq 1$%
, fix\'{e} pour tout $j\geq 1,$

\begin{equation*}
\text{\ }F_{n(j)}(x(h_{p}))\leq \varepsilon .
\end{equation*}

\noindent En faisant $j\rightarrow \infty$, nous aurons

\begin{equation*}
\text{\ }F(x(h_{p}))\leq \varepsilon .
\end{equation*}

\noindent Maintenant, par continuit\'{e} \`{a} droite, nous avons quand $p\uparrow
+\infty ,$%
\begin{equation*}
\text{\ }F(x)\leq \varepsilon .
\end{equation*}

\noindent Nous concluons alors : pour tout $\varepsilon >0,$%
\begin{equation*}
\exists (1\leq i\leq k),x_{i}<-M\Longrightarrow F(x)\leq \varepsilon .
\end{equation*}

\noindent ce qui prouve bien (\ref{condFRMINF}). Ceci avec (\ref{condFRPINF}) montre que $F$
est bien une fonction de r\'{e}partition.
\section{Applications}

\subsection{Th\'{e}or\`{e}me de continuit\'{e} de L\'{e}vy}

\begin{theorem} \label{cv.levy}
Soit une suite de fonctions caract\'{e}ristiques $\psi _{n}$ sur $\mathbb{R}$
convergente vers une fonction $\psi $ continue en zero, alors $\psi $ est
une fonction caract\'{e}ristique.
\end{theorem}

\bigskip \noindent \textbf{PREUVE}. Il faut d'abord remarquer que $\psi (0)=0$ puisque $\psi _{n}(n)$
pour tout $n\geq 1.$\ On peut bien supposer les $\psi _{n}$ sont les
fonctions caract\'{e}ristiques de suites $X_{n},$ c'est-\`{a}-dire que%
\begin{equation*}
\psi _{n}(t)=E(e^{itX_{n}}),t\in \mathbb{R}.
\end{equation*}

\noindent D'apr\`es le lemme \ref{funct.facts.lem01} de la section \ref{funct.facts} du chapitre \ref{funct}, nous avons $\left\vert \sin a\right\vert \leq a$
pour $\left\vert a\right\vert \geq 2.$ Et donc%
\begin{equation*}
1_{(\left\vert \delta x\right\vert >2)}\leq 2(1-\frac{\sin \delta x}{\delta x%
})
\end{equation*}

\noindent et, par l'\'{e}galit\'{e} de droite facile \`{a} montrer, nous obtenons
\begin{equation*}
1_{(\left\vert \delta x\right\vert >2)}\leq 2(1-\frac{\sin \delta x}{\delta x%
})=\frac{1}{\delta }\int_{-\delta }^{\delta }(1-\cos tx)dt.
\end{equation*}

\noindent Appliquons cette formule \`{a} $X_{n}$ pour avoir%
\begin{equation*}
1_{(\left\vert X_{n}\right\vert >2/\delta )}\leq \frac{1}{\delta }%
\int_{-\delta }^{\delta }(1-\cos tX_{n})dt
\end{equation*}

\noindent et en prenant les esp\'{e}rances math\'{e}matiques et en appliquant le th\'eor`eme de Fubini
pour des fonctions int\'{e}grables, nous arrivons \`a
\begin{equation*}
\mathbb{P}(\left\vert X_{n}\right\vert >\frac{2}{\delta })\leq \frac{1}{\delta }%
\int_{-\delta }^{\delta }R_{e}(1-Ee^{itX_{n}})dt.
\end{equation*}

\noindent En appliquant le th\'{e}or\`{e}me de convergence domin\'{e} \`{a} $%
R_{e}(1-Ee^{itX_{n}})\rightarrow R_{e}(1-\psi (t)),$ on a 
\begin{equation}
\liminf_{n\rightarrow \infty }P(\left\vert X_{n}\right\vert >\frac{2}{%
\delta })\leq \frac{1}{\delta }\int_{-\delta }^{\delta }R_{e}(1-\psi (t))dt.
\label{tensionLevy1}
\end{equation}

\noindent La fonction \textit{partie r\'elle}, $R_{e}(\cdot )$ est continue et donc, par hypoth\`{e}se, $%
R_{e}(1-\psi (t))\rightarrow 0$ lorsque $\delta \rightarrow 0.$ D\`{e}s lors
lorsque $\delta \rightarrow 0$ dans (\ref{tensionLevy1})$,$ nous obtenons%
\begin{equation*}
\liminf_{n\rightarrow \infty }P(\left\vert X_{n}\right\vert >\frac{2}{%
\delta })=0.
\end{equation*}

\noindent Cela implique que la suite est uniform\'{e}ment tendue. Donc elle contient
une sous-suite $X_{n_{k}}$ qui converge vers une variable $X.$ Par le th\'{e}%
or\`{e}me de convergence domin\'{e}e, pour tout $t\in R,$%
\begin{equation*}
\psi _{n_{k}}(t)=E(\exp (itX_{n_{k}}))\rightarrow E(\exp (itX))=\psi _{0}(t).
\end{equation*}

\noindent Par l'unicit\'{e} des limites dans $\mathbb{R}$,%
\begin{equation*}
\psi =\psi _{0}
\end{equation*}

\noindent et donc $\psi $ est bien une fonction caract\'{e}ristique.\\

\noindent Passons \`{a} la caract\'{e}risation de la convergence vague par les
fonctions caract\'{e}ristiques.\\

\subsection{Une autre preuve de la caract\'{e}risation de la convergence vague par les fonctions caract\'{e}ristiques.}

\begin{theorem} \label{cv.tension.ConvFC}
Une suite variables al\'{e}atoires $X_{n}$ \`{a} valeurs dans $\mathbb{R}^{k}$
converge vaguement vers le vecueur al\'{e}atoire $X\in \mathbb{R}^{k}$ si et
seulement si pour tout $u\in \mathbb{R}^{k},$ $E(\exp (i<u,X_{n}>)\rightarrow E(\exp
(i<u,X>).$
\end{theorem}

\noindent \textbf{PREUVE}. Le sens direct d\'{e}coule du th\'{e}or\`{e}me de convergence domin\'{e}e. Prouvons le sens indirect et supposons que pour tout $u\in \mathbb{R}^{k}$ 
\begin{equation*}
\psi _{n}(u)=E(\exp (i<u,X_{n}>)\rightarrow E(\exp (i<u,X>)=\psi _{n}(u).
\end{equation*}

\noindent Pour tout $i$ fix\'{e}, $1\leq i\leq k,$ la suite des $i$-i\`{e}me
composantes $X_{n}^{(i)}$ verifie, pour tout $t\in \mathbb{R}$, 
\begin{equation*}
\psi _{n^{(i)}}(t)=\psi _{n}(\underset{i-\text{i\`{e}me place}}{\underbrace{%
0,..,t,..0}})=E(\exp (itX_{n}^{(i)})\rightarrow \psi (\underset{i-\text{i%
\`{e}me place}}{\underbrace{0,..,t,..0}})=\psi ^{(i)}(t).
\end{equation*}

\noindent La fonction $\psi $ est continue en z\'{e}ro puisqu'elle est une fonction caract%
\'{e}ristique. D\`{e}s lors, les fonctions partielles sont aussi continues
en z\'{e}ro. La formule ci-dessus dit que la suite des fonctions caract\'{e}%
ristiques $\psi _{n^{(i)}}(t)=E(\exp (itX_{n}^{(i)}))$ converge vers une
fonction $\psi ^{(i)}(t)$\ continue en z\'{e}ro. Le th\'{e}or\`{e}me pr\'{e}c%
\'{e}dent nous fait conclure que la suite $X_{n}^{(i)}$ est tendue. Ainsi
les suites des composantes de $X_{n}$\ sont tendues. Par la proposition \ref%
{tensTensprop6}, la suite $X_{n}$ est tendue.\\

\noindent Nous pouvons conclure en deux \'etapes.\\

\noindent Etape 1 : Chaque sous-suite de of $X_{n}$ contiens une sous-suite qui converge vaguement vers une certaine mesure de probabilit\'e  $L$. Par hypoth\`ese, $X_{n}$ converge vers $\mathbb{P}_{X}$. Par la caract\'erisation des lois de probabilit\'es dans et par unicit\'e de la limite vague en loi, nous avons donc $L=\mathbb{P}_{X}$. Ainsi, il existe une mesure de probabilit\'e (qui est $L_0=\mathbb{P}_{X}$) telle que toute 
sous-suite  $X_{n}$ contient une sous-suite qui converge vers $\mathbb{P}_{X}$.\\

\noindent Etape 2 : Soit $f : \mathbb{R}^k \mapsto \mathbb{R}$ une fonction continue et born\'ee. Soit une sous-suite $\mathbb{E}f(X_{n_{j}})$, $j\geq 1$, de la suite  $\mathbb{E}f(X_n)$, $n\geq 1$.\\

\noindent Cette sous-suite $X_{n_{j}}$, $j\geq 1$, contient une sous-suite $X_{n_{j_{\ell}}}$, $\ell \geq 1$, qui converge vers $L_0$ quand $\ell \rightarrow +\infty$. Donc $\mathbb{E}f(X_{n_{j_{\ell}}})$ converge vers $\int f dL_0$. Ainsi $A=\int f dL_0$ est un nombre r\'eel tel que toute sous-suite $\mathbb{E}f(X_n)$, $n\geq 1$, contient une sous-suite convergente vers $A$.\\

\noindent Par le crir\`ere de Prohorov (Voir Exercice 4 de la section \ref{funct.sec.1} du chapitre  \ref{funct}), $\mathbb{E}f(X_n)$ converge 
vers $\int f dL_0$. Donc $X_n$ converge vaguement vers $L_0=P_{X}$

 

%% file: asymptotics_cv_03_fr.tex
\chapter{Outils Particuliers pour la Convergence Vague dans $\mathbb{R}$} \label{cv.R}

Ce chapitre se focalise sur des outils sp\'ecifiques \`a la convergence vague de suites de variables al\'eatoires r\'eelles. Pour de telles variables, nous pouvons utiliser les repr\'esentations dites de Renyi par le biais de variables al\'eatoires exponentielles standard ou uniformes standard. De telles representations sont bas\'ees sur les fonctions inverses g\'en\'eralis\'ees sur les quelles le premier chapitre se concentre.\\

\noindent En plus, en relation avec les r\'esultats de la section \ref{cv.CvCp} et du th\'eor\`eme \ref{cv.skorohodWichura} du chapitre \ref{cv},  le traitement de la convergence vague sur le m\^eme espace de probabilit\'e peut devenir une affaire de calculs directs et plus ou moins automatique. Ce chapitre donne des outils dans ce sens.

\section{Inverses g\'en\'eralis\'ees des fonctions monotones} \label{cv.sec2}

\bigskip
\noindent Cette th\'eorie est faite pour les fonctions croissantes et  continues \`a droite. Elle peut aussi \^etre faite pour les fonctions 
d\'ecroissantes et continues \`a gauche. Mais il suffit de la faire pour un des cas de monotonie et d'adpater les d\'efinitions et les r\'esultats \`a l'autre cas.\\

 \noindent Soit $F$ be une fonction croissante (au sens large) et  continue \`a droite  d\'efinie de 
$\mathbb{R}$ vers $\mathbb{R}$. D\'efinissons l'inverse g\'en\'eralis\'ee de $F$  comme suit :

\begin{equation*}
F^{-1}(u)=\inf \{x\in \mathbb{R},F(x)\geq u\}, u\in \mathbb{R}.
\end{equation*}

\bigskip
\noindent En raison de l'importance de cette transformation, dite des quantiles, pour la th\'eorie univari\'ee des valeurs extr\^emes, nous allons exposer une liste de quelque de ses propri\'et\'es importantes. Puisque nous voulons que ces propri\'et\'es soient retenues une fois pour toute, nous allons les donner sous formes de points et de pourvoir les preuves \`a la fin de la liste.\\

\bigskip \noindent \textbf{Point (1)} Pour tout $u\in\mathbb{R}$ et Pour tout $t\in\mathbb{R}$ 
\begin{equation}
F(F^{-1}(u))\geq u \tag{A}
\end{equation}

\noindent et

\begin{equation}
F^{-1}(F(x))\leq x. \tag{B}
\end{equation}

\bigskip \noindent \textbf{Point (2)} Pour tout $(u,t)\in R^{2},$%

\begin{equation}
(F^{-1}(u)\leq t)\Longleftrightarrow (u\leq F(t)) \tag{A}
\end{equation}%

\noindent et

\begin{equation}
(F^{-1}(u)>t)\Longleftrightarrow (u>F(t)) \tag{B}
\end{equation}

\bigskip \noindent \textbf{Point (3)} $F^{-1}$  est croissante et continue \`a gauche.\\

\bigskip \noindent \textbf{Point (4)} D\'efinissons la convergenge vague d'une suite de fonction monotones $((F_n)_{n\geq1}$ vers $F$ quand 
$n\rightarrow +\infty$, not\'e par $F_n \rightsquigarrow F$, si et seulement si :

$$
\forall(x\in C(F)), \text{ } F_n(x)\rightarrow F(x) \text{ as } n\rightarrow +\infty,
$$

\noindent o\`u $C(F)$ d\'esigne l'ensemble des points de continuit\'e de $F$. Nous avons : quand $n\rightarrow +\infty$.

$$
(F_n \rightsquigarrow F) \Rightarrow (F^{-1}_n \rightsquigarrow F^{-1})
$$ 

\bigskip \noindent \textbf{Point (5)} Supposons que $F_n$ et $F$ soient \textbf{les fonctions de distribution de variables al\'eatoires r\'eelles} et 
$F_n \rightsquigarrow F$. Si $F$ est \textbf{continue}, alors nous avons la convergence uniforme

$$
\sup_{ x\in R} |F_n(x)-F(x)| \rightarrow 0 \text{ as } n \rightarrow +\infty
$$ 

\bigskip \noindent \textbf{Point (6)} Une fonction de distribution $F$ sur $\mathbb{R}$ a, au plus, un nombre d\'enombrable de points de discontinuit\'e.\\

\bigskip \noindent \textbf{Point (7)} Soit $\mathbb{P}$ une probabilit\'e quelconque sur $\mathbb{R}$ avec un support $(a,b)$, ce qui signifie que $\mathbb{P}((a,b)^c)=0$ ou ce qui signifie que 
$$
a=\inf \{x, \text{ }\mathbb{P}(]-\infty ,x])>0\} \text{ and } b=\inf \{x,\text{ }\mathbb{P}(]-\infty ,x])=1\}.
$$ 

\noindent Alors, pour $0<\varepsilon <1$, il existe un nombre fini de partition de $(a,b)$, 
$$
a=t_{0}<t_{1}<t_{2}<...<t_{k}<t_{k+1}=b
$$ 

\noindent tels que pour $0<i<k$,

$$
\mathbb{P}(]t_{i},t_{i+1}[)\leq \varepsilon.
$$ 

\noindent Nous pouvons toujours \'etendre les limites \`a

$$
-\infty \leq t_{0}<t_{1}<t_{2}<...<t_{k}<t_{k+1}\leq +\infty
$$ 

\bigskip \noindent puisque $\mathbb{P}(]-\infty ,a[)=0$ et  $\mathbb{P}(]b, +\infty[)=0$.

\bigskip \noindent \textbf{Point (8)} Soient $F$ et $G$ deux fonctions de distributions \`a la fois d\'ecroissante ou \`a la fois croissante. Si aucune d'entre elles n'est d\'eg\'en\'er\'e, alors il existe deux points de continuit\'e de toutes les deux $F$, not\'es $G$ $x_1$ et $x_2$ tels que $x_1 < x_2$ et
$$
F(x_1) < F(x_2) \text{ and } G(x_1) < G(x_2).
$$

\bigskip \noindent \textbf{Point (9)} Soit $F$ une fonction croissante de $\mathbb{R}$ to $[a,b]$, sans hypoth\`ese
de continuit\'e \`a gauche ou \`a droite. Alors, pour $y\in ]a,b[$,
\begin{equation*}
F(F^{-1}(y)-0)\leq y\leq F(F^{-1}(y)+0),
\end{equation*}

o\`u $F(\cdot+0)$ and $F(\cdot-0)$ d\'esignent respectivement la limte \`a droite et la limite \`a gauche de $x$.\\

\noindent Si la fonction est d\'ecroissante, l'inverse g\'en\'eralis\'ee sera d\'efinie par 
\begin{equation*}
F^{-1}(y)=\inf \{x\in R,F(x)\leq y\},y\in (a,b).
\end{equation*}

\noindent Dans ce cas, nous avons pour tout  $y\in (a,b)$

\begin{equation*}
F(F^{-1}(y)+)\leq y\leq F(F^{-1}(y)-)
\end{equation*}

\bigskip \noindent \textbf{B - Preuves des diff\'erents points}.\\

\bigskip \noindent \textbf{Preuves du Point 1}. Partie (A). Posons  

\begin{equation*}
A_{u}=\left\{ x\in \mathbb{R},F(x)\geq u\right\}, u\in \mathbb{R}.
\end{equation*}

\noindent Comme $F^{-1}(u)=\inf A_{n}$, il existe une suite $(x_{n})_{n\geq 0}\in A_{u}$ tels que
\begin{equation*}
\left\{ 
\begin{array}{ccc}
F(x_{n}) & \geq  & u \\ 
&  &  \\ 
x_{n} & \downarrow  & F^{-1}(u)%
\end{array}%
\right. 
\end{equation*}

\noindent Par la continuit\'e \`a droite de $F$, nous avons 
\begin{equation*}
F(F^{-1}(u))\geq u.
\end{equation*}
\bigskip 

\noindent Cela prouve la premi\`ere formule (A). En ce qui concerne la formule (B), considerons $x\in \mathbb{R}$ et posons 
\begin{equation*}
F^{-1}(F(x))=\inf A_{F(x)}
\end{equation*}

\noindent En divisant  $A_{F(x)}$ dans
\begin{equation*}
A_{F(x)}=\left[ -\infty ,x\right[ \cap A_{F(x)}+\left[ x,+\infty \right]
\cap A_{F(x)}
\end{equation*}

\begin{equation*}
=:A_{F(x)}(1)+A_{F(x)}(2)
\end{equation*}%

\bigskip \noindent Par le \textbf{fait 1} \`a la fin de ce paragraphe, nous avons
\begin{equation*}
\inf A_{F(x)}=\min (\inf A_{F(x)}(1),\inf A_{F(x)}(2))
\end{equation*}

\noindent Mais, nous avons

\begin{equation*}
y\in A_{F(x)}(1)\Longrightarrow y\leq x,\text{ then }\inf A_{F(x)}(1)\leq x
\end{equation*}%

\noindent Ensuite, nous avons \'evidemment
\begin{equation*}
\inf  A_{F(x)}(2)=\{y\geq x, F(y)\geq F(x) \}=x.
\end{equation*}

\noindent Ainsi, il vient que

\begin{equation*}
\inf  A_{F(x)}\leq x.
\end{equation*}

\noindent C'est-\`a-dire:

\begin{equation*}
F^{-1}(F(x))=\inf A_{F(x)}\leq x.
\end{equation*}

\noindent Cela met fin \`a la preuve du point 1 .\\

\bigskip \noindent \textbf{Preuves du Point 2}. Il est \'evident que chacune des formules (A) et (B) peut \^etre d\'eriv\'ee l'une de l'autre en prenant les compl\'ementaires. Donc nous pouvons seulement prouver l'une d'elles, par exemple (B). Supposons que $\left( u>F(t)\right)$. Par continuit\'e \`a droite de $F$ en $t$, nous pouvons trouver $\epsilon$ tels que  
\begin{equation*}
u>F(t+\epsilon ).
\end{equation*}

\noindent Maintenant, pour $x\in A_{u}$ nous avons s\^urement 
\begin{equation*}
x>t+\epsilon.
\end{equation*}

\noindent Sinon, nous aurions
\begin{equation*}
x\leq t+\epsilon \Longrightarrow F(x)\leq F(t+\epsilon )<u,
\end{equation*}

\noindent et cela aurait conduit \`a la conclusion que $x \notin A_{u}$, qui est en contradiction avec l'hypoth\`ese. Alors, $x>t+\epsilon$ pour tout $x\in A_{u}$. Cela implique que
\begin{equation*}
\inf A_{u}=F^{-1}(u)\geq t+\epsilon >t.
\end{equation*}%

\noindent Nous avons prouv\'e les sens direct de la premi\`ere formule. Pour prouver les sens indirect, considerons $F^{-1}(u)>t$. Ensuite, supposons que $u>F(t)$ n'est pas v\'erif\'e. Cela implique que $F(t) \geq u$, et par suite $t\in A_{u}$. Puisque 
\begin{equation*}
\inf A_{u}=F^{-1}(u)\leq t, 
\end{equation*}
 
\noindent ce qui est contraire \`a l'hypoth\`ese. Donc $u>F(t)$. 

\bigskip 
\noindent \textbf{Preuves du Point 3}. Nous commen\c{c}ons par \'etablir que $F^{-1}$ est non d\'ecroissante. Nous avons
\begin{equation*}
\forall u\leq u^{\prime },A_{u^{\prime }}\leq A_{u}\Longrightarrow \inf A_{u^{\prime }}\leq \inf A_{u}.
\end{equation*}

\noindent Ceci implique que

\begin{equation*}
F^{-1}(u^{\prime })\geq F^{-1}(u)
\end{equation*}

\noindent Ensuite, nous devons prouver que $F^{-1}$ est continue \`a gauche.  Soit $u\in \mathbb{R}$. Nous avons pour tout $h\geq 0$,
\begin{equation*}
F^{-1}(u-h)\leq F^{-1}(u).
\end{equation*}

\noindent Ainsi
\begin{equation*}
\underset{h\downarrow 0}{\lim }F^{-1}(u-h)\leq F^{-1}(u).
\end{equation*}

\noindent Supposons que 
\begin{equation*}
\underset{h\downarrow 0}{\lim }F^{-1}(u-h)=\alpha <F^{-1}(u).
\end{equation*}

\noindent Nous pouvons trouver $\epsilon >0$ tels que $\alpha +\epsilon <F^{-1}(u)$. Maintenant, pour tout $h\geq 0$,

\begin{equation*}
F^{-1}(u-h)<\alpha +\epsilon.
\end{equation*}

\noindent Par d\'efinition de l'infinimum, il existe $x$  tel que 
\begin{equation*}
F(x)\geq u-h\text{ and }F^{-1}(u-h)<\alpha +\epsilon.
\end{equation*}

\noindent Par la formule (A) du Point 1,

\begin{equation*}
F^{-1}(u-h)<\alpha +\epsilon \Longrightarrow u-h\leq F(\alpha +\epsilon ).
\end{equation*}%

\noindent Ensuite, nous obtenons comme $h\downarrow 0$%
\begin{equation*}
u\leq F(\alpha +\epsilon).
\end{equation*}
\noindent Etant donn\'e que cela est vrai pour tout $\epsilon >0$, nous pouvons faire tendre $\epsilon \downarrow 0$ pour obtenir

\begin{equation*}
u\leq F(\alpha).
\end{equation*}
\noindent Mais, par la formule (B) du Point (2) et en utilisant l'hypoth\`ese, nous arrivons \`a
\begin{equation*}
\left( \alpha <F^{1}(u)\right) \Leftrightarrow \left( F^{-1}(u)>\alpha \right)
\Leftrightarrow \left( u>F(\alpha )\right) .
\end{equation*}

\noindent Ceci est clairement une contradiction. Nous concluons que 
\begin{equation*}
\underset{h\downarrow 0}{\lim }F^{-1}(u-h)=F^{-1}(u).
\end{equation*}

\noindent Et ensuite $F^{-1}$ est  continue \`a gauche.\newline

\bigskip \noindent \textbf{Preuves du Point 4}

\noindent Supposons que $F_{n}\rightsquigarrow F$. Soit $y\in \mathbb{R}$ et soit $%
\varepsilon >0 $. Comme le nombre de points de discontinuit\'e de $F$ est au plus d\'enombrable, nous pouvons trouver un point de continuit\'e $x$ de $F$ dans un intervalle ouvert  $(F^{-1}(y)-\varepsilon, F^{-1}(y))$. Par le Point 2, $(F^{-1}(y)>x)$ est is equivalent \`a $(F(x) < y)$. Comme $x\in C(F)$, $F_{n}(x)\rightarrow F(x)$. Alors pour toute valeur de $n$ assez grande,  nous avons $F_{n}(x)<y$ et alors $x<F_{n}^{-1}(y)$. Nous avons
\begin{equation*}
F^{-1}(y)-\varepsilon \leq x<F_{n}^{-1}(y)
\end{equation*}
\bigskip

\noindent c'est-\`a-dire pour tout $\varepsilon >0$,
\begin{equation*}
F_{n}^{-1}(y)>F^{-1}(y)-\varepsilon .
\end{equation*}
\bigskip

\noindent Faisons $n$ tendre vers $+\infty$ et $\varepsilon$ d\'ecroitre vers $0$, nous obtenons pour tout $y\in \mathbb{R}$,
\begin{equation*}
\liminf_{n\rightarrow \infty }F_{n}^{1}(y)\geq F^{-1}(y).
\end{equation*}

\noindent Maintenant , soit $y$  un point de continuit\'e de $F^{-1}$. Pour tout $y^{\prime}>y$, nous pouvons trouver un point de continuit\'e point $x$ de $F$ tels que 
\begin{equation}
F^{-1}(y^{\prime })<x<F^{-1}(y^{\prime })+\varepsilon.  \label{cv.d2}
\end{equation}

\noindent Par le point 1, $x>F^{-1}(y^{\prime})\Longrightarrow F(x)\geq F(F^{-1}(y))\geq y^{\prime }$. Alors

\begin{equation*}
y<y^{\prime }\leq F(x).
\end{equation*}

\bigskip
\noindent Puisque  $x\in C^{\prime }F),$ $F_{n}(x)\rightarrow F(x)$,nous avons pour les grandes valeurs de $n$, $y<F_{n}(x)$  et par la Formule (A) du Point 2, 
$F_{n}^{-1}(y)\leq x$. En combinant  cela avec (\ref{cv.d2}), nous obtenons
\begin{equation*}
F^{-1}(y^{\prime })\geq x\geq F_{n}^{-1}(y).
\end{equation*}

\noindent Maintenant, soit $n\rightarrow +\infty$ pour obtenir 
\begin{equation*}
\lim \sup_{n\rightarrow \infty }F_{n}^{-1}(y)\leq F^{-1}(y^{\prime }).
\end{equation*}
\bigskip
 
\noindent Maintenant, soit $y^{\prime} \downarrow y$, et nous avons $F^{-1}(y^{\prime
})\downarrow F^{-1}(y)$ par la continuit\'e de $F^{-1}$ en $y$. Nous arrivons \`a  
\begin{equation*}
\lim \sup_{n\rightarrow \infty }F_{n}^{-1}(y)\leq F^{-1}(y).\leq \lim
\inf_{n\rightarrow \infty }F_{n}^{1}(y).
\end{equation*}
\bigskip

\noindent Nous concluons enfin
\begin{equation*}
F^{-1}(y)=\lim \sup_{n\rightarrow \infty }F_{n}^{-1}(y)=\lim
\inf_{n\rightarrow \infty }F_{n}^{-1}(t).
\end{equation*}

\bigskip \noindent \textbf{Preuve du Point 5}. Comme $F$ est croissante, $x$ est un point de discontinuit\'e de $F$ si et seulement si le saut de discontinuit\'e $F(x+)-F(x-)$ est strictement positif. Notons par $D$ l'ensemble de tous les points de discontinuit\'e de  $F$, et pour tout $k\geq 1$, notons par $D_{k}$ l'ensemble des points de discontinuit\'e tels que $F(x+)-F(x-) > 1/k$ et par $D_{k,n}$ l'ensemble des points de discontinuit\'e dans l'intervalle $[-n,n]$ tels que  $F(x+)-F(x-) > 1/k$. Nous allons montrer que $D_{k,n}$ est fini.\\

\noindent Supposons que nous pouvons trouver $m$ points $x_1$, ..., $x_m$ dans $D_{k,n}$. Comme $F$ est croissante, on peut voir que la somme des sauts de discontinuit\'e est inf\'erieure \`a $F(n)-F(-n)$. On peut faire un dessin simple pour $m=3$ et projeter les sauts sur l'axe des ordonn\'ees pour voir cela facilement. Ainsi,
$$
\sum_{1\leq j \leq m} F(x+)-F(x-) \leq F(n)-F(-n).
$$ 
\noindent Puisque chacun de ces sauts d\'epasse $1/k$, nous avons
$$
\sum_{1\leq j \leq m} (1/k) \leq \sum_{1\leq j \leq m} F(x+)-F(x-) \leq F(n)-F(-n).
$$ 
\noindent
$$
m/k  \leq F(n)-F(-n)$$

\noindent c'est-\`a-dire

$$
m \leq k(F(n)-F(-n))
$$
\bigskip

\noindent Nous concluons en disant que nous ne pouvons pas avoir plus que $[k(F(n)-F(-n))]$ points dans $D_{k,n}$, donc $D_{k,n}$ est fini. Comme
$$
D=\cup_{n\geq 1} \cup_{k\geq 1} D(k,n)
$$

\noindent Nous voyons que  $D$ est d\'enombrable. Ceci met fin \`a la preuve.

\bigskip \noindent \textbf{Preuve du Point 6}. Soit $0<\varepsilon <1$. Let $F(t)=\mathbb{P}(]-\infty ,t])$. Ceci est une fonction de r\'epartition, donc v\'erifiant $F(\infty)=0$ et $F(+\infty)=1$. Fixons $0<\varepsilon<1$. Posons $k=[1/\varepsilon]$, o\`u $[t]$ repr\'esente la partie enti\`ere de $t$. Nous avons alors
\begin{equation*}
k\varepsilon \leq 1\leq k\varepsilon +\varepsilon.
\end{equation*}

\noindent Notons 
\begin{equation*}
s_{i}=i\varepsilon \text{, for i=1,...,k and }s_{k+1}=1.
\end{equation*}

\noindent et
 
\begin{equation*}
t_{i}=F^{-1}(s_{i})=\inf \{u,\text{ G(u)}\geq s_{i}\}.
\end{equation*}

\noindent Par le Point 1,
\begin{equation}
F(t_{i})\geq s_{i}. (\label{tens_a}) 
\end{equation}

\noindent Ensuite,pour tout $1\leq i<k,$ 
\begin{equation*}
F(t_{i+1}-)=\lim_{h\downarrow 0}F(t_{i+1}-h).
\end{equation*}

\noindent Par d\'efiniton de $t_{i+1}$, qui est l'infimum des valeurs de $u$ tels que  $F(u) \geq (i+i)\varepsilon$, nous avons s\^urement, 
\begin{equation*}
F(t_{i+1}-h)<(i+1)\varepsilon .
\end{equation*}

\noindent  En allant \`a la limite , nous obtenons
\begin{equation}
B(t_{i+1}-)\leq (i+1)\varepsilon.  \label{tens_b}
\end{equation}

\noindent  En mettant ensemble (\ref{tens_a}) et (\ref{tens_b}), nous avons  
\begin{equation*}
\mathbb{P}(]t_{i},t_{i+1}[)=F(t_{i+1}-)-F(t_{i})\leq (i+1)\varepsilon
-i\varepsilon =\varepsilon ,
\end{equation*}

\noindent pour  $i=1,..,k$. Enfin $i=k$, nous avons $F(t_{k+1})=1$ et   
\begin{equation*}
\mathbb{P}(]t_{k},t_{k+1}[)=1-F(t_{k})\leq 1-k\varepsilon \leq \varepsilon .
\end{equation*}

\noindent  Pour $i=0$, comme $F(t_{0})\geq 0$, nous avons
\begin{equation*}
\mathbb{P}(]t_{i},t_{i+1}[)=F(t_{1}-)-F(t_{0})\leq F(t_{1}-)\leq \varepsilon
\end{equation*}

\bigskip
\noindent Nous venons de prouver que $0\leq i\leq k,$%
\begin{equation*}
\mathbb{P}(]t_{i},t_{i+1}[)=F(t_{i+1}+)-F(t_{i})\leq (i+1)\varepsilon-i\varepsilon=\varepsilon .
\end{equation*}

\bigskip \noindent \textbf{Preuve du Point 7}. Nous allons appliquer le Point 6. Consid\'erons la mesure de probabilit\'e de Lebesgues-Stieljes g\'en\'er\'ee par $F$ et caract\'eris\'ee par  
$$
\mathbb{P}(]u,v])=F(v)-F(u), u  \leq v.
$$

\noindent En particulier, nous avons $\mathbb{P}(]-\infty,v])=F(v)-F(u)$. Fixons $\varepsilon >0$ et consid\'erons une subdivision 
$$
-\infty=t_{0}<t_{1}<t_{2}<...<t_{k}<t_{k+1}=+\infty
$$ 

\noindent tel que pour tout $0\leq j \leq k$,
\begin{equation*}
F(t_{j+1}-)-F(t_{j})=\mathbb{P}_{X}(]t_{i},t_{i+1}[)\leq \varepsilon .
\end{equation*}
\bigskip
\noindent Maintenant, nous voulons prouver la convergence uniforme. Soit $x$ un des  $t_{j}$.  Nous avons 
\begin{equation*}
F_{n}(x)-F(x)\leq \sup_{0\leq j\leq k+1}\left| F_{n}(t_{j})-F(t_{j})\right| .
\end{equation*}
\bigskip

\noindent Tout autre $x$ est dans l'un des intervalles $]t_{j},t_{j+1}[$. Utilisons la croissance de $F$ et $F_{n}$ pour avoir   

\begin{eqnarray*}
F_{n}(x)-F(x)&\leq& F_{n}(t_{j+1}-)-F(x)\\
&\leq &F_{n}(t_{j+1}-)-F(t_{j+1}-)+F(t_{j+1}-)-F(x)\\
&\leq & F_{n}(t_{j+1}-)-F(t_{j+1}-)+F(t_{j+1}-)-F(t_{j})\\
&\leq & \sup_{0\leq j\leq k+1}\left| F_{n}(t_{j}-)-F(t_{j}-)\right| +\varepsilon
\end{eqnarray*}

\noindent et

\begin{eqnarray*}
F(x)-F_{n}(x)&\leq & F(x)-F_{n}(t_{j})\\
&\leq & F(x)-F(t_{j})+F(t_{j})-F_{n}(t_{j})\\
&\leq & F(t_{j+1}-)-F(t_{j})+F(t_{j})-F_{n}(t_{j})\\
&\leq & \sup_{0\leq j\leq k+1}\left|F_{n}(t_{j})-F(t_{j})\right| +\varepsilon
\end{eqnarray*}

\noindent A l'arriv\'ee, nous avons pour tout $x$ diff\'erent des $t_j$, 
\begin{equation*}
\left| F(x)-F_{n}(x)\right| \leq \max (\sup_{0\leq j\leq k+1}\left|
F_{n}(t_{j}-)-F(t_{j}-)\right| ,\sup_{0\leq j\leq k+1}\left|
F_{n}(t_{j})-F(t_{j})\right| )+\varepsilon .
\end{equation*}

\noindent Alors, en notant,

\begin{equation*}
\sup_{x\in \mathbb{R} }\left| F_{n}(x)-F(x)\right|=\left\| F_{n}-F\right\|_{\infty },
\end{equation*}

\noindent nous avons obtenu
\begin{equation*}
\left\| F_{n}-F\right\|_{\infty }\leq  \max (\sup_{0\leq j\leq k+1}\left|
F_{n}(t_{j}-)-F(t_{j}-)\right| ,\sup_{0\leq j\leq k+1}\left|
F_{n}(t_{j})-F(t_{j})\right| )+\varepsilon.
\end{equation*}

\noindent A ce stade , nous avons la conclusion plus g\'en\'erale. Si pour tous les r\'eels $x$, $F_n(x) \rightarrow F(x)$ et $F_n(x-) \rightarrow F(x-)$, alors nous pouvons conclure que

$$
\limsup  \left\| F_{n}-F\right\|_{\infty } \leq  \varepsilon,
$$

\noindent pour un $\varepsilon>0$ arbitraire. Ainsi, nous avons

$$
\lim \left\| F_{n}-F\right\|_{\infty } =0.
$$

\bigskip

\noindent Pour \'etendre cette conclusion pour le cas o\`u $F$ est continue et $F_n(x) \rightarrow F(x)$, nous allons prouver que $F_n(x-) \rightarrow F(x)$ pour tout $x$.\\
\bigskip

\noindent Pour le prouver cela pour un quelconque $x$ fix\'e et avec $0\leq h_p \downarrow 0$ comme $p \uparrow +\infty$. Nous avons pour chaque $n$,

$$
F_n(x-)-F(x) \leq F_n(x)-F(x) \leq |F_n(x-)-F(x)| 
$$

\noindent et

$$
F(x)-F_n(x-) \leq F(x)-F_n(x-h_p) \leq |F(x)-F(x-h_p)|+|F(x-h_p)-F_n(x-h_p)| 
$$

\bigskip \noindent En combinant ces deux points, nous avons 
$$
|F(x)-F_n(x-)| \leq max(|F_n(x-)-F(x)|, |F(x)-F(x-h_p)|+|F(x-h_p)-F_n(x-h_p)|). 
$$

\bigskip

\noindent Maintenant fixons $p$ et soit  $n\rightarrow +\infty$ pour avoir

$$
\limsup_{n\rightarrow +\infty}|F(x)-F_n(x-)| \leq |F(x)-F(x-h_p)|. 
$$

\noindent Finalement, faisons $p\rightarrow +\infty$ pour avoir la conclusion par le fait de la continuit\'e de $F$.\\

\bigskip \noindent \textbf{Preuve du Point 8}. $F$ est d\'eg\'en\'er\'ee si et seulement si il y a  un unique point de croissance, par exemple $a$, en lequel elle pr\'esente un saut de discontinuit\'e.  Il est sous la forme : $F(x)=c_1$ pour $x<a$ et $F(x)=c_2$ pour $x \geq a$, 
avec $c_1 < c_2$. Donc , si $F$ est non d\'eg\'en\'er\'ee, il comporte au moins deux points de croissance.  Par cons\'equent, nous pouvons trouver trois points de continuit\'e de $F$ : $a_1 < a_1 < a_3$ tels que  $F(a_1) < F(a_1) < F(a_3)$. Si $G$ est \'egalement non d\'eg\'en\'er\'ee, on peut trouver aussi trois points de continuit\'e de $G$ : $b_1 < b_1 < b_3$ tels que $F(b_1) < F(b_1) < F(b_3)$. Nous consid\'erons deux cas.\\
\bigskip

\noindent Cas 1. L'intervalle $[a_1, a_3]$ et $[b_1, b_3]$ sont disjoints ou, ont une intersection r\'eduite \`a un point qui forc\'ement un point de continuit\'e. Supposons par exemple que $a_3 \leq b_1$. Prenons $x_1=a_1$ et $x_2=b_3$. Nous avons
$$
F(x_1)<F(a_3)\leq F(b_1) \leq F(b_3)=F(x_2)
$$ 

\noindent et
$$
G(x_1) \leq <G(a_3) \leq G(b_1) < F(b_3)=G(x_2).
$$ 
\bigskip

\noindent Cas 2. L'intervalle $[a_1, a_3]$ et $[b_1, b_3]$ ont une intersection non vide et non r\'eduite \`a un point. Cette intersection poss`ede un int\'erieur non vide. On peut donc y choisir un point $t$ de continuit\'e de $F$ et de $G$. Nous avons alors s\^urment $F(a_1)<F(t)$ ou $F(t)<F(a_3)$. Sinon , nous aurions $F(a_1)=F(a_3)$, ce qui violerait ce qui est ci-dessus. De fa\c{c}on similaire $G(b_1)<G(t)$ ou $G(t)<G(a_3)$.

\noindent Maitenant, prenons $x_1=min(a_1,b_1)$ et $x_x=min(a_3,b_3)$. Discutons :\\

\noindent Si $F(a_1)<F(t)$, nous avons  

$$
F(x_1)\leq F(a_1) < F(t) \leq F(a_3) \leq F(x_2)
$$

\noindent Si $F(t)<F(a_3)$, nous avons  

$$
F(x_1)\leq F(a_1) \leq < F(t) < \leq F(a_3) \leq F(x_2)
$$ 

\bigskip \noindent Nous concluons que 

$$
F(x_1) < F(x_2)
$$ 

\noindent Nous montrons de fa\c{c}on similaire que

$$
G(x_1) < G(x_2)
$$ 
\bigskip

\noindent Par le Point 5, nous savons que les points de discontinuit\'e de $F$ et $G$ sont au plus d\'enombrable. Nous, nous pouvons ajuster $x_1$ et $x_2$ de sorte qu'ils soient \`a la fois des points de continuit\'e de $F$ et de $G$.\\

\bigskip \noindent \textbf{Preuve du Point (9)}. Commen\c{c}ons par le premier dans lequel la fonction $F$ est
croissante. Par d\'{e}finition de l'inverse g\'{e}n\'{e}ralis\'{e}e, nous
avons pour tout $h>0$

\begin{equation*}
F(F^{-1}(y)+h)\geq y
\end{equation*}

et and
\begin{equation*}
F(F^{-1}(y)-h)<y.
\end{equation*}

\noindent En passant \`a la limite quand $h\downarrow 0$, nous obtenons
\begin{equation*}
F(F^{-1}(y)-)\leq y\leq F(F^{-1}(y)+).
\end{equation*}

\noindent Similarement, si $F$ est d\'ecroissante, nous avons pour tout  $h>0,$%
\begin{equation*}
F(F^{-1}(y)+h)\leq y
\end{equation*}

\noindent et 
\begin{equation*}
F(F^{-1}(y)-h)>y.
\end{equation*}

\noindent En passant \`a la limite quand $h\downarrow 0$, nous obtenons

\begin{equation*}
F(F^{-1}(y)+)\leq y\leq F(F^{-1}(y)-), \ \  y \in (a,b).
\end{equation*}

\bigskip \bigskip 

\noindent \textbf{Fact 1}. Soit $A$ et  $B$ deux sous-ensembles disjoint de $\mathbb{R}$. Nous avons
have 
\begin{equation*}
\inf A \cup B = min(\inf A, \inf B). 
\end{equation*}

\noindent En effet, de mani\`ere \'evidente, $\inf A \cup B$ est inf\'erieur \`a $\inf A$ et \`a
$\inf B$, et par suite, nous avons $\inf A \cup B \leq min(\inf A, \inf B)$. A l'inverse, supposons ue nous n'avons pas l'\'egalit\'e, i.e.,

\begin{equation*}
\inf A \cup B < min(\inf A, \inf B). 
\end{equation*}

\bigskip \noindent Dans ce cas, il existerait une suite $(z_n)_{n\leq 0}$ de points $A \cup B$ d\'ecroissante vers 
$\inf (A \cup B)$. Alorsn s\^urement, lorsque $n$ est assez grand, $z_n$ sera strictement plus petit que
chacun des nombres $\inf A$ et $\inf B$. Et pourtant, il soit soit dans $A$ soit dans $B$. Ceci est absurde. Nous concluons que nous avons l'\'egalit\'e.\\

\bigskip \bigskip

\subsection{Applications des fonctions inverses g\'en\'eralis\'ees}

\bigskip
\noindent La premi\`ere application est la repr\'esentation d'une variable al\'eatoire r\'eelle par une variable al\'eatoire uniforme standard $ U \sim \mathcal{U}(0,1)$ associ\'ee \`a la fonction de r\'epartition $H(x)=0$ pour $x <0$, $H(x)=x$ pour $x \in (0,1)$ et $H(x)=1$ pour $x>0$. Nous avons :

\begin{lemma} \label{lemmatooldf10} Soit $F$  une fonction de distribution tels que $F(-\infty)=0$ et $F(+\infty)=1$. Soit $U \sim \mathcal{U}(0,1)$, definie sur un espace de probabilit\'e $(\Omega,\mathcal{A},\mathbb{P})$. Alors $F$ est la fonction de r\'epartition de $X=F^{-1}(U)$.
\end{lemma}

\noindent  \textbf{Preuve}. Nous avons par la formule (A) du point 2 ci-dessus que
$$
\mathbb{P}(X \leq x)=P(X \leq x)=\mathbb{P}(F^{-1}(U) \leq x)=\mathbb{P}(U \leq F(x))=F(x).
$$ 

\noindent Une deuxi\`eme application est cette forme simple du Th\'eor\`eme de Skorohod.
\bigskip

\begin{theorem} \label{lemmatooldf11} Soit $F_n \rightsquigarrow F$, o\`u $F_n$ et $F$ sont des fonction de distribution  tels que $F_n(-\infty)=0$ et $F_n(+\infty)=1$, pour $n\geq 0$, $F(-\infty=0)$ et $F(+\infty)=1$. Alors, il existe un espace de probabilit\'e $(\Omega,\mathcal{A},\mathbb{P})$ portant une suite variables al\'eatoires r\'eelles $X_n$ et une variable al\'eatoire $X$ telles que pour tout $n\geq 0$, $F_n$ est la fonction de r\'epartition de $X_n$, c'est-\`a-dire  
$F_n(.)=\mathbb{P}(X_n\leq .)$, et $F$ est la r\'epartition $X$, i.e. $F(.)=\mathbb{P}(X \leq .)$ et

$$
X_n \rightarrow \ \ p.s \ \ as \ \ n\rightarrow +\infty.
$$  
\end{theorem}
\bigskip

\noindent \textbf{Preuve}. Consid\'erons $([0,1],\mathcal{B}([0,1]),\lambda)$ o\`u $\lambda$ est la mesure de Lebesgues sur  $[0,1]$, qui est une mesure de probabilit\'e. Consid\'erons la fonction identit\'e : $U$ : $([0,1],\mathcal{B}([0,1]),\lambda) \mapsto ([0,1],\mathcal{B}([0,1]),\lambda)$. Alors $U$ suit une loi uniforme standard puisque pour tout $x\in (0,1)$, 

$$
\lambda(U \leq x)=\lambda(U^{-1}(]-\infty,x]))
$$

\noindent o\`u $U^{-1}$ est l'inverse de $U$ et alors $U^{-1}(]-\infty,x])=]-\infty,x]$ . Ainsi
$$
\lambda(U \leq x)=\lambda(U^{-1}(]-\infty,x]))=\lambda(]-\infty,x])=x.
$$
\bigskip

\noindent Considerons $X_n=F_n^{-1}(U)$ et $X=F^{-1}(U)$. Alors chaque $F_n$ est une fonction de r\'epartition de $X_n$ et $F$ est une fonction de r\'epartition de  $X$. Montrons que $X_n$ converge vers $X$ presque s\^urement. En utilisant le point 4 ci-dessus, nous avons $F_n^{-1} \rightsquigarrow F^{-1}$. \\
\bigskip

\noindent Nous avons
$$
1 \geq \lambda(X_n \rightarrow X) = \lambda(\{u\in [0,1], X_n(u) \rightarrow X(u) \})
$$

$$
=\lambda(\{u\in [0,1], F_n^{-1}(u) \rightarrow F^{-1}(u)\})
$$

$$ 
\geq \lambda(\{u \in [0,1], \text{u est un point de continuit\'e de } F \})=1,
$$

\noindent \'etant donn\'e que le compl\'ementaire de $\{u \in [0,1], \text{u est un point de continuit\'e de} F \}$ est d\'enombable et le mesure d'un ensemble par la mesure de Lebesgue est nulle.\\

\section{Repr\'esentation uniform et exponentilles de Renyi} \label{cv.R.renyi}

\noindent Dans cette section, nous nous int\'eressons aux repr\'esentations de statistique d'ordre $X_{1,n}\leq ...\leq X_{n,n}$ , $n\geq 1,$ de $n$
variables al\'eatoires ind\'ependantes $X_{n},...,X_{n}$ de fonction de r\'epartition commune $F$ par rapport \`a des variables 
uniformes ou exponentielles standard.\\

\noindent Dans toute la section, les variables sont d\'efinies sur le m\^eme espace de probabilit\'e $\left( \Omega ,\mathcal{A},\mathbb{P}\right).$\\

\noindent Nous commen\c{c}ons par rappeler la densit\'e de probabilit\'e (\textit{dp}) de la statistique d'ordre issue d'un \'echantillon de \textit{dp} $h$.

\subsection{Densit\'e de probabilit\'e de la statistique d'ordre}

\begin{proposition} \label{cv.R.malmquist01} Soit $Z_{1},Z_{2},...,Z_{n}$ $n$ copies ind\'ependantes d'une variable al\'eatoire $Z$
de loi absol\^ument continue et de \textit{dp} $h$. Alors la statistique d'ordre associ\'ee $(Z_{1,n},Z_{2,n},...,Z_{n,n}$ est de loi absol\^ument continue et sa \textit{dp} jointe est donn\'ee par 
\begin{equation*}
h_{\left( Z_{1,n},...,Z_{n,n}\right) }\left( z_{1,...,}z_{n}\right)
=n!\prod\limits_{i=1}^{n}h\left( z_{i}\right) 1_{\left( z_{1}\leq ...\leq
z_{n}\right) .}
\end{equation*} 
\end{proposition}

\bigskip \noindent \textbf{Preuve}. \\
Supposons que les hypoth\`eses de la proposition soient vraies. Trouvons la densit\'e conjointe des $r$ statistiques d'ordre $Z_{n_{1},n}\leq
Z_{n_{2},n}\leq ...\leq Z_{n_{r},n}$, avec $1\leq r\leq n$, $1\leq n_{1}<n_{2}<...<n_{r}$.\newline

\noindent Puisque la variable al\'eatoire $Z$ est absol\^ument continues, ses observations sont distinctes presque s\^urement et nous avons
$Z_{n_{1},n}<Z_{n_{2},n}<_{.}..<Z_{n_{r},n}$, $p.s$. Soit $dz_{i}$, $1 \leq i \leq r$, assez petits. Alors pour $z_{1}<z_{2}<...<z_{r}$, l'\'ev\`enement 
\begin{equation}
(Z_{n_{i},n}\in ]z_{i}-dz_{i}/2,z_{i}+dz_{i}/2[,\ \ 1\leq i\leq r). \label{cvR.oderStatPos}
\end{equation}

\noindent se r\'ealise avec $n_{1}-1$ observations de l'\'echantillon $Z_{1},...,Z_{n}$ tombant strictement \`a gauche de $z_{1}$, une
dans $]z_{1}-dz_{1}/2,z_{1}+dz_{1}/2[$, $n_{2}-n_{1}-1$ entre $z_{1}+dz_{1}/2$ et  $z_{2}-dz_{1}/2$, une dans  $]z_{2}-dz_{2}/2,z_{1}+dz_{2}/2[$,
etc., et $n-k_{k}$ observations strictement \`a droite de $z_{r}$.\\

\begin{figure}[htbp]
	\centering
		\includegraphics[width=1.00\textwidth]{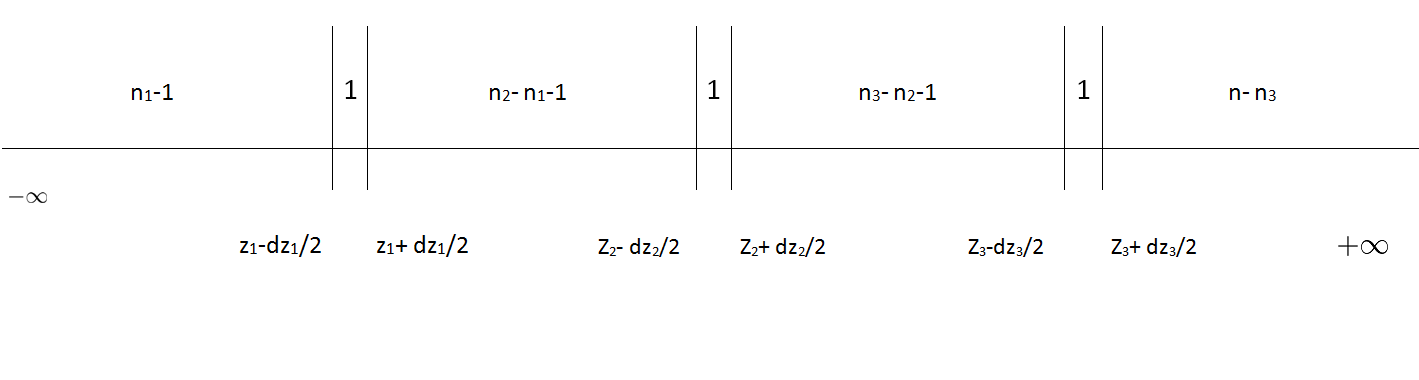}
	\caption{Placement des observations aux points $z_1<...<z_r$ lorsque l'\'ev\`enement (\ref{cvR.oderStatPos}) a lieu}
	\label{cv.R.FigOS}
	\end{figure}

\noindent Cette description est repr\'esent\'ee dans la figure \ref{cv.R.FigOS} pour $r=3$.\newline

\noindent Par d\'efinition, la densit\'e de probabilit\'e $f_{(Z_{n_{1},n},...,Z_{n_{r},n})}$, si elle 
existe, doit v\'erifier

\begin{equation}
\frac{\mathbb{P}(Z_{n_{i},n}\in ]z_{1}-dz_{i}/2,z_{i}+dz_{i}/2[,\
\ 1\leq i\leq r)}{dz_{1}\times ...\times dz_{r}}
=f_{(Z_{n_{1},n},...,Z_{n_{r},n})}(z_{1},...,z_{r})(1+\varepsilon(dz_{1},...,dz_{r})),  \label{cv.densityExtreme02}
\end{equation}

\noindent o\`u $\varepsilon (dz_{1},...,dz_{r})\rightarrow 0$ quand si chaque 
$\Delta _{i}\rightarrow 0$ $(1\leq i\leq r)$. Maintenant, en utilisant l'id\'ependance de l'\'echantillon $Z_{1},...,Z_{n}$, la probabilit\'e 
 $\mathbb{P}(Z_{n_{i},n}\in ]z_{1}-dz_{i}/2,z_{i}+dz_{i}/2[,\ \ 1\leq i\leq r)$ est ibtenue avec par un distribution multinomiale. En utilisant
enfin le fait que $h$ est la densit\'e de probability conjointe des observations, nous obtenons

\begin{eqnarray*}
&&\frac{\mathbb{P}(Z_{n_{i},n}\in ]z_{1}-dz_{i}/2,z_{i}+dz_{i}/2[,\ \ 1\leq i\leq r)}{dz_{1}\times ...\times dz_{r}} \\
&=&n!\times \frac{h(z_{1})^{n_{1}-1}}{(n_{1}-1)!}\times \frac{%
(F(z_{2})-F(z_{1}))^{n_{2}-n_{1}-1}}{(n_{2}-n_{1}-1)!}\\
&\times& ....\times \frac{(F(z_{j})-F(z_{j-1}))^{n_{j}-n_{j-1}-1}}{((n_{j}-n_{j-1}-1)!)}\\
&\times&...\times \frac{(F(z_{r})-F(z_{r-1}))^{n_{r}-n_{r-1}-1}}{(n_{r}-n_{r-1}-1)!}\\
&\times& \frac{(1-F(z_{r}))^{n-n_{r}}}{(n-n_{r})!} \\
&\times& \prod \frac{\mathbb{P}(Z_{n_{i},n}\in ]z_{i}-dz_{i}/2,z_{i}+dz_{i}/2[)}{1!\Delta _{i}}.
\end{eqnarray*}

\noindent Le dernier facteur du produit est donn\'e par

$$
\prod_{i=1}^{r} h(z_{i})(1+dz_{i}).
$$
 
\noindent En posant $n_{0}=0$ et $n_{r}=n+1$ et pour tout  $-\infty
=z_{0}<z_{1}<...<z_{r}<z_{r+1}=+\infty ,$%
\begin{equation*}
f_{(Z_{n_{1},n},...,Z_{n_{r},n})}(z_{1},...,z_{r})=n!\prod_{j=1}^{r+1}\frac{%
h(z_{j})(F(z_{j})-F(z_{j-1}))^{n_{j}-n_{j-1}-1}}{(n_{j}-n_{j-1}-1)!},
\end{equation*}

\noindent nous voyons que $f_{(Z_{n_{1},n},...,Z_{n_{r},n})}$ satisfait (\ref{cv.densityExtreme02}). D\`es lors, nous avons la conclusion 
partielle :

\bigskip 

\begin{lemma} \label{cv.R.malmquist01P} Supposons que $Z_{1},Z_{2},...,Z_{n}$ soient $n$ observations ind\'ependantes
d'une variable al\'eatoire $Z$ de loi absol\^ument continue par rapport \`a la mesure de Lebesgues et densit\'e de probabilit\'e $h$ par
rapport \`a la mesure de Lebesgues et de fonction de r\'epartition $H$, d\'efinies sur le m\^eme espace de probabilit\'e $\left( \Omega ,\mathcal{A},\mathbb{P}\right)$. Soit $1\leq r\leq n$, $1\leq n_{1}<n_{2}<...<n_{r}$. Alors les $r$ statistiques d'ordre $Z_{n_{1},n}<Z_{n_{2},n}<_{.}..<Z_{n_{r},n}$ admettent pour densit\'e de probabilit\'e conjointe par rapport \`a la mesure de Lebesgues, la fonction ci-dessous d\'efinie en  $(z_{1},...,z_{r})$ par

\begin{equation}
n!\prod_{j=1}^{r+1}\frac{h(z_{j})(F(z_{j})-F(z_{j-1}))^{n_{j}-n_{j-1}-1}}{%
(n_{j}-n_{j-1}-1)!}1_{(z_{1}<...<z_{r}),}  \label{cv.densityExtreme03}
\end{equation}

\noindent o\`u, par convention, $n_{0}=0$ et  $n_{r}=n+1,z_{0}=-\infty =z_{0}$ et $z_{r+1}=+\infty .$
\end{lemma}

\bigskip \noindent Pour finir la preuve, posons $r=n$ et $n_{1}=1,n_{2}=2,...,n_{n}=n.$
Puisque les nombres $n_{j}-n_{j-1}-1$ s'annullent dans (\ref{cv.densityExtreme03}), il vient que les statistiques $Z_{1,n}<Z_{2,n}<_{.}..<Z_{n,n}$ ont la densit\'e de probabilit\'e conjointe \'egale \`a 

\begin{equation*}
n!\prod_{j=1}^{n}h(z_{j})1_{(z_{1}<...<z_{r}),}
\end{equation*}

\noindent Dans le reste de la section, nous nous focalisons sur les relations entre statistics d'ordre issues de loi uniforme ou exponentielles standard.

\begin{proposition} \label{cv.R.malmquist02} Soit $n\geq 1$ fix\'e et $U_{1,n}\leq
U_{2,n}\leq ...\leq U_{n,n}$ la statistique d'ordre associ\'ee \`a $%
U_{1},U_{2},...,U_{n}$, $n$ variables al\'eatoires ind\'ependantes uniform\'ement distribu\'ees sur  $\left( 0,1\right)$. Soit $%
E_{1},E_{2},...,E_{n},E_{n+1},$ $\left( n+1\right)$ une suite de $n+1$ variables al\'eatoites independentes suivant une loi exponentielle stantadrd, i.e., \begin{equation*}
\forall x\in \mathbb{R},\text{ \ }\mathbb{P}\left( E_{i}\leq n\right)
=\left( 1-e^{-x}\right) 1_{\left( x\geq 0\right) },i=1,...,n+1.
\end{equation*}

\noindent Soit $S_{j}=E_{1}+...+E_{j},$ $1\leq j\leq n+1$. Alors nous avons l'\'egalit\'e en loi suivante :
\begin{equation*}
\left( U_{1;n},U_{2;n},...,U_{n;n}\right) =^{d}\left( \frac{S_{1}}{S_{n+1}}%
,...,\frac{S_{n}}{S_{n+1}}\right) 
\end{equation*}
\end{proposition}

\bigskip \noindent \textbf{Preuve}. D'une part, d'apr\`es (\ref{cv.R.malmquist01}), la \textit{dp} conjointe de  
$U=(U_{1,n}, U_{2,n}, U_{n,n}$ est donn\'ee par 
\begin{equation*}
\forall \left( u_{1},...,u_{n}\right) \in \mathbb{R}^{n},f_{U}\left(
u_{1},...,u_{n}\right) =n!1_{\left( 0\leq u_{1}\leq ...\leq u_{n}\leq
1\right) }.
\end{equation*}

\noindent Nous allons trouver la distribution de \ $Z_{n+1}^{\ast }=\left(
S_{1},S_{2},...,S_{n},S_{n+1}\right) $ conditionnellement \`a $S_{n+1}=t$, $t>0.$ Nous avons
pour $y=\left( y_{1},y_{2},...,y_{n}\right) \in \mathbb{R}^{n},$%
\begin{eqnarray}
f_{Z_{n}^{\ast }}^{S_{n+1}=t}\left( y\right)  &=&\frac{f_{\left( Z_{n}^{\ast
},S_{n+1}\right) }\left( y,t\right) }{f_{S_{n+1}}\left( t\right) } \label{cv.m01} \\
&=&\frac{f_{Z_{n+1}^{\ast }}\left( y,t\right) }{f_{S_{n+1}}\left( t\right) }%
\boldsymbol{1}_{\left( 0\leq y_{1}\leq ...\leq y_{n}\leq t\right) }. \notag
\end{eqnarray}

\bigskip \noindent Mais $S_{n+1}$ suit une loi gamma de param\`etres  $n+1$ et $1,$ i.e., $S_{n+1}\sim \gamma \left( n+1,1\right) ,$ et sa \textit{dp} est :  
\begin{equation}
f_{S_{n+1}}\left( t\right) =\frac{t^{n}e^{-t}}{\Gamma \left( n+1\right) }%
1_{\left( t\geq 0\right) }=\frac{t^{n}}{n!}e^{-t}\boldsymbol{1}_{\left(
t\geq 0\right) }.  \label{cv.m02}
\end{equation}

\bigskip \noindent Trouvons maintenant la distribution de $\left( S_{1},...,S_{n+1}\right)$. Elle vientde celle de la  transformation
\begin{equation*}
\left( 
\begin{array}{c}
E_{1} \\ 
. \\ 
. \\ 
. \\ 
. \\ 
E_{n+1}%
\end{array}%
\right) =\left( 
\begin{array}{cccccc}
1 &  &  &  &  &  \\ 
-1 & 1 &  &  &  &  \\ 
0 & -1 & 1 &  &  &  \\ 
. &  & . & . &  &  \\ 
. &  &  & . & . &  \\ 
. & . & . & . & -1 & 1%
\end{array}%
\right) \left( 
\begin{array}{c}
S_{1} \\ 
. \\ 
. \\ 
. \\ 
. \\ 
S_{n+1}%
\end{array}%
\right) 
\end{equation*}

\bigskip \noindent 
Soit $B$ la matrice de la formule ci-dessus. Le ajcobien de cette transformation est :  $\left\vert B\right\vert =1$
and 
\begin{equation*}
B\left( 
\begin{array}{c}
y_{1} \\ 
. \\ 
. \\ 
y_{n+1}%
\end{array}%
\right) =\left( y_{1},y_{2}-y_{1},...,y_{n+1}-y_{n}\right) 
\end{equation*}

\bigskip \noindent Donc la \textit{dp} de $\left( S_{1},...,S_{n+1}\right) $ est la suivante

\begin{eqnarray*}
f_{Z_{n+1}^{n}}(y_{1},...,y_{n+1}) &=&f_{\left( E_{1},...,E_{n+1}\right)
}\left( B\left( y_{1},...,y_{n+1}\right) \right) \boldsymbol{1}_{\left(
0\leq y_{1}\leq ...\leq y_{n}\leq y_{n+1}\right) } \\
&=&\prod\limits_{i=1}^{n+1}e^{-\left( y_{i}-y_{i-1}\right) }\boldsymbol{1}%
_{\left( 0\leq y_{1}\leq ...\leq y_{n}\leq y_{n+1}\right) } \\
&=&e^{-y_{n+1}}\boldsymbol{1}_{\left( 0\leq y_{1}\leq ...\leq y_{n}\leq
y_{n+1}\right) }.
\end{eqnarray*}

\noindent avec $y_{0}=0$ par convention. Retournant \`a  (\ref{cv.m01}) et \`a (\ref{cv.m02}), nous ibtenons avec  $y=\left(
y_{1},y_{2},...,y_{n}\right) ,$%
\begin{equation}
f_{Z_{n}^{\ast }}^{S_{n+1}=t}\left( y\right) =\frac{n!}{t^{n}}\boldsymbol{1}%
_{\left( 0\leq y_{1}\leq ...\leq y_{n}\leq t\right) }.  \label{cv.m03}
\end{equation}

\noindent Posons  $u=\left( u_{1},u_{2},...,u_{n}\right) \in \mathbb{R}^{n},$, nous avons
\begin{equation*}
f_{\left( \frac{S_{1}}{S_{n+1}},...,\frac{S_{n}}{S_{n+1}}\right)
}^{S_{n+1}=t}\left( u\right) =f_{\left( \frac{S_{1}}{t},...,\frac{S_{n}}{t}%
\right) }^{S_{n+1}=t}\left( u_{1},u_{2},...,u_{n}\right) .
\end{equation*}

\noindent Ceci est la \textit{dp} obtenue en  (\ref{cv.m03}) par la transformation 
\begin{equation*}
\left( y_{1},y_{2},...,y_{n}\right) =t\left( u_{1},u_{2},...,u_{n}\right)
\Longleftrightarrow \left( u_{1},u_{2},...,u_{n}\right) =\frac{1}{t}\left(
y_{1},y_{2},...,y_{n}\right) 
\end{equation*}

\noindent de jacobien $t^{n}.$ D\`es lors, 
\begin{eqnarray*}
f_{\left( \frac{S_{1}}{S_{n+1}},...,\frac{S_{n}}{S_{n+1}}\right)
}^{S_{n+1}=t}\left( u_{1},u_{2},...,u_{n}\right)  &=&f_{Z_{n}^{\ast
}}^{S_{n+1}=t}\left( t\left( u_{1},u_{2},...,u_{n}\right) \right) t^{n}%
\boldsymbol{1}_{\left( 0\leq tu_{1}\leq ...\leq tu_{n}\leq t\right) } \\
&=&n!\boldsymbol{1}_{\left( 0\leq u_{1}\leq ...\leq u_{n}\leq 1\right) }.
\end{eqnarray*}

\bigskip \noindent Ceci est exactement (\ref{cv.m01}). Nous concluons que la distribution de  $%
Z=\left( \frac{S_{1}}{S_{n+1}},...,\frac{S_{n}}{S_{n+1}}\right)$ conditionnellement \`a $%
S_{n+1}=t$ ne d\'epend pas de $t$. Il s'en suit que cette distribution conditionnelle est la distribution inconditionnelle si bien que 

$$\left( \frac{S_{1}}{S_{n+1}},...,\frac{S_{n}}{S_{n+1}}\right)
$$ 

\noindent poss\`ede la m\^eme loi que  $U=(U_{1,n}, U_{2,n}, ..., U_{n,n})$ et de plus, elle est independante $S_{n+1}$. Ceci finit la preuve.\\

\bigskip \noindent Nous formalisons la derni\`ere conclusion par
\begin{lemma} \label{cv.R.malmquist03} Soit \ $E_{1},E_{2},...,E_{n},E_{n+1},n\geq 1$ $n$ variables al\'eatoires ind\'ependantes de loi
exponentielle standard d\'efinies sur le m\^eme espace de probabilit\'e. Soit \ $S_{i}=E_{1}+E+...+E_{i},1\leq i\leq n+1,$. Alors 
$$
\left( 
\frac{S_{1}}{S_{n+1}},...,\frac{S_{n}}{S_{n+1}}\right)
$$

\noindent est ind\'ependant de $S_{n+1}.$
\end{lemma}

\bigskip \noindent La proposition pr\'ec\'edente exprime une repr\'esentation de la statistique d'ordre issue d'une loi uniforme standard en finction
de celles ussues d'une loi exponentielle standard. La proposition qui suit propose une voie inverse.

\begin{proposition} \label{cv.R.malmquist03} Adoptons les notations de la proposition \ref{cv.R.malmquist02}. Alors pour tout $n\geq 1$, nous avons 
\begin{equation*}
\left( -\log U_{1,n},...,-\log U_{n,n}\right) =^{d}\left(
E_{1,n},...,E_{n,n}\right) ,
\end{equation*}

\noindent o\`u $E_{1,n}\leq ...\leq E_{n,n}$ est la statistique d'ordre associ\'ee \`a $E_{1},E_{2},...,E_{n}$, $n$ variables al\'eatoires ind\'ependentes suivant la loi exponentielle standard.
\end{proposition}

\bigskip \noindent \textbf{Preuve}. Par la proposition \ref{cv.R.malmquist01}, la \textit{dp} de $E_{1,n}\leq
...\leq E_{n,n}$ est
 
\begin{equation}
f_{Z}\left( z\right) =n!e^{-\sum_{i=1}^{n}z_{i}}\boldsymbol{1}_{\left( 0\leq
z_{1}\leq ...\leq z_{n}\right) },%
\text{ }z=\left( z_{1},...,z_{n}\right) \in \mathbb{R}^{n}.  \label{cv.m05}
\end{equation}

\noindent o\`u $Z=\left( E_{1,n},...,E_{n,n}\right)$. La  distribution de $Z^{\ast
}=\left( -\log U_{1,n},...,-\log U_{n,n}\right) $ r\'esulte de celle de $U=\left( U_{1,n},...,U_{n,n}\right)$ par la transformation diff\'eomorphe
$\left(z_{1},...,z_{n}\right) =\left( -\log u_{1},...,-\log u_{n}\right)$ qui pr\'eserve l'ordre des arguments. Le jacobien de cette transformation est donn\'e par : 

\begin{eqnarray*}
\left\vert \frac{\partial U_{i}}{\partial z_{j}}\right\vert  &=&\left\vert
\partial \frac{\partial e^{-z_{i}}}{\partial z_{j}}\right\vert =\left\vert
diag\left( -e^{-z_{1}},...,-e^{-z_{n}}\right) \right\vert  \\
&=&e^{-\sum_{i=1}^{n}z_{i}}.
\end{eqnarray*}%

\noindent D\`es lors, la \textit{pd} de $Z^{\ast }$ est  
\begin{eqnarray*}
f_{Z^{\ast }}\left( z_{1},...,z_{n}\right)  &=&f_{U}\left(
-e^{-z_{1}},...,-e^{-z_{n}}\right) e^{-\sum_{i=1}^{n}z_{i}}\boldsymbol{1}%
_{\left( 0\leq z_{1}\leq ...\leq z_{n}\right) } \\
&=&n!e^{-\sum_{i=1}^{n}z_{i}}\boldsymbol{1}_{\left( 0\leq z_{1}\leq ...\leq
z_{n}\right) }.
\end{eqnarray*}

\noindent Cette \textit{dpf} est exactement celle de $\left( E_{1,n},...,E_{n,n}\right)$ en vertu de (\ref{cv.m05}). Ceci finit la preuve.\\

\bigskip \noindent Donnons une autre version de ce r\'esultat. Il est \'evident que pour toute variable al\'eatoire $U$ suivant la loi
uniforme standard, nous avons $U=^{d}1-U$. D\`es lors, pour tout $n\geq 1$, nous avons
\begin{equation*}
\left( U_{1,n},...,U_{n,n}\right) =^{d}\left(
1-U_{1,n},...,U_{n-i+1,n},...,1-U_{n,n}\right) .
\end{equation*}

\noindent L'\'egalit\'e en loi dans la proposition \ref{cv.R.malmquist03} devient :
pour tout $n\geq 1$, 
$$
\left( -\log \left( 1-U_{n,n}\right) ,...,-\log \left(1-U_{1,n}\right) \right) =^{d}\left( E_{1,n},...,E_{n,n}\right)
$$ 

\bigskip \noindent Allons plus loi et posons 
\begin{equation*}
\alpha _{i,n}=-\log \left( 1-U_{i,n}\right) ,1\leq i\leq n.
\end{equation*}

\bigskip \noindent Pour tout  $n\geq1$, soit la transformation
\begin{equation*}
\left( 
\begin{array}{l}
n\alpha _{1,n} \\ 
\left( n-1\right) \left( \alpha _{2,n}-\alpha _{1,n}\right)  \\ 
. \\ 
. \\ 
\left( n-i+1\right) \left( \alpha _{i,n}-\alpha _{i-1,n}\right)  \\ 
. \\ 
. \\ 
1\left( \alpha _{n,n}-\alpha _{n-1,n}\right) 
\end{array}%
\right) =\left( 
\begin{array}{l}
V_{1} \\ 
V_{2} \\ 
. \\ 
. \\ 
V_{i} \\ 
. \\ 
\\ 
V_{n}%
\end{array}%
\right) .
\end{equation*}

\noindent We have 
\begin{equation*}
\left( 
\begin{array}{l}
\alpha _{1,n} \\ 
\alpha _{2,n} \\ 
. \\ 
. \\ 
. \\ 
. \\ 
\\ 
\alpha _{n,n}%
\end{array}%
\right) =\left( 
\begin{array}{l}
V_{1}/n \\ 
V_{1}/n+V_{2}/\left( n-1\right)  \\ 
. \\ 
. \\ 
. \\ 
. \\ 
\\ 
V_{1}/n+V_{2}/\left( n-1\right) +...+V_{n-1}/2+V_{1}/1%
\end{array}%
\right) 
\end{equation*}

\noindent Par la formule de changement de variable classique, la \textit{dp} de $\left( V_{1},...,V_{n}\right)$ est donn\'ee par  
\begin{eqnarray*}
f_{V}\left( v_{1},...,v_{n}\right)  &=&f_{(\alpha _{1,n},...,\alpha
_{n,n})}\left( v_{1}/n,v_{1}/n+v2/(n-1),...,v_{1}/n+v_{2}/\left( n-1\right)
+...+v_{n}\right)  \\
&&\times \left\vert J\left( v\right) \right\vert \times \boldsymbol{1}%
_{D_{V}}\left( v\right) .
\end{eqnarray*}

\noindent La jacobien de cette transformation est donn\'ee par
\begin{equation*}
\left\vert J\left( v\right) \right\vert =\frac{1}{n!}
\end{equation*}

\noindent et le support de $V$ est 
\begin{equation*}
D_{V}=\mathbb{R}_{+}^{n}
\end{equation*}

\noindent Nous pouvons conclure en utilisant (\ref{cv.m05}) qui donne la \textit{dp} de $(\alpha
_{1,n},...,\alpha _{n,n})$, et en posant $s_{i}=v_{1}/n+v2/(n-1)+...+v_{i}/(n-i+1),$ $i=1,...,n$. Nous obtenons
 
\begin{eqnarray*}
f_{V}\left( v_{1},...,v_{n}\right)  &=&\frac{1}{n!}\times
n!e^{-\sum_{i=1}^{n}s_{i}}\boldsymbol{1}_{\left( v_{1}\geq 0,..,v_{n}\geq
0\right) } \\
&=&e^{-\sum_{i=1}^{n}s_{i}}\boldsymbol{1}_{\left( v_{1}\geq 0,..,v_{n}\geq
0\right) }
\end{eqnarray*}

\noindent Nous pouvons v\'erifier facilement que $s_{1}+...+s_{n}=v_{1}+...+v_{n}$. Tout ceci nous am\`ene \`a 
\begin{equation*}
f_{V}\left( v_{1},...,v_{n}\right) =\prod\limits_{i=1}^{n}e^{-v_{i}}%
\boldsymbol{1}_{\left( v_{i}\geq 0\right) }.
\end{equation*}

\noindent Ceci nous dit que le vecteur $\left( V_{1},...,V_{n}\right)$ est \`a composantes ind\'ependantes suivant chacune une loi exponentielle standard. R\'esumons ce fait par :\\

\begin{proposition} \label{cv.R.malmquist04} Adoptons les notations pr\'ec\'edentes. Soit $\alpha _{i,n}=-\log \left( 1-U_{i,n}\right) ,$ 
$i=1,...,n$, pour $n\geq 1$. Alors, les variables al\'eatoires $n\alpha _{1,n}$, $\left( n-1\right)
\left( \alpha _{2,n}-\alpha _{1,n}\right)$ ,..., $\left( n-i+1\right) \left(
\alpha _{i,n}-\alpha _{i-1,n}\right)$,..., $\left( \alpha _{n,n}-\alpha
_{n-1,n}\right)$ sont ind\'eprendantes et suivent la loi exponentielle standard.
\end{proposition}

\bigskip \noindent Allons plus loin et posons  $1\leq i\leq n$, 
\begin{equation*}
\left( n-i+1\right) \left( \alpha _{i,n}-\alpha _{i-1,n}\right) =\left(
n-i+1\right) \log \left( 
\frac{1-U_{i-1,n}}{1-U_{i,n}}\right) .
\end{equation*}

\noindent De ce qui pr\'ec\`ede, nous avons que les variables al\'eatoires 
\begin{equation*}
E_{n-i+1}^{\ast }=\left( n-i+1\right) \left( \alpha _{i,n}-\alpha
_{i-1,n}\right) =\log \left( \frac{1-U_{n-i,n}}{1-U_{n-i+1,n}}\right)
^{\left( n-i+1\right) }
\end{equation*}

\noindent sont ind\'ependantes et suivent la loi exponentielle standard. Utilisons la repr\'esentation distributionnelle des  $U_{n-i,n}$ par les $U_{i+1,n}$ pour arriver \`a :

\begin{proposition} \label{cv.R.malmquist05} (Malmquist representation). Soit \ $U_{1},U_{2},...,U_{n}$ une suite de $n\geq 1$ variables al\'eatoires ind\'ependantes et suivant toute une loi uniforme standard. Soit $0\leq U_{1,n}<,U_{2,n}<,...,<U_{n,n}\leq 1$ la statistique d'ordre associ\'ee. Alors , les variables al\'eatoires  
\begin{equation*}
\log \left( \frac{U_{i+1,n}}{U_{i,n}}\right) ^{i},i=1,...,n
\end{equation*}

\noindent sont ind\'ependantes et suivent la loi exponentielle standard.
\end{proposition}

 

%% file: asymptotics_cv_04_fr.tex
\chapter[Processus empirique uniforme]{Le processus empirique fonctionnel comme outil g\'{e}n\'{e}ral en statistique asymptotique}

\label{cv.empTool}

\section{Utilisation du petit $o$ et du grand $O$}

Dans ce chapitre, nous montrerons comment ombiner tous les concepts que
nous avons \'{e}tudi\'{e} aussi loin pour obtenir des outils encore simples
mais puissants qui peuvent \^{e}tre systematiquement utilis\'{e}s pour
trouver des lois asymptotiques normales dans une grande vari\'{e}t\'{e} de
probl\`{e}mes, m\^{e}me dans de probl\`{e}mes de recherche d'au jour
d'ajourd'hui. Nous \'{e}tudierons d'abord les manipulations des symboles $o_{%
\mathbb{P}}$ et $O_{\mathbb{P}}$ concernant les limites en probabilit\'{e}.
Ensuite, nous pr\'{e}sentons le processus empirique fonctionnel qui est
utilis\'{e} ici seulement dans le contexe du cas des distributions finies.
Et, nous donnerons quelques cas comme illustrations.\newline

\noindent Il est important de noter pour une fois que la m\'{e}thode donn\'{e}e est
valide pour les suites de variables al\'{e}atoires et la limite de variables
al\'{e}atoires d\'{e}finie sur le m\^{e}me espace de probabilit\'{e}. En
consequence, nous traitons des suites de variables al\'{e}atoires $%
(X_{n})_{n\geq 1},$ $(Y_{n})_{n\geq 1},(Z_{n})_{n\geq 1}$,... d\'{e}finies
sur le m\^{e}me espace de probabilit\'{e} $(\Omega ,\mathcal{A},\mathbb{P})$ 
\`{a} valeurs dans $\mathbb{R}^{k}$, $k\geq 1$ et $(a_{n})_{n\geq 1},$ $%
(b_{n})_{n\geq 1},(c_{n})_{n\geq 1}$ sont des nombres al\'{e}atoires
strictement positifs.\newline

\noindent \textbf{I - Grand $O$ et petit $o$ presque s\^{u}rement}.\newline

\noindent \textbf{DEFINITIONS}.\newline

\noindent \textbf{(a)} La suite de variables al\'{e}atoires r\'{e}elles
\noindent $(X_{n})_{n\geq 1}$ est dite \^{e}tre un o (lire le nom de la
lettre o) de $a_{n}$ presque s\^{u}rement quand n$\rightarrow +\infty $, not%
\'{e} par 
\begin{equation*}
X_{n}=o(a_{n}),a.s.\text{ quand }n\rightarrow +\infty ,
\end{equation*}

\noindent si et seulement si%
\begin{equation}
\lim_{n\rightarrow +\infty }X_{n}/a_{n}=0\text{ }a.s.  \label{eqso}
\end{equation}

\bigskip \noindent \textbf{(b)} La suite de variables al\'{e}atoires r\'{e}%
elles $(X_{n})_{n\geq 1}$ est dite \^{e}tre un grand O de $a_{n}$ presque s%
\^{u}rement quand n$\rightarrow +\infty ,$ not\'{e} par 
\begin{equation*}
X_{n}=O(a_{n}),a.s.\text{ quand }n\rightarrow +\infty ,
\end{equation*}

\noindent si et seulement si la suite $\{\left\vert X_{n}\right\vert
/a_{n},n\geq 1\}$ est presque s\^{u}rement born\'{e}e, c'est \`{a} dire

\begin{equation}
\lim_{n\rightarrow +\infty }\sup \left\vert X_{n}\right\vert /a_{n}<+\infty ,%
\text{ }a.s.  \label{eqbO}
\end{equation}

\bigskip \noindent \textbf{ATTENTION}. Les signes d'\'{e}galit\'{e} utilis%
\'{e}s dans (\ref{eqso}) et (\ref{eqbO}) doivent \^{e}tre lus dans une
direction seulement,  dans le sens suivant : Le membre de gauche est un petit o de $%
a_{n}$\ ou un grand O de $a_{n}$. N'inversons jamais l'\'{e}galit\'{e} de la
gauche \`{a} la droite. Par exemple, si $X_{n}$ est un $o(n)$, c'est aussi
un $o(n^{2})$ et nous pouvons \'{e}crire $o(n)=o(n^{2})$ $p.s.$, mais nous  ne
pouvons pas \'{e}crire $o(n^{2})=o(n)$ $p.s.$ Un exemple : $X_{n}=n^{3/2}$
est un $o(n^{2})$ mais n'est pas un $o(n).$ Cette remarque s'\'{e}tendra aux
notations des petits $o$ et des grands $O$ en probabilit\'{e} qui sont d\'{e}%
finie ci-dessous.\newline

\bigskip \noindent \textbf{Cas particuliers concernant les constantes}. Si $%
a_{n}=C>0$ pour tout $n\geq 1,$ not\'{e} $a_{n}\equiv C,$ nous avons :%
\newline

\bigskip \noindent \textbf{(i)} $X_{n}=O(C)$ p.s. si et seulement si $X_{n}/C
$ est born\'{e} $p.s.$ si et seulement si $X_{n}$ est born\'{e} $p.s..$ et
nous \'{e}crivons  
\begin{equation*}
X_{n}=O(1)\text{ }p.s.
\end{equation*}

\bigskip \noindent \textbf{(ii)} $X_{n}=o(C)$ $p.s.$ si et seulement si $%
X_{n}/C$ $\rightarrow 0$ $p.s.$ si et seulement si $X_{n}/C\rightarrow 0$ $%
p.s.\ $et nous \'{e}crivons 
\begin{equation*}
X_{n}=o(1)\text{ }p.s.
\end{equation*}

\bigskip \noindent \textbf{(iii)} Pour toute constante $C>0$, nous \'{e}%
crivons $C=O(1)$.\newline

\bigskip \noindent \textbf{PROPRIETES}.\newline

\noindent Les propri\'{e}t\'{e}s sont tr\`{e}s nombreuses et les
utilisateurs doivent souvent v\'{e}rifier celle d\'{e}pendant de ses
travaux. Mais certaines d'entre elles doivent \^{e}tre connues et pr\^{e}tes 
\`{a} \^{e}tre utilis\'{e}es. Classons les en trois groupes.\newline

\bigskip \noindent \textbf{Groupe A}. Propri\'{e}ti\'{e}s du petit o.\newline

\noindent \textbf{(1)} $o(a_{n})o(b_{n})=o(a_{n}b_{n})$ $p.s.$.\newline

\noindent \textbf{(2)} (1) $o(o(a_{n}))=o(a_{n})$ $p.s.$\newline

\noindent \textbf{(3)} Si $b_{n}\geq a_{n}$ pour tout $n\geq
1,o(a_{n})=o(b_{n})$ $p.s.$\newline

\noindent \textbf{(4)} $o(a_{n})+o(a_{n})=o(a_{n})$ $p.s.$.\newline

\noindent \textbf{(5)} $o(a_{n})+o(b_{n})=o(a_{n}+b_{n})$ $a.s.$ et $%
o(a_{n})+o(b_{n})=o(a_{n}\vee b_{n})$ $p.s.$ o\`{u} $a_{n}\vee b_{n}=\max
(a_{n},b_{n}).$\newline

\noindent \textbf{(6)} $o(a_{n})=a_{n}o(1)$ $a.s.$ et $a_{n}o(1)=O(a_{n})$ $%
a.s.$\newline

\bigskip \noindent \textbf{PREUVES}. Chacune de ces propri\'{e}ti\'{e}s est
rapidement prouv\'{e}e dans :\newline

\noindent \textbf{(1)} Si \bigskip $X_{n}=o(a_{n})$ et $Y_{n}=o(b_{n}),$
alors 
\begin{equation*}
\lim_{n\rightarrow +\infty }\frac{\left\vert X_{n}Y_{n}\right\vert }{%
a_{n}b_{n}}=\lim_{n\rightarrow +\infty }\frac{\left\vert X_{n}\right\vert }{%
a_{n}}\times \lim_{n\rightarrow +\infty }\frac{\left\vert Y_{n}\right\vert }{%
b_{n}}=0\text{ }p.s
\end{equation*}

\noindent et alors $X_{n}Y_{n}=o(a_{n}b_{n})$ $p.s.$\newline

\noindent \textbf{(2)} Si $Y_{n}=o(a_{n}),$ $p.s.$ et $X_{n}=o(Y_{n}),$ $%
p.s.,$%
\begin{equation*}
\lim_{n\rightarrow +\infty }\frac{\left\vert X_{n}\right\vert }{a_{n}}%
=\lim_{n\rightarrow +\infty }\left\vert \frac{X_{n}}{Y_{n}}\right\vert
\times \frac{\left\vert Y_{n}\right\vert }{a_{n}}=\lim_{n\rightarrow +\infty
}\left\vert \frac{X_{n}}{Y_{n}}\right\vert \times \lim_{n\rightarrow +\infty
}\frac{\left\vert Y_{n}\right\vert }{a_{n}}=0,
\end{equation*}

\noindent c'est \`{a} dire $X_{n}=o(a_{n})$ $p.s.$\newline

\noindent \textbf{(3)} Si \bigskip $X_{n}=o(a_{n})$ et $b_{n}\geq a_{n}$
pour tout n$\geq 1,$ alors 
\begin{equation*}
0\leq \lim_{n\rightarrow +\infty }\sup \frac{\left\vert X_{n}\right\vert }{%
b_{n}}=\lim_{n\rightarrow +\infty }\sup \frac{\left\vert X_{n}\right\vert }{%
a_{n}}\frac{a_{n}}{b_{n}}\leq \lim_{n\rightarrow +\infty }\sup \frac{%
\left\vert X_{n}\right\vert }{a_{n}}=0\text{ }p.s
\end{equation*}

\noindent et $\left\vert X_{n}/b_{n}\right\vert \rightarrow 0$ $p.s.,$ c'est 
\`{a} dire $X_{n}=o(b_{n}),$ $p.s.$\newline

\noindent \textbf{(4)} Si \bigskip $X_{n}=o(a_{n})$ et $Y_{n}=o(a_{n}),$
alors%
\begin{equation*}
\lim_{n\rightarrow +\infty }\sup \frac{\left\vert X_{n}+Y_{n}\right\vert }{%
a_{n}}\leq \lim_{n\rightarrow +\infty }\frac{\left\vert X_{n}\right\vert }{%
a_{n}}+\lim_{n\rightarrow +\infty }\frac{\left\vert Y_{n}\right\vert }{b_{n}}%
=0\text{ }p.s
\end{equation*}

\noindent et donc $X_{n}+Y_{n}=o(a_{n})$ $p.s.$\newline

\noindent \textbf{(5)} Pour prouver que $o(a_{n})+o(b_{n})=o(a_{n}+b_{n})$ $%
p.s.$, utiliser le Point (3) pour voir que $o(a_{n})=o(a_{n}+b_{n})$ $p.s.$
puisque $a_{n}+b_{n}\geq a_{n}$ pour tout $n\geq 1$ et aussi bien que$\
o(b_{n})=o(a_{n}+b_{n})$ $p.s.$et alors utiliser le Point (3) pour conclure.
Nous prouvons que $o(a_{n})+o(b_{n})=o(a_{n}\vee b_{n})$ $p.s.$ de la m\^{e}%
me mami\`{e}re.\newline

\noindent \textbf{(6)} C'est une simple utilisation de d\'{e}finition.%
\newline

\bigskip \noindent \textbf{Groupe B}. Propri\'{e}ti\'{e}s du grand O.\newline

\noindent \textbf{(1)} $O(a_{n})O(b_{n})=O(a_{n}b_{n})$ $p.s.$\newline

\noindent \textbf{(2)} $O(O(a_{n}))=O(a_{n})$ $p.s.$\newline

\noindent \textbf{(3)} Si $b_{n}\geq a_{n}$ pour tout $n\geq
1,O(a_{n})=O(b_{n})$ $p.s.$\newline

\noindent \textbf{(4)} $O(a_{n})+O(a_{n})=O(a_{n})$ $p.s.$\newline

\noindent \textbf{(5)} $O(a_{n})+O(b_{n})=O(a_{n}+b_{n})$ $p.s.$ et $%
O(a_{n})+O(b_{n})=O(a_{n}\vee b_{n})$ $p.s.$ o\`{u} $a_{n}\vee b_{n}=\max
(a_{n},b_{n}).$\newline

\noindent \textbf{(6)} $O(a_{n})=a_{n}O(1)$ $p.s.,$ et $a_{n}O(1)=O(a_{n})$ $%
p.s.$\newline

\bigskip \noindent \textbf{PREUVES}. Ces propri\'{e}t\'{e}s sont prouv\'{e}%
es exactement comme celles du \textbf{Groupe A}, o\`{u} les limites sont
utilis\'{e}es \`{a} la place des limites sup\'{e}rieures.\newline

\noindent \textbf{Groupe C}. Propri\'{e}ti\'{e}s des combinations des petits
o's et grands O's.\newline

\noindent \textbf{(1)} $o(a_{n})O(b_{n})=o(a_{n}b_{n})$ $p.s.$\newline

\noindent \textbf{(2)} $o(O(a_{n}))=o(a_{n})$ $p.s.$ et $%
O(o(a_{n}))=o(a_{n}),$ $p.s.$\newline

\noindent \textbf{(3a)} Si $a_{n}=O(b_{n}),$ $p.s.,$ alors $%
o(a_{n})+O(b_{n})=O(b_{n})$ $p.s.$\newline

\noindent \textbf{(3b)} Si $b_{n}=O(a_{n}),$ $p.s.,$ alors $%
o(a_{n})+O(b_{n})=O(a_{n})$ $p.s.$\newline

\noindent \textbf{(3c)} Si $b_{n}=o(a_{n}),$ $p.s.,$ alors $%
o(a_{n})+O(b_{n})=o(b_{n})$ $p.s.$\newline

\noindent \textbf{(4)} $(1+o(a_{n}))^{-1}-1=o(a_{n})$, $p.s.$\newline

\bigskip \noindent \textbf{PREUVES}.\newline

\noindent \textbf{(1)} Si $X_{n}=o(a_{n})$ et $Y_{n}=O(b_{n}),$ alors%
\begin{equation*}
\lim_{n\rightarrow +\infty }\sup \frac{\left\vert Y_{n}\right\vert }{b_{n}}%
=C<+\infty \text{ }p.s.
\end{equation*}

\noindent et 
\begin{eqnarray*}
\lim_{n\rightarrow +\infty }\sup \frac{\left\vert X_{n}Y_{n}\right\vert }{%
a_{n}b_{n}} &=&\lim_{n\rightarrow +\infty }\left( \left\vert \frac{X_{n}}{%
a_{n}}\right\vert \times \frac{\left\vert Y_{n}\right\vert }{b_{n}}\right)
=\lim_{n\rightarrow +\infty }\sup \left\vert \frac{X_{n}}{a_{n}}\right\vert
\times \lim_{n\rightarrow +\infty }\sup \frac{\left\vert Y_{n}\right\vert }{%
b_{n}} \\
&\leq &C\lim_{n\rightarrow +\infty }\sup \left\vert \frac{X_{n}}{a_{n}}%
\right\vert =0\text{ }p.s.
\end{eqnarray*}

\noindent \textbf{(2)} Utiliser les Points (6) des Groupes A et B pour dire 
\begin{equation*}
o(O(a_{n}))=o(1)\times O(a_{n})=a_{n}\times o(1)\times O(1)=a_{n}\times
o(1)=o(a_{n})
\end{equation*}

\noindent et 
\begin{equation*}
O(o(a_{n}))=o(a_{n})O(1)=a_{n}\times o(1)\times O(1)=a_{n}\times
o(1)=o(a_{n}).
\end{equation*}

\noindent \textbf{(3a-b-c)} Ces trois points sont prouv\'{e}s de mani\`{e}%
res similaires. Donnons les details de (3b) par exemple. Soient $%
X_{n}=o(a_{n})$ et $Y_{n}=O(b_{n})$ et$b_{n}=O(a_{n}).$ On a 
\begin{eqnarray*}
o(a_{n})+O(b_{n}) &=&o(a_{n})+O(O(a_{n}))=o(a_{n})+O(a_{n}) \\
&=&a_{n}(o(1)+O(1))=a_{n}\times O(1)=O(a_{n}).
\end{eqnarray*}

\noindent \textbf{(4)} Nous avons 
\begin{eqnarray*}
(1+o(a_{n}))^{-1}-1 &=&\frac{o(a_{n})}{1+o(a_{n})}=\frac{o(a_{n})}{1+o(a_{n})%
}=o(a_{n})O(1) \\
&=&a_{n}o(1)O(1)=a_{n}o(1)=a_{n}o(a_{n}).
\end{eqnarray*}

\bigskip \noindent \textbf{II - Grand $O$ et petit $o$ en probabilit\'{e}}.%
\newline

\noindent \textbf{DEFINITIONS}.\newline

\noindent \textbf{(a)} La suite de variables al\'{e}atoires r\'{e}elles
\noindent $(X_{n})_{n\geq 1}$ est dite \^{e}tre un o (lire le nom de la
lettre o) de $a_{n}$ en probabilit\'{e} quand n$\rightarrow +\infty ,$ not%
\'{e} par%
\begin{equation*}
X_{n}=o_{\mathbb{P}}(a_{n}),\text{ quand }n\rightarrow +\infty ,
\end{equation*}

\noindent si et seulement si%
\begin{equation*}
\lim_{n\rightarrow +\infty }X_{n}/a_{n}=0
\end{equation*}

\noindent en probabilit\'{e}$,$

\noindent c'est \`{a} dire pour tout $\lambda >0$ 
\begin{equation*}
\lim_{n\rightarrow +\infty }\mathbb{P}(\left\vert X_{n}\right\vert >\lambda
a_{n})=0.
\end{equation*}

\noindent \textbf{(b)} La suite de variables al\'{e}atoires r\'{e}elles $%
(X_{n})_{n\geq 1}$ est dite \^{e}tre un grand O de $a_{n}$ en probabilit\'{e}
quand n$\rightarrow +\infty ,$ not\'{e} par

\begin{equation*}
X_{n}=O_{\mathbb{P}}(a_{n}),\text{ quand }n\rightarrow +\infty ,
\end{equation*}

\noindent si et seulement si la suite $\{\left\vert X_{n}\right\vert
/a_{n},n\geq 1\}$ est born\'{e}e en probabilit\'{e}, c'est \`{a}dire : Pour
tout $\varepsilon >0,$ il existe une constante $\lambda >0$, telle que%
\begin{equation}
\inf_{n\geq 1}\mathbb{P}(\left\vert X_{n}\right\vert \leq \lambda a_{n})\geq
1-\varepsilon   \label{big-O1}
\end{equation}

\noindent ce qui est \'{e}quivalente \`{a}

\begin{equation}
\liminf_{\lambda \uparrow +\infty }\limsup_{n\rightarrow +\infty }\mathbb{P}%
(\left\vert X_{n}\right\vert >\lambda a_{n})=0.  \label{big-O2}
\end{equation}

\noindent Avant d'aller plus loin, prouvons ceci :

\begin{lemma}
\label{oO.01} Chacune des (\ref{big-O1}) et (\ref{big-O2}) est \'{e}%
quivalente \`{a} : Pour tout $\varepsilon >0,$ il existe un entier $N\geq 1$
une constante $\lambda >0$, telle que%
\begin{equation}
\inf_{n\geq N}\mathbb{P}(\left\vert X_{n}\right\vert \leq \lambda a_{n})\geq
1-\varepsilon .  \label{big-O3}
\end{equation}
\end{lemma}

\bigskip \noindent \textbf{PREUVES}. Pour prouver cela pour (\ref{big-O1})
et (\ref{big-O3}), il sufira de prouver que (\ref{big-O3}) $\Longrightarrow $
(\ref{big-O1}) puisque l'implication inverse est \'{e}vidente. Supposons que
(\ref{big-O3}), c'est \`{a} dire pour $\varepsilon >0,$ ils existent $N\geq 1
$ et un nombre r\'{e}el $\lambda _{0}>0$ tels que

\begin{equation}
\forall (n\geq N),\text{ }\mathbb{P}(\left\vert X_{n}/a_{n}\right\vert \leq
\lambda _{0})\geq 1-\varepsilon .
\end{equation}

\noindent Si $N=1,$ alors (\ref{big-O1})\ est v\'{e}rifi\'{e}e. Si non, nous
avons $j\in \{1,...,N-1\}$ fix\'{e}$,$ $(\left\vert X_{j}/a_{j}\right\vert
\leq \lambda )\uparrow \Omega $ quand $\lambda \uparrow +\infty .$ Donc par
le Th\'{e}or\`{e}me de la Convergence Monotone, il existe pour chaque $j\in
\{1,...,N-1\}$ un nombre r\'{e}al $\lambda _{j}>0$ tel que $\mathbb{P}%
(\left\vert X_{j}/a_{j}\right\vert \leq \lambda _{j})>1-\varepsilon $. Nous
pr\'{e}nons $\lambda =\max (\lambda _{0},\lambda _{1},...,\lambda _{N-1})$
pour obtenir%
\begin{equation*}
\forall (n\geq 1),\text{ }\mathbb{P}(\left\vert X_{n}/a_{n}\right\vert \leq
\lambda )\geq 1-\varepsilon ,
\end{equation*}

\noindent ce qui est (\ref{big-O1}). Maintenant, prouvons que (\ref{big-O2})$%
\Longleftrightarrow $(\ref{big-O3}). D'abord (\ref{big-O2}) veut dire%
\begin{equation*}
\lim_{\lambda \uparrow +\infty }\limsup_{n\rightarrow +\infty }\mathbb{P}%
(\left\vert X_{n}\right\vert >\lambda a_{n})=0,
\end{equation*}

\noindent puisque $\limsup_{n\rightarrow +\infty }\mathbb{P}(\left\vert
X_{n}\right\vert >\lambda a_{n})$ est d\'{e}croissante quand $\lambda
\uparrow +\infty $ sur $[0,1].$ Nous obtenons pour tout $\varepsilon >0,$ il
existe un nombre r\'{e}el $\lambda >0$ tel que 
\begin{equation*}
\limsup_{n\rightarrow +\infty }\mathbb{P}(\left\vert X_{n}\right\vert
>\lambda a_{n})=\lim_{N\uparrow +\infty }\sup_{n\geq N}\mathbb{P}(\left\vert
X_{n}\right\vert >\lambda a_{n})\leq \varepsilon /2.
\end{equation*}

\noindent Alors pour un quelconque $N>0$, 
\begin{equation*}
\sup_{n\geq N}\mathbb{P}(\left\vert X_{n}\right\vert >\lambda a_{n})\leq
\varepsilon ,
\end{equation*}

\noindent c'est \`{a} dire 
\begin{equation*}
\inf_{n\geq N}\mathbb{P}(\left\vert X_{n}\right\vert \leq \lambda a_{n})\geq
1-\varepsilon ,
\end{equation*}

\noindent ce qui est (\ref{big-O3}). Qui \`{a} son tour, la reformulation de
ceci donne : pour tout $\varepsilon >0,$ il existe $N_{0}>0$ et un nombre r%
\'{e}el $\lambda _{0}>0$ tel que

\begin{equation}
\inf_{n\geq N}\mathbb{P}(\left\vert X_{n}\right\vert \leq \lambda
_{0}a_{n})\geq 1-\varepsilon ,
\end{equation}

\noindent c'est \`{a} dire

\begin{equation*}
\sup_{n\geq N}\mathbb{P}(\left\vert X_{n}\right\vert >\lambda
_{0}a_{n})<\varepsilon ,
\end{equation*}

\noindent ce qui conduit \`{a} 
\begin{equation*}
\inf_{N\geq 1}\sup_{n\geq N}\mathbb{P}(\left\vert X_{n}\right\vert >\lambda
_{0}a_{n})<\varepsilon ,
\end{equation*}

\noindent et ensuite 
\begin{equation*}
\inf_{\lambda >0}\inf_{N\geq 1}\sup_{n\geq N}\mathbb{P}(\left\vert
X_{n}\right\vert >\lambda a_{n})<\varepsilon ,
\end{equation*}

\noindent ce qui est (\ref{big-O2}).\newline

\bigskip

\bigskip \noindent \textbf{COMMENTAIRES, NOTATIONS ET QUELQUES LEMMES}.%
\newline

\noindent \textbf{(a)} A partir du Chapitre \ref{cv.tensRk}, un $O_{\mathbb{P%
}}(1)$ is simplement une suite de variables al\'{e}atoires tendues. A partir
du Theorem \ref{tensTheo2} du Chapitre \ref{cv.tensRk}, nous avons que
n'importe quelle suite $X_{n}=O_{\mathbb{P}}(a_{n})$ contient une sous suite 
$(X_{n_{k}})_{k\geq 1}$ telle que $(X_{n_{k}}/a_{n_{k}})_{k\geq 1}$ converge
vaguement dans $\mathbb{R}$.\newline

\noindent \textbf{(b)} Il peut \^{e}tre convenant de reformuler (\ref{big-O1}%
) \`{a} la phrase suivante. Pour tout $\varepsilon >0,$ il existe un nombre r%
\'{e}el $\lambda >0$ tel que $\left\vert X_{n}\right\vert \leq \lambda a_{n}$
Avec une Probabilit\'{e} Au Moins Egale \`{a} $1-\varepsilon $ - not\'{e} $%
APAME(1-\varepsilon )$ pour tout $n\geq 1,$ ou $\left\vert X_{n}\right\vert
>\lambda a_{n}$ Avec une Probabilit\'{e} Au Plus Egale \`{a}  $\varepsilon $
- denoted $APAPE(\varepsilon )$ pour tout $n\geq 1.$(\ref{big-O3}), nous
rempla\c{c}ons \textbf{(pour }$n\geq 1)$ par \textbf{(pour des grandes
valeures de }$n)$ ou par \textbf{(pour }$n$\textbf{\ plus grand qu'un
quelconque nombre }$n_{0}>1)$.\newline

\noindent Pour des d\'{e}monstrations tr\`{e}s longues, l'utilisation des
types de phrases d\'{e}crits ci dessus peut \^{e}tre utile.\newline

\noindent Maintenant, nous pouvons donner quelques propri\'{e}t\'{e}s
importantes des petits $o^{\prime }s$ et des grands $O^{\prime }s$ en
probabilit\'{e}.\newline

\bigskip \noindent Le lemme suivant est aussi utile.

\begin{lemma}
\label{oO-02} Nous avons les propri\'{e}t\'{e}s suivantes

\noindent (a) Si $X_{n}$ est une suite de $k-$vecteurs al\'{e}atoires
convergeant en probabilit\'{e} vers $k$-vecteur al\'{e}atoire $X,$ alors $%
\left\Vert X_{n}\right\Vert =O_{\mathbb{P}}(1)$.\newline

(b) Si $X_{n}$ est une suite de vecteurs al\'{e}atoires \`{a} valeurs dans
un espace m\'{e}trique $(S,d)$ convergeant en probabilit\'{e} vers une
constante $C\in S$ et si $g$ est une application m\'{e}surable de $(S,d)$
vers un espace m\'{e}trique $(E,r).$ Ensuite si $g$ est continue en $C,$
alors $g(X_{n})$ converge en probabilit\'{e} vers $C.$

\bigskip \noindent (c) Alors pour toute suite de $k-$vecteurs al\'{e}atoires 
$(X_{n})_{n\geq 1}$ convergeant vers zero en probabilit\'{e}. Soient $R(x)$
une fonction r\'{e}elle de $x\in \mathbb{R}^{k}$ telle que $R(0)=0.$ Si $%
R(x)=o(\Vert x\Vert ^{p})$ quand $x\rightarrow 0$ pour un certain $p>0$,
alors $R(X_{n})=o_{\mathbb{P}}(\Vert X_{n}\Vert ^{p}).$ Si $R(x)=O(\Vert
x\Vert ^{p})$ quand $x\rightarrow 0$ pour un certain $p>0$, alors $%
R(X_{n})=O_{\mathbb{P}}(\Vert X_{n}\Vert ^{p})$
\end{lemma}

\bigskip \noindent \textbf{Preuve}.\newline

\noindent \textbf{Preuve du Point (a)}. Si $X_{n}\rightarrow _{\mathbb{P}}X$%
, alors par la Proposition \ref{cv.CvCp.01}, $X_{n}\rightarrow _{w}X$ et par
le Th\'{e}or\`{e}me \ref{cv.mappingTh} de l'application continue du Chapitre %
\ref{cv} $\left\Vert X_{n}\right\Vert \rightarrow _{w}\left\Vert
X\right\Vert .$ Alors par le Th\'{e}or\`{e}me \ref{cv.theo.portmanteau.rk},
nous avons pour tout point de continuit\'{e} de $F_{\left\Vert X\right\Vert
}(\lambda )=P(\left\Vert X\right\Vert \leq \lambda ),$

\begin{equation*}
\lim \mathbb{P}\left( \left\vert X_{n}\right\vert >\lambda \right) =\mathbb{P%
}\left( \left\Vert X\right\Vert >\lambda \right) =F_{\left\Vert X\right\Vert
}(\lambda ).
\end{equation*}

\noindent Puisque l'ensemble des points de continuit\'{e} de $F_{\left\Vert
X\right\Vert }$ est au plus denombrable (see \ ). Alors appliquer la formule
ci dessus pour $\lambda \rightarrow +\infty $ lorsque les $\lambda $ sont
des points de continuit\'{e}. Puisque $1-F_{\left\Vert X\right\Vert
}(\lambda )\rightarrow 0$ quand $\lambda \rightarrow +\infty ,$ alors pour
tout $\varepsilon >0,$ nous sommes capable de choisir une valeur de $\lambda
(\varepsilon )$ qui est un point de continuit\'{e} de $F_{\left\Vert
X\right\Vert }$ satisfaisant $1-F_{\left\Vert X\right\Vert }(\lambda
)<\varepsilon .$ Pour tout $\ \varepsilon >0,$ nous avons trouv\'{e} $%
\lambda (\varepsilon )>0$ tel que

\begin{equation*}
\limsup_{n\rightarrow +\infty }\mathbb{P}\left( \left\vert X_{n}\right\vert
>\lambda (\varepsilon )\right) \leq \varepsilon ,
\end{equation*}

\noindent ce qui implique  
\begin{equation*}
\lim_{\lambda \rightarrow +\infty }\limsup_{n\rightarrow +\infty }\mathbb{P}%
\left( \left\vert X_{n}\right\vert >\lambda \right) =0.
\end{equation*}

\noindent Le Point (a) est prouv\'{e}.\newline

\noindent \textbf{Preuve du Point (b)}. Supposons les notations de ce point
et supposons que $g$ est continue en $C$. Soit $\varepsilon >0.$ Par la
continuit\'{e} de $g$ en $C$, il existe $\eta >0$ tel que

\begin{equation*}
d(x,C)<\eta \Longrightarrow r(g(x),g(C))<\varepsilon /2.
\end{equation*}

\noindent Maintenant

\begin{eqnarray*}
\mathbb{P}(r(g(X_{n}),g(C) &>&\varepsilon )=\mathbb{P}(\left\{
r(g(X_{n}),g(C)>\varepsilon \right\} \cap \left\{ d(X_{n},C)\geq \eta
\right\} ) \\
&+&\mathbb{P}(\left\{ r(g(X_{n}),g(C)>\varepsilon \right\} \cap \left\{
d(X_{n},C)<\eta \right\} ) \\
&\leq &\mathbb{P}(d(X_{n},C)\geq \eta ),
\end{eqnarray*}

\noindent puisque $(\left\{ r(g(X_{n}),g(C)>\lambda \right\} \cap \left\{
d(X_{n},C)<\eta \right\} )\subset $ $(\left\{ r(g(X_{n}),g(C)>\varepsilon
\right\} \cap \left\{ d(X_{n},C)<\eta \right\} \cap \left\{ \left\{
r(g(X_{n}),g(C)<\varepsilon /2\right\} \right\} )=\emptyset .$ Alors,
puisque $X_{n}\rightarrow _{\mathbb{P}}C,$ nous avons

\begin{equation*}
\limsup_{n\rightarrow +\infty }\mathbb{P}(r(g(X_{n}),g(C)>\varepsilon )\leq
\limsup_{n\rightarrow +\infty }\mathbb{P}(d(X_{n},C)\geq \eta )=0.
\end{equation*}

\noindent Ainsi le Point (b) est vrai.\newline

\noindent \textbf{Preuve du Point (c-1)}. Soit $R(x)=o(\Vert x\Vert ^{p})$
quand $x\rightarrow 0.$ Alors $g(x)=\left\vert R(x)/\Vert x\Vert
^{p}\right\vert \rightarrow 0$ $x\rightarrow g(0)=0.$ Par la continuit\'{e}
de $g$ en zero et par le Point (a), $g(X_{n})=\left\vert R(X_{n})\right\vert
/\Vert X_{n}\Vert ^{p}\rightarrow _{\mathbb{P}}0$ quand $n\rightarrow
+\infty $ lorsque $X_{n}\rightarrow _{\mathbb{P}}X$ quand $n\rightarrow
+\infty .$ D'o\`{u} $R(X_{n})=o_{\mathbb{P}}(\Vert X_{n}\Vert ^{p}).$.%
\newline

\noindent \textbf{Preuve du Point (c-2)}. Soit $R(x)=O(\Vert x\Vert ^{p})$
quand $x\rightarrow 0.$ Alors pour tout $\varepsilon >0,$\ ils existent $%
\eta >0$ et $C>0$ tel que $\left\vert R(x)\right\vert /\Vert x\Vert ^{p}\leq
C$ pour tout $\left\Vert x\right\Vert <\eta .$ Alors pour $\lambda >C,$

\begin{eqnarray*}
\mathbb{P}(\left\vert R(X_{n})\right\vert /\Vert X_{n}\Vert ^{p} &>&\lambda
)=\mathbb{P}(\left\{ \left\vert R(X_{n})\right\vert /\Vert X_{n}\Vert
^{p}>\lambda \right\} \cap \left\{ \left\Vert X_{n}\right\Vert \geq \eta
\right\} ) \\
&+&\mathbb{P}(\left\{ \left\vert R(X_{n})\right\vert /\Vert X_{n}\Vert
^{p}>\lambda \right\} \cap \left\{ \left\Vert X_{n}\right\Vert <\eta
\right\} ) \\
&\leq &\mathbb{P}(\left\Vert X_{n}\right\Vert \geq \eta ),
\end{eqnarray*}

\noindent puisque $(\left\{ \left\vert R(X_{n})\right\vert /\Vert X_{n}\Vert
^{p}>\lambda \right\} \cap \left\{ \left\Vert X_{n}\right\Vert <\eta
\right\} )\subset (\left\{ \left\vert R(X_{n})\right\vert /\Vert X_{n}\Vert
^{p}>\lambda \right\} \cap \left\{ \left\Vert X_{n}\right\Vert <\eta
\right\} \cap \left\{ \left\{ \left\vert R(X_{n})\right\vert /\Vert
X_{n}\Vert ^{p}<C\right\} \right\} )=\emptyset .$ Alors pour tout $\lambda
>C,$

\begin{equation*}
\limsup_{n\rightarrow +\infty }\mathbb{P}(\left\vert R(X_{n})\right\vert
/\Vert X_{n}\Vert ^{p}>\lambda )=0
\end{equation*}

\noindent et ensuite

\begin{equation*}
\lim_{\lambda \rightarrow +\infty }\limsup_{n\rightarrow +\infty }\mathbb{P}%
(\left\vert R(X_{n})\right\vert /\Vert X_{n}\Vert ^{p}>\lambda )=0.
\end{equation*}

\noindent Nous concluons que $\left\vert R(X_{n})\right\vert /\Vert
X_{n}\Vert ^{p}=O_{\mathbb{P}}(1).$.\newline

\bigskip \noindent \textbf{PROPRIETES PRINCIPALES}.\newline

\noindent \textbf{(1)} Si $X_{n}=o(1)$ $p.s.,$ alors $X_{n}=o_{\mathbb{P}}(1)
$.\newline

\noindent \textbf{(2)} $o_{\mathbb{P}}(a_{n})=a_{n}o_{\mathbb{P}}(1)$ et $%
a_{n}o_{\mathbb{P}}(1)=O_{\mathbb{P}}(a_{n})$.\newline

\noindent \textbf{(3)} $o_{\mathbb{P}}(a_{n})o_{\mathbb{P}}(b_{n})=o_{%
\mathbb{P}}(a_{n}b_{n})$.\newline

\noindent \textbf{(4)} $o_{\mathbb{P}}(o_{\mathbb{P}}(a_{n}))=o_{\mathbb{P}%
}(a_{n})$.\newline

\noindent \textbf{(5)} Si $b_{n}\geq a_{n}$ pour tout $n\geq 1,o_{\mathbb{P}%
}(a_{n})=o_{\mathbb{P}}(b_{n})$.\newline

\noindent \textbf{(6)} $o_{\mathbb{P}}(a_{n})+o_{\mathbb{P}}(a_{n})=o_{%
\mathbb{P}}(a_{n})$.\newline

\noindent \textbf{(7)} $o_{\mathbb{P}}(a_{n})+o_{\mathbb{P}}(b_{n})=o_{%
\mathbb{P}}(a_{n}+b_{n})$ et $o_{\mathbb{P}}(a_{n})+o_{\mathbb{P}}(b_{n})=o_{%
\mathbb{P}}(a_{n}\vee b_{n})$ $p.s.$ o\`{u} $a_{n}\vee b_{n}=\max
(a_{n},b_{n})$.\newline

\noindent \textbf{(8)} Si $X_{n}=O(1)$ $p.s.,$ alors $X_{n}=O_{\mathbb{P}}(1)
$.\newline

\noindent \textbf{(9)} $O_{\mathbb{P}}(a_{n})=a_{n}O_{\mathbb{P}}(1)$

\noindent \textbf{(10)} $O_{\mathbb{P}}(a_{n})O_{\mathbb{P}}(b_{n})=O_{%
\mathbb{P}}(a_{n}b_{n})$.\newline

\noindent \textbf{(11)} $O_{\mathbb{P}}(O_{\mathbb{P}}(a_{n}))=O_{\mathbb{P}%
}(a_{n})$ $a.s.$.\newline

\noindent \textbf{(12)} Si $b_{n}\geq a_{n}$ pour tout $n\geq 1,O_{\mathbb{P}%
}(a_{n})=O_{\mathbb{P}}(b_{n})$.\newline

\noindent \textbf{(13)} $O_{\mathbb{P}}(a_{n})+O_{\mathbb{P}}(a_{n})=O_{%
\mathbb{P}}(a_{n})$.\newline

\noindent \textbf{(14)} $O_{\mathbb{P}}(a_{n})+O_{\mathbb{P}}(b_{n})=O_{%
\mathbb{P}}(a_{n}+b_{n})$ et $O_{\mathbb{P}}(a_{n})+O_{\mathbb{P}}(b_{n})=O_{%
\mathbb{P}}(a_{n}\vee b_{n})$ o\`{u} $a_{n}\vee b_{n}=\max (a_{n},b_{n})$.%
\newline

\noindent \textbf{(15)} $o_{\mathbb{P}}(a_{n})O_{\mathbb{P}}(b_{n})=o_{%
\mathbb{P}}(a_{n}b_{n})$.\newline

\noindent \textbf{(16)} $o_{\mathbb{P}}(O_{\mathbb{P}}(a_{n}))=o_{\mathbb{P}%
}(a_{n})$ et $O_{\mathbb{P}}(o_{\mathbb{P}}(a_{n}))=o_{\mathbb{P}}(a_{n})$.%
\newline

\noindent \textbf{(17a)} Si $a_{n}=O_{\mathbb{P}}(b_{n})$ alors $o(a_{n})+O_{%
\mathbb{P}}(b_{n})=O_{\mathbb{P}}(b_{n})$ $p.s.$\newline

\noindent \textbf{(17b)} Si $b_{n}=O_{\mathbb{P}}(a_{n})$ alors $o_{<\mathbb{%
P}}(a_{n})+O_{\mathbb{P}}(b_{n})=O_{\mathbb{P}}(a_{n})$.\newline

\noindent \textbf{(17c)} Si $b_{n}=o_{\mathbb{P}}(a_{n}),$ $p.s.,$ alors $%
o(a_{n})+O_{\mathbb{P}}(b_{n})=o(b_{n})$.\newline

\noindent \textbf{(18)} $(1+o_{\mathbb{P}}(a_{n}))^{-1}-1=o_{\mathbb{P}%
}(a_{n})$.\newline

\noindent \textbf{(19)} Un $o_{\mathbb{P}}(1)$ est un $O_{\mathbb{P}}(1)$

\bigskip \noindent \textbf{PREUVES}.\newline

\bigskip \noindent (1) Ceci vient de l'implication : $X_{n}\rightarrow 0$
p.s. $\Longrightarrow X_{n}\rightarrow _{P}0$ (Voir Proposition \ref%
{cv.CvCp.01}).\newline

\bigskip \noindent (2) Si $X_{n}=o_{\mathbb{P}}(a_{n})\Longleftrightarrow
\left\vert X_{n}/a_{n}\right\vert \longrightarrow _{P}0\Longleftrightarrow
X_{n}/a_{n}=o_{\mathbb{P}}(1)\Longleftrightarrow X_{n}=a_{n}o_{\mathbb{P}%
}(1).$\newline

\bigskip \noindent (3) Par le Point (2) ci dessus, $o_{\mathbb{P}}(a_{n})o_{%
\mathbb{P}}(b_{n})=a_{n}b_{n}\times o_{\mathbb{P}}(1)o_{\mathbb{P}%
}(1)=a_{n}b_{n}\times o_{\mathbb{P}}(1)=o_{\mathbb{P}}(a_{n}b_{n})$ (V\'{e}%
rifions $o_{\mathbb{P}}(1)o_{\mathbb{P}}(1)=o_{\mathbb{P}}(1)$ en Propri\'{e}%
t\'{e} (A1) dans la sous-section Appendice ci dessous).\newline

\bigskip \noindent (4) $o_{\mathbb{P}}(o_{\mathbb{P}}(a_{n}))=o_{\mathbb{P}%
}(a_{n})o_{\mathbb{P}}(1)=a_{n}\times o_{\mathbb{P}}(1)o_{\mathbb{P}%
}(1)=a_{n}\times o_{\mathbb{P}}(1)=o_{\mathbb{P}}(a_{n})$ (Utilisons encore
la Propri\'{e}t\'{e} (A1) de la sous-section Appendice ci dessous).\newline

\bigskip \noindent (5) Soit $b_{n}\geq a_{n}$ pour tout $n\geq 1,X_{n}=o_{%
\mathbb{P}}(a_{n}).$ Pour tout $\eta >0,0\leq \lim_{n\rightarrow +\infty
}\sup P(\left\vert X_{n}/b_{n}\right\vert >\eta )\leq \lim_{n\rightarrow
+\infty }P(\left\vert X_{n}/a_{n}\right\vert >\eta )=0.$\newline

\bigskip \noindent (6) Soient $X_{n}=o_{\mathbb{P}}(a_{n})$ et $Y_{n}=o_{%
\mathbb{P}}(a_{n}).$ Utiliser l'outil classique, pour $\eta >0,$%
\begin{equation*}
\left( \frac{\left\vert X_{n}\right\vert }{a_{n}}>\eta /2\right) \cap \left( 
\frac{\left\vert Y_{n}\right\vert }{a_{n}}>\eta /2\right) \subset \left( 
\frac{\left\vert X_{n}+Y_{n}\right\vert }{a_{n}}>\eta \right) .
\end{equation*}

\noindent Alors pour $\eta >0,$
\begin{eqnarray}
\limsup_{n\rightarrow +\infty }\mathbb{P}\left( \frac{\left\vert
X_{n}+Y_{n}\right\vert }{a_{n}}>\eta \right)  &\leq &\limsup_{n\rightarrow
+\infty }\mathbb{P}\left( \frac{\left\vert X_{n}\right\vert }{a_{n}}>\eta
/2\right)  \notag \\
&+&\limsup_{n\rightarrow +\infty }\mathbb{P}\left( \frac{\left\vert
Y_{n}\right\vert }{a_{n}}>\eta /2\right) =0. \label{decompSum}
\end{eqnarray}

\bigskip \noindent (7) Pour prouver ce point, il suffit de combiner les Points (5) et (6)
ci-dessus.\newline

\bigskip \noindent (8) $X_{n}=O(1)$ $a.s.$ quand $n\rightarrow +\infty $
signifie qu'il existe $\Omega _{0}$ m\'{e}surable telle que $\mathbb{P}%
(\Omega _{0})=1$ et pour tout $\omega \in \Omega _{0},$ 
\begin{equation*}
\limsup_{n\rightarrow +\infty }\left\vert X_{n}(\omega )\right\vert
=\inf_{n\geq 1}\sup_{p\geq n}\left\vert X_{p}\right\vert =M(\omega )<+\infty
.
\end{equation*}

\noindent Nous avons pour tout $n\geq 1,$%
\begin{equation*}
\mathbb{P}\left( \left\vert X_{n}\right\vert >\lambda \right) \leq \mathbb{P}%
\left( \sup_{p\geq n}\left\vert X_{p}\right\vert >\lambda \right) 
\end{equation*}

\noindent Nous avons $Y_{n}=\sup_{p\geq n}\left\vert X_{p}\right\vert
\downarrow M$ $p.s.$ Alors $Y_{n}1_{\Omega _{0}}\rightarrow _{\mathbb{P}%
}M1_{\Omega _{0}}$. Nous travaillons avec les variables al\'{e}atoires et
nous pouvons appliquer les r\'{e}sultats de la convergence vague. Nous
obtenons $Y_{n}1_{\Omega _{0}}\rightarrow _{w}M1_{\Omega _{0}}$ par la
Proposition \ref{cv.CvCp.01}$.$ Par le Th\'{e}or\`{e}me \ref%
{cv.theo.portmanteau.rk}, nous avons pour tout point de $F_{M}(\lambda
)=P(M1_{\Omega _{0}}\leq \lambda )$. Utilisons le Th\'{e}or\`{e}me de la
Convergence Monotone pour obtenir

\begin{eqnarray*}
\limsup_{n\rightarrow +\infty }\mathbb{P}\left( \left\vert X_{n}\right\vert
>\lambda \right) &=&\limsup_{n\rightarrow +\infty }\mathbb{P}\left(
\sup_{p\geq n}\left\vert X_{p}\right\vert 1_{\Omega _{0}}>\lambda \right) \\
&=&\limsup_{n\rightarrow +\infty }\mathbb{P}\left( \sup_{p\geq n}\left\vert
X_{p}\right\vert 1_{\Omega _{0}}>\lambda \right) =\mathbb{P}\left(
M1_{\Omega _{0}}>\lambda \right) =1-F_{M}(\lambda ).
\end{eqnarray*}

\noindent Puisque l'ensemble des points de discontinuit\'{e} de $F_{M}$ est
au plus d\'{e}nombrable (Voir \ ). Nous appliquons alors la formule ci
dessus pour $\lambda \rightarrow +\infty $ lorsque les $\lambda $ \ sont des
points de continuit\'{e}. Puisque $1-F_{M}(\lambda )\rightarrow 0$ quand $%
\lambda \rightarrow +\infty ,$ alors pour tout $\varepsilon >0,$ nous
pouvons choisir une valeur de $\lambda (\varepsilon )$ qui est un point de
continuit\'{e} de$F_{M}$ satisfaisant $1-F_{M}(\lambda )<\varepsilon .$ Pour
tout $\varepsilon >0,$ nous trouver $\lambda (\varepsilon )>0$ tel que

\begin{equation*}
\limsup_{n\rightarrow +\infty }\mathbb{P}\left( \left\vert X_{n}\right\vert
>\lambda (\varepsilon )\right) \leq \varepsilon ,
\end{equation*}

\noindent ce qui implique 
\begin{equation*}
\lim_{\lambda \rightarrow +\infty }\limsup_{n\rightarrow +\infty }\mathbb{P}%
\left( \left\vert X_{n}\right\vert >\lambda \right) =0.
\end{equation*}

\noindent Donc $X_{n}=O_{\mathbb{P}}(1).$\newline

\bigskip \noindent (9) Soit $X_{n}=O_{\mathbb{P}}(a_{n}).$ On a

\begin{equation*}
\lim_{\lambda \rightarrow +\infty }\limsup_{n\rightarrow +\infty }\mathbb{P}%
\left( \left\vert X_{n}/a_{n}\right\vert >\lambda \right) =0.
\end{equation*}

\noindent C'est la d\'{e}finition que $X_{n}/a_{n}=O_{\mathbb{P}}(1)$ et
alors $X_{n}=a_{n}O_{\mathbb{P}}(1).$\newline

\bigskip \noindent (10) $O_{\mathbb{P}}(a_{n})O_{\mathbb{P}%
}(b_{n})=a_{n}b_{n}O_{\mathbb{P}}(1)O_{\mathbb{P}}(1)=a_{n}b_{n}O_{\mathbb{P}%
}(1)$ (V\'{e}rifions la Propri\'{e}t\'{e} ci dessous $O_{\mathbb{P}}(1)O_{%
\mathbb{P}}(1)=O_{\mathbb{P}}(1)$ dans (A2) dans la sous-section Appendice).%
\newline

\bigskip \noindent \textbf{(11)} $O_{\mathbb{P}}(O_{\mathbb{P}}(a_{n}))=O_{%
\mathbb{P}}(a_{n})O_{\mathbb{P}}(1)=a_{n}O_{\mathbb{P}}(1)O_{\mathbb{P}%
}(1)=O_{\mathbb{P}}(a_{n}).$.\newline

\bigskip \noindent \textbf{(12)} Soit $b_{n}\geq a_{n}$ pour tout $n\geq 1$
et $X_{n}=O_{\mathbb{P}}(a_{n})$. On a%
\begin{equation*}
\lim_{\lambda \rightarrow +\infty }\limsup_{n\rightarrow +\infty }\mathbb{P}%
\left( \left\vert X_{n}/b_{n}\right\vert >\lambda \right) \leq \lim_{\lambda
\rightarrow +\infty }\limsup_{n\rightarrow +\infty }\mathbb{P}\left(
\left\vert X_{n}/a_{n}\right\vert >\lambda \right) =0.
\end{equation*}

\noindent Donc $X_{n}=O_{\mathbb{P}}(b_{n})$.\newline

\bigskip \noindent \textbf{(13)} Soit $X_{n}=O_{\mathbb{P}}(a_{n})$ et $%
X_{n}=O_{\mathbb{P}}(a_{n})$. Utilisons la m\^{e}me technique comme dans la
formule \ref{decompSum} ci-dessous pour trouver 

\begin{eqnarray}
\lim_{\lambda \rightarrow +\infty }\limsup_{n\rightarrow +\infty }\mathbb{P%
}\left( \frac{\left\vert X_{n}+Y_{n}\right\vert }{a_{n}}>\lambda \right)
&\leq& \lim_{\lambda \rightarrow +\infty }\limsup_{n\rightarrow +\infty }%
\mathbb{P}\left( \frac{\left\vert X_{n}\right\vert }{a_{n}}>\lambda
/2\right) \notag\\
&+&\lim_{\lambda \rightarrow +\infty }\limsup_{n\rightarrow
+\infty }\mathbb{P}\left( \frac{\left\vert Y_{n}\right\vert }{a_{n}}>\lambda
/2\right) =0.
\end{eqnarray}

\bigskip \noindent \textbf{(14)} Combinons les Points (12) et (13) pour
trouver celui-ci.\newline

\bigskip \noindent \textbf{(15)} $o_{\mathbb{P}}(a_{n})O_{\mathbb{P}%
}(b_{n})=a_{n}b_{n}o_{\mathbb{P}}(1)O_{\mathbb{P}}(1)=a_{n}b_{n}o_{\mathbb{P}%
}(1)=o_{\mathbb{P}}(a_{n}b_{n})$. (V\'{e}rifions que $o_{\mathbb{P}}(1)O_{%
\mathbb{P}}(1)$ dans la proprit\'e (A3) dans la sous-section Appendice ci
dessous)$.$\newline

\bigskip \noindent \textbf{(16)} $o_{\mathbb{P}}(O_{\mathbb{P}}(a_{n}))=O_{%
\mathbb{P}}(a_{n})o_{\mathbb{P}}(1)=a_{n}O_{\mathbb{P}}(1)o_{\mathbb{P}%
}(1)=a_{n}o_{\mathbb{P}}(1)=o_{\mathbb{P}}(a_{n})$ et $O_{\mathbb{P}}(o_{%
\mathbb{P}}(a_{n}))=o_{\mathbb{P}}(a_{n})O_{\mathbb{P}}(1)=o_{\mathbb{P}%
}(a_{n})$.\newline

\bigskip \noindent \textbf{(17a)} Soient $a_{n}=O(b_{n})$ et $X_{n}=o_{%
\mathbb{P}}(a_{n})$ et $Y_{n}=O_{\mathbb{P}}(b_{n}).$ Il existe $C>0$ telle
que $a_{n}\leq Cb_{n}$ pour tout $n\geq 1.$ Alors $X_{n}=o_{\mathbb{P}%
}(a_{n})=o_{\mathbb{P}}(Cb_{n})$ par le Point (5). Mais \'{e}videmment $%
X_{n}=o_{\mathbb{P}}(Cb_{n})=o_{\mathbb{P}}(b_{n})$ et alors $X_{n}=O_{%
\mathbb{P}}(b_{n})$ par le Point (19) ci dessous$.$ Finalement $%
X_{n}+Y_{n}=O_{\mathbb{P}}(b_{n})+O_{\mathbb{P}}(b_{n})=O_{\mathbb{P}}(b_{n})
$.\newline

\bigskip \noindent \textbf{(17b)} Soient $b_{n}=O(a_{n})$ et $X_{n}=o_{%
\mathbb{P}}(a_{n})$ et $Y_{n}=O_{\mathbb{P}}(b_{n}).$ Nous \'{e}changeons
les r\^{o}les de $a_{n}$ et $b_{n}$ pour obtenir $b_{n}\leq Ca_{n}$ et $%
Y_{n}=O_{\mathbb{P}}(Ca_{n})=O_{\mathbb{P}}(a_{n})$ par le Point (1) et
finalement $X_{n}+Y_{n}=o_{\mathbb{P}}(a_{n})+O_{\mathbb{P}}(a_{n})=O_{%
\mathbb{P}}(a_{n})+O_{\mathbb{P}}(a_{n})=O_{\mathbb{P}}(a_{n}).$\newline

\bigskip \noindent \textbf{(17c)} Soient $b_{n}=o_{\mathbb{P}}(a_{n})$ et $%
X_{n}=o_{\mathbb{P}}(a_{n})$ et $Y_{n}=O_{\mathbb{P}}(b_{n}).$ Nous avons $%
X_{n}+Y_{n}=o_{\mathbb{P}}(a_{n})+O_{\mathbb{P}}(o_{\mathbb{P}}(a_{n}))=o_{%
\mathbb{P}}(a_{n})+o_{\mathbb{P}}(a_{n})$ par le Point (16). Finalement $%
X_{n}+Y_{n}=o_{\mathbb{P}}(a_{n}).$\newline

\bigskip \noindent \textbf{(18)} Nous avons 
\begin{equation*}
(1+o_{\mathbb{P}}(a_{n}))^{-1}-1=\frac{o_{\mathbb{P}}(a_{n})}{1+o_{\mathbb{P}%
}(a_{n})}.
\end{equation*}

\noindent Par le Point (b) , $(1+o_{\mathbb{P}}(a_{n}))^{-1}\rightarrow _{P}1
$ et par le Point (a) du m\^{e}me lemme, $(1+o_{\mathbb{P}}(a_{n}))^{-1}=O_{%
\mathbb{P}}(1).$ Ainsi 
\begin{equation*}
(1+o_{\mathbb{P}}(a_{n}))^{-1}-1=O_{\mathbb{P}}(1)o_{\mathbb{P}}(a_{n})=o_{%
\mathbb{P}}(a_{n}),
\end{equation*}

\noindent par le Point (15).\newline

\bigskip \noindent $(19)$ Par le Lemme \ref{oO-02}, un $o_{\mathbb{P}}(1)$
converge vers $O$ en probabilit\'{e} et est alors un $O_{\mathbb{P}}(1)$.%
\newline

\subsection{Extensions}

Les concepts de petit\textit{\ o} et grand\textit{\ O} dans $\mathbb{R}^{k}$
de la fa\c{c}on suivante :\newline

\noindent \textbf{(1)} La suite de vecteurs al\'{e}atoires $(X_{n})_{n\geq 1}
$ dans $\mathbb{R}^{k},$ est un $o(a_{n})$ $p.s.$ si et seulement si$%
\left\Vert X_{n}\right\Vert /a_{n}=o(1)$ $p.s.$, et est un $o_{\mathbb{P}%
}(a_{n})$ si et seulement si $\left\Vert X_{n}\right\Vert /a_{n}=o_{\mathbb{P%
}}(1)$.\newline

\noindent \textbf{(b)} La suite de vecteurs al\'{e}atoires $(X_{n})_{n\geq 1}
$ dans $\mathbb{R}^{k},$ est un $O(a_{n})$ $p.s.$ si seulement si $%
\left\Vert X_{n}\right\Vert /a_{n}=O(1)$ $p.s.$, et est un $O_{\mathbb{P}%
}(a_{n})$ si et seulement si $\left\Vert X_{n}\right\Vert /a_{n}=O_{\mathbb{P%
}}(1)$.\newline

\noindent A partir de l\`{a}, traiter ces concepts est facile en combinant
leurs propri\'{e}t\'{e}s dans $\mathbb{R}$ et celles des normes dans $%
\mathbb{R}^{k}.$

\subsection{\protect\bigskip Suites \'{e}quilibr\'{e}es}

Il peut aider dans certains cas d'avoir des suites $X_{n}$ telles que $%
\left\Vert X_{n}\right\Vert $ et $1/\left\Vert X_{n}\right\Vert $ soient born%
\'{e}es en probabilit\'{e}. Donnons quelques notations pour les suites r\'{e}%
elles.\newline

\noindent \textbf{(1)} Pour $0\leq a<b<+\infty ,$ nous notons par $X_{n}=O_{%
\mathbb{P}}(a,b,a_{n},b_{n})$ la propri\'{e}t\'{e} que pour tout $%
\varepsilon >0,$ il existe $\lambda >0$ tel que nous avons $(a+\lambda \leq
\left\vert X_{n}\right\vert /a_{n},\left\vert X_{n}\right\vert /a_{n}\leq
b-\lambda )$ $APAME(1-\varepsilon ),$ pour des grandes valeurs de $n.$ Si $%
a_{n}=b_{n}$ pour tout $n\geq 1,$ nous \'{e}crivons simplement $X_{n}=O_{%
\mathbb{P}}(a,b,a_{n}).$\newline

\noindent \textbf{(1)} Pour $0\leq a,$ nous notons par $X_{n}=O_{\mathbb{P}%
}(a,+\infty ,a_{n},b_{n})$ la propri\'{e}t\'{e} que pour tout $\varepsilon
>0,$ il existe $\lambda >0$ telle que nous avons $(a+\lambda ^{-1}\leq
\left\vert X_{n}\right\vert /a_{n},\left\vert X_{n}\right\vert /a_{n}\leq
\lambda )$ $APAME(1-\varepsilon ),$ pour des grandes valeurs de $n.$ If $%
a_{n}=b_{n}$ pour tout $n\geq 1,$ nous \'{e}crivons simplement $X_{n}=O_{%
\mathbb{P}}(a,+\infty ,a_{n}).$.\newline

\noindent Un exemple d'un $X_{n}=O_{\mathbb{P}}(0,+\infty ,1)$ est une suite 
$X_{n}$ convergeant vaguement vers $X>0$ $p.s.$ Dans ce cas $%
1/X_{n}\rightarrow 1/X$ est finie $p.s.$ et alors $X_{n}=O_{\mathbb{P}}(1)$
et $1/X_{n}=O_{\mathbb{P}}(1).$ Combinant ces deux points conduit \`{a} $%
X_{n}=O_{\mathbb{P}}(0,+\infty ,1).$\\

\newpage

\subsection{Appendice}

\noindent \textbf{(A1)} Si $X_{n}\rightarrow _{\mathbb{P}}a\in \mathbb{R}$
et $X_{n}\rightarrow _{\mathbb{P}}b\in \mathbb{R},$ alors $%
X_{n}Y_{n}\rightarrow _{\mathbb{P}}ab.$\newline

\noindent \textbf{Preuve}. Nous avons $(\eta +\left\vert b\right\vert )\eta
+\left\vert a\right\vert \eta \rightarrow 0$ $\ $quand $\eta \rightarrow 0.$
pour tout $\varepsilon >0,$ pour tout $\delta >0,$ choisissons une valeur de 
$\eta >0$ telle que $(\eta +\left\vert b\right\vert )\eta +\left\vert
a\right\vert \eta <\delta .$ Nous appliquons la d\'{e}finition des
convergences $X_{n}\rightarrow _{\mathbb{P}}a$ et $X_{n}\rightarrow _{%
\mathbb{P}}b$ pour obtenir une valeur $N_{0}\geq 1$ telle que pour tout $%
n\geq N_{0},$

\begin{equation*}
\mathbb{P}(\left\vert X_{n}-a\right\vert \geq \eta )\leq \varepsilon /2\text{
et }\mathbb{P}(\left\vert Y_{n}-b\right\vert \geq \eta )\leq \varepsilon /2.
\end{equation*}

\noindent Mais 
\begin{eqnarray*}
\left\vert X_{n}Y_{n}-ab\right\vert  &=&\left\vert
X_{n}Y_{n}-aY_{n}+aY_{n}-ab\right\vert  \\
&\leq &\left\vert Y_{n}\right\vert \left\vert X_{n}-a\right\vert +\left\vert
a\right\vert \left\vert Y_{n}-b\right\vert  \\
&\leq &(\left\vert Y_{n}-b\right\vert +\left\vert b\right\vert )\text{ }%
\left\vert X_{n}-a\right\vert +\left\vert a\right\vert \left\vert
Y_{n}-b\right\vert 
\end{eqnarray*}

\noindent Sur $(\left\vert X_{n}-a\right\vert \geq \eta )^{c}\cap
(\left\vert Y_{n}-b\right\vert \geq \eta )^{c},$

\begin{equation*}
\left\vert X_{n}Y_{n}-ab\right\vert \leq (\eta +\left\vert b\right\vert
)\eta +\left\vert a\right\vert \eta \leq \delta .
\end{equation*}

\noindent D'o\`{u} pour $n\geq N_{0},$%
\begin{equation*}
(\left\vert X_{n}-a\right\vert \geq \eta )^{c}\cap (\left\vert
Y_{n}-b\right\vert \geq \eta )^{c}\subset (\left\vert
X_{n}Y_{n}-ab\right\vert \leq \delta ),
\end{equation*}

\noindent c'est \`{a} dire 
\begin{equation*}
\mathbb{P}(\left\vert X_{n}-a\right\vert \geq \eta )^{c}\cap (\left\vert
Y_{n}-b\right\vert \geq \eta )^{c}\leq \mathbb{P}(\left\vert
X_{n}Y_{n}-ab\right\vert \leq \delta ),
\end{equation*}

\noindent et en pr\'{e}nant les compl\'{e}mentaires, 
\begin{equation*}
\mathbb{P}(\left\vert X_{n}Y_{n}-ab\right\vert >\delta )\leq \mathbb{P}%
((\left\vert X_{n}-a\right\vert \geq \eta )\cup (\left\vert
Y_{n}-b\right\vert \geq \eta ))\leq \varepsilon /2+\varepsilon
/2=\varepsilon .
\end{equation*}

\noindent Par suite 
\begin{equation*}
X_{n}Y_{n}\rightarrow _{\mathbb{P}}ab.
\end{equation*}

\bigskip \noindent \textbf{Propri\'{e}t\'{e} (A2)}. $X_{n}=O_{\mathbb{P}}(1)$
et $X_{n}=O_{\mathbb{P}}(1)$ alors $X_{n}Y_{n}=O_{\mathbb{P}}(1).$\newline

\noindent \textbf{Preuves}. En appliquant la d\'{e}finition d'un $O_{\mathbb{%
P}}(1)$, nous pouvons trouver pour tout $\varepsilon >0,$ deux nombres
entiers $N_{1}$ et $N_{2}$ et deux nombres positifs $\lambda _{1}>0$ et $%
\lambda _{2}>0$ tels que

\begin{equation*}
\forall (n\geq N_{1}),\mathbb{P(}\left\vert X_{n}\right\vert \leq \lambda
_{1})\geq 1-\varepsilon /2\text{ et }\forall (n\geq N_{2}),\mathbb{P(}%
\left\vert Y_{n}\right\vert \leq \lambda _{2})\geq 1-\varepsilon /2.
\end{equation*}

\noindent Pour $n\geq \max (N_{1,}N_{2}),$ 
\begin{equation*}
\mathbb{(}\left\vert X_{n}\right\vert \leq \lambda _{1})\cap \mathbb{(}%
\left\vert Y_{n}\right\vert \leq \lambda _{2})\subset \mathbb{(}\left\vert
X_{n}Y_{n}\right\vert \leq \lambda _{1}\lambda _{2})
\end{equation*}

\noindent ce qui est \'{e}quivalent \`{a} 
\begin{equation*}
\mathbb{(}\left\vert X_{n}Y_{n}\right\vert >\lambda _{1}\lambda _{2})\subset 
\mathbb{(}\left\vert X_{n}\right\vert >\lambda _{1})\cup \mathbb{(}%
\left\vert Y_{n}\right\vert >\lambda _{2})
\end{equation*}

\noindent ce qui implique pour $n\geq \max (N_{1,}N_{2}),$%
\begin{equation*}
\mathbb{P(}\left\vert X_{n}Y_{n}\right\vert >\lambda _{1}\lambda _{2})\leq 
\mathbb{P(}\left\vert X_{n}\right\vert >\lambda _{1})+\mathbb{P(}\left\vert
Y_{n}\right\vert >\lambda _{2})\leq \varepsilon /2+\varepsilon
/2=\varepsilon .
\end{equation*}

\noindent Ainsi pour tout $\varepsilon >0,$ il existe un nombre positif $N$ $%
(=\max (N_{1,}N_{2})),$ il existe $\lambda >0$ $(=\lambda _{1}\lambda _{2})$
et pour tout $n\geq N,$ 
\begin{equation*}
\mathbb{P(}\left\vert X_{n}Y_{n}\right\vert \leq \lambda )\geq 1-\varepsilon
..
\end{equation*}

\noindent D'o\`{u} $X_{n}Y_{n}=O_{\mathbb{P}}(1)$.\newline

\noindent

\begin{equation*}
X_{n}Y_{n}\rightarrow _{\mathbb{P}}ab.
\end{equation*}

\bigskip \textbf{Propri\'{e}t\'{e} (A3)}. Si $X_{n}=o_{\mathbb{P}}(1)$ et $%
Y_{n}=O_{\mathbb{P}}(1)$ alors $X_{n}Y_{n}=o_{\mathbb{P}}(1).$\newline

\noindent \textbf{Preuve}. Fixons $\varepsilon >0.$ En appliquant la d\'{e}%
finition d'un $O_{\mathbb{P}}(1)$ ils existent un nombre entier $N_{1}$ et
un nombre positif $\lambda >0$ tels que 
\begin{equation*}
\mathbb{P(}\left\vert Y_{n}\right\vert \leq \lambda )\geq 1-\varepsilon /2.
\end{equation*}

\noindent Maintenant soit $\eta >0.$ Appliquons la d\'{e}finition de $%
X_{n}=o_{\mathbb{P}}(1)$ pour obtenir qu'il existe un entier positif $N_{2}$
tel que

\begin{equation*}
\forall (n\geq N_{2}),\text{ }\mathbb{P(}\left\vert X_{n}\right\vert >\eta
/\lambda )\leq \varepsilon /2.
\end{equation*}

\noindent donc pour $n\geq \max (N_{1,}N_{2}),$ 
\begin{equation*}
\mathbb{(}\left\vert X_{n}\right\vert \leq \eta /\lambda )\cap \mathbb{(}%
\left\vert Y_{n}\right\vert \leq \lambda )\subset \mathbb{(}\left\vert
X_{n}Y_{n}\right\vert \leq \eta )
\end{equation*}

\noindent ce qui est \'{e}quivalent \`{a} 
\begin{equation*}
\mathbb{(}\left\vert X_{n}Y_{n}\right\vert >\eta )\subset \mathbb{(}%
\left\vert X_{n}\right\vert >\eta /\lambda )\cup \mathbb{(}\left\vert
Y_{n}\right\vert >\lambda )
\end{equation*}

\noindent ce qui implique pour $n\geq \max (N_{1,}N_{2}),$%
\begin{equation*}
\mathbb{P(}\left\vert X_{n}Y_{n}\right\vert >\eta )\leq \mathbb{P(}%
\left\vert X_{n}\right\vert >\eta /\lambda )+\mathbb{P(}\left\vert
Y_{n}\right\vert >\lambda )\leq \varepsilon /2+\varepsilon /2=\varepsilon .
\end{equation*}

\noindent Ainsi pour tout $\varepsilon >0,$ pour tout $\eta >0,$ il existe
un nombre positif $N$ $(=\max (N_{1,}N_{2})),$ pour tout $n\geq N,$%
\begin{equation*}
\mathbb{P(}\left\vert X_{n}Y_{n}\right\vert >\eta )\leq \varepsilon .
\end{equation*}

\noindent D'o\`{u} $X_{n}Y_{n}=o_{\mathbb{P}}(1)$.

\newpage

\section{M\'ethodes Delta} \label{cv.empTool.sec2}

\bigskip La m\'ethode Delta, dans ses diff\'erentes versions, constitue un moyen rapide pour obtenir de nouvels lois asymptotiques
à partir d'une suited de de variables al\'eatoires d\'ejà convergente vaguement et d\'efinies sur un m\^eme espace de probabilit\'e 
$(\Omega ,\mathcal{A},\mathbb{P})$. Nous pr\'esenterons ici les versions univari\'es, vectorielles et matricielles. Nous verrons ici encore l'utilit\'e
des r\'esultats de la section \ref{cv.CvCp} du chapitre \ref{cv} lorsqu'ils sont combin\'es avec les r\`egles de manipulations des petits $o^{\prime }s$ et des grands $O^{\prime }s$ en probabilit\'e expos\'es dans la premi\`ere section de ce chapitre.\\

\noindent Commen\c{c}ons avec la m\'ethode Delta univari\'ee, sur $\mathbb{R}$.\\ 

\subsection{Version univari\'e de la m\'ethode Delta}

\begin{proposition} \label{delta01} Soit $(X_{n})_{n\geq 1}$ une suite de variables al\'eatoires d\'efinies sur les m\^eme espace de probabilit\'e $(\Omega ,A,\mathbb{P})$ et soit $\theta$ un nombre r\'eel  et $(a_{n}>0)_{n\geq 1}$ une suite de nombres r\'eels v\'erifiant $a_{n}\rightarrow +\infty$ as $n\rightarrow +\infty$.\\

\noindent Soit $g:D\rightarrow \mathbb{R}$ une fonction de classe $C^{1}$, telle que $D$ soit un domaine de $\mathbb{R}$, 
que $\theta$ soit dans l'int\'erieur $\overset{o}{D}$ de $D$, et que l'on ait $\{X_{n},n\geq 1\}\subset \overset{o}{D}$.\\

\noindent  Si $a_{n}(X_{n}-\theta )$ converge vaguement vers une variable al\'eatoire $Z$ quand $n\rightarrow +\infty$, alors  
$$
a_{n}(g(X_{n})-g(\theta ))\rightsquigarrow g^{\prime }(\theta )Z \ \ quand \ \  n\rightarrow +\infty,
$$

\noindent o\`u $\nabla g(a)=g^{\prime}(\theta)$ est la d\'eriv\'ee  $g$ en $\theta$.
\end{proposition}

\bigskip \noindent \textbf{Preuve de la proposition \ref{delta01}}. Supposons que toutes les hypoth\`eses de la proposition proposition soient vraies. Par le point (a) du lemme \ref{oO-02}, nous avons $a_{n}(X_{n}-\theta )=O_{P}(1)$ et alors,
\begin{equation*}
X_{n}=\theta +O_{P}(1)a_{n}^{-1}\rightarrow _{\mathbb{P}}\theta. 
\end{equation*}

\noindent Ceci, en vertu de la proposition \ref{cv.CvECp} de la section \ref{cv.CvCp} du chapter \ref{cv}, est \'equivalent \`a la convergence vague suivante
\begin{equation*}
X_{n}\rightsquigarrow \theta.
\end{equation*}

\noindent Maintenant, par the th\'eor\`eme des accroissements finis, nous avons
\begin{equation}
g(X_{n})-g(\theta )=g^{\prime }(Y_{n})(X_{n}-\theta ),  \label{appliMVT}
\end{equation}

\noindent o\`u
\begin{equation*}
\min (X_{n},\theta )\leq Y_{n}\leq \max (X_{n},\theta ),
\end{equation*}

\noindent i.e.,

\begin{equation*}
\left\vert Y_{n}-\theta \right\vert \leq \left\vert X_{n}-\theta \right\vert.
\end{equation*}

\bigskip
\noindent Il vient que $Y_{n}\rightarrow _{\mathbb{P}}\theta $ et puisque  $g^{\prime }$
est continu, nous avons $g^{\prime }(Y_{n})\rightarrow _{\mathbb{P}}g(\theta )
$ en vertu du point (b) du lemme \ref{oO-02}. D\`es lors, en utilisant la proposition \ref{cv.CvECp} de la section \ref{cv.CvCp} du chapitre ref{cv}, nous voyons que ceci est \'equivalent \`a 
\begin{equation*}
g^{\prime }(Y_{n})\rightsquigarrow g^{\prime }(\theta ).
\end{equation*}

\noindent Par la propri\'et\'e de Slutsky donn\'ee dans la formule \ref{cv.slutsky} de la section section 
\ref{cv.CvCp} du chapitre \ref{cv}, nous avons la convergence vague suivante
\begin{equation*}
(g^{\prime }(Y_{n}),a_{n}(X_{n}-\theta ))\rightsquigarrow (g^{\prime
}(\theta),Z)
\end{equation*}

\noindent et par suite, en vertu du th\'eor\`eme de la transformation continue \ref{cv.mappingTh} du chapitre \ref{cv}, combin\'ee avec (\ref{appliMVT}), nous obtenons le r\'esultat final
 
\begin{equation*}
a_{n}(g(X_{n})-g(\theta ))=(g^{\prime }(Y_{n})\times a_{n}(X_{n}-\theta
))\rightsquigarrow g^{\prime }(\theta )Z.
\end{equation*}

\bigskip \bigskip \bigskip \textbf{Remarque}. Utilisons, non pas le nombre d\'eriv\'e, mais plut\^ot  la fonction (lin\'eaire) d\'eriv\'ee d\'efinie par

\begin{equation*}
h\rightarrow g_{\theta }^{\prime }(h)=g^{\prime }(\theta )h,
\end{equation*}

dans la proposition  \ref{delta01}. Nous pouvons \'ecrire alors le r\'esultat sous la forme
\begin{equation*}
a_{n}(g(X_{n})-g(\theta ))=g_{\theta }^{\prime }(Z).
\end{equation*}

\noindent Cette expression sugg\`ere la possibilit\'e d'avoir ce r\'esultats en dimension sup\'erieurs and plus tard dans des espaces fonctionnels.
Pour l'instant, passons aux formes multivari\'ees.

\subsection{Versions Multivari\'ees}

Le premier r\'esultat concerne la transformation de la convergence vague d'une suite suites variables al\'eatoires de $k$ composantes par
une fonction r\'eelle de $k$ arguments. Ceci donne la version vectorielle de la m\'ethode Delta.

\begin{proposition} \label{delta02} Soit $(X_{n})_{n\geq 1}$ une suite de vecteurs al\'eatoires de dimenstion $k$, $k\geq 1$, d\'efinis sur
le m\^eme espace de probabilit\'e $(\Omega ,A,\mathbb{P})$ et soit $\theta \in 
\mathbb{R}^{k}$ et  $(a_{n}>0)_{n\geq 1}$ une suite de nombres al\'eatoires tels que $a_{n}\rightarrow +\infty$ quand $n\rightarrow +\infty$.\\

\noindent Soit $g:D\rightarrow \mathbb{R}$ une fonction de classe $C^{1}$, telle que $D$ soit un domain de $\mathbb{R}^{k}$, que $\theta$ soit
$\overset{o}{D}$ soit l'int\'erieur de  $D$, et que l'on ait $\{X_{n},n\geq 1\}\subset \overset{o}{D}$.\\

\noindent Si $a_{n}(X_{n}-\theta)$ converge vaguement vers un vecteur al\'eatoire $Z$ de dimenstion $k$, lorsque $n\rightarrow +\infty$, alors nous avons

\begin{equation*}
a_{n}(g(X_{n})-g(\theta ))\rightsquigarrow \text{ }^{t}\nabla g(\theta
)Z=<\nabla g(\theta ),Z>\text{ as }n\rightarrow +\infty ,
\end{equation*}

\noindent o\`u 
\begin{equation*}
^{t}\nabla g(\theta )=(\frac{\partial g(\theta )}{\partial \theta _{1}},...,%
\frac{\partial g(\theta )}{\partial \theta _{k}})
\end{equation*}

\noindent est le vecteur gradient de la fonction $g$ en $\theta$.
\end{proposition}

\bigskip \noindent Le deuxi\`eme r\'esultat concerne la transformation de la convergence vague d'une suite suites variables al\'eatoires de $k$ composantes par une fonction de $m$ composantes de $k$ arguments. Ceci donne la version matricielle de la m\'ethode Delta.

\begin{proposition} \label{delta03} Soit $(X_{n})_{n\geq 1}$ une suite de vecteurs al\'eatoires de dimenstion $k$, $k\geq 1$, d\'efinis sur
le m\^eme espace de probabilit\'e $(\Omega ,A,\mathbb{P})$ et soit $\theta \in 
\mathbb{R}^{k}$ et  $(a_{n}>0)_{n\geq 1}$ une suite de nombres al\'eatoires tels que $a_{n}\rightarrow +\infty$ quand $n\rightarrow +\infty$.\\

\noindent Soit $g:D\rightarrow \mathbb{R}^{m}$ une fonction de classe $C^{1}$, telle que $D$ soit un domain de $\mathbb{R}^{k}$, que $\theta$ soit
$\overset{o}{D}$ soit l'int\'erieur de  $D$, et que l'on ait $\{X_{n},n\geq 1\}\subset \overset{o}{D}$. Notons par $g_j$, $1\leq j \leq m$, les composantes de la fonction $g$.\\

\noindent Si $a_{n}(X_{n}-\theta)$ converge vaguement vers un vecteur al\'eatoire $Z$ de dimension $k$, lorsque $n\rightarrow +\infty$, alors nous avons

\begin{equation*}
a_{n}(g(X_{n})-g(\theta ))\rightsquigarrow \text{ }g_{\theta }^{\prime }Z=%
\text{ as }n\rightarrow +\infty ,
\end{equation*}

\noindent o\`u $g_{\theta }^{\prime }$ est la matrice des d\'eriv\'ees partielles de premier ordre

\begin{equation*}
g_{\theta }^{\prime }=\left( 
\begin{tabular}{lllll}
$\frac{\partial g_{1}}{\partial \theta _{1}}$ & ... & $\frac{\partial g_{1}}{%
\partial \theta _{j}}$ & .. & $\frac{\partial g_{1}}{\partial \theta _{k}}$
\\ 
...  & ... & ... & ... & ... \\ 
$\frac{\partial g_{i}}{\partial \theta _{1}}$ & ... & $\frac{\partial g_{i}}{%
\partial \theta _{j}}$ & ... & $\frac{\partial g_{j}}{\partial \theta _{k}}$
\\ 
... & ... & ... & .... & ... \\ 
$\frac{\partial g_{m}}{\partial \theta _{1}}$ & ... & $\frac{\partial g_{m}}{%
\partial \theta _{j}}$ & ... & $\frac{\partial g_{m}}{\partial \theta _{k}}$%
\end{tabular}%
\right) ,
\end{equation*}
\end{proposition}

\bigskip \noindent \textbf{Preuve de la proposition \ref{delta02}}. Suppsons que les hypoth\`eses de la proposition aient lieu.\\

 \noindent Utilisons le d\'eveloppement de la fonction $g$ au premier ordre en  $\theta =$ $^{t}(\theta
_{1},...,\theta _{k})$ au point $x=$ $^{t}(x_{1},...,x_{k})$ :

\begin{equation}
g(x)-g(\theta )=(x_{1}-\theta _{1})\frac{\partial g}{\partial \theta _{1}}%
(\theta )+...+(x_{1}-\theta _{k})\frac{\partial g}{\partial \theta _{k}}%
(\theta )+o(\left\Vert x-\theta \right\Vert ).  \label{develop01}
\end{equation}

 \noindent Puisque $a_{n}(T_{n}-\theta )=$ $a_{n}^{t}((T_{1,n},...,T_{k,n})-^{t}(\theta
_{1},...,\theta _{k}))\rightsquigarrow Z=$ $^{t}(Z_{1},...,Z_{k})$, il suit de l'application du th\'eor\`eme de la transformation continue
\ref{cv.mappingTh} du chapitre \ref{cv}, que pour tout $1\leq i\leq k$,  $a_{n}(T_{j,n}-\theta _{j})$ converge vaguement vers $Z_{j}$ et alors par le point (a) du lemme \ref{oO-02}, nous avons
\begin{equation*}
1\leq j\leq k,(T_{j,n}-\theta _{j})=O_{\mathbb{P}}(a_{n}^{-1}).
\end{equation*}

 \noindent Par les points (10) and (13) des propri\'et\'es principales de la partie II de la section pr\'ec\'edente, nous obtenons

\begin{equation}
\left\Vert T_{n}-\theta \right\Vert =\left\{ \sum_{j=1}^{k}(T_{j,n}-\theta
_{j})^{2}\right\} ^{1/2}=O_{\mathbb{P}}(a_{n}^{-1})=o_{\mathbb{P}}(1).
\label{develop02}
\end{equation}

\noindent Le terme $o(\left\Vert x-\theta \right\Vert )$ dans \ref{develop01} est continu en tant que 
diff\'erence de deux functions continues et prend la valeur $0$ pour $\left\Vert x-\theta \right\Vert =0$. En utilisant la partie (c) du lemme 
\ref{oO-02}, et en combinant cela avec (\ref{develop01}) et (\ref{develop02}), nous arrivons \`a
 
\begin{eqnarray*}
a_{n}(g(x)-g(\theta )) &=&a_{n}(T_{1,n}-\theta _{1})\frac{\partial g}{%
\partial \theta _{1}}(\theta )+...+a_{n}(T_{k,n}-\theta _{k})\frac{\partial g%
}{\partial \theta _{k}}(\theta )+a_{n}o_{\mathbb{P}}(O_{\mathbb{P}%
}(a_{n}^{-1})) \\
&=&\text{ }^{t}\nabla g(\theta )(a_{n}(T_{n}-\theta )+o_{\mathbb{P}}(1).
\end{eqnarray*}

 \noindent Cela implique que $a_{n}(g(x)-g(\theta ))$ et $^{t}\nabla g(\theta
)(a_{n}(T_{n}-\theta )$ sont \'equivalents en probabilit\'e. Puisque  $^{t}\nabla
g(\theta )(a_{n}(T_{n}-\theta )\rightsquigarrow $ $^{t}\nabla g(\theta )Z$
en vertu du th\'eor\`eme de la transformation \ref{cv.mappingTh} of Chapter \ref{cv}, l'application de la proposition \ref{cvcp.prop4} de la Section \ref{cv.CvCp} du chapitre \ref{cv}, m\`ene \`a 

\begin{equation*}
a_{n}(g(x)-g(\theta ))\rightsquigarrow ^{t}\nabla g(\theta )Z.
\end{equation*}

\noindent Ceci finit la preuve de la proposition \ref{delta03}.\\

\bigskip \noindent \textbf{Preuve de la proposition \ref{delta03}}. Assume that the hypthesese of the proposition hold.\\

 \noindent La fonction $g$ poss\`ede  $m$ composantes  $g_{j} \in \mathbb{R}^{m}$ et nous pouvons \'ecrire $g=^{t}(g_{1},...,g_{m})$. Chaque composante est de classe $C^{1}$. Utilisons la conclusion de la proposition \ref{delta02} pour chaque composante en  $\theta =$ $%
^{t}(\theta _{1},...,\theta _{k})$ for $x=$ $^{t}(x_{1},...,x_{k})$ pour obtenir

\begin{equation}
g_{j}(x)-g_{j}(\theta )=(x_{1}-\theta _{1})\frac{\partial g_{j}}{\partial
\theta _{1}}(\theta )+...+(x_{1}-\theta _{k})\frac{\partial g_{j}}{\partial
\theta _{k}}(\theta )+o(\left\Vert x-\theta \right\Vert ).
\end{equation}

 \noindent Ecrivons cela sous forme matricielle ainsi qu'il suit 
\begin{equation*}
g(x)-g(\theta )=g_{\theta }^{\prime }(x-\theta )+o^{(m)}(\left\Vert x-\theta
\right\Vert ),
\end{equation*}

 \noindent \`u $o^{(m)}(\left\Vert x-\theta \right\Vert )$ est un vecteur de $m$
coordonn\'ees tel que chaque composante est une fonction continue, et aussi un  $
o(\left\Vert x-\theta \right\Vert )$. Une notation similaire est aussi utilis\'ee pour 
$o_{\mathbb{P}}(\circ)$. En appliquant la m\'ethode d\'ej\`a utilis\'ee dans la proposition \ref{delta02}, nous obtenons
 
\begin{equation*}
g(T_{n})-g(\theta )=g_{\theta }^{\prime }(T_{n}-\theta )+o_{\mathbb{P}%
}^{(m)}(a_{n}^{-1})
\end{equation*}

 \noindent et
\begin{equation*}
a_{n}(g(T_{n})-g(\theta ))=g_{\theta }^{\prime }a_{n}(T_{n}-\theta )+o_{%
\mathbb{P}}^{(m)}(1).
\end{equation*}

 \noindent Nous avons donc
\begin{equation*}
\left\Vert a_{n}(g(T_{n})-g(\theta ))-g_{\theta }^{\prime
}a_{n}(T_{n}-\theta )\right\Vert _{\mathbb{R}^{m}}=\left\Vert o_{\mathbb{P}%
}^{(m)}(1)\right\Vert _{\mathbb{R}^{m}}=o_{\mathbb{P}}(1).
\end{equation*}

 \noindent D\`es lors $a_{n}(g(T_{n})-g(\theta ))$ has the same weak limit as $g_{\theta
}^{\prime }a_{n}(T_{n}-\theta )$ which is $g_{\theta }^{\prime }Z$ by the
continuous mapping.\\

\newpage

\section{Utilisation du Processus Empirique Fonctionnel en Statistique Asymptotique} \label{cv.empTool.sec3}

\subsection{Processus Empirique Fonctionnel}

Le Processus Empirique Fonctionnel (\textit{PEF}) est un outil puissant qui peut \^etre utilis\'e pour obtenir des distributions limites. Il est similaire \`a la m\'ethode Delta. Cependant, la m\'ethode du \textit{PEF} poss\`ede un avantage que nous d\'ecrivons ci-bas.\\

\noindent ETant donn\'ee une suite de variables al\'eatoires $Z_{1}$, $Z_{2}$, ..., d\'efinies sur le m\^eme espace de probabilit\'e, ind\'ependantes et identiquement distribu\'ees,  de loi de probabilit\'e  commune $\mathbb{P}_{0}$ \`a valeurs dans un espace m\'etrique $S$, que nous prenons \'egal à $\mathbb{R}^k$, $k\geq 1$ ici. Nous serons en mesure \\

\noindent \textbf{(1)} de trouver un processus gaussien $\mathbb{G}_{\mathbb{P}_{0}}$ index\'e par des fonctions $f$ de $S$ dans $\mathbb{R}$,\\

\noindent et \\

\noindent \textbf{(2)} d'exprimer les distribution asymptotiques de statistiques qui sont des fonctions of  $Z_{1},Z_{2},...,Z_{n}$ par rapport au processus gaussien $\mathbb{G}_{0}$. \\

\bigskip \noindent Cette m\'ethode a l'avantage \'enorme  de pouvoir travailler sur toutes les statistiques fonctions de $Z_{1},Z_{2},...,Z_{n}$ de mani\`ere modulaire. En effet, on peut travailler de mani\`ere s\'epar\'ee pour chacune de ces statistiques comme and un catalogue. Et \`a tout moment, on peut prendre un certain nombre de ces statistiques et obtenir imm\'eduatement leur loi conjointe limite, gr\^ace au champ gaussien 
$\mathbb{G}_{\mathbb{P}_{0}}$. \textbf{Nous dirons que ces dsitributions limites sont exprim\'ees dans le champ gaussien de $\mathbb{G}_{\mathbb{P}_{0}}$}. Au contraire la m\'ethode Delta ne poss\`ede pas cette souplesse modulaire. Dans chaque cas, il faut faire des calculs sp\'ecifiques. \\

\noindent Un autre fait \`a signaler est que les lois limites conjointes obtenues par la m\'ethode du \textit{PEF} utilisent des variances et covariances exprim\'ees sous des formes fonctionnelles. Cela rend les calculs directs et fluides. A l'arriv\'ee, ces variances peuvent sembler compliqu\'ees et non calculables \`a la main. Mais nous ne nous en occupons pas. Les ordinateurs puissants de notre \`ere sont là pour les calculer en des temps tr\`es courts. Cela implique que l'utilisation de cette m\'ethode, dans beaucoup de cas, doit s'accompagner avec l'\'ecriture de codes dans des logiciels courants.\\

\noindent  Avant de pr\'esenter le Processus Empirique Fonctionnel, nous voudrions rassurer le lecteur que nous n'utiliserons ici que les limites en distributions finies de ce processus, restant ainsi dans le cadre vectoriel de ce monograph. Nous ne ferons pas appel aux outils puissants de convergence uniforme ni aux classes de Vapnick-Cervonenkis.\\

\noindent Soit $Z_{1}$, $Z_{2}$, ... une suite de copies ind\'ependantes d'une variable al\'eatoire $Z$ d\'efinies sur le m\^eme espace de 
probabilit\'e  $(\Omega ,\mathcal{A},\mathbb{P})$ \`a valeurs dans un espace m\'etrique $(S,d)$. D\'efinissons pour chaque $n\geq 1$, le Processus Empirique Fonctionnel (\textit{PEF}) par  
\begin{equation*}
\mathbb{G}_{n}(f)=\frac{1}{\sqrt{n}}\sum_{i=1}^{n}(f(Z_{i})-\mathbb{E}
f(Z_{i})),
\end{equation*}

\bigskip \noindent o\`u  $f$ est une fonction mesurable d\'efinie de $\mathbb{R}$ dans 

\begin{equation}
\mathbb{V}_{Z}(f)=\int \left( f(x)-\mathbb{P}_{Z}(f)\right)
^{2}dP_{Z}(x)<\infty ,  \label{fep.var}
\end{equation}

\noindent ce qui implique que 

\begin{equation}
\mathbb{P}_{Z}(\left\vert f\right\vert )=\int \left\vert f(x)\right\vert
dP_{Z}(x)<\infty \text{.}  \label{fep.esp}
\end{equation}

\bigskip \noindent Notons par  $\mathcal{F}(S)$ - $\mathcal{F}$ en court - la class des fonctions r\'eelles d\'efinies sur $S$ v\'erifiant 
(\ref{fep.var}). Cette espace $\mathcal{F}$, muni de l'addition et de la multiplication externe par des scalaires de $\mathbb{R}$, est un espace linaire.\\

\noindent L'une des propri\'et\'es les plus importantes de $\mathbb{G}_{n}$ est qu'elle est un op\'erateur lin\'eaire sur $\mathcal{F}$, i.e.,
pour tout $f$ et $g$ \'el\'ements de $\mathcal{F}$ et pour tout $(a,b)\in \mathbb{R}{^{2}}$, nous avons

\begin{equation*}
a\mathbb{G}_{n}(f)+b\mathbb{G}_{n}(g)=\mathbb{G}_{n}(af+bg).
\end{equation*}

\bigskip \noindent La convergence vague en distribution est la suivante

\begin{lemma} \label{fep.lemma.tool.1} Given the notation above, then for any finite
number of elements $f_{1},...,f_{p}$ of $\mathcal{S}$, $k\geq 1$, we have

\begin{equation*}
^{t}(\mathbb{G}_{n}(f_{1}),...,\mathbb{G}_{n}(f_{k}))\rightsquigarrow 
\mathcal{N}_{k}(0
\end{equation*}

\bigskip \noindent where 
\begin{equation*}
\Gamma (f_{i},f_{j})=\int \left( f_{i}-\mathbb{P}_{Z}(f_{i})\right) \left(
f_{j}-\mathbb{P}_{Z}(f_{j})\right) d\mathbb{P}_{Z}(x),1\leq ,j\leq k.
\end{equation*}
\end{lemma}

\bigskip \noindent Ce lemme nous dit que la limite vague de la suite de vecteurs $^{t}$ $(\mathbb{G}_{n}(f_1), \mathbb{G}_{n}(f_2), ...,\mathbb{G}_{n}(f_k))$ a la m\^eme loi que le vecteur $^{t}$ $(\mathbb{G}(f_1), \mathbb{G}(f_2), ...,\mathbb{G}(f_k))$, o\`u $\{\mathbb{G}(f),f\in \mathcal{F}\}$ est un process gaussien de fonction variance-covariance

\begin{equation}
\Gamma (f,g)=\int \left( f-\mathbb{P}_{Z}(f)\right) \left(
g-\mathbb{P}_{Z}(g)\right) d\mathbb{P}_{Z}(x), \ \ (f,g)\in \mathcal{F}^2. \label{cv.pefTools.cov01}
\end{equation}

\noindent En appliquant le th\'eor\`eme de Skorohod-Wichura \ref{cv.skorohodWichura} (Voir Chapitre \ref{cv}), nous pouvons supposer que nous sommes sur un espace de probabilit\'e et que nous avons l'approximation suivante

\begin{equation}
\mathbb{G}_{n}(f_{1})=\mathbb{G}_{n}(f_{1})+o_{\mathbb{P}}(1). 1\leq i \leq p. \label{cv.pefTools.rep01}
\end{equation}

\noindent Nous reviendrons sur l'application de cette formule.\\

\bigskip \noindent \textbf{Preuve du lemme \ref{fep.lemma.tool.1}}. Il suffira d'appliquer le crit\`ere de Cram\'er-Wold 
(voir Proposition \ref{cv.wold} in \ref{ChapRevCvRk}), en montrant que pour tout $a=^{t}(a_{1},...,a_{k})\in \mathbb{R}^{k}$, nous avons 

$$
<a,T_{n}>\rightsquigarrow <a,T>
$$

\noindent o\`u nous avons not\'e $T_{n}=^{t}(\mathbb{G}_{n}(f_{1}),...,\mathbb{G}_{n}(f_{k}))$ et o\`u $T$ suit une loi $\mathcal{N}
_{k}(0,\Gamma (f_{i},f_{j})_{1\leq i,j\leq k})$\  et  $<\circ ,\circ >$
d\'esigne le produit scalaire usuel de $\mathbb{R}^{k}$.\\

\noindent Mais, par le th\'eor\`eme central limite dans $\mathbb{R}$, nous avons

\begin{equation*}
<a,T_{n}>=\mathbb{G}_{n}\left( \sum\limits_{i=1}^{k}a_{i}f_{i}\right)
\rightsquigarrow \mathcal{N}(0,\sigma _{\infty }^{2}),
\end{equation*}

\bigskip \noindent o\`u, avecla notation $g=\sum_{1\leq i\leq k}a_{i}f_{i}$, nous avons
\begin{equation*}
\sigma _{\infty }^{2}=\int \left( g(x)-\mathbb{P}_{Z}(g)\right) ^{2}dP_{Z}(x)
\end{equation*}

\bigskip \noindent et ceci donne ais\'ement
\begin{equation*}
\sigma _{\infty }^{2}=\sum\limits_{1\leq i,j\leq k}a_{i}a_{j}\Gamma
(f_{i},f_{j}),
\end{equation*}

\noindent si bien que $N(0,\sigma _{\infty }^{2})$ est la loi de $<a,T>$. La preuve est finie.\\

\subsection{Comment utiliser l'outil du \textit{PEF}?}

\bigskip \noindent En g\'en\'eral, les statistiques avec lesquelles nous travaillons en g\'en\'eral utilisent des donn\'ees univari\'ees ou multivari\'ees ou tics, c'est-\`a-dire que nous travaillons souvent dans $\mathbb{R}^{k}$. Une fois que nous avons notre \'echatillon  $Z_{1},Z_{2},...$, constitu\'e de vecteurs d\'efinis sur le m\^eme espace de probabilit\'e et \`a valeurs dans $\mathbb{R}^{k}$, l'essentiel des statistiques \'etudi\'ees 
sont de la forme 
\begin{equation*}
H_{n}=\frac{1}{n}\sum\limits_{i=1}^{k}H(Z_{i})
\end{equation*}

\bigskip \noindent o\`ur $H\in \mathcal{F}$. Dans un tel cas, nous utilisons Lemma \ref{fep.lemma.tool.1} et le point point (a) du lemme 
\ref{oO-02}, pour avoir ce d\'ev\'elopment tr\`es simple 
$\mu (H)=\mathbb{E}H(Z),$ 
\begin{equation}
H_{n}=\mu (H)+n^{-1/2}\mathbb{G}_{n}(H).  \label{fep.expan}
\end{equation}

\bigskip \noindent Nous avons que $\mathbb{G}_{n}(H)$ est asymptotiquement born\'e en probabilit\'e puisque  $\mathbb{G}_{n}(H)$ converge vaguement, disons vers $M(H)$, et par la suite, par le th\'eor\`eme de la transformation continue (see Theorem \ref{cv.mappingTh}, Chapitre \ref{cv}), nous aurons  $\left\Vert \mathbb{G}_{n}(H)\right\Vert \rightsquigarrow \left\Vert M(H)\right\Vert$. Puisque toutes variables $\mathbb{G}_{n}(H)$ sont d\'efinies sur le m\^eme espace de probabilit\'e, nous aurons en vertu de l'assertion du th\'eor\`eme Portmanteau (Theorem \ref{cv.theo.portmanteau}, Chapitre \ref{cv}) relative aux ensembles ouverts, que pour tout $\lambda >0$,

\begin{equation*}
\lim \sup_{n\rightarrow \infty }P(\left\Vert \mathbb{G}_{n}(H)\right\Vert
>\lambda )\leq P(\left\Vert M(H)\right\Vert >\lambda).
\end{equation*}

\bigskip \noindent Par la suite, nous avons
\begin{equation*}
\liminf_{\lambda \rightarrow \infty}\limsup_{n\rightarrow \infty
}\mathbb{P}(\left\Vert \mathbb{G}_{n}(H)\right\Vert >\lambda )\leq \liminf_{\lambda \rightarrow \infty } \mathbb{P}(\left\Vert M(H)\right\Vert >\lambda )=0.
\end{equation*}

\bigskip \noindent A partir d'ici, nous faisons appel aux notation $O_{\mathbb{P}}$, en disant que nous avons : $\mathbb{G}_{n}(H)=O_{\mathbb{P}}(1)$. La formule (\ref{fep.expan}) devient 
\begin{equation*}
H_{n}=\mu (H)+n^{-1/2}\mathbb{G}_{n}(H)=\mu (H)+O_{\mathbb{P}}(n^{-1/2})
\end{equation*}

\bigskip \noindent A pr\'esent, utilisons la m\'ethode delta. En effet, soit 
$g:\mathbb{R}\longmapsto \mathbb{R}$ une function contin\^ument diff\'erentiable dans un voisinage de $\mu (H)$. The th\'eor\`eme des accroissements finis m\`ene à
\begin{equation}
g(H_{n})=g(\mu (H))+g^{\prime }(\mu _{n}(H))\text{ }n^{-1/2}\mathbb{G}_{n}(H)
\label{fep.expan01}
\end{equation}

\bigskip \noindent o\`u 
\begin{equation*}
\mu _{n}(H)\in \lbrack (\mu (H)+n^{-1/2}\mathbb{G}_{n}(H))\wedge \mu
(H),(\mu (H)+n^{-1/2}\mathbb{G}_{n}(H))\vee \mu (H)]
\end{equation*}

\bigskip \noindent si bien que 
\begin{equation*}
\left\vert \mu _{n}(H)-\mu (H)\right\vert \leq n^{-1/2}\mathbb{G}_{n}(H)=O_{%
\mathbb{P}}(n^{-1/2}).
\end{equation*}

\bigskip \noindent D\`es lors $\mu _{n}(H)$ converge vers $\mu _{n}(H)$ en probabilit\'e (not\'e $\mu _{n}(H)$ $\rightarrow _{\mathbb{P}}\mu (H))$. Mais la convergence en probabilit\'e vers une constante \'equivaut \`a la convergence vague. Donc $\mu _{n}(H)$ $\rightsquigarrow \mu (H)$. En utilisant encore  le th\'eor\`eme de la transformation continue (see Theorem \ref{cv.mappingTh}, Chapitre \ref{cv}), nous obtenons $g^{\prime }(\mu _{n}(H))\rightsquigarrow g^{\prime }(\mu (H))$. En r\'e-utilisant l'\'equivalence entre convergence vague et convergence en probabilit\'e vers une constante, il vient que  $g^{\prime }(\mu_{n}(H))\rightarrow _{\mathbb{P}}g^{\prime }(\mu (H))$. Maintenant, (\ref{fep.expan01}) devient

\begin{eqnarray*}
g(H_{n}) &=&g(\mu (H))+(g^{\prime }(\mu (H)+o_{P}(1))\text{ }n^{-1/2}\mathbb{%
G}_{n}(H) \\
&=&g(\mu (H))+g^{\prime }(\mu (H)\times \text{ }n^{-1/2}\mathbb{G}%
_{n}(H)+o_{P}(1))\text{ }n^{-1/2}\mathbb{G}_{n}(H) \\
&=&g(\mu (H))+\text{ }n^{-1/2}\mathbb{G}_{n}(g^{\prime }(\mu
(H)H)+o_{P}(n^{-1/2})
\end{eqnarray*}

\bigskip \noindent Nous arrivons au d\'eveloppement final
\begin{equation}
g(H_{n})=g(\mu (H))+\text{ }n^{-1/2}\mathbb{G}_{n}(g^{\prime }(\mu
(H)H)+o_{P}(n^{-1/2}).  \label{fep.expanFinal}
\end{equation}

\noindent En utilisant le repr\'esentation de Skorohod-Wichura, nous obtenons \`a travers la formule \label{cv.pefTools.rep01}, que

\begin{equation}
g(H_{n})=g(\mu (H))+\text{ }n^{-1/2}\mathbb{G}(g^{\prime }(\mu
(H)H)+o_{P}(n^{-1/2}).  \label{fep.expanFinal01}
\end{equation}

\bigskip \noindent La m\'ethode consiste \`a utiliser le d\'eveloppement (\ref{fep.expanFinal}) autant de fois que possible, et de faire certains calculs alg\'ebriques suir eux.\\

\noindent Ces calculs alg\'ebriques mentionn\'es ci-haut consistent par ailleurs \`a appliquer les r\'esultats du lemma ci-dessous.

\begin{lemma}
\label{fep.lemma.tool.2} Let ($A_{n})$ and ($B_{n})$ be two sequences of real
valued random variables defined on the same probability space holding the
sequence $Z_{1}$, $Z_{2}$, ...

\noindent Let $A$ and $B$ be two real numbers and Let $L(z)$
and $H(z)$ be two real-valued functions$\ of$ $z\in S$, with $(L,H)\in \mathcal{F}^2$.\\

\noindent Suppose that 
$$
A_{n}=A+n^{-1/2}\mathbb{G}_{n}(L)+o_{P}(n^{-1/2})
$$ 

\noindent and 

$$
A_{n}=B+n^{-1/2}\mathbb{G}_{n}(H)+o_{P}(n^{-1/2}).
$$ 

\noindent Then

\begin{equation*}
A_{n}+B_{n}=A+B+n^{-1/2}\mathbb{G}_{n}(L+H)+o_{P}(n^{-1/2}),
\end{equation*}

\noindent and

\begin{equation*}
A_{n}B_{n}=AB+n^{-1/2}\mathbb{G}_{n}(BL+AH)
\end{equation*}

\noindent and if $B\neq 0$,
\begin{equation*}
\frac{A_{n}}{B_{n}}=\frac{A}{B}+n^{-1/2}\mathbb{G}_{n}(\frac{1}{B}L-\frac{A}{%
B^{2}}H)+o_{P}(n^{-1/2})
\end{equation*}
\end{lemma}

\bigskip \noindent \textbf{Preuve}. Cette preuve est donn\'ee au lecteur en guise d'exercices.\\

\bigskip \noindent En mettant ensemble toutes les \'etapes pr\'ec\'edente de mani\`ere intelligente, la m\'ethodologie aboutit à un r\'esultat de la forme
\begin{eqnarray*}
T_{n}&=&t+n^{-1/2}\mathbb{G}_{n}(h)+o_{P}(n^{-1/2})\\
&=& t+n^{-1/2}\mathbb{G}(h)+o_{P}(n^{-1/2})
\end{eqnarray*}

\noindent qui entra\^ine la convergence vague

\begin{eqnarray*}
\sqrt{n}(T_{n}-t)&=&\mathbb{G}_{n}(h)+o_{P}(1)\rightsquigarrow \mathbb{N}(0,\Gamma(h,h))\\
&=&\mathbb{G}(h)+o_{P}(1).
\end{eqnarray*}

\bigskip \noindent Dans cette derni\`ere partie, nous allons montrer comment appliquer cet outil \`a un cas non relatif au coefficient de 
corr\'elation lin\'eaire.

\subsection{Un example} \label{fep.subsec4}

\noindent Appliquons la m\'ethode \`a l'estimateur plug-in coefficient de 
corr\'elation lin\'eaire. Un estimateur plug-in d'une statistique est sa forme empirique. L'estimateur plug-in
du coefficient de corr\'elation lin\'eaire de deux variables al\'eatoires r\'eelles $(X,Y)$, telle que ni $X$ et ni $Y$ ne sont d\'eg\'en\'er\'ees, est d\'efini par
\begin{equation*}
\rho =\frac{\sigma _{xy}}{\sigma _{x}^{2}\sigma _{y}^{2}}
\end{equation*}

\noindent o\`u

\begin{equation*}
\mu _{x}=\int x\text{ }dP_{X}(x),\text{ }\mu _{y}=\int x\text{ }dP_{X}(x),%
\text{ }\sigma _{xy}=\int (x-\mu _{x})(y-\mu _{y})dP_{(X,Y)}(x,y)
\end{equation*}

\noindent et

\begin{equation*}
\sigma _{x}^{2}=\int (x-\mu _{x})^{2}dP_{X}(x),\text{ }\sigma _{y}^{2}=\int
(x-\mu _{x})(y-\mu _{y})dP_{X}(y).
\end{equation*}

\noindent Nous ecartons aussi le cas o\`u nous avons $\left\vert \rho \right\vert =1$, cas dans lequel une des variables
$X$ et $Y$ est une fonction affine de l'autre, c'est-\`a-dire par exemple que nous avons $X=aY+b$ pour $(a,b)\mathbb{R}^2$.\\

\noindent Il est clair que le centrage des variables  $X$ et  $Y$ en leur moyennes et leur normalisation par les \'ecart-types  $\sigma _{x}$ et $\sigma _{y}$ n'affecte pas le coefficient de corr\'elation lin\'eaire $\rho$. Ainsi nous pouvons centrer et normaliser $X$ et $Y$ et par la suite, assumer que nous avons

\begin{equation*}
\mu _{x}=\text{ }\mu _{y}=0,\text{ }\sigma _{x}=\sigma _{y}=1.
\end{equation*}

\bigskip \noindent Cependant, nous laisserons figurer les moyennes et \'ecart-types. C'est seulement à la conclusion que nous prendrons les valeurs particuli\`eres.
\newline

\bigskip \noindent Construisons l'estimateur plug-in estimator de $\rho $. Pour cela, soit $(X_{1},Y_{1}),$ $(X_{2},Y_{2}),...$ 
une suite d'observations ind\'ependantes de observations of $(X,Y)$. Pour tout $n\geq 1$, l'estimateur plug-in est 

\begin{equation*}
\rho _{n}=\left\{ \frac{1}{n}\sum_{i=1}^{n}(X_{i}-\overline{X})(Y_{i}-%
\overline{Y})\right\} \left\{ \frac{1}{n^{2}}\sum_{i=1}^{n}(X_{i}-\overline{X%
})^{2}\times \sum_{i=1}^{n}(X_{i}-\overline{X})^{2}\right\} ^{-1/2}.
\end{equation*}

\bigskip \noindent Nous allons donner les propri\'et\'es asymptotiques de $\rho _{n}$
en tant qu'estimateur de  $\rho$.  Introduisons les notations :

\begin{equation*}
\mu _{(p,x),(q,y)}=E((X-\mu _{x})^{p}(Y-\mu _{y})^{q}),\mu _{4,x}=E(X-\mu
_{x})^{4}\text{, }\mu _{4,x}=E(X-\mu _{x})^{4})
\end{equation*}

\bigskip \noindent L'application de la m\'ethode donne le r\'esultat suivant.

\begin{theorem}
\label{fep.theo1}

\bigskip \noindent Supposons ni $X$ ni $Y$ ne soit d\'eg\'en\'er\'ee et ont toutes leus deux des moments finis jusqu'à l'ordre 4, que  $X^{3}Y$ et $XY^{3}$ soit de moyenne finies. Alors, quand  $n\rightarrow \infty ,$%
\begin{equation*}
\sqrt{n}(\rho _{n}-\rho )\rightsquigarrow N(0,\sigma ^{2}),
\end{equation*}

\bigskip \noindent o\`u
\begin{eqnarray*}
\sigma ^{2} &=&\sigma _{x}^{-2}\sigma _{y}^{-2}(1+\rho ^{2}/2)\mu
_{(2,x),(2,y)}+\rho ^{2}(\sigma _{x}^{-4}\mu _{4,x}+\sigma _{y}^{-4}\mu
_{4,y})/4 \\
&&-\rho (\sigma _{x}^{-3}\sigma _{y}^{-1}\mu _{(3,x),(1,y)}+\sigma
_{x}^{-1}\sigma _{y}^{-3}\mu _{(1,x),(3,y)})
\end{eqnarray*}
\end{theorem}

\bigskip \noindent Ce r\'esultat permet de tester l'ind\'ependance entre $X$ et  $Y$, ou de tester l'absence de correlation
lin\'eaire dans le sens suivant.

\bigskip

\begin{theorem}
\label{fep.theo2} Supposons que les hypoth\`eses du th\'eor\`eme \ref{fep.theo1} soitent vraies. Alors, nous avons les assertions suivantes :
\newline
\bigskip \noindent \textbf{(1)} Si $X$ et  $Y$ ne sont pas lin\'eairement correl\'ees, i.e.  $\rho=0$, nous avons 
\begin{equation*}
\sqrt{n}\rho _{n}\rightsquigarrow \mathcal{N}(0,\sigma _{1}^{2}),
\end{equation*}

\bigskip \noindent o\`u
\begin{equation*}
\sigma _{1}^{2}=\sigma _{x}^{-2}\sigma _{y}^{-2}\mu _{(2,x),(2,y)}.
\end{equation*}

\bigskip \noindent \textbf{(2)} Si $X$ et  $Y$ sont ind\'ependantes entre elles, alors $\rho
=0$, et \bigskip 
\begin{equation*}
\sqrt{n}\rho _{n}\rightsquigarrow \mathcal{N}(0,1)
\end{equation*}
\end{theorem}

\bigskip \noindent \textbf{Preuves}. Nous allons utiliser le \textit{PEF} bas\'e sur les observations  $(X_{i},Y_{i}),i=1,2,...$, qui sont des copies i d\'ependantes de $(X,Y)$. Ecrivons

\begin{equation*}
\rho _{n}^{2}=\frac{\frac{1}{n}\sum_{i=1}^{n}X_{i}Y_{i}-\overline{X}\text{ }%
\overline{Y}}{\left\{ \frac{1}{n}\sum_{i=1}^{n}X_{i}^{2}-\overline{X}%
^{2}\right\} ^{1/2}\left\{ \frac{1}{n}\sum_{i=1}^{n}Y_{i}^{2}-\overline{Y}%
^{2}\right\} ^{1/2}}=\frac{A_{n}}{B_{n}}.
\end{equation*}

\bigskip \noindent Nous disons une fois pour toute que toutes les fonctions de $Z=(X,Y)$
qui apparaissent dans les calculs sont mesurables et ont des seconds moments finis. Traitons s\'epar\'ement
le num\'erateur et le d\'enominateur. l'utilisation du \textit{PEF} pour $A_{n}$ gives

\begin{equation}
\left\{ 
\begin{tabular}{l}
$\frac{1}{n}\sum_{i=1}^{n}X_{i}Y_{i}=\mu _{xy}+n^{-1/2}G_{n}(p),$ \\ 
$\overline{X}=\mu _{x}+n^{-1/2}G_{n}(\pi _{1}),$ \\ 
$\overline{Y}=\mu _{y}+n^{-1/2}G_{n}(\pi _{2}),$%
\end{tabular}%
\right.  \label{fep.casAN}
\end{equation}

\bigskip \noindent o\`u $p(x,y)=xy$, $\pi _{1}(x,y)=x,$ $\pi _{2}(x,y)=y$. A partir de l\`a, nous utilisons le fait que
 $G_{n}(g)=O_{P}(1)$ pour $\mathbb{E}(g(X,Y)^{2})<+\infty$ et obtenons

\begin{equation}
A_{n}=\mu _{xy}+n^{-1/2}G_{n}(p)-(\mu _{x}+n^{-1/2}G_{n}(\pi _{1}))(\mu
_{y}+n^{-1/2}G_{n}(\pi _{2})).  \label{fep.haut}
\end{equation}

\bigskip \noindent Cela m\`ene \`a 
\begin{equation*}
A_{n}=\sigma _{xy}+n^{-1/2}G_{n}(H_{1})+o_{P}(n^{-1/2})
\end{equation*}

\bigskip \noindent avec
\begin{equation*}
H_{1}(x,y)=p(x,y)-\mu _{x}\pi _{2}-\mu _{y}\pi _{1.}
\end{equation*}

\bigskip \noindent Maintenant, nous devons traiter $B_{n}$. Par sym\'etrie des roles de $\left\{ \frac{1}{n}\sum_{i=1}^{n}X_{i}^{2}-\overline{X}^{2}\right\} ^{1/2}$ et de $\left\{ \frac{1}{n}\sum_{i=1}^{n}Y_{i}^{2}-\overline{Y}^{2}\right\}
^{1/2}$, nous traitons l'une de ces expressions et \'etendons les r\'esultats à l'autre. Travaillons avec $\left\{ \frac{1}{n}\sum_{i=1}^{n}X_{i}^{2}-\overline{X}^{2}\right\} ^{1/2}$. La combinaison de (\ref{fep.casAN}) et de la m\'ethode Delta conduit \`a
\begin{equation*}
\overline{X}^{2}=\left( \mu _{x}+n^{-1/2}G_{n}(\pi _{1})\right) ^{2}=\mu
_{x}^{2}+2\mu _{x}n^{-1/2}G_{n}(\pi _{1})+o_{P}(n^{-1/2}),
\end{equation*}

\bigskip \noindent c'est-\`a-dire,
\bigskip
\begin{equation*}
\overline{X}^{2}=\left( \mu _{x}+n^{-1/2}G_{n}(\pi _{1})\right) ^{2}=\mu
_{x}^{2}+n^{-1/2}G_{n}(2\mu _{x}\pi _{1})+o_{P}(n^{-1/2}).
\end{equation*}

\bigskip \noindent Nous obtenons
\begin{eqnarray*}
\frac{1}{n}\sum_{i=1}^{n}X_{i}^{2}-\overline{X}^{2}
&=&m_{2,x}+n^{-1/2}G_{n}(\pi _{1}^{2})-\overline{X}^{2} \\
&=&m_{2,x}-\mu _{x}^{2}+n^{-1/2}G_{n}(\pi _{1}^{2}-2\mu _{x}\pi
_{1})+o_{P}(n^{-1/2}) \\
&=&\sigma _{x}^{2}+n^{-1/2}G_{n}(\pi _{1}^{2}-2\mu _{x}\pi
_{1})+o_{P}(n^{-1/2}).
\end{eqnarray*}

\bigskip \noindent Une nouvelle application de la methode Delta donne 
\begin{equation*}
\left\{ \frac{1}{n}\sum_{i=1}^{n}X_{i}^{2}-\overline{X}^{2}\right\}
^{1/2}=\sigma _{x}+n^{-1/2}G_{n}(\frac{1}{2\sigma _{x}}\left\{ \pi
_{1}^{2}-2\mu _{x}\pi _{1}\right\} )+o_{P}(n^{-1/2}).
\end{equation*}

\bigskip \noindent Etendons les r\'esultats sur les $Y_i$ pour avoir
\begin{equation*}
\left\{ \frac{1}{n}\sum_{i=1}^{n}Y_{i}^{2}-\overline{Y}^{2}\right\}
^{1/2}=\sigma _{y}+n^{-1/2}G_{n}(\frac{1}{2\sigma _{y}}\left\{ \pi
_{2}^{2}-2\mu _{y}\pi _{2}\right\} )+o_{P}(n^{-1/2}).
\end{equation*}

\bigskip \noindent Nous arrivons \`a
\begin{eqnarray*}
B_{n} &=&\left\{ \frac{1}{n}\sum_{i=1}^{n}X_{i}^{2}-\overline{X}^{2}\right\}
^{1/2}\left\{ \frac{1}{n}\sum_{i=1}^{n}Y_{i}^{2}-\overline{Y}^{2}\right\}
^{1/2} \\
&=&\sigma _{x}\sigma _{y}+n^{-1/2}G_{n}(\frac{\sigma _{y}}{2\sigma _{x}}%
\left\{ \pi _{1}^{2}-2\mu _{x}\pi _{1}\right\} +\frac{\sigma _{x}}{2\sigma
_{y}}\left\{ \pi _{2}^{2}-2\mu _{y}\pi _{2}\right\} )+o_{P}(n^{-1/2}).
\end{eqnarray*}

\bigskip \noindent En mettant
\begin{equation*}
H_{2}(x,y)=\frac{\sigma _{y}}{2\sigma _{x}}\left\{ \pi _{1}^{2}-2\mu _{x}\pi
_{1}\right\} +\frac{\sigma _{x}}{2\sigma _{y}}\left\{ \pi _{2}^{2}-2\mu
_{y}\pi _{2}\right\},
\end{equation*}

\bigskip \noindent nous avons

\begin{equation}
B_{n}=\sigma _{x}\sigma _{y}+n^{-1/2}G_{n}(H_{2})+n^{-1/2}.  \label{fep.bas}
\end{equation}

\bigskip \noindent Maintenant, la combinaison de (\ref{fep.haut}) et de (\ref{fep.bas}), et l'utilisation du lemme \ref{fep.lemma.tool.2} donne

\begin{equation*}
\sqrt{n}(\rho _{n}^{2}-\rho ^{2})=n^{-1/2}G_{n}(\frac{1}{\sigma _{x}\sigma
_{y}}H_{1}-\frac{\sigma _{xy}}{\sigma _{x}^{2}\sigma _{y}^{2}}%
H_{2})+o_{P}(1).
\end{equation*}

\bigskip \noindent Posons
\begin{equation*}
H=\frac{1}{\sigma _{x}\sigma _{y}}(p(x,y)-\mu _{x}\pi _{2}-\mu _{y}\pi _{1})-%
\frac{\rho }{\sigma _{x}\sigma _{y}}\left\{ \frac{1}{2\sigma _{x}^{2}}%
\left\{ \pi _{1}^{2}-2\mu _{x}\pi _{1}\right\} +\frac{1}{2\sigma _{y}^{2}}%
\left\{ \pi _{2}^{2}-2\mu _{y}\pi _{2}\right\} \right\}.
\end{equation*}

\bigskip \noindent Nous continuons avec les variables centr\'ees et normalis\'ees. Nous avons

\begin{equation*}
H(x,y)=p(x,y)-\frac{\rho }{2}(\pi _{1}^{2}+\pi _{2}^{2})
\end{equation*}

\bigskip \noindent et
\begin{equation*}
H(X,Y)=XY-\frac{\rho }{2}(X^{2}+Y^{2}).
\end{equation*}

\bigskip \noindent Adoptons la notation

\begin{equation*}
\mu _{(p,x),(q,y)}=E((X-\mu _{x})^{p}(Y-\mu _{y})^{q}).
\end{equation*}

\bigskip \noindent Nous obtenons
\begin{equation*}
\mathbb{E}H(X,Y)=\sigma _{xy}-\rho =0
\end{equation*}

\bigskip \noindent avec la remarque que $varH(X,Y)$ is \'egale \`a 
\begin{equation*}
\mu _{(2,x),(2,y)}+\rho ^{2}(\mu _{4,x}+\mu _{4,y})/4-\rho (\mu
_{(3,x),(1,y)}+\mu _{(1,x),(3,y)})+\rho ^{2}\mu _{(2,x),(2,y)}/2
\end{equation*}

\noindent et finalemeny 

$$
varH(X,Y)=\sigma _{0}^{2}
$$ 

\noindent avec
\begin{equation*}
\sigma _{0}^{2}=(1+\rho ^{2}/2)\mu _{(2,x),(2,y)}+\rho ^{2}(\mu _{4,x}+\mu
_{4,y})/4-\rho (\mu _{(3,x),(1,y)}+\mu _{(1,x),(3,y)}).
\end{equation*}

\bigskip \noindent Ceci donne la conclusion pour des variables
$X$ et $Y$ normalis\'ees :

\begin{equation*}
\sqrt{n}(\rho _{n}-\rho )\rightsquigarrow N(0,\sigma _{0}^{2}).
\end{equation*}

\bigskip \noindent Enfin, lorsque nous utilisons les coefficients de normalisation dans in $\sigma
_{0}$, nous obtenons

\begin{eqnarray*}
\sigma ^{2} &=&\sigma _{x}^{2}\sigma _{y}^{2}(1+\rho ^{2}/2)\mu
_{(2,x),(2,y)}+\rho ^{2}(\sigma _{x}^{4}\mu _{4,x}+\sigma _{y}^{4}\mu
_{4,y})/4 \\
&&-\rho (\sigma _{x}^{3}\sigma _{y}\mu _{(3,x),(1,y)}+\sigma _{x}\sigma
_{y}^{3}\mu _{(1,x),(3,y)})
\end{eqnarray*}

\bigskip \noindent pour conclure dans le cas g\'en\'eral que 

\begin{equation*}
\sqrt{n}(\rho _{n}-\rho )\rightsquigarrow N(0,\sigma ^{2})
\end{equation*}

\bigskip \noindent La preuve du th\'eor\`eme \ref{fep.theo2} d\'ecoule à la suite calculs directs sous les conditions particuli\`eres de $\rho $ et en pr\'esence de l'ind\'ependance.

%% file: asymptotics_math_01_fr.tex
\chapter{Th\'eorie des fonctions et \'el\'ements d'analyse r\'eelle \`a travers des exercices}
\label{funct}


\section{Revue des limites dans $\overline{\mathbb{R}}$. Ce que  nous ne devons pas ignirer des limites.} \label{funct.sec.1}

\noindent \textbf{Definition}:Un nombre r\'eel $\ell \in \overline{\mathbb{R}}$ est un point d'accumulation d'une suite $(x_{n})_{n\geq 0}$ finie ou infinie de r\'eel  , dans $\overline{\mathbb{R}}$,si et seulement si il existe  une sous-suite $(x_{n(k)})_{k\geq 0}$ de $(x_{n})_{n\geq 0}$ telle que $x_{n(k)}$ converge vers $\ell $, quand $k\rightarrow +\infty $.\newline
\bigskip

\noindent \textbf{Exercise 1 : } Soient les ensembles $y_{n}=\inf_{p\geq n}x_{p}$ et $z_{n}=\sup_{p\geq n}x_{p} $ d\'efinis pour tout $n \geq 0$. Montrer que :\newline
\bigskip

\noindent \textbf{(1)} $\forall n\geq 0,y_{n}\leq x_{n}\leq z_{n}$\newline
\bigskip

\noindent \textbf{(2)} Justifier l'existence de la limite de $y_{n}$ appel\'ee limite inf\'erieure de la suite $(x_{n})_{n\geq 0}$, d\'efinie par $\lim inf x_{n}$ ou $\underline{\lim }$ $x_{n},$ et \'egale \`a l'expression ci-dessous
\begin{equation*}
\underline{\lim }\text{ }x_{n}=\lim \inf x_{n}=\sup_{n\geq 0}\inf_{p\geq
n}x_{p}
\end{equation*}
\bigskip

\noindent \textbf{(3)} Justifier l'existence de la limite de $z_{n}$ appel\'ee limite sup\'erieure de la suite $(x_{n})_{n\geq 0}$ d\'efinie par $
\lim \sup x_{n}$ ou $\overline{\lim }$ $x_{n},$  et elle est \'egale \`a 
\begin{equation*}
\overline{\lim }\text{ }x_{n}=\lim \sup x_{n}=\inf_{n\geq 0}\sup_{p\geq
n}x_{p}
\end{equation*}
\bigskip

\noindent \textbf{(4)} Establir que  
\begin{equation*}
-\liminf x_{n}=\limsup (-x_{n})\noindent \text{ \ \ et \ }-\limsup
x_{n}=\liminf (-x_{n}).
\end{equation*}
\bigskip

\newpage 
\noindent \textbf{(5)} 

Montrer que la limite sup\'erieure est sous- additive et la limite inf\'erieure est sur- additive, i.e. :  pour les deux suites
$(s_{n})_{n\geq 0}$ et  $(t_{n})_{n\geq 0}$ 
\bigskip

\begin{equation*}
\limsup (s_{n}+t_{n})\leq \limsup s_{n}+\limsup t_{n}
\end{equation*}
et
\begin{equation*}
\lim \inf (s_{n}+t_{n})\leq \lim \inf s_{n}+\lim \inf t_{n}
\end{equation*}

\noindent \textbf{(6)} D\'eduire de (1) que si
\begin{equation*}
\lim \inf x_{n}=\lim \sup x_{n},
\end{equation*}%
alors  $(x_{n})_{n\geq 0}$ a une limite et 
\begin{equation*}
\lim x_{n}=\lim \inf x_{n}=\lim \sup x_{n}
\end{equation*}

\bigskip

\noindent \textbf{Exercise 2.} Points d'accumulation de la suite $(x_{n})_{n\geq 0}$.\newline

\noindent \textbf{(a)} Montrer que si $\ell _{1}$=$\lim \inf x_{n}$ et $\ell
_{2}=\lim \sup x_{n}$ sont des points d' accumulation $(x_{n})_{n\geq 0}.
$ Montrer un cas et d'en d\'eduire le second en utilisant le point(3) de l'exercice 1.\newline
\bigskip
\noindent \textbf{(b)} Montrer que $\ell _{1}$ est le point d'accumulation plus petit $(x_{n})_{n\geq 0}$ et $\ell _{2}$ est le plus grand.
(De m\^eme, montrer un cas et d'en d\'eduire le second en utilisant le point(3)de l'exercice 1).\newline
\bigskip
\noindent \textbf{(c)} D\'eduire de (a) que si $(x_{n})_{n\geq 0}$ a une limite $\ell ,$ alors elle est \'egal au point d'accumulation unique et donc,
\begin{equation*}
\ell =\overline{\lim }\text{ }x_{n}=\lim \sup x_{n}=\inf_{n\geq
0}\sup_{p\geq n}x_{p}.
\end{equation*}
\bigskip
\noindent \textbf{(d)} Combinez ce r\'esultat au point\textbf{(6)} de l'exercice 1 pour montrer qu'une suite $(x_{n})_{n\geq 0}$ de $\overline{\mathbb{R}}
$ a une limite $\ell $ in $\overline{\mathbb{R}}$ si et seulement si \ $\lim \inf
x_{n}=\lim \sup x_{n}$ et donc
\begin{equation*}
\ell =\lim x_{n}=\lim \inf x_{n}=\lim \sup x_{n}
\end{equation*}
\bigskip
\newpage

\noindent \textbf{Exercice 3. }  Soit $(x_{n})_{n\geq 0}$  une suite  croissante de $\overline{\mathbb{R}}$. Etudier sa limite sup\'erieure et la limite inf\'erieure et en d\'eduire que
\begin{equation*}
\lim x_{n}=\sup_{n\geq 0}x_{n}.
\end{equation*}
\bigskip
\noindent En d\'eduire que pour une s\'equence non croissante $(x_{n})_{n\geq 0}$
de $\overline{\mathbb{R}},$%
\begin{equation*}
\lim x_{n}=\inf_{n\geq 0}x_{n}.
\end{equation*}

\bigskip

\noindent \textbf{Exercice 4.} (Les crit\`eres de convergence)\newline

\noindent \textbf{Crit\`ere 1 .} Soit $(x_{n})_{n\geq 0}$ une suite r\'eelle de $\overline{\mathbb{R}}$ et un nombre r\'eel $\ell \in \overline{\mathbb{R}}$ tel que: Pour toute sous-suite de $(x_{n})_{n\geq 0}$ a \'egalement une sous-suite( c'est une sous-suite de $(x_{n})_{n\geq 0}$ )  qui converge vers $\ell .$
Alors, la limite de $(x_{n})_{n\geq 0}$ existe et est \'egale \`a $\ell .$\newline
\bigskip
\noindent \textbf{Crit\`ere 2. } Intersections maximales et intersections minimales. \newline
\bigskip
\noindent Soit $(x_{n})_{n\geq 0}$ une suite dans $\overline{\mathbb{R}}$ et  deux nombres r\'eels $a$ et $b$ tel que  $a<b.$
Nous d\'efinissons%
\bigskip
\begin{equation*}
\nu _{1}=\left\{ 
\begin{array}{cc}
\inf  & \{n\geq 0,x_{n}<a\} \\ 
+\infty  & \text{if (}\forall n\geq 0,x_{n}\geq a\text{)}%
\end{array}%
\right. .
\end{equation*}%
\bigskip
Si $\nu _{1}$ est fini , soit
\begin{equation*}
\nu _{2}=\left\{ 
\begin{array}{cc}
\inf  & \{n>\nu _{1},x_{n}>b\} \\ 
+\infty  & \text{if (}n>\nu _{1},x_{n}\leq b\text{)}%
\end{array}%
\right. .
\end{equation*}%
.
\bigskip
\noindent Tant que les $\nu _{j}'s$ sont finis,  nous pouvons d\'efinir pour $\nu
_{2k-2}(k\geq 2)$
\bigskip
\begin{equation*}
\nu _{2k-1}=\left\{ 
\begin{array}{cc}
\inf  & \{n>\nu _{2k-2},x_{n}<a\} \\ 
+\infty  & \text{if (}\forall n>\nu _{2k-2},x_{n}\geq a\text{)}%
\end{array}%
\right. .
\end{equation*}%
et pour $\nu _{2k-1}$ fini, 
\begin{equation*}
\nu _{2k}=\left\{ 
\begin{array}{cc}
\inf  & \{n>\nu _{2k-1},x_{n}>b\} \\ 
+\infty  & \text{si(}n>\nu _{2k-1},x_{n}\leq b\text{)}%
\end{array}%
\right. .
\end{equation*}
\bigskip
\noindent Nous nous arr\^etons une fois qu'un $\nu _{j}$ est $+\infty$. Si $\nu
_{2j}$ est fini, alors 
\begin{equation*}
x_{\nu _{2j}}-x_{\nu _{2j-1}}>b-a. 
\end{equation*}
\bigskip
\noindent On dit alors : par ce passage de $x_{\nu _{2j-1}}$ \`a $x_{\nu
_{2j}},$ nous avons r\'ealis\'e un passage (vers le haut ) du segment $[a,b]$
appel\'e \textit{up-crossings}. De m\^eme, si l'un des $\nu _{2j+1}$
est fini, alors  le segment $[x_{\nu _{2j}},x_{\nu _{2j+1}}]$ est une intersection de point minimal (downcrossing) pour le segment $[a,b].$ Soit
\begin{equation*}
D(a,b)=\text{ nombre d'intersection de points maximums }[a,b]\text{.}
\end{equation*}

\bigskip

\noindent \textbf{(a)} Quelle est la valeur de $D(a,b)$ si \ $\nu _{2k}$ est fini et $\nu
_{2k+1}$ infini.\newline

\noindent \textbf{(b)} Quelle est la valeur de $D(a,b)$ si \ $\nu _{2k+1}$ est fini et  $\nu
_{2k+2}$ infini.\newline

\noindent \textbf{(c)} Quelle est la valeur de $D(a,b)$ si tous les $\nu _{j}'s$ sont finis.%
\newline

\noindent \textbf{(d)} Montrer que $(x_{n})_{n\geq 0}$ admet une limite si pour tout $a<b,$ $D(a,b)<\infty.$\newline

\noindent \textbf{(e)} Montrer que $(x_{n})_{n\geq 0}$ admet une limite si pour tout $a<b,$ $(a,b)\in \mathbb{Q}^{2},D(a,b)<\infty .$\newline

\bigskip

\noindent \textbf{Exercice 5. } (Crit\`ere de Cauchy). Soit $(x_{n})_{n\geq 0}$ $\mathbb{R}$ une suite de (\textbf{nombres r\'eels}).\newline

\noindent \textbf{(a)} montrer que si $(x_{n})_{n\geq 0}$ est de Cauchy,
alors elle admet un point accumulation unique $\ell \in \mathbb{R}$ qui est sa limite.\newline

\noindent \textbf{(b)} Montrer que si une suite  $(x_{n})_{n\geq 0}\subset 
\mathbb{R}$ \ converge vers $\ell \in \mathbb{R},$ alors, elle est de Cauchy.%
\newline

\noindent \textbf{(c)} D\'eduire le crit\`ere de Cauchy pour les suites de nombres r\'eels .
\newpage

\begin{center}
\textbf{SOLUTIONS}
\end{center}

\noindent \textbf{Exercice 1}.\newline

\noindent \textbf{Question (1) :}. Il est \'evident que :%
\begin{equation*}
\underset{p\geq n}{\inf }x_{p}\leq x_{n}\leq \underset{p\geq n}{\sup }x_{p},
\end{equation*}

\noindent \'etant donn\'e que $x_{n}$ est un \'el\'ement de $\left\{ x_{n},x_{n+1},...\right\} $ sur laquelle nous prenons le supremum ou infinimum.
\newline

\noindent \textbf{Question (2) :}. Soit $y_{n}=\underset{p\geq 0}{\inf }
x_{p}=\underset{p\geq n}{\inf }A_{n},$ où $A_{n}=\left\{
x_{n},x_{n+1},...\right\} $ est une suite non croissante  de l'ensemble : $\forall n\geq 0$,
\begin{equation*}
A_{n+1}\subset A_{n}.
\end{equation*}

\noindent Ainsi l'infinimum dans $A_{n}$ est croissant. Si $y_{n}$ est croissante dans $%
\overline{\mathbb{R}},$ sa limite est au dessus des bornes finie ou infinie. Ainsi
\begin{equation*}
y_{n}\nearrow \underline{\lim }\text{ }x_{n},
\end{equation*}%
est un nombre fini ou infini.\newline

\noindent \textbf{Question (3) :}. Nous d\'emontrons aussi que $z_{n}=\sup A_{n}$ d\'ecroit et $z_{n}\downarrow \overline{\lim }$ $x_{n}$.\newline

\noindent \textbf{Question (4) \label{qst4}:}. Nous rappelons que  
\begin{equation*}
-\sup \left\{ x,x\in A\right\} =\inf \left\{ -x,x\in A\right\}. 
\end{equation*}

\noindent lequel , nous permet d'\'ecrire  
\begin{equation*}
-\sup A=\inf -A.
\end{equation*}

\noindent Ainsi,

\begin{equation*}
-z_{n}=-\sup A_{n}=\inf -A_{n} = \inf \left\{-x_{p},p\geq n\right\}..
\end{equation*}

\noindent Le terme de la droite tend vers $-\overline{\lim}\ x_{n}$ et celui de la gauche vers $\underline{\lim} \ -x_{n}$ et donc 

\begin{equation*}
-\overline{\lim}\ x_{n}=\underline{\lim }\ (-x_{n}).
\end{equation*}

\bigskip \noindent De m\^eme, nous montrons que:
\begin{equation*}
-\underline{\lim } \ (x_{n})=\overline{\lim} \ (-x_{n}).
\end{equation*}

\noindent 

\noindent \textbf{Question (5)}. Ces propri\'et\'es viennent des formules suivante ,où $A\subseteq \mathbb{R},B\subseteq \mathbb{R}$ :%
\begin{equation*}
\sup \left\{ x+y,A\subseteq \mathbb{R},B\subseteq \mathbb{R}\right\} \leq
\sup A+\sup B.
\end{equation*}

\noindent De ce fait : 
\begin{equation*}
\forall x\in \mathbb{R},x\leq \sup A
\end{equation*}

\noindent et
\begin{equation*}
\forall y\in \mathbb{R},y\leq \sup B.
\end{equation*}

\noindent Ainsi,
\begin{equation*}
x+y\leq \sup A+\sup B,
\end{equation*}

\noindent où 
\begin{equation*}
\underset{x\in A,y\in B}{\sup }x+y\leq \sup A+\sup B.
\end{equation*}%
De m\^eme,%
\begin{equation*}
\inf (A+B\geq \inf A+\inf B.
\end{equation*}

\noindent De ce fait:

\begin{equation*}
\forall (x,y)\in A\times B,x\geq \inf A\text{ and }y\geq \inf B.
\end{equation*}

\noindent Ainsi, 
\begin{equation*}
x+y\geq \inf A+\inf B.
\end{equation*}

\noindent Ains,
\begin{equation*}
\underset{x\in A,y\in B}{\inf }(x+y)\geq \inf A+\inf B
\end{equation*}

\noindent \textbf{Application}.\newline

\begin{equation*}
\underset{p\geq n}{\sup } \ (x_{p}+y_{p})\leq \underset{p\geq n}{\sup } \ x_{p}+\underset{p\geq n}{\sup } \ y_{p}.
\end{equation*}

\noindent Toutes ces suites sont non croissantes. Prenant l'infinimum, nous obtenons la limite superieur:

\begin{equation*}
\overline{\lim }\text{ }(x_{n}+y_{n})\leq \overline{\lim }\text{ }x_{n}+%
\overline{\lim }\text{ }x_{n}.
\end{equation*}

\bigskip

\noindent \textbf{Question (6) :} Mettons

\begin{equation*}
\underline{\lim } \ x_{n}=\overline{\lim } \ x_{n},
\end{equation*}

\noindent Comme : 
\begin{equation*}
\forall x\geq 1,\text{ }y_{n}\leq x_{n}\leq z_{n},
\end{equation*}%

\begin{equation*}
y_{n}\rightarrow \underline{\lim} \ x_{n}
\end{equation*}%

\noindent et

\begin{equation*}
z_{n}\rightarrow \overline{\lim } \ x_{n},
\end{equation*}

\noindent Nous appliquons le th\'eor\`eme de Sandwich  pour conclure que la limite de $x_{n}$ existe et :

\begin{equation*}
\lim \text{ }x_{n}=\underline{\lim }\text{ }x_{n}=\overline{\lim }\text{ }%
x_{n}.
\end{equation*}

\bigskip 
\noindent \textbf{Exercice 2}.\newline

\noindent \textbf{Question (a).}\\

\noindent La question (4) de l'exercice 1 \'etant prouv\'ee, il suffit maintenant de montrer cette propri\'et\'e pour l'une des limites.Considerons la limite sup\'erieur et les trois cas d'\'etude suivants:\\

\noindent \textbf{le cas de la limite sup\'erieure finie} :

\begin{equation*}
\underline{\lim \text{ }}x_{n}=\ell \text{ finite.}
\end{equation*}

\noindent Par d\'efinition, 
\begin{equation*}
z_{n}=\underset{p\geq n}{\sup }x_{p}\downarrow \ell .
\end{equation*}

\noindent Donc: 
\begin{equation*}
\forall \varepsilon >0,\exists (N(\varepsilon )\geq 1),\forall p\geq
N(\varepsilon ),\ell -\varepsilon <x_{p}\leq \ell +\varepsilon .
\end{equation*}

\noindent Prenons le moins  que celui l\`a :

\begin{equation*}
\forall \varepsilon >0,\exists n_{\varepsilon }\geq 1:\ell -\varepsilon
<x_{n_{\varepsilon }}\leq \ell +\varepsilon.
\end{equation*}

\noindent Nous pouvons construire une sous-suite  convergeant vers $\ell$.\\

\noindent Soit $\varepsilon =1:$%
\begin{equation*}
\exists N_{1}:\ell -1<x_{N_{1}}=\underset{p\geq n}{\sup }x_{p}\leq \ell +1.
\end{equation*}

\noindent Mais si 
\begin{equation}
z_{N_{1}}=\underset{p\geq n}{\sup }x_{p}>\ell -1, \label{cc}
\end{equation}

\noindent il existe s\^urement un $n_{1}\geq N_{1}$ tel que 
\begin{equation*}
x_{n_{1}}>\ell -1.
\end{equation*}

\noindent sinon nous nous pourrons  avoir  
\begin{equation*}
( \forall p\geq N_{1},x_{p}\leq \ell -1\ ) \Longrightarrow \sup \left\{
x_{p},p\geq N_{1}\right\} =z_{N_{1}}\geq \ell -1,
\end{equation*}
lequel est contraire \`a (\ref{cc}). Donc, il existe $n_{1}\geq N_{1}$ tel que 
\begin{equation*}
\ell -1<x_{n_{1}}\leq \underset{p\geq N_{1}}{\sup }x_{p}\leq \ell -1.
\end{equation*}

\noindent i.e.

\begin{equation*}
\ell -1<x_{n_{1}}\leq \ell +1.
\end{equation*}

\noindent Nous progressons avec l'\'etape $\varepsilon =\frac{1}{2}$ et nous consid\'erons la suite 
 $(z_{n})_{n\geq n_{1}}$dont la limite reste \'egale \`a $\ell$. Donc, il existe $N_{2}>n_{1}:$
\begin{equation*}
\ell -\frac{1}{2}<z_{N_{2}}\leq \ell -\frac{1}{2}.
\end{equation*}

\noindent Nous d\'eduisons comme pr\'ec\'edemment que $n_{2}\geq N_{2}$ tel que
\begin{equation*}
\ell -\frac{1}{2}<x_{n_{2}}\leq \ell +\frac{1}{2}
\end{equation*}

\noindent Avec $n_{2}\geq N_{1}>n_{1}$.\\

\noindent Ensuite, nous mettons $\varepsilon =1/3,$ il existera $N_{3}>n_{2}$ tel que 
\begin{equation*}
\ell -\frac{1}{3}<z_{N_{3}}\leq \ell -\frac{1}{3}
\end{equation*}

\noindent et nous pourrons voir un $n_{3}\geq N_{3}$ tel que

\begin{equation*}
\ell -\frac{1}{3}<x_{n_{3}}\leq \ell -\frac{1}{3}.
\end{equation*}

\noindent Pas \`a pas, nous d\'eduisons l'existence de $ x_{n_{1}},x_{n_{2}},x_{n_{3}},...,x_{n_{k}},...$ avec $n_{1}<n_{2}<n_{3}%
\,<...<n_{k}<n_{k+1}<...$ tel que

$$
\forall k\geq 1, \ell -\frac{1}{k}<x_{n_{k}}\leq \ell -\frac{1}{k},
$$

\noindent i.e.

\begin{equation*}
\left\vert \ell -x_{n_{k}}\right\vert \leq \frac{1}{k}.
\end{equation*}

\noindent Ce qui implique que: 
\begin{equation*}
x_{n_{k}}\rightarrow \ell 
\end{equation*}

\noindent Conclusion : $(x_{n_{k}})_{k\geq 1}$ est bien une sous-suite tel que $n_{k}<n_{k+1}$ pour tout $k \geq 1$ 
et il converge vers $\ell$, qui est alors un point d'accumulation.\\

\noindent \textbf{Le cas de la limite sup\'erieur \'egale \`a $+\infty$} : 
$$
\overline{\lim} \text{ } x_{n}=+\infty.
$$
\noindent Lorsque $z_{n}\uparrow +\infty ,$ nous avons : $\forall k\geq 1,\exists
N_{k}\geq 1,$ 
\begin{equation*}
z_{N_{k}}\geq k+1.
\end{equation*}

\noindent Pour $k=1$, soit $z_{N_{1}}=\underset{p\geq N_{1}}{\inf }%
x_{p}\geq 1+1=2.$ Donc il existe 
\begin{equation*}
n_{1}\geq N_{1}
\end{equation*}%
tel que :%
\begin{equation*}
x_{n_{1}}\geq 1.
\end{equation*}

\noindent Pour $k=2:$ considerons les suites $(z_{n})_{n\geq n_{1}+1}.$
De la m\^eme mani\`ere, nous trouvons  
\begin{equation*}
n_2 \geq n_{1}+1
\end{equation*}%
\noindent et
\begin{equation*}
x_{n_{2}}\geq 2.
\end{equation*}

\noindent Etape par \'etape, nous trouvons pour tout $k\geq 3$, et $n_{k}\geq n_{k-1}+1$ tel que 
\begin{equation*}
x_{n_{k}}\geq k.
\end{equation*}

\noindent Ce qui conduit \`a $x_{n_{k}}\rightarrow +\infty $ comme $k\rightarrow +\infty $.\\

\noindent \textbf{ Cas de la limit  superieur \'egale \`a $-\infty$} : 

$$
\overline{\lim }x_{n}=-\infty.
$$

\noindent Cela implique : $\forall k\geq 1,\exists N_{k}\geq 1,$ tel que %
\begin{equation*}
z_{n_{k}}\leq -k.
\end{equation*}

\noindent Pour $k=1,\exists n_{1}$ tel que%
\begin{equation*}
z_{n_{1}}\leq -1.
\end{equation*}
Mais
\begin{equation*}
x_{n_{1}}\leq z_{n_{1}}\leq -1
\end{equation*}

\noindent Soit $k=2$. Considerons $\left( z_{n}\right) _{n\geq
n_{1}+1}\downarrow -\infty .$ Il existera $n_{2}\geq n_{1}+1:$%
\begin{equation*}
x_{n_{2}}\leq z_{n_{2}}\leq -2
\end{equation*}

\noindent Etape par \'etape, nous trouvons $n_{k1}<n_{k+1}$ de telle sorte que $x_{n_{k}}<-k$ pour tout $k$ plus grand que $1$. Donc,
\begin{equation*}
x_{n_{k}}\rightarrow +\infty 
\end{equation*}

\bigskip

\noindent \textbf{Question (b).}\\

\noindent Soit $\ell$ un point d'accumulation de $(x_n)_{n \geq 1}$, la limite d'un de ses suites $(x_{n_{k}})_{k \geq 1}$. Nous avons 
$$
y_{n_{k}}=\inf_{p\geq n_k} \ x_p \leq x_{n_{k}} \leq  \sup_{p\geq n_k} \ x_p=z_{n_{k}}
$$

\noindent Le terme de gauche est une sous suite de $(y_n)$ tendant vers la limite inf\'erieure et le c\^ot\'e droit est une suous suite de $(z_n)$ tendant vers la limitesuperieur. Donc, nous obtenons :

$$
\underline{\lim} \ x_{n} \leq \ell \leq \overline{\lim } \ x_{n},
$$

\noindent ce qui montre que $\underline{\lim} \ x_{n}$ est le plus petit point d'accumulation et $\overline{\lim } \ x_{n}$ est le plus grand.\\

\noindent \textbf{Question (c).} Si la suite $(x_n)_{n \geq 1}$  admet une limite $\ell$, alors cette limite est celle de toutes ses sous-suites, dons les sous-suite  tend vers les limites sup\'erieur ou inf\'erieu. Ce qui r\'epond \`a la question (b). \\

\noindent \textbf{Question (d).} NNous r\'epondons \`a cette question en combinant le point (d) de cette exercice et le point (\textbf{6}) de l'exercice \textbf{1}.\\ 

\noindent \textbf{Exercice 3}. Soit $(x_{n})_{n\geq 0}$ une suite non-d\'ecroissante, nous avons:
\begin{equation*}
z_{n}=\underset{p\geq n}{\sup} \ x_{p}=\underset{p\geq 0}{\sup} \ x_{p},\forall
n\geq 0.
\end{equation*}

\noindent Pourquoi? Parce que par la croissance nous avons,
\begin{equation*}
\left\{ x_{p},p\geq 0\right\} =\left\{ x_{p},0\leq p\leq n-1\right\} \cup
\left\{ x_{p},p\geq n\right\}
\end{equation*}

\bigskip

\noindent Etant donn\'e que tous les \'el\'ements de $\left\{ x_{p},0\leq p\leq
n-1\right\} $ sont plus petits que ceux de $\left\{ x_{p},p\geq n\right\} ,$
le supremum est atteint en $\left\{ x_{p},p\geq n\right\} $ et donc
\begin{equation*}
\ell =\underset{p\geq 0}{\sup } \ x_{p}=\underset{p\geq n}{\sup }x_{p}=z_{n}
\end{equation*}%
Ainsi
\begin{equation*}
z_{n}=\ell \rightarrow \ell .
\end{equation*}

\noindent Nous avons aussi $y_n=\inf \left\{ x_{p},0\leq p\leq n\right\}=x_n$ qui est une suite non d\'ecroissante et ainsi converge vers
$\ell =\underset{p\geq 0}{\sup } \ x_{p}$. \\

\bigskip

\noindent \textbf{Exercice 4}.\\

\noindent Soit $\ell \in \overline{\mathbb{R}}$ ayant la propri\'et\'e indiqu\'ee. Soit $\ell ^{\prime }$ un point d'accumulation donn\'e.%
\begin{equation*}
 \left( x_{n_{k}}\right)_{k\geq 1} \subseteq \left( x_{n}\right) _{n\geq 0}%
\text{ tel que }x_{n_{K}}\rightarrow \ell ^{\prime}.
\end{equation*}

\noindent Par hypoth\`ese cette suite $\left( x_{n_{K}}\right) $
a \`a son tour une sou-suite $\left( x_{n_{\left( k(p)\right) }}\right)_{p\geq 1} $ tel que $x_{n_{\left( k(p)\right) }}\rightarrow
\ell $ as $p\rightarrow +\infty $.\newline

\noindent Mais comme une sous-suite de $\left( x_{n_{\left( k\right)
}}\right) ,$ 
\begin{equation*}
x_{n_{\left( k(\ell )\right) }}\rightarrow \ell ^{\prime }.
\end{equation*}%
Ainsi
\begin{equation*}
\ell =\ell ^{\prime}.
\end{equation*}

\noindent En appliquant cela \`a la limite sup\'erieure et \`a la limite inf\'erieure, nous avons:%
\begin{equation*}
\overline{\lim} \ x_{n}=\underline{\lim}\ x_{n}=\ell.
\end{equation*}

\noindent Et donc $\lim x_{n}$ existe et est \'egal \`a $\ell$.\\

\noindent \textbf{Exercice 5}.\\

\noindent \textbf{Question (a)}. Si $\nu _{2k}$ fini et $\nu _{2k+1}$ infini, il a alors exactement $k$ un point d'intersection  sup\'erieur: 
$[x_{\nu_{2j-1}},x_{\nu _{2j}}]$, $j=1,...,k$ : $D(a,b)=k$.\\

\noindent \textbf{Question (b)}. Si$\nu _{2k+1}$ fini et $\nu _{2k+2}$ infini, il a alors exactement $k$ un point d'intersection  sup\'erieur:
$[x_{\nu_{2j-1}},x_{\nu_{2j}}]$, $j=1,...,k$ : $D(a,b)=k$.\\

\noindent \textbf{Question (c)}. SI tous les $\nu_{j}'s$ sont finis,alors , il existent  un nombre inifini de points d'intersection  sup\'erieur: 
$[x_{\nu_{2j-1}},x_{\nu_{2j}}]$, $j\geq 1k$ : $D(a,b)=+\infty$.\\

\noindent \textbf{Question (d)}. Supposons qu'il existe $a < b$  des nombres rationnels tels que $D(a,b)=+\infty$. 
Alors tous les $\nu _{j}'s$ sont finis. La sous-suite $x_{\nu_{2j-1}}$ est strictement inf\'erieure $a$. 
Donc, sa limite inf\'erieure est inf\'erieure \`a $a$. Cette limite inf\'erieure est un point d'accumulation de la suite $(x_n)_{n\geq 1}$, 
il en est de plus $\underline{\lim}\ x_{n}$, qui est inf\'erieur \`a $a$.\\

\noindent De m\^eme, la sous-suite $x_{\nu_{2j}}$ est strictement inf\'erieure $b$. Donc, la limite sup\'erieure est au-dessus de $a$. 
Cette limite sup\'erieure est un point d'accumulation de la suite $(x_n)_{n\geq 1}$, donc elle est inf\'erieure \`a $\overline{\lim}\ x_{n}$, 
qui est directement au-dessus $b$. Qui conduit \`a:

$$
\underline{\lim}\ x_{n} \leq a < b \leq \overline{\lim}\ x_{n}. 
$$

\noindent Cela implique que la limite de $(x_n)$ n'existe pas. En revanche, nous venons de prouver que la limite de $(x_n)$ existe, 
Entre-temps pour tous les nombres r\'eels $a$ et $b$ tel que $a<b$, $D(a,b)$ est fini.\\

\noindent Maintenant , supposons que la limite de $(x_n)$ n'existe pas. Alors,
$$
\underline{\lim}\ x_{n} < \overline{\lim}\ x_{n}. 
$$

\noindent Nous pouvons alors trouver deux rationnels $a$ et $b$ tel que $a<b$ et un nombre $\epsilon$ tel que $0<\epsilon$, tous les

$$
\underline{\lim}\ x_{n} < a-\epsilon < a < b < b+\epsilon <  \overline{\lim}\ x_{n}. 
$$

\noindent Si $\underline{\lim}\ x_{n} < a-\epsilon$, nous pouvons revenir de la question \textbf{(a)} de exercice \textbf{2} et construire une suite $(x_n)$
ce qui tend vers $\underline{\lim}\ x_{n}$ tout en restant en dessous de $a-\epsilon$. De m\^eme, si $b+\epsilon < \overline{\lim}\ x_{n}$, 
nous pouvons cr\'eer une suite de $(x_n)$ ce qui tend vers $\overline{\lim}\ x_{n}$ tout en restant au-dessus $b+\epsilon$. 
Il est \'evident avec ces deux suites que nous pourrions d\'efinir avec ces deux suites tous $\nu_{j}$ fini et donc $D(a,b)=+\infty$.\\

\noindent Nous venons de montrer par l'absurde que si tous les $D(a,b)$ sont finis pour tous les rationnels $a$ et $b$  tel que $a<b$, 
alors,la limite de $(x)n)$ existe.\\

\noindent \textbf{Exercice 5}. crit\`ere de Cauchy dans $\mathbb{R}$.\\

\noindent Supposons que la suite est Cauchy, $i.e.$,
$$
\lim_{(p,q)\rightarrow (+\infty,+\infty)} \ (x_p-x_q)=0.
$$

\noindent Ensuite, soient $x_{n_{k,1}}$ et $x_{n_{k,2}}$  deux sous-suites  convergeant respectivement vers $\ell_1=\underline{\lim}\ x_{n}$ et $\ell_2=\overline{\lim}\ x_{n}$. Alors

$$
\lim_{(p,q)\rightarrow (+\infty,+\infty)} \ (x_{n_{p,1}}-x_{n_{q,2}})=0.
$$

\noindent En premi\`ere location $p\rightarrow +\infty$,nous avons

$$
\lim_{q\rightarrow +\infty} \ \ell_1-x_{n_{q,2}}=0.
$$

\noindent Ce qui montre que la limite $\ell_1$ est finie, sinon $\ell_1-x_{n_{q,2}}$ demeurerait infinie et ne tendent pas vers $0$. 
En \'echangeant les r\^oles de $p$ et $q$, nous avons aussi que $\ell_2$ est fini.\\

\noindent Enfin, en laissant $q\rightarrow +\infty$, dans la derni\`ere \'equation , nous obtenons
$$
\ell_1=\underline{\lim}\ x_{n}=\overline{\lim}\ x_{n}=\ell_2.
$$

\noindent ce qui prouve l'existence d'une limite finie de la suite $(x_n)$.\\

\noindent Supposons maintenant que la limite finie $\ell$ de $(x_n)$ existe. Alors

$$
\lim_{(p,q)\rightarrow (+\infty,+\infty)} \ (x_p-x_q)=\ell-\ell=0.
$$0
\noindent Ce qui montre que la suite est de Cauchy.\\

\bigskip
\section{Miscelleanuous facts} \label{funct.facts}

\bigskip \noindent 
\begin{lemma} \label{funct.facts.lem01}. For any $a\in \mathbb{R},$ 
\begin{equation*}
\left\vert e^{ia}-1\right\vert =\sqrt{2(1-\cos a)}\leq 2\left\vert \sin
(a/2)\right\vert \leq 2\left\vert a/2\right\vert ^{\delta }.
\end{equation*}
\end{lemma}

\bigskip \noindent \textbf{Proof}. This is easy for $\left\vert a/2\right\vert >1.$ Indeed
for $\delta >0,\left\vert a/2\right\vert ^{\delta }>0$ and%
\begin{equation*}
2\left\vert \sin (a/2)\right\vert \leq 2\leq 2\left\vert a/2\right\vert
^{\delta }
\end{equation*}

\bigskip \noindent Now for $\left\vert a/2\right\vert >1,$ we have the
expansion

\begin{eqnarray*}
2(1-\cos a) &=&a^{2}-\sum\limits_{k=2}^{\infty }(-1)^{2}\frac{a^{2k}}{(2k)!}%
=x^{2}-2\sum\limits_{k\geq 2,k\text{ }even}^{\infty }\frac{a^{2k}}{(2k)!}-%
\frac{a^{2(k+1)}}{(2(k+1))!} \\
&=&a^{2}-2x^{2(k+1)}\sum\limits_{k\geq 2,k\text{ }even}^{\infty }\frac{1}{%
(2k)!}\left\{ \frac{1}{a^{2}}-\frac{1}{(2k+1)((2k+2)...(2k+k)}\right\} .
\end{eqnarray*}

\bigskip \noindent For each $k\geq 2,$ for $\left\vert a/2\right\vert <1,$%
\begin{equation*}
\left\{ \frac{1}{a^{2}}-\frac{1}{(2k+1)((2k+2)...(2k+k)}\right\} \geq
\left\{ \frac{1}{4}-\frac{1}{(2k+1)((2k+2)...(2k+k)}\right\} \geq 0.
\end{equation*}

\bigskip \noindent Hence 
\begin{equation*}
2(1-\cos a)\leq a^{2}.
\end{equation*}

\bigskip \noindent But for $\left\vert a/2\right\vert ,$ the function $%
\delta \hookrightarrow \left\vert a/2\right\vert ^{\delta }$ is
non-increasing $\delta ,0\leq \delta \leq 1$. Then%
\begin{equation*}
\sqrt{2(1-\cos a)}\leq \left\vert a\right\vert =2\left\vert a/2\right\vert
^{1}\leq 2\left\vert a/2\right\vert ^{\delta }.
\end{equation*}